\documentclass[10pt]{amsart}
\usepackage{amsmath}
\usepackage{amsxtra}
\usepackage{amscd}
\usepackage{amsthm}
\usepackage{amsfonts}
\usepackage{amssymb}
\usepackage{eucal} 
\newtheorem{lem}{Lemma}[section]%%%%%%%

\theoremstyle{definition}

\theoremstyle{remark}

\newcommand{\thmref}[1]{Theorem~\ref{#1}}
\newcommand{\secref}[1]{\S\ref{#1}}
\newcommand{\lemref}[1]{{ Lemma~\ref{#1}}}
\newcommand{\defref}[1]{Definition~\ref{#1}}
\newcommand{\propref}[1]{{ Proposition~\ref{#1}}}
\newcommand{\corref}[1]{Corollary~\ref{#1}}

\newcommand{\remref}[1]{Remark~\ref{#1}}
\newcommand{\exref}[1]{Example~\ref{#1}}

\newcommand{\aref}[1]{Assumption~\ref{#1}}
\newcommand{\nc}{\newcommand}
\nc{\renc}{\renewcommand}
\nc{\ssec}{\subsection}
\nc{\sssec}{\subsubsection}
\nc{\on}{\operatorname}
\nc{\remm}[1]{\<{remark} \ \lbl{#1} \>{remark}}

\renc\k{\mathbf{k}}
\nc\ol{\overline} 
\nc\wh{\widehat}
\nc\tboxtimes{\wt{\boxtimes}}

\nc{\Aa}{{\mathbb{A}}}
 \nc{\Gg}{{\mathbb{G}}}  

\def\g{\gamma}
\def\te{{\tilde{e}}}
\def\Ee{{\mathcal E}}
\nc{\Hh}{{\mathbb{H}}}

 \nc{\Nn}{{\mathbb{N}}}
\nc{\Pp}{{\mathbb{P}}}
\nc{\Rr}{{\mathbb{R}}}
\nc{\BV}{{\mathbb{V}}}
\nc{\BW}{{\mathbb{W}}}
\nc{\Zz}{{\mathbb{Z}}}
\nc{\Qq}{{\mathbb{Q}}}
\nc{\Ss}{{\mathbb{S}}}
\nc{\Cc}{{\mathbb{C}}}
\nc{\Ff}{{\mathbb{F}}}
\nc{\Oo}{{\mathcal{O}}}
\nc{\Mm}{{\mathcal{M}}}

\nc{\eo}{{\mathbf{\tau}}}
\nc{\dU}{{\overset{\bullet}{\bigcup}}{}}
\nc{\du}{\, {\overset{.}{\cup}}\, {}}
\nc{\dual}[1]{{\overset{\vee}{#1}{}}}

\nc{\cM}{{\check{\mathcal M}}{}}
 \nc{\oM}{{\overset{\circ}{\mathcal M}}{}}

 \nc{\fB}{{\mathfrak{B}}}
\nc{\tT}{{\widetilde{T}}}  

\nc{\fF}{{\mathcal{F}}}

\nc{\bb}{{\mathbf{b}}}
\nc{\bc}{{\mathbf{c}}}
\nc{\bd}{{\mathbf{d}}}
\nc{\be}{{\mathbf{e}}}
\nc{\bj}{{\mathbf{j}}}
\nc{\bn}{{\mathbf{n}}}

\nc{\bph}{{\mathbf{\phi}}}
\nc{\bp}{{\mathbf{p}}}
\nc{\bq}{{\mathbf{q}}}
\nc{\bF}{{\mathbf{F}}}
\nc{\bu}{{\mathbf{u}}}
\nc{\bv}{{\mathbf{v}}}
\nc{\bx}{{\mathbf{x}}}
\nc{\bh}{{\mathbf{h}}}
\nc{\bs}{{\mathbf{s}}}
\nc{\by}{{\mathbf{y}}}
\nc{\bw}{{\mathbf{w}}}
\nc{\bA}{{\mathbf{A}}}
\nc{\bK}{{\mathbf{K}}}
\nc{\bI}{{\mathbf{I}}}
\nc{\bB}{{\mathbf{B}}}
\nc{\bG}{{\mathbf{G}}}
\nc{\bC}{{\mathbf{C}}}
\nc{\bD}{{\mathbf{D}}}
\nc{\bP}{{\mathbf{P}}}
\nc{\bH}{{\mathbf{H}}}
\nc{\bM}{{\mathbf{M}}}
\nc{\bN}{{\mathbf{N}}}
\nc{\bV}{{\mathbf{V}}}
\nc{\bU}{{\mathbf{U}}}
\nc{\bL}{{\mathbf{L}}}
\nc{\bT}{{\mathbf{T}}}
\nc{\bW}{{\mathbf{W}}}
\nc{\bX}{{\mathbf{X}}}
\nc{\bY}{{\mathbf{Y}}}
\nc{\bZ}{{\mathbf{Z}}}
\nc{\bS}{{\mathbf{S}}}
\nc{\ba}{{\mathbf{a}}}

\nc{\sA}{{\mathsf{A}}}
\nc{\sB}{{\mathsf{B}}}
\nc{\sC}{{\mathsf{C}}} 
\nc{\sF}{{\mathsf{F}}}
\nc{\sG}{{\mathsf{G}}}
\nc{\sK}{{\mathsf{K}}}
\nc{\sM}{{\mathsf{M}}}
\nc{\sO}{{\mathsf{O}}}
\nc{\sQ}{{\mathsf{Q}}}
\nc{\sP}{{\mathsf{P}}}
\nc{\sZ}{{\mathsf{Z}}}
\nc{\sfp}{{\mathsf{p}}}
\nc{\sr}{{\mathsf{r}}}
\nc{\sg}{{\mathsf{g}}}
\nc{\sff}{{\mathsf{f}}}
\nc{\sfb}{{\mathsf{b}}}
\nc{\sfc}{{\mathsf{c}}} 

\nc{\tA}{{\widetilde{{A}}}}
\nc{\tD}{{\widetilde{{A}}}}
\nc{\tH}{{\widetilde{{A}}}}
\nc{\tB}{{\widetilde{{B}}}}
\nc{\tg}{{\widetilde{\mathfrak{g}}}}
\nc{\tG}{{\widetilde{G}}}
\nc{\TM}{{\widetilde{\mathbb{M}}}{}}
\nc{\tO}{{\widetilde{\mathsf{O}}}{}} 
\nc{\TZ}{{\tilde{Z}}}
\nc{\tx}{{\tilde{x}}}
\nc{\tf}{{\tilde{f}}}
\nc{\tz}{{\tilde{\zeta}}}
\nc{\tmu}{{\tilde{\mu}}}
\nc{\td}{{\tilde{d}}}
\nc{\tX}{{\widetilde{X}}}
  
   \nc{\E}{{\mathop{\operatorname{\rm E }}}}
 \nc{\Mor}{{\mathop{\operatorname{\rm Mor \,}}}}
\nc{\Ob}{{\mathop{\operatorname{\rm Ob \,}}}}
  \nc{\Sym}{{\mathop{\operatorname{\rm Sym}}}}
   \nc{\Aut}{{\mathop{\operatorname{\rm Aut}}}}
 \nc{\Spec}{{\mathop{\operatorname{\rm Spec}}}}
\nc{\Ker}{{\mathop{\operatorname{\rm Ker}}}}
 \nc{\dom}{{\mathop{\operatorname{\rm dom}}}}
\nc{\End}{{\mathop{\operatorname{\rm End}}}}
 \nc{\Hom}{\on{\Hom}} 
 \nc{\GL}{{\mathop{\operatorname{\rm GL}}}}
 \nc{\Id}{{\mathop{\operatorname{\rm Id}}}}
 \nc{\rk}{{\mathop{\operatorname{\rm rk}}}}
\nc{\irk}{{\mathop{\operatorname{\rm i-rk}}}}
 \nc{\length}{{\mathop{\operatorname{\rm length}}}}
\nc{\supp}{{\mathop{\operatorname{\rm supp}}}}
\nc{\val}{{\rm val}}
\def\abs#1{  | #1 | }
\nc{\valr}{{\rm val_{rv}}}
\nc{\valrv}{\valr}
\def\bdd{^{\rm bdd}}
 \nc{\Res}{{\rm res}}
 \nc{\Aff}{{\rm Aff}}
\nc{\res}{{\mathop{\operatorname{\rm res}}}}
\nc{\rad}{{\mathop{\operatorname{\rm rad}}}}

\def\tensor{{\otimes}}

\def\meet{\cap}
\def\union{\cup}

\def\si{\sigma}
\def\Sum{\Sigma}
\def\G{\Gamma}
\def\<{\begin}
 \def\>{\end}

\nc{\tV}{{\widetilde{{V}}}}
\nc{\hb}[1]{\hbox{#1}}
\def\rv{{\rm rv}}
 
\nc{\seq}[1]{\stackrel{#1}{\sim}}

\nc{\oeq}[1]{\underset{#1}{=}} 

\def\inv {{^{-1}}}
 \def\beq#1{\begin{equation} \label{#1}}   
\def\eeq{\end{equation}}

\def\Uu{\mathbb U}  \def\Vv{\mathbb V}

\def\uu{u}
 
 \def\e{\epsilon} 
 
\def\iso{\simeq}

\def\prf{\begin{proof}}
\def\eprf{\end{proof} }

\def\acl{\mathop{\rm acl}\nolimits}
 \def\dcl{\mathop{\rm dcl}\nolimits}

 \def\lbl#1{  \label{#1}  }

\author{Ehud Hrushovski, David Kazhdan}

\address{\newline Institute of Mathematics, the Hebrew
University of Jerusalem, Givat Ram, Jerusalem, 91904, Israel.} 
 
\email{ehud@math.huji.ac.il,kazhdan@math.huji.ac.il}

\nc{\Claim}[1]{{\noindent \bf Claim{ #1 }}}
  \nc{\pr}{{\mathop{\operatorname{\rm pr}}}}

\nc{\Mmm}{{(1+\Mm)}}
  \def\RVi{{\RV_\infty}}
   \def\RVp{{\RV^{>0}}}  
   \def\RVpi{{\RV^{>0}_\infty}}
 \def\km{{\k^*}}     
     
\def\Hf{{H_{\rm fin}}}

 \title{Integration in valued fields}

\begin{document}

\<{abstract}  
We develop a theory of  integration over valued fields of residue characteristic zero.   In particular we obtain new and base-field independent foundations for integration
over local fields of large residue characteristic, extending results of Denef,Loeser, Cluckers.
The method depends on an  analysis of definable sets up to definable bijections.  We obtain a precise description of the Grothendieck semigroup of such sets in terms of related groups over the residue
field and value group.  This yields new invariants of all definable bijections, as well as 
invariants of measure preserving bijections.  
  
\>{abstract}

\maketitle

 \def\Vvs{ { \Vv ^s}}
 \def\tri{{\mathsf B}}
 \def\tres{{\mathbf{t}}} 
 \def\fG{{\mathcal G}}
 \def\fC{{\mathcal C}}
\def\fD{{\mathcal D}}
\def\fGV{{\fG}_{\Vv}}
\def\ObG{{\rm Ob}_{\fG}}  
\def\ObGV{{\rm Ob}_{\fG_{\Vv}}}
\def\MorG{{\rm Mor}_{\fG}}
 %---------------------------------------------------------------------------------------
 \def\cB{{\mathfrak B}_{\rm cl}}
  \def\oB{{\mathfrak B}_{\rm open}}
 \nc{\av}{{\mathop{\operatorname{\rm av}}}}
 \def\RV{{\rm RV}}

\def\st{^{*}}

 \nc{\JG}{{\mathop{\operatorname{{ Jcb_{\G}  }}}}}
 \nc{\JRV}{{\mathop{\operatorname{{ Jcb_{\RV}  }}}}}
 \nc{\JVF}{{\mathop{\operatorname{{ Jcb  }}}}}

\def\rvi{ \rv^{-1}}
\def\A{\alpha}
\def\fCg{{\fC}^{RV}_\G}
 \def\fCgb{{\fC}^{RV}_{\G^{bdd}}}
 
\def\fCga{ {\G_A}}
 \def\fCgab{{\fC}^{RV}_{\G_A^{bdd}}}
  \def\emp{\emptyset}
  \renc\b{\beta}
 \renc\a{\alpha}
 
  \nc{\tF}{{\widetilde{F}}}
 % \nc{\tV}{{\widetilde{V}}}
    \nc{\tU}{{\widetilde{U}}}
 \def\m{\setminus}

\nc{\conjrv}{\underset{rv}{\sim}}
\def\rvinv{{$\conjrv$-invariant }}

\def\rvi{ \rv^{-1}}
\def\A{\alpha} 
 
%----------------t4
\def\tchi{\widetilde{\chi}}
\def\tC{\widetilde{C}}
\def\d{\delta}
\def\wX{\overline{X}}
\def\wY{\overline{Y}}
 
\def\Gpos{{\G^{\geq 0}}}
  %---------------------------------------------------------------------------------------
 \def\Om{{\Omega}}
 \def\om{{\omega}}

  \def\AA{{\mathcal A}}
 
\def\IA{{I \AA}} 
\def\IAk{I Fn(\k,\K(\RES))} 
\def\RR{{\mathcal R}}
\def\FF{{\mathcal F}}
\def\DD{{\mathcal D}}
\def\CC{{\mathcal C}}
 \newcommand{\isomto}{\overset{\rightarrow}{\cong}}

\def\fnr{{Fn^{\RV}}}

\def\rh{{\vartheta}}
 \nc{\SG}{{\mathop{\operatorname{{  K_+ }}}}} 
   \nc{\SGe}{{\mathop{\operatorname{{  K_+^{eff} }}}}}
 \def\SK{\SG}
 \def\SKe{\SGe}
 \nc{\K}{{\mathop{\operatorname{ { K} }}}}
 \nc{\Ke}{{\mathop{\operatorname{ { K^{eff}} }}}}
 \def\KM{\mu \K}
 
 \def\Isp{{\rm I_{sp}}}
  \def\Ispm{{\rm {\mu I_{sp}}}}
 \def\Ispmg {{\rm {\mu_\G I_{sp}}}}
 \def\VFdn{ \VF[n]/ (\dim < n)}
\def\Ispd{ \Isp'[n]}
  \def\Var{{\rm Var}}  \def\VarK{\Var_{Kz}}
  \def\stq{\subseteq}
  \def\Kz{K}

\def\RVlmn{{  \mu _l \RVni}}

 \def\Ivf {{\rm I_{vf}}}
 \def\Ivfm{{\rm I_{vf_\mu}}}
 
  \def\XXint#1#2#3{{ \setbox0=\hbox{$#1{#2#3}{\int}$} \vcenter{\hbox{$#2#3$}} \kern- .5\wd0}}
 
 \def\Xint#1{\mathchoice {\XXint \displaystyle \textstyle{#1}}%
 {\XXint\textstyle\scriptstyle{#1}}%
  {\XXint\scriptstyle\scriptscriptstyle{#1}}%
  {\XXint\scriptscriptstyle\scriptscriptstyle{#1}}%
 \!\int}
 
  \def\ints{{\Xint {\rm D}}}
   \def\inte{\Xint {\rm e}} 
 
 \def\pRV{{\rm pRV }}
 \def\uX{{\underline{X}}}     
 \def\L{{\mathbb L}}
    \def\uL{{\underline{\L}}}       
   \def\LY{{\L Y}}
   \def\LX{{\L X}}
      \def\uY{{\underline{Y}}}

\def\Ks{{\rm ! K}}
\def\fh{{\mathfrak h}}

    \def\lemm#1{    \begin{lem} \lbl{#1} }
\def\T{{\bf T}}
\def\bt{{\mathbf t}}
\def\V{{\rm V}}

\def\Tdcl{{\dcl^T}}
\def\Ttp{{{\rm tp}^T}}

\def\tP{{\widetilde{P}}}
\def\tQ{{\widetilde{Q}}}
\def\tY{{\widetilde{Y}}}
\def\tW{{\widetilde{W}}}

\def\fCgb{{\rm \G^{bdd}_A}} 

\def\RES{{\rm RES}}
\def\sd{\ltimes}
\def\vol{{\rm vol}}

\def\RCF{{\rm RCF}}
\def\RCVF{{\rm RCVF}} 
\def\DOAG{{\rm DOAG}}
\def\ACVF{{\rm ACVF}}

\def\ACVFR{{\rm ACVF^{R}}}
\def\TRV{{\rm TRV}}
\def\RES{{\rm RES}}
\def\VAL{{\rm VAL}}
\def\pCF{{\rm pCF}}

 \def\intm{\int}

 \tableofcontents

\def\VF{{ \mathbb{\rm VF}}}
\def\VFR{{\rm VFR}} 
\def\VFRm{{\rm VFR_\mu}} 
\def\VFn{{ \mathbb{\rm  VF [n]}}}
\def\idot{\cdot}
\def\VFni{{ \mathbb{\rm  VF [n,\idot]}}}
\def\RVn{{ \mathbb{\rm  RV[n]}}}
\def\RVni{{  \RV[n,\idot]}}
\def\mudot{ \overset{.}{\mu}}
 \def\VFm{  \mu {\VF}}
  \def\VFmg{  \mu_\G {\VF}}
  \def\VFma{  |\mu| {\VF}}
 \def\VFvol{ {\VF}_{\rm vol}}
  \def\VFvolg{ {\VF}_{\rm vol_\G}}
 \def\volVF{{\VF}_{\rm vol}}
 \def\RVm{  \mu {\RV}}      \def\RVmg{ \mu_\G {\RV}}
    
    \def\RVma{  |\mudot| {\RV}} 
 \def\RVvol{{\RV}_{\rm vol}} 
  \def\RVvolg{{\RV}_{\Gamma {\rm -vol}}}
 \def\RESvolg{{\RES}_{\Gamma {\rm -vol}}}
 \def\RVleqn{{\RV[\leq n,\idot]}}
\def\RVnm{{ \mu \RV[\leq n,\idot]}}
   \def\VFmn{\mu \VF [n]}
   \def\RESmg{{\mu_\G RES}}
   \def\VFgn{\mu_\G \VF [n]}
   \def\mG{{\mu  \Gamma }} 
\def\volG{{{\rm vol } \G}}
\def\sggf{\SG(\G^{fin})}
\def\sggfm{\SG({\mG}^{fin})}
  
     \def\RVmn{{\mu  \RV [n,\idot]}}
 
    \def\VFmb{ \VFm _{\rm bdd}}
 \def\VFmbg{ \VFm _{\G;\rm bdd}}
 \def\mVFr{ {\VFm} {\rv}}
 
 \def\RVmb {\RVm _{\rm bdd}}

\nc{\ev}{{\mathop{\operatorname{\rm ev }}}}

\pagebreak
\<{section}{Introduction}

Since Weil's {\em Foundations}, algebraic varieties have been understood independently of a particular base field; thus an algebraic group $G$ exists prior to the abstract or topological groups of points $G(F)$,
taken over various fields $F$.  For Hecke algebras, or other geometric objects whose definition requires integration, no comparable viewpoint exists.  One uses the  topology and measure theory of each local field separately; since a field $F$ has measure zero from the point of 
view of any nontrivial finite extension, at the foundational level there is no direct connection between
the objects obtained over different fields.   The main thrust of this paper is the development of 
a theory of integration over valued fields, that is geometric in the sense of Weil.   At present the theory covers local fields of   residue characteristic zero or, in applications, large
positive residue characteristic.  

Our approach to integration continues a line traced by Kontsevich, Denef-Loeser,  Loeser-Cluckers 
{   (cf. \cite{cluckersloeser})}.  In integration over non-archimedean local fields there are two sources
for the numerical values.  The first  is counting points of varieties over the residue field. 
Kontsevich explained that these numerical values can be replaced, with a gain of geometric information, by the 
isomorphism classes of the varieties themselves up to appropriate
transformations, or more precisely by their classes in a a certain Grothendieck ring.   This makes it possible to understand geometrically the changes in integrals upon unramified
base change.  In this aspect our approach is very similar.  The main difference
is a  slight generalization of the 
the notion of variety over the residue field, that allows us to avoid what amounted to a choice
of uniformizer in the previous theory.

The second source of the numbers  is the piecewise linear geometry of the value group.  We geometrize this ingredient too, obtaining a theory of integration taking
values in an entirely geometric ring, a tensor product of a Grothendieck ring of generalized
varieties over the residue field, and a Grothendieck ring of piecewise linear varieties over the value group.   

Viewed in this way, the integral is an invariant of measure preserving definable bijections.
We actually find all such invariants.   In addition we consider and determine all possible
invariants of definable bijections; we obtain in particular two Euler characteristics on 
definable sets, with values in the Grothendieck group of generalized  varieties over the residue
field.

At the level of foundations, until an additive character is introduced, we 
  are able to work   with Grothendieck semigroups rather than with classes in  Grothendieck groups.

\ssec{The logical setting}   Let $L$ be a valued field, with valuation ring $\Oo_L$.  $\Mm$ denotes the
  maximal ideal.
  We let $\VF^n(L) = L^n$.  The notation $\VF^n$
is analogous to the symbol $\Aa^n$ of algebraic geometry, denoting affine $n$-space.
Let $\RV^m (L) = L^* / (1+\Mm)$, $\G(L) =  L^*/ \Oo^*_L$, $\k(L) = \Oo_L / \Mm_L$.  
Let $\rv: \VF \to \RV$ and $\val: \VF \to \G$ be the natural maps.  The natural map $\RV \to \G$
is denoted $\valr$.  The
exact sequence 
$$0 \to \k^* \to \RV \to \G \to 0$$
shows that $\RV$ is, at first approximation, just a way to wrap together the residue field
and value group.   

We consider expressions
of the form $h(x)=0$ and 
$\val f(x) \geq \val g(x)$ where $f,g,h \in L[X], X = (X_1,\ldots,X_n)$.  A {\em semi-algebraic formula} is a finite Boolean combination of such basic expressions.  A semi-algebraic
formula $\phi$  clearly defines a subset $D(L)$ of $\VF^n(L)$.  Moreover if $f,g,h \in L_0[X]$, we obtain
a functor $L \mapsto D(L)$ from valued field extensions of $L_0$ to sets.   We will later
describe more general definable sets; but for the time being take a {\em definable subset of $\VF^n$}
to be a  functor $D=D_\phi$ of this form.   

An intrinsic
 description of definable subsets of   $\RV^m$ is given in \secref{language}.
  In particular,
 definable
 of $(\k^*)^m$ coincide with constructible sets in the usual Zariski sense; while modulo
 $(\k^*)^m$, a definable set is a piecewise linear subset of $\G^m$.  The structure of 
 arbitrary definable subsets of $\RV^m$ is   analyzed
 in \secref{groupext}.  
 
The advantages of this approach are identical to the benefits in algebraic geometry of
working with arbitrary algebraically closed fields, over arbitrary base fields.  One can use
Galois theory to describe rational points over subfields.  Since function fields are treated
on the same footing, one has a mechanism to inductively reduce higher dimensional
geometry to questions in dimension one, and often in fact to  dimension zero (as in algebraic geometry,
statements about fields, applied to generic points, can imply birational statements about varieties.)

\ssec{Model theory}  Since topological tools are no longer available, it is necessary to define
notions such as dimension in a different way. The basic framework comes from
\cite{Cmin}; we recall and develop it further in \S 2 and \S 4.  It is in many respects analogous
to the o-minimal framework of \cite{vddries}, that has become well-accepted in real algebraic geometry.

  In addition, whereas in geometry all varieties
are made as it were of the same material,  here a number of rather different types of objects
co-exist, and the interaction between them must be clarified.  In particular, the residue
field and the value group are orthogonal in a sense that will be defined below;  
definable subsets of one can never be isomorphic to subsets of the other, unless both are finite.  This orthogonality
has an effect on definable subsets of $\VF^n$ in general; for example, closed disks 
behave very differently from open ones.    Here we follow and further develop \cite{acvf1}.

Note that the set of rational points
of closed and open disks over discrete valuation rings, for instance, cannot be distinguished;
as in rigid geometry, the geometric setting is required to make sense of the notions.  Nevertheless they have immediate consequences to local fields.  As an example,   we define the notion of a definable distribution;
these are defined as a function on the space of polydisks with certain properties.  Making use of  model
theoretic properties of the space of polydisks, we show that any definable distribution 
agrees outside a proper subvariety with one 
obtained by integrating a function.  This is valid over any valued field of sufficiently large 
residue characteristic.  
 In particular  for large $p$, the $p$-adic Fourier transform of a rational polynomial is a locally constant function away from an exceptional subvariety, in the usual sense
(  \corref{lL1}.)
The analogue for $\Rr$ and $\Cc$ was proved by  Bernstein   
using D-modules.  {   For an individual $\Qq_p$, the same result can be shown  using 
Denef integration and a similar analysis of definable sets over $\Qq_p$.  These results were obtained
independently by Cluckers and Loeser, cf. \cite{clu-loe2}.}

\ssec{More general definable sets}   Throughout the paper we discuss not semi-algebraic sets, but definable subsets of a theory with the requisite geometric properties (called $\V$-minimality.)
This includes also the rigid analytic structures of \cite{lipshitz}.   The adjective ``geometrically''
can be take to mean here:  in the sense of the $\V$-minimal theory.  

While we work geometrically throughout the paper, the isomorphisms we obtain are canonical
and so specialize to rational points over substructures.  Thus a posteriori our results
apply to definable sets over any Hensel field of large residue characteristic.  See \secref{exp-rat}.  

For model theorists, this systematic use of  of algebraically closed valued fields to 
apply to other Hensel fields 
is only beginning to be familar.  As an illustration, see \propref{henselQE}, where it is shown that after a little analysis of definable sets over algebraically closed valued fields, quantifier elimination for Henselian fields of residue characteristic zero  becomes a consequence
of Robinson's earlier quantifier elimination in the algebraically closed case.

A third kind of generalization is an a-posteriori expansion of the language in the $\RV$ sort.
Such an expansion involves loss of information
in the integration theory, but is sometimes useful.  For instance one may want to use the Denef-Pas language, splitting the exact sequence
into a product of residue field and value group.  Another example occurs in   \thmref{transi} where
it is explained, given a valued field whose residue field is also a valued field, what happens when one integrates twice.  To discuss this, the residue field is expanded so as to become
itself a valued field.

\ssec{Generalized algebraic varieties}  We now describe the basic ingredients  in more detail.
Let $L_0$ be a valued field with residue field $\k_0$
and value group $A$.
For each point $\g \in \Qq \tensor A$ we have one-dimensional  
 $\k$-vector space 
$$V_\g = \{0 \} \union \{x \in K: \val(x)=\g \} / (1+ \Mm)$$ 
As discussed above, $V_\g$ should be viewed as a functor $L \mapsto V_\g(L)$ 
on valued field extensions $L$ of $L_0$, giving a vector space over
the residue field functor.  If $\g-\g' \in A$ then $V_\g,\V_\g'$ are definably isomorphic,
so one essentially has $V_\g$ for $\g \in (\Qq \tensor A) / A$.

Fix ${\bar{\g}} = (\g_1,\ldots,\g_n)$, and $V_i = V_{\g_i}$, $V_{\bar{\g}} =  \Pi_i V_{\g_i}$.    A {\em $\bar{\g}$-  polynomial}
is an polynomial $H(X) = \sum a_\nu X^{\nu}$ with $\val_p(a_\nu) + \sum_i \nu(i) \g_i = 0 $
for each nonzero term $a_\nu X^{\nu}$.   The coefficients $a_\nu$ are described
in  \secref{weighting}; for the purposes of the introduction, and of  \thmref{d+} below, it suffices to think of   integer coefficients.  
Such a polynomial clearly defines
a function $H: V_{\bar{\g}} \to \k$.  
In particular one has the set of zeroes $Z(H)$.  The {\em generalized residue structure} $\RES_{L_0}$ 
is the residue field, together with the collection of one-dimensional vector spaces $V_\g (\g \in \Qq \tensor A)$
over it, and the functions $H: V_{\bar{\g}} \to \k$ associated to each ${\bar{\g}}$-polynomial.  

The intersection $W$ of finitely many zero sets $Z(H)$ is called a {\em generalized algebraic variety over the residue field}.   Given a valued field extensions $L$ of $L_0$, we
have the set of points $W(L) \subseteq V_{\bar{\g}}(L)$.  
When $L$ is a local field, $W(L)$ is finite.

We will systematically use the Grothendieck group of   generalized varieties  residue field,
rather than the usual Grothendieck group of varieties.    
They are fundamentally of a similar nature:  base change to an algebraically closed value field makes them isomorphic.   But the generalized residue
field makes it possible to to see canonically objects that are only visible after base
change in the usual approach.   One application is  \thmref{d+} below.

$\SG \RES_{L_0} [n]$ denotes the Grothendieck group of 
generalized varieties of dimension $\leq n$; in the paper we will omit $L_0$ from the notation.

\ssec{Rational polyhedra over ordered Abelian groups}  

Let $A$ be an ordered Abelian group.  A {\em   rational polyhedron} $\Delta$ {\em over $A$}
is given by an expression 
$$\Delta = \{x:  Mx \geq b\}$$
with $x =(x_1,\ldots,x_n)$, $M$ a $k \times n$ - matrix with rational coefficients, and $b \in A^k$.
We view this as a functor $B \mapsto \Delta(B)$ on ordered Abelian group extensions
$B$ of $A$.  This functor is already determined by its value at $B = \Qq \tensor A$.  In particular
when $A \leq \Qq$, $\Delta$ is an ordinary rational polyhedron.

$\SG \G_A[n]$ is the semigroup generated by such polyhedra, up to piecewise $GL_n(\Zz)$-transformations and $A$-translations; see \secref{exgam}.  When $A$ is fixed it is omitted from the notation.    

In our applications, $A$ will be the value group of a valued field $L_0$.  If $B$ is the value
group of a valued field extension $L$, write $\Delta(L)$ for $\Delta(B)$.

 \ssec{The Grothendieck semiring of definable sets}  

Fix a base field $L_0$.      The word ``definable'' will mean:  $\T_{L_0}$-definable, with $\T$ a fixed
$\V$-minimal theory.  To have an example in mind one can read ``semi-algebraic over $L_0$''
in place of ``definable''.

Let $ \VF[n]$ be the category of definable subsets $X$ of $n$-dimensional algebraic varieties over  $L_0$; a morphism $X \to X'$ is a 
  definable bijection $X \to X'$ (see \defref{VFcat} for equivalent definitions.)  
$\SG \VF[n]$ denotes the Grothendieck semi-group, i.e. the set of isomorphism classes
of $\VF[n]$ with the disjoint sum operation.  $[X]$ denotes the class of $X$ in the 
Grothendieck semi-group.

We explain how an isomorphism class of $\VF[n]$ is determined precisely by isomorphism
classes of generalized algebraic varieties, and rational polyhedra, whose dimensions add up to $n$.

If $X \subseteq \RES^m$ and $f: X \to \RES^n$ is a finite-to-one map, let
$$\L (X,f) = \VF^n \times _{\rv,f} X = \{(v_1,\ldots,v_n,x): v_i \in \VF, x \in X, \rv(v_i) = f_i(x)\}$$
The $\VF[n]$-isomorphism class $[\L(X,f)]$ does not depend on $f$, and
   is also denoted $[\L X]$.    

 When $S$ is a smooth 
scheme over $\Oo$, $X$ a definable subset of $S(\k)$, $\pi: S(\Oo) \to S(\k)$ the natural
reduction map, we have $[\L X] = [\pi \inv X]$,

We let $\RES[n]$ be the category of pairs $(X,f)$ as above; a morphism $(X,f) \to (X',f')$
is just a definable bijection $X \to X'$.   
 Let $\SG \RES[*]$ be the direct sum of the Grothendieck semigroups $\SG \RES[n]$.

On the other hand, we have already defined $\SG \G[n]$.  Let $\SG \G[*]$ be the
direct sum of the $\SG \G[n]$.  
An element of $\SG \G[n]$
is represented by a definable $X \subseteq  \G[n]$.  Let $\L X = \val \inv (X)$,
$\L [X] = [\L X]$.  
 
It is shown in \propref{tensor} that the Grothendieck semiring of $\RV$ is
the tensor product $\SG \RES[*] \tensor \SG \G [*]$ over the semiring
$K_+\Gamma ^{fin}$ of classes of finite subsets of $\G$; see \secref{exgam}.  

 Note that $\L ([1]_1) = \L([1]_0) + \L([(0,\infty)]_1)$,
where $[1]_1 \in \SG \RES[1]$, $[1]_0 \in \SG \RES[0]$ are the classes
of the singleton set $1$, and $ [(0,\infty)]_1$ is the class in $\SG \G[1]$ of the 
semi-infinite segment $(0,\infty)$.   Indeed $\L ([1]_1)$ is the unit open ball around $1$,
$ \L([1]_0)$ is the point $\{1\}$, while $\L([(0,\infty)]_1)$ is the unit open ball around $0$,
isomorphic by a shift to the unit open ball around $1$.  This is the one relation
that cannot be understood in terms of the Grothendieck semi-ring of $\RV$; it will be seen
to correspond to the analytic summation of geoemtric series in the Denef theory.
Let $\Isp$ be the congruence on the ring  $\SG \RES[*] \tensor \SG \G [*]$ generated
by:  $[1]_1  \sim  [1]_0 +  [(0,\infty)]_1$.

The following theorem summarizes the relation between definable sets in $\VF$ and in $\RV$; 
it follows from \thmref{summ} together with \propref{tensor} in the text.

\<{thm}\lbl{1} $\L$ induces an  surjective homomorphism of filtered semirings
 $$\SG \RES[*] \tensor \SG \G [*]  \to   \SG(\VF) $$
 The kernel is precisely the congruence $\Isp$. 
   \>{thm}
 
The inverse   isomorphism $\SG(\VF) \to \SG \RES[*] \tensor \SG \G [*]  / \Isp$ 
 can be viewed as a kind of Euler characteristic, 
respecting products and disjoint sums, and functorial in various other ways.  

The values of this Euler characteristc are themselves geometric objects, both on the algebraic-geometry
side ($\RES$) and the combinatorial-analytic side ($\G$).  This is valuable for some purposes; 
in particular it becomes clear that 
the isomorphism is compatible with taking rational points over Henselian subfields (cf. \propref{rationalpoints})

For other applications, however, it would be useful to obtain more manageable numerical invariants;
for this purpose one needs to analyze the structure of $\SG(\G[*])$.  We do not fully do this here, 
but using a number of homomorphisms on $\SG(\G[*])$ we obtain a number of invariants.
In particular using the  $\Zz$-valued Euler characteristics on $\K ( \G [*])$ (cf. \secref{exgam} and \cite{marikova},
\cite{kageyama-fujita}), we obtain   two homomorphisms on $\SG \VF[n]$ essentially to  $\K (\RES [n])$.  The reason there
are two rather than one has to do with  Poincar\'e duality,
see   \thmref{retract2}.   

For instance,  when $F$ is a field of characteristic $0$ we obtain an
invariant of  rigid analytic varieties 
over $F((t))$, with values in  the Grothendieck ring $\K(\Var_F)$ of algebraic varieties over $F$;
and another in $\K(\Var_F)[[\Aa_1] \inv]$ (\propref{ls1}).
It is instructive to compare this with the invariant of \cite{loeser-sebag}, with values   in $\K(\RES[n]) / [G_m]$.
\footnote{The setting is somewhat different:  Loeser-Sebag can handle positive characteristic too, but
assume smoothness.}   Since any two closed balls are isomorphic, via additive translation and  multiplicative contractions,      all closed balls must have the same invariant.   
Working with a discrete  value group tends to force $[G_m]=0$, since
it appears that a closed ball $B_0$ of valuation 
radius $0$ equals $G_m$ times a closed ball $B_1$ of valuation radius $1$.  Since our
technology is based on divisible value groups, the ``equation'' $[B_0]=[B_1][G_m]$ is 
replaced for us by $[B_0]=[B_0^o][G_m]$, where $B_0^o$ is the open ball of valuation radius
$0$.  Though $B_1$ and $B_0^0$ have the same $F((t))$-rational points, 
they are geometrically distinct (cf. \lemref{image1}) and so no collapse takes place.  

See  also subsections \ref{tim} and \ref{cl-ha} for two previously known cases.
 
 By such Euler characteristic  methods we can prove a statement purely concerning
  algebraic varieties, partially answering a question of 
Gromov   and Kontsevich (\cite{gromov}, p. 121.).  In particular,   two elliptic curves
with isomorphic complements in projective space were previously known  to be isogenous,
by zeta function methods; we show that they are isomorphic. 
{    This also follows from 
\cite{larsen-luntz}; the method there requires strong forms of resolution of singularities. }
 See \thmref{complement}.

\ssec{Integration of forms up to absolute value}   \lbl{intro-forms}

Over local fields, data for integration consists of a triple $(X,V,\omega)$, with 
$X$ a   definable subset  of a smooth variety $V$, and $\omega$ volume form  
on $V$.  We are interested an integral of the form $\int_X |\omega|$, so that 
multiplication of $\omega$ by a function with norm one does not count as a change,
nor does removing a subvariety of $V$ of smaller dimension.  
Using an equivalent description of $\VF[n]$, where the objects come with a distinguished finite-to-one map into affine space,
we can represent an integrand 
as a pair $(X,\om)$ with $X \in \Ob \VF[n]$ and $\om$ a function on $X$ into 
$\G$.    Isomorphisms are
essential bijections, preserving the form up to a function of norm 1.  
See \defref{VFmcat} for a precise definition of this category, the  category $\VFmg[n]$.
 
Integration is intended to be an invariant of isomorphisms in this category.  Thus
we can find the integral if we determine all invariants.  We do this in complete analogy
with \thmref{1}. 
 
For $n \geq 0$ let  $\G[n]$ be the category whose objects are  finite unions of rational polyhedra over the group $A$  of definable points of $\G$.  A morphism $f: X \to Y$ of $\G[n]$ is a bijection such that for some partition $X = \union_{i=1}^k X_i$ into rational polyhedra, $f|X_i$ is given by an element of $GL_n(\Zz) \sd A^n$.   Let $\mG [n]$ be the category of  pairs $(X,\om)$, 
with $X$ an object of $\G[n]$, 
and $\om: X \to \G$ a piecewise affine map.  A morphism $f: (X,\om)  \to (X',\om')$ is a morphism  $f: X \to X'$ of $\G[n]$
such that
$\sum _{i=1}^l x_i + \om(x) = \sum_{i=1}^l x'_i + \om'(x')$
  whenever $(x_1',\ldots,x_n')=f(x_1,\ldots,x_n)$. 
 Given $(X,\om) \in \Ob \mG [n]$, define $\L X$ as above, and adjoint the pullback of $\om$
 to obtain an object of $\VFmg [n]$.   This gives a homomorphism
 $\SG \mG [n] \to \SG \VFmg [n]$.   
  
 \<{thm}\lbl{2}   $\L$ induces an  surjective homomorphism of filtered semirings
 $$\SG \RES[*] \tensor_{\Nn}  \SG {\mG} [*]  \to   \SG(\VFmg)[*] $$
 The kernel is generated by the relations $p \tensor 1 = 1 \tensor [(\valr(p),\infty)]$ 
 and $1 \tensor a = \valr \inv(a) \tensor 1$.
   \>{thm} 
 
 In the statement of the theorem, $p$ ranges over definable points of $\RES$
 (actually one value suffices), and $a$ ranges over definable points of $\G$.  

 This can also be written as $(\SG \RES[*] \tensor_{\sggfm}  \SG {\mG} [*]) /\Ispm \iso \SG(\VFmg)[*]$, where $\sggfm$ is the subsemiringof   subsets  of $\mG$ with finite support, and
 $\Ispm$ is a semiring congruence defined similarly to $\Isp$.  The base of the tensor
 leads to the identification of a point of $\G$ with   with a coset of $\k^*$ in $\RES$, while 
$\Ispm$ identifiese a point
 of $\RES$ with an infinite interval of $\G$. 
   The inverse isomorphism 
can be viewed as an integral. 

We introduce neither additive nor multiplicative inverses in $\SG \RES[*] $ formally, so that
the target of integration is completely geometric.  

We proceed to give an application of the first part of the theorem (the surjectivity) 
  in terms of ordinary $p$-adic integration.
 
\ssec{Integrals over local fields:  uniformity over ramified extensions}

Let $L$ be a local field, finite extension of $\Qq_p$ or $\Ff_p((t))$.  
We   normalize the Haar measure  $\mu$ in such a way that the maximal ideal has measure 1,
the norm by $|a| = \mu \{x: |x| < |a| \}$.  Let $\RES_L$ be the generalized residue field, and 
$\G_L$ be the value group.  We assume $\Qq_p$ or $\Ff_p((t))$ have value group $\Zz$,
and identify $\G_L$ with a subgroup of $\Qq$.

Given $c=(c_1,\ldots,c_k) \in L^k$ and $s=(s_1,\ldots,s_k) \in \Rr^k$ with $s_i \geq 1$, 
let $|c|^s = \Pi_{i=1}^k |c_i|^{s_i}$.  

Let  $\lambda$ be a multiplicative character $\Rr^n \to \Rr^*$.  Define:
$$ \ev_\lambda(\Delta(B)) = \sum_{b \in \Delta(B)} \lambda(b) $$
provided this sum is absolutely convergent. 
Given linear functions $h_0,\ldots,h_k$ on $\Rr^n$
and $s_1,\ldots,s_k \in \Rr$, let 
 $\ev_{h,s,Q} =  \ev_\lambda$ where $\lambda(x) = Q^{h_0(x) +  \sum s_i h_i(x)}$  
 
 \<{thm}\lbl{d+}  Fix $n,d,k \in \Nn$.  Let $p$ be a large prime compared to $n,d,k$, 
and let $f \in \Qq_p[X_1,\ldots,X_n]^k$ have degrees $\leq d$.  
Then there exist finitely many generalized varieties $X_i$ over $\RES(\Qq_p)$,  
   rational  polyhedra $\Delta_i$,   $\g(i) \in \Qq^{\geq 0}$, $n_i \in \Nn$,
 and linear functions $h^i_0,\ldots,h^i_k$   with rational coefficients,  
   such that for any finite extension
$L$ of $\Qq_p$ with residue field $GF(q)$ and $\val(L^*)=(1/r) \Zz, \val(p)=1$,   and any    $s \in \Rr_{\geq 1}^k$, 
$$\int_{\Oo_L^n} |f|^s = \sum_i q^{r\g(i)} (q-1)^{n_i} |X_i(L)|   \ev_{h^i,s,q^r}(\Delta_i(L))$$  
\>{thm}
 
 Note:
 
 1)$\Delta_i(L)$ depends only on the ramification degree $r$ of $L$ over $\Qq_p$.

2) The formula is a sum of non-negative terms.

3)  $\ev_{h,s,q^r}(\Delta_i((1/r) \Zz))$ can be written  in closed form as a 
rational function of $q^{rs}$.    This follows from Denef, who shows it for more
general sets $\Delta_i$ definable in  Pressburger arithmetic;   such  analytic summation is an essential component of his integration theory.   Since  it plays no role in our approach we leave the statement in geometric form.

4)   The generalized varieties $X_i$ and   polyhedra $\Delta_i$ are simple  
 functions of the coefficients $f$.   
 Here we wish to emphasize not this,
 but the uniformity of the expression over  ramified extensions of $\Qq_p$. 
 
 The proof follows \propref{tensor-m} (it uses only the easy surjectivity 
 in this proposition,  and \propref{tr2}.)

 \ssec{Bounded and unbounded sets}   
 
 The isomorphism of semirings of \thmref{2} obviously induces an isomorphism of 
rings.  However introducing additive inverses loses information on the $\G$ side;
the class of the interval $[0,1)$ becomes $0$, since $[0,\infty)$ and $[1,\infty)$
are isomorphic.  The classical remedy is to cut down to bounded sets before
groupifying.  This presents no difficulty, since the isomorphism respects boundedness.

In higher dimensional local fields, stronger notions of boundedness may be useful,
such as those introduced by Fesenko.  Since these questions are not entangled
with the theory of integration, and can be handled a posteriori, we will deal with them in
a sequel.

Here we mention only that even if one insists on integrating all definable integrands, with no 
boundedness condition, into a ring, some but
not all information is lost.  This is due to the existence of Euler characteristics
on $\G$, and thus again to the fact that we work geometrically, with divisible groups, even
if the base field has a discrete group.  We will see (\lemref{gamma-volume}) that
$\SG({\mG}[n])$ can be identified with the group of definable functions $\G \to  \SG(\G[n])$.  Applying an appropriate Euler characteristic  reduces to the group of piecewise constant
functions on $\G$ into $\Zz$.  Recombining with $\RES$  we obtain 
a consistent definition of an integral on unbounded integrands, compatible with measure preserving maps, sums and products, with values in $\K(\RES)[A] / [\Aa_1]_1 \K(\RES)[A]$, where $A$ is the group of definable points of $\G$, and  $[\Aa_1]_1$ is the class
of the affine line.  See \thmref{retract-m}.

 \ssec{Finer volumes}  
 
 We also consider a finer category of definable sets with $\RV$-volume forms.  This means
 that a volume form $\om$ is identified with $g \om$ only when $g-1 \in \Mm$; $\val(g)=0$
 does not suffice.  We obtain an integral whose values themselves are definable sets
 with volume forms; in particular including  algebraic varieties with volume forms over the residue field.  
 
  \<{thm}\lbl{3}   $\L$ induces an  surjective homomorphism of graded semirings
 $$\SG \RVm[*]   \to   \SG(\VFm)[*] $$
 The kernel is precisely the congruence $\Ispm$. 
   \>{thm} 
 $  {\RVm}$ is the category of definable subsets of $\RVm^*$ enriched with volume
 forms;   see \defref{RVm}.  Again, an isomorphism is induced in the opposite
 direction, that can be viewed as a motivic integral
 $$\int:  \SG(\VFm)[*] \to  \SG \RVm[*] / \Ispm$$
 This allows an iteration of the integration theory, either with an integral of the same nature
 if the residue field is a valued field, or with a different kind of integral if for instance
 the residue field is $\Rr$.

\ssec{Hopes}  We mention three.  Until now, a deep obstacle existed to extending Denef's theory to positive 
characteristic; namely the theory was based on quantifier elimination for Hensel fields
of residue characteristic $0$, or for finitely ramified extensions of $\Qq_p$, and it is known
that no similar quantifier elimination is possible for $\Ff_p((t))$, if any is.   On the other hand
Robinson's quantifier elimination is perfectly valid in positive characteristic.  This raises
hopes of progress in this direction, though other obstacles remain.

It is natural to think that the theory can be applied to higher dimensional local fields;
we will consider this in a sequel.

Another important target is asymptotic integration over $\Rr$.  Nonstandard
extensions of $\Rr$ admit natural valued field structures.  This is the basis
of Robinson's nonstandard analysis.  These valued fields have divisible value groups,
and so previous theories of definable integration do not apply.  The theory of this paper
applies however, and we expect that it will yield connections between $p$-adic integration
and asymptotics of real integrals.   

\ssec{Course of the paper}

After recalling some basic model theory in \S 2, we proceed in \S 3 to $\V$-minimal theories.  

In  \secref{descent-objects} we show that any definable   subset  of $\VF^n$ admits a constructible bijection
 with some $\L(X,f)$.  In fact only a very limited class  bijections is needed; a typical one has
 the form: $(x_1,x_2) \mapsto (x_1, x_2+f(x_1,x_2))$, so it is clearly measure preserving.
The proof is simple and brief, and uses only a little
 of the preceding material.  We note here that for many applications this 
  statement is already sufficient; in particular it suffices to give the surjectivity
  in Theorems \ref{1} and \ref{2}, and hence the application \thmref{d+}.  

In \S 5 we return to the geometry of $\V$-minimal structures, developing a theory
of differentiation.  We show the compatibility between differentiation in $\RV$
and in $\VF$.  This is needed for \thmref{3}.  Differentiation in $\VF$
involves much finer scales than in $\RV$; in effect $\RV$ can only see
distances measured by valuation $0$, while the derivative in $\VF$ involves
distances of arbitrarily large valuation.  The proof uses a continuity argument on scales.
It fails in  positive characteristic, in its present form.

\S 6 is devoted to showing that $\L$ yields a well-defined map $\SG(\RV) \to \SG(\VF)$;
in other words not only objects, but also isomorphism can be lifted.  

Sections \S 7 and \S 8 investigate  the kernel of $\L$ in \thmref{1}.  This is the most technical part of the paper, and we have not been able to give a  proof as functorial as we would have liked.
See Question \ref{q}.

In \S 9 we study the piecewise linear Grothendieck group; see the introduction to this section.

\S 10 decomposes the Grothendieck group of $\RV$ into the components $\RES$ and $\G$,
used througout  this introduction.  

\S 11 introduces an additive character, and hence the Fourier transform.   The isomorphism
of volumes given by \thmref{3} suffices for this extension; it is not necessay to redo the 
theory from scratch, but merely to follow through the functoriality.   

\S 12 contains the extension to definable sets over Hensel fields mentioned above, 
and \S 13 the application to the Grothendieck group of varieties.

Thanks to Aviv Tatarsky and Moshe Kaminsky, and to Lou Van den Dries, Clifton Ealy,
and Jana Ma\v{r}\'{i}kov\'{a},
 for useful comments and corrections.
 
\>{section}%{Introduction}

\<{section}{First order theories}

The bulk of this paper uses no deep results from logic, beyond Robinson's  quantifier elimination
the theory of algebraically closed valued fields (\cite{robinson}).   However,   it is imbued with a model theoretic viewpoint.  We will not explain the most basic notions of logic: language, theory, model.  Let us just mention that a language consists of basic relation and function symbol, and formulas are built out of these, using symbols for Boolean operations and quantifiers; cf. e.g. \cite{enderton} or \cite{johnstone}, or the first section of \cite{CK}); but we attempt  in this section to bridge the gap between these and the model-theoretic language used in the paper.  

A language $L$ consists of a family of ``sorts'' $S_i$, a collection of variables ranging over each sort,
a set of relation symbols $R_j$, each intended to denote a subset of a finite product of sorts,
and a set of function symbols $F_k$ intended to denote functions from a given finite product of sorts
to a given sort.  From these, and the logical symbols $\&,\neg,\forall,\exists$ one forms {\em formulas}. 
A sentence is a formula with no free variables (cf. \cite{enderton}.)
A theory $T$ is a set of sentences of $L$.   A theory is called {\em complete} if for every sentence $\phi$
of $L$, either $\phi$ or its negation $\neg \phi$ is in $T$.

 A {\em universe} $M$ for the language $L$ consists, by definition, 
of a set $S(M)$ for each sort $S$ of $L$.  An $L$-structure consists of such a universe,
together with an interpretation of each relation and function symbol of $L$.   One can define the truth
value of a sentence in a structure $M$; more generally,     if $\phi(x_1,\ldots,x_n)$ is 
a formula, with $x_i$ a variable of sort $S_i$,   then one defines
the interpretation $\phi(M)$ of $\phi$ in $M$, as the set of all $d \in S_1(M) \times \ldots \times S_n(M)$    of which $\phi$ is true.  If
 every sentence in $T$ is true in $M$, one says that $M$ is a model
of $T$ ($M \models T$).   The fundamental theorem here is a consequence of G\"odel's completeness theorem called 
  the {\em compactness theorem:} 
a theory $T$ has a model   if every finite subset of $T$ has a model.  

The language $L_{rings}$ of rings, for example, has one sort, three function symbols $+,\cdots,-$, two constants $0,1$;
any ring is an $L_{rings}$-structure; one can obviously write down a theory $T_{fields}$
in this language
 whose models are precisely the fields.

\ssec{Basic examples of theories}  \lbl{language}

   We will work with a number of theories
associated with valued fields.

\begin{enumerate}   
  \item ACF , the theory of algebraically closed fields.   The language is the language of rings $\{+,\cdot,-,0,1\}$,
 mentioned earlier.  The theory states that the model is a field, and for each $n$, that every monic
 polynomial of degree $n$ has a root.  For instance for $n=2$:
 $$(\forall u_1)(\forall u_0)(\exists x)(x^2+u_1x+u_0=0)$$
 ACF(0) includes in addition the sentences:  $1+1 \neq 0, 1+1+1 \neq 0, \ldots$.  This theory is complete
    (Tarski-Chevalley.)  It will arise 
  as the theory of the residue field of our valued fields.   
   
  \item Divisible ordered Abelian groups ($\DOAG$).   The language consists of a single sort, a binary relation 
  symbol $<$, a binary 
  function symbol $+$, a unary function symbol $-$ and a constant symbol $0$.   The theory states
  that a model is an ordered Abelian group.  In addition axioms asserting divisibility by $n$ for each $n$,
  for instance:  $(\forall x)(\exists y) (y+y=x)$.

  This is the theory  of the
  value group of a model of $\ACVF$.  
  
\item  The $\RV$ sort.
(Extension  of (2) by (1)).   
  The language has one official sort, denoted $\RV$, and includes:  
 Abelian group operations $\cdot,/$ on $\RV$,  
a unary predicate $\k^*$ for a subgroup, and an operation $+: \k^2 \to \k$,
where $\k$ is $\k^*$ augmented by a constant $0$.     Finally,
there is a partial ordering; the theory states that $\k^*$ is the equivalence class of $1$; that $\leq $  is a total ordering
on $\k^*$-cosets, making $\RV/\k^* =: \Gamma$ a divisible ordered
Abelian group, and that $(\k,+,\cdot)$ is an     algebraically closed field.  
(We thus have an  exact sequence  $0 \to \k^*  \to \RV \to \Gamma \to 0$, but we treat $\Gamma$ as an imaginary sort.)    This theory  $\TRV$ is complete too.

We will sometimes view $\RV$ as an autonomous structure; but it
will arise from an algebraically closed valued field, as in (5) below.

\item  Let $M \models \TRV$, and let $A$ be a subgroup of $\G(M)$.  Within $\TRV_A$ we see 
an interpretation of ACF, namely the algebraically closed field $\k$.  In addition for each $a \in A$
 we have a 1-dimensional $\k$-space, 
the fiber of $\RV$ lying over $\G$ augmented by $0$. Collectively the field $\k$ with this collection of vector spaces will be denoted $\RES$.

  \item $\ACVF$, the theory of  algebraically closed valued fields.  According to Robinson, 
   the completions, denoted 
$\ACVF(q,p )$, are obtained by specifying the characteristic $q$  and
  residue characteristic $p$.   We will   be concerned with $\ACVF(0,0)$ in this paper.  However
 since any sentence of $\ACVF(0,0)$ lies in $\ACVF(0,p)$ for almost all primes $p$, the results
 will a posteriori apply also to valued fields of characteristic zero and large residue characteristic.

        We will take  $\ACVF(0,0)$ to have
  two sorts, $\VF$ and $\RV = \VF^* /(1+ \Mm)$.  The language includes the language of rings (1) 
  on the $\VF$ sort, the language (3) on the $\RV$ sort, and a function symbol $\rv$ for a function
     $\VF^* \to \RV$.   Denote $\rv \inv (\RV^{\geq 0}) = \Oo, \rv \inv (0) = \Mm$.  
     
     The theory states that $\VF$ is a valued field, with valuation ring $\Oo$ and maximal ideal $\Mm$;
     that $\rv: \VF^* \to \RV$ is a surjective group homomorphism, and the restriction to
     $\Oo$ (augmented by $0 \mapsto 0$) is a surjective ring homomorphism.

The structure $\ACVF_A$ induces on $\G$ is of a uniquely
divisible Abelian group, with constants for the elements of $\G(A)$.
Thus every definable subset of $\G$ is a finite union of points
and open intervals (possibly infinite.)

\item{Rigid analytic expansions (Lipshitz)}  The theory $\ACVFR$ of algebraically closed valued
fields expanded by a family $R$ of analytic functions.  See    \cite{lipshitz}
and  \cite{lipshitz-robinson}.
  Our theory
of definable sets will be carried out axiomatically, and are thus also valid for these
   rigid analytic expansions.

\end{enumerate}

 A {\em definable set}  $D$
is not really a set, but a functor from the category of models of $T$  to the category of sets of the form $M \mapsto \phi(M)$, where $\phi$ is a formula of $L$.  
Model theorists do not really distinguish between the definable set $D$ and the formula $\phi$
defining it;   we will usually refer to definable sets rather than to formulas.    If $R \subseteq D \times D'$ and for any model $M \models T$,  $R(M)$ is the graph of a  function $D(M) \to D(M')$, we say $R$ is a {\em definable function of $T$.}
Similarly we say $D$ is {\em finite} if $D(M)$ is finite for any $M \models T$, etc.  It follows from the compactness
theorem   that if $D$ is finite, then for some integer $m$ we have $|D(M)| \leq m$
for any  $M \models T$.
We sometimes   write $S^*$ to denote $S^n$ for some unspecified $n$.

   By a {\em map} between $L$-structures $A,B$ we mean a 
 family $f=(f_S)$ indexed by the sorts of $L$,    with $f_S: S(A) \to S(B)$;
   one extends $f$ to products of sorts by setting $f((x_1,\ldots,x_n)) = (f(x_1),\ldots,f(x_n))$.
 $f$ is {\em an embedding of structures} if 
 $f ^{-1} R(B) = R(A)$ for any atomic formula $R$ of $L$.  
 Taking $R$ to be the equality relation, this includes in particular the statement that
 each $f_S$ is injective.

On occasion we will use $\infty$-definable sets. An $\infty$-definable set is a functor of the form
$M \mapsto \meet {\mathcal D}$, where $ {\mathcal D} $ is a given collection of definable sets.  
In a complete theory
a definable set is determined by the value it has at a single model, this is of course false for $\infty$-definable sets.

We write:  $a \in D$ to mean:  $a \in D(M)$ for some $M \models T$.   It is customary, since Shelah,
to choose a single universal domain $\Uu$ embedding all ``small'' models, and let $a \in D$
mean:  $a \in D(\Uu)$; we will not require this interpretation, but the reader is welcome to take it.

We will sometimes consider {\em imaginary sorts}.  If $D$ is a definable
set, and $E$ a  definable equivalence relation   on $D$, then $D/E$ may be considered
to be an {\em imaginary sort}; as a definable set it is just the functor $M \mapsto D(M) / E(M)$.
 A definable subset of a product $\Pi_{i=1}^n D_i/E_i$ of imaginary sorts (and ordinary sorts)
 is taken to be a subset whose preimage in $\Pi_{i=1}^n D_i$ is definable; the notion of a definable function
 is thus also defined.  In this way the imaginary sorts can be treated on the same footing as the others. 
 The set of all elements of all imaginary sorts of a structure $M$ is denoted $M^{eq}$.   It is easy to construct a theory
 $T^{eq}$ in a language $L^{eq}$ whose category of models is (essentially) $\{M^{eq}: M \models T \}$.  See \cite{shelah}, \cite{poizat}, \S 16d.

 Given a definable set $D \subseteq  S \times X$, where $S,X$ are definable sets, and
given $s \in S$, let $D(s) = \{x \in X: (s,x) \in D\}$.   
Thus $D$ is viewed as
a {\em family} of definable subsets of $X$, namely $\{D(s): s \in S\}$.  If 
$s \neq s'$ implies $D(s) \neq D(s')$, we say that the parameters are {\em canonical},
or that $s$ is a {\em code} for $D(s)$.  In particular, if $E$ is a definable equivalence relation,
the imaginary elements $a/E$ can be considered as codes for the classes of $E$.

$T$ is said to {\em eliminate imaginaries} if every imaginary sort admits a definable
injection into a product of some of the sorts of $L$.  For instance, the theory of algebraically closed fields eliminates imaginaries.  See \cite{poizat-imag} for an excellent exposition of these issues.  We note 
that $T$ admits elimination of imaginaries iff 
 for any family $D \subseteq S \times X$ there exists a  family $D' \subseteq S' \times X$ 
such that for any $t \in S$ there exists a {\em unique} $t' \in S'$ with $D(t)=D'(t')$. 

  (Recall that  $t \in S$ means:
$t \in S(M)$ for some $M \models T$.  The uniqueness of $t'$ implies in this case that one can choose $t' \in S'(M)$
too.)     In this case we also say that $t'$ is called a {\em canonical parameter} or {\em code} for $D(t)$.

\<{example} \lbl{exprim} \rm let $b$ be  a nondegenerate closed ball in a model  the theory   $\ACVF$  of algebraically closed valued
 fields.  Then $b= \{x: \val(x-c) \geq \val(c-c') \}$ for some elements $c \neq c'$ of the field .  $b$ is coded by $\bar{b}= (c,c')/E$, where $(c,c') E (d,d')$ iff $\val(c-c')=\val(d-d') \leq \val(c-d)$.  However
 we often fail to distinguish notationally between $b$ and $\bar{b}$, and in particular
 we write $A(b) = A(\bar{b})$.   \>{example} 

The only imaginary sorts  that will really be essential for us are the sorts
$\fB$ of closed and open balls.  The closed balls around $0$ can be identified with their radius,
hence the   valuation group $\G(M) = \VF ^*(M) / \Oo ^*(M)$
of a valued field $M$ is embedded as part of $\fB$. 

{\bf Notation}    Let $\fB = \fB^o \union \fB^{cl}$, the sorts of open and closed sub-balls  of $\VF$.
  Let    ${\Gamma^+} = \{\gamma \in \Gamma:  \gamma \geq 0 \}$.  
   $$\fB^{cl} = \dU_{\gamma \in \Gamma} \fB^{cl}_\gamma, \ \ \fB^{cl}_\gamma = \VF / \gamma \Oo$$
  $$\fB^{o}= \dU_{ \gamma \in \Gamma } \fB^{o}_\gamma, \ \ \fB^{o}_\gamma = \VF / \gamma \Mm$$ 
Here $\gamma \Mm = \{x \in \VF: \val(x) > \gamma \}$,
$\gamma \Oo = \{x \in \RES: \val(x) \geq \gamma \}$.  
The elements of $\fB^{cl}_\gamma$ , $fB^{o}_\gamma$ will
be referred to as {\em closed and open balls of valuative radius $\gamma$};
though this valuative definition of radius means that bigger balls
have smaller radius.  The word ``distance'' will be used similarly.

By a {\em thin annulus} we will mean:  a closed ball of valuative radius $\gamma$, 
with an open ball of valuative radius $\gamma$ removed.

 Fix a model $M$ of $T$.  A {\em substructure} $A$ of  $M$  (written:  $A \leq M$)
consists of a subset $A_S$ of $S(M)$, for each sort $S$ of $L$, closed under all
definable functions of $T$.   For example,   the substructures of models of $T_{fields}$ are the integral domains.

 In general, the   {\em definable closure} of a set $A_0 \subset M$ is
the smallest substructure containing $A_0$; it is denoted $\dcl(A_0)$ or $<A_0>$. 
An element of $<A_0>$ can be written as $g(a_1,\ldots,a_n)$ with $a_i \in A_0$ and $g$ a definable
function; i.e., it is an element satisfying a formula $\phi(x,a_1,\ldots,a_n)$ of $L_{A_0}$ in one variable  that has exactly one solution in $M$.  
 If $A$ is a substructure,     $\dcl(A \union \{c\})$ is also denoted $A(c)$.  These notions apply equally 
 when $A,c$ contain elements of the imaginary sorts.
If $B$ is contained in sorts $S_1,\ldots,S_n$, then 
$\dcl(B)$ is said to be an $S_1,\ldots,S_n$-generated substructure.  In the special case
of valued fields, where one of the sorts $\VF$ is the ``main'' valued field sort, 
a $\VF$-generated structure will be said to be field-generated, or sometimes just ``a field''.

For any definable set $D$, we let $D(A)$ be the set of points of $D(M)$ with coordinates in $A$.
If $S = D/E$ is an imaginary sort, $S(A)$ is the set of $a \in S$ whose preimage is defined over $A$.
We have  $D(A) / E(A) \subseteq S(A)$.     $D(A) / E(A)$ is of course
closed under definable functions $S^m \to S$ that lift to definable functions $D^m \to D$,
  but it is not 
  necessarily closed under arbitrary definable functions,
  i.e functions whose graph   is the image of a definable subset of $D^m \times D$.
  For example $x \mapsto (1/n)x $ is a definable function on
  the value group of a model of $\ACVF$, but if  $A \leq M \models \ACVF$,
  $\G(A)$ need not be divisible.
  
   When $A \leq M, B \leq N$ with $M,N \models T$,   a function $f: A \to B$ 
  is called a  {\em (partial) elementary embedding} $(A,M) \to (B,N)$ if for any definable set $D$ of $L$,  $f \inv D(B) = D(A)$.  In particular, when $A=M, B=N$, one say that $M$ is an
  {\em elementary submodel} of $N$.  
  
    By a {\em constructible set} over    $A$, we mean the functor
$L \mapsto \phi(L)$ on models $M \models T_A$,  
where  $\phi  = \phi(x_1,\ldots,x_n,a_1,\ldots,a_m)$
is a  quantifier-free formula 
 with parameters from $A$.

 We say that $T$ {\em admits
   quantifier-elimination}  if every  definable set  coincides
   with a constructible set.   It follows in this case that for any $A$, any $A$-definable set is
   $A$-constructible.  
   When $T$ admits quantifier elimination,$f:A \to B$ is a partial elementary embedding
   iff it is an 
      embedding of structures.  
      
      The theories   (1-5) of \secref{language} admit    quantifier elimination in their natural algebraic  languages
(theorems of Tarski-Chevalley and Robinson; cf.  \cite{acvf1}.)  The sixth admits quantifier elimination
in a language that needs to be formulated with more care, see \cite{lipshitz}.

      In all of this 
  paper except for \secref{expansionsS}, \secref{rational}, we will use only structural properties
  of definable sets, and not explicit formulas.  In this situation quantifier elimination can be assumed
  softly, by merely increasing the language by definition so that all definable sets become
 equivalent to quantifier-free ones.  The above distinctions will only directly come into play in the 
 two sections \secref{expansionsS}, \secref{rational}.   
  
If $A \leq M \models T$, $L_A$ is the language $L$ expanded by a constant $c_a$
for each element $a$ of $A$, so that an $L_A$-structure is the same as an $L$-structure $M$
together with a function $A_S \to S(M)$, for each sort $S$.   $T_A$ is the set of $L_A$
sentences true in $M$ when the constant symbol   $c_a$ is interpreted as $a$; the models of $T_A$ 
 are models $M$ of $T$, together with an isomorphic embedding of $A$ as
a substructure of $M$. 
In particular, $M$ with the inclusion of $A$ in $M$ is an $L_A$-structure
denoted $M_A$.   For any subset $A_0 \subseteq M$, we write
$T_{A_0}$ for $T_{<A_0>}$, where $<A_0>$ is the substructure generated by $A_0$.

A definable set of $T_A$ will also be referred to as $A$-definable; similarlly for other notions such as those defined just below.  

   A {\em parametrically definable set} of $T$ is by definition a $T_A$ - definable set for some $A$.
 
  An {\em almost definable} set is the union of classes of a definable  equivalence
relation with finitely many classes.  An element $e$ is called {\em algebraic} (resp. definable)
if the singleton set $\{e\}$ is almost definable (resp. definable).  When $T$ is a complete theory,
the set of algebraic (definable) elements of a model $M$ of $T$ forms a substructure
that does not depend on $M$, up to (a unique) isomorphism.  
 
 Let $A_0 \subseteq M \models T$; 
the set of  $e \in M $  almost  definable over $A_0$ 
is called the {\em algebraic closure } of $A_0$, $\acl(A_0)$.  If $A_0$ is contained in sorts
$S_1,\ldots, S_n$, any substructure of $\acl(A_0)$ containing $\dcl(A_0)$ is said to be {\em almost 
$S_1,\ldots,S_n$-generated.}

\<{example} \lbl{def-lin}  \rm If a definable set $D$ carries a definable linear ordering, then every
algebraic element of $D$ is definable.  This is because the {\em least} 
element of a finite  definable set $F$ is clearly definable; the rest are contained in  a
smaller finite definable subset of $D$, so are definable by induction.  \rm 

If in addition $D$ has elimination of imaginaries, and $Y$ is almost definable and definable 
with parameters from $D$, then $Y$ is definable.   Indeed using   elimination of imaginaries in $D$, the set $Y$ can be defined using 
canonical parameters.  These are algebraic elements of $D$, hence definable.  
 \>{example}
 
Two definable functions $f: X \to Y,f': X \to Y'$ will be called {\em isogenous} if
for all $x \in X$, $\acl(f(x))=\acl(f'(x))$.

\sssec{Compactness}  \lbl{varconst}

Compactness often allows to replace arguments in relative dimension one over a definable set,
by arguments in dimension one over a different base structure.   Here is an example:

\<{lem} \lbl{collect}  Let $f_i: X_i \to Y$ be   definable maps between definable sets of $T$
($i=1,2$.)    Assume that for any $M \models T$ and
$b \in Y(M)$, $X_1(b) := f_1 \inv(b) $ is $T_b$-definably isomorphic to $X_2(b)= f_2 \inv (b)$.  Then $X_1,X_2$
are definably isomorphic.  \>{lem}

\<{proof}  Let  $\fF$ be the family of pairs $(U,h)$, where $U$ is a definable subset of $Y$,
and $h: f_1 \inv U \to f_2 \inv U$ is a definable bijection.  

\Claim{} For any $b \in Y(M), M \models T$, there exists $(U,h) \in \fF$ with $b \in U$.  

\<{proof}  Let $b \in Y(M)$.  There exists a $T_b$-definable bijection
$X_1(b) \to X_2(b)$.  This bijection can be written:  $x \mapsto g(x,b)$, where $g$
is a definable function.  Let $U = \{y \in Y: (x \mapsto g(x,y)) \text{ is a bijection } X_1(y) \to X_2(y)\}$.  Then $(U, g(x,f_1(x))) \in \fF$, and $b \in U$.   \>{proof}

Now by compactness, there exist a finite number of definable subsets $U_1,\ldots,U_k$
of $Y$, with $Y = \union_i U_i$, and with $(U_i,h_i) \in \fF$ for some $h_i$.
Let $U_i' = U_i \m (U_1 \union \ldots \union U_{i-1})$, and define $h = \union_i h_i | U_i'$.  
Then $h: X_1 \to X_2$ is the required bijection.   \>{proof}

 Here is another example of the use of compactness:
 
 \<{example}\lbl{alldepend}  If $D$ is a definable set, and for any $a,b \in D$, $a \in \acl(b)$, then
 $D$ is finite.   More generally, if $a \in \acl(b)$ for any $b \in D$, then $a \in \acl(\emptyset)$\>{example}
 
 \prf We prove the first statement, the second being similar.   For any model $M$, pick $a \in M$; then $D(M) \subseteq \acl(a)$.   For $b \in \acl(a)$, let $\phi_b$ be the formula $x \neq b \& D(x)$.  
 So the set of formulas $Th(M)_M \union \{\phi_b\}$ has no common 
 solution.  
By compactness, some finite subset already has no solution;   this is only possible if $D(M)$ is finite.
 \eprf

 \sssec{Transitivity, orthogonality}  

  A  definable set $D$
is  {\em transitive}    if it has no proper, nonempty  definable subsets.  
(The usual word is ``atomic''.  One also says that {\em $D$  generates a complete type}.)  
 It is {\em (finitely) primitive}   if it admits no nontrivial definable equivalence relation (with finitely many classes).

 \<{remark} \lbl{exprim1} Let $A$ be a $\VF$-generated substructure of a model of $\ACVF$.  
  When $A$ is $\VF$-generated, we will see that an $\ACVF_A$-definable ball $b$ is never transitive in $\ACVF_A$; indeed it always contains an $A$-definable finite set.  But  
 $b$ is always $\ACVF_{A(b)}$-definable, and quite often it is transitive; cf. \lemref{transitive}.
  \>{remark}

 Two definable sets $D,D'$ are said to be {\em orthogonal} if 
 any definable subset of $D^m \times D^l$
  is a finite union of rectangles $E \times F$, $E \subseteq D^m$, $F \subseteq D^l$.
In this case, the rectangles $E,F$ can be taken to be almost definable. 
If the rectangles can actually be taken definable, we say $D,D'$ are {\em strongly
orthogonal}.

  \sssec{Types}  \lbl{Types}
  Let $S$ be a product of sorts, and let $M \models T$, $a \in S(M)$.  
 We write $tp(a)=tp(a;M)$ (the type of $a$) for the set of definable sets $D$ with $a \in D$;
 when $p = tp(a)$ we write $a \models p$.
    A {\em complete type}
 is the type of some element in some model.  If $q=tp(a)$, we say $a$ is a {\em realization} of $q$.
 The set ${\mathcal Tp}_S$ of complete types belonging to $S$ can be topologized:  a basic open set is the set of types including a given definable set $D$.
 The {\em compactness theorem} of model theory implies that this is a compact topological space:  if $\{D_i\}$ is any collection of definable sets with nonempty finite intersections, the compactness
 theorem asserts the existence of $M \models T$ with $\meet_i D_i(M) \neq \emptyset$.
 
 The compactness theorem is often used by way of a construction called {\em saturated models}; cf. \cite{CK}.
 These are models where all types over ``small'' sets are realized.  They enjoy excellent Galois-theoretic properties:  in particular if $M$ is saturated, then $\dcl(A_0) = Fix Aut(M/A_0)$ for any finite $A_0 \subseteq M$.  
 If $D$ is $\acl(A_0)$-definable, then there exists an $A_0$-definable $D'$ which is a finite union of
 $Aut(M/A_0)$-conjugates of $D$.   \lbl{types}

 A type $p$ can also be identified with the functor $P$ from models of $T$ (under elementary embeddings) into sets; $P(M) = \{a \in M: a \models p \}$.  As with definable sets, we speak as if $P$ is simply a set.  Unlike definable sets, the value of $P(M)$ at a single model does not determine $P$ (it could be empty; but it does determine $P$ if $M$ is sufficiently saturated.)

   Any definable map 
 $f: S \to S'$ induces a map $f_*: {\mathcal Tp}_S \to {\mathcal Tp}_{S'}$; as another consequence of the compactness theorem, $f_*$ is continuous.  We also have
 a restriction map from types of $T_A$ to types of $T$, $tp_{T(A)}(a) \mapsto tp_T(a)$.

 If $L \subseteq L'$ and $T \subseteq T'$, we say that $T'$ is an {\em expansion} of $T$.
In this case any $T'$-type $p'$ restricts to a $T$-type $p$.  If $p'$ is the {\em unique}
type of $T'$ extending $p$, we say that $p$ implies $p'$.  

The simplest kind of expansion is an expansion by constants, i.e. a theory $T_A$
(where  $A \leq M \models T$.)  If $c \in M^n$, or more generally if $c \in M^{eq}$, 
 the type of $c$ for $M_A$ is denoted
$tp(c/A)$.  It is rare for $tp(c)$ to imply $tp(c/A)$, but significant when it happens.

An instance of this is strong orthogonality:  it is easy to see
that strong orthogonality of   two definable sets $D,D'$
is equivalent to the following condition:

(*)   If   $A'$ is generated by elements of $D'$, then
 any type of elements of
$D$ generates a complete type over $A'$. 

The asymmetry in (*) is therefore only apparent.

Similarly, we have:

\<{lem} \lbl{tp-f} Let $D,D'$ be definable sets.   Then $(1) \iff (2), (3) \iff (4)$.
\begin{enumerate}
  \item Every definable function $f:D \to D'$ is piecewise constant, i.e. there exists a partition $D = \union_{i=1}^n D_i$ of $D$ into definable sets, with $f$ constant on $D_i$.
  \item If $d \in D, d' \in D'$, $d' \in \dcl(d)$ then $d' \in \dcl(\emptyset)$.
  \item  If $f:E \to D$ is a definable finite-to-one map, and $g: E \to D'$ is definable,
  then $g(E)$ is finite.
  \item   If $d \in D, d' \in D'$, $d' \in \acl(d)$ then $d' \in \acl(\emptyset)$.
\end{enumerate}
\>{lem}

\<{proof}  Let us show that  (3) implies (4).  Let  $M \models T$, $d \in D(M), d' \in D'(M)$, $d' \in \acl(d)$.  Then 
$d'$ lies in some finite $T_{d}$-definable set $D'(d) \subseteq D'$.  Since $T_d$ is obtained from $T$
by adding a constant symbol for $d$, there exists a formula $\phi(x,y)$ of the language
of $T$ and some $m$ such that $M \models \phi(d,d')$ and $M \models (\exists ^{\leq m} z) \phi(d,z)$.  Let $X_0 = \{( x:   (\exists ^{\leq m} y) \phi(x,y)\}$, $E = \{(x,y): x \in X_0, \phi(x,y)\}$,
$f(x,y)=x$, $g(x,y)=y$.   Then by (3) $g(E)$ is finite, but $d' \in g(E)$, so $d' \in \acl(\emptyset)$.

Next, (4) implies (3):  let $f,E,g$ be as in (3) , and suppose $g(E)$ is infinite.  
In particular, for any finite $F \subseteq \acl(\emptyset)$ there exists $d' \in g(E) \m F$.
Thus the family consisting of $g(E)$ and the complement of all finite definable sets has nonempty intersections
of finite subfamilies, so by the compactness theorem, in some $M \models T$ there exists $d' \in g(E) \m \acl(\emptyset)$.

Let $d \in E(M)$ be such that $d' =g(d)$.  Then $d' \in \acl(f(d))$, but $f(d) \in D$, contradicting (4).
Thus (4) implies (3).  

The equivalence of (1),(2) is similar.  \>{proof}

\<{example} \lbl{11nd}  Let $P$ be a complete type, and $f$ a definable function.
Then $f(P)$ is a complete type $P'$.  If $f$ is injective on $P$ then 
 there exist definable $D \supseteq P$,
$D' \supseteq P'$ such that $f$ restricts to a bijection of $D$ with $D'$.  \>{example}

\prf  For any definable $D'$, $f \inv D'$ is definable, so $P \subseteq f \inv D'$ or
$P \meet f \inv D' = \emptyset$.  Thus $P' \subseteq D'$ or $P' \meet D' = \emptyset$.
So $P'$ is complete.

Let $\{D_i\}$ be the family of definable sets containing $P$.  Let
 $R_i=\{(x,y) \in D_i^2: x \neq y, f(x)=f(y) \}$.  Then $\meet _i R_i = \emptyset$.   Since
the family of $\{D_i\}$ is closed under finite intersections, it follows from the compactness theorem that
for some $i$, $R_i = \emptyset$.  Let $D=D_i, D' = f(D)$.  \eprf

\sssec{Naming almost definable sets} \lbl{naming}
As special case of an expansion by constants, we can  move from a complete  theory $T$ to the theory $T_A$, where
$A=acl(\emptyset)$ is the set of all algebraic elements of a model $M$ of $T$, including imaginaries. 
The effect is a theory where  
  each class
of any definable equivalence relation $E$ with   finitely many classes is definable.   Since $T$ is complete, the  isomorphism type of 
$acl(\emptyset)$ in a model $M$ does not depend on the choice of model;  so the theory
$T_A$ is determined.
  A definable set in this theory corresponds to an almost definable set in $T$.

    When $D$ is a constructible set, $T|D$ denotes the
theory induced on $D$:  If $T$ eliminates quantifiers, 
the language is just the restriction to $D$ of the relations and functions of $L$.
If the language is countable, the countable models of $D_A$
are of the form $D(M)$, where $M$ is a countable model of $T_A$.)

\sssec{Stable embeddedness}  \lbl{stab-emb}  A definable subset $D$ of any product
of sorts (possibly imaginary) is called
  {\em stably embedded} (in $T$) if for any $A$, any $T_A$-definable
subset of $D^m$ is $T_B$-definable for some $B \subset D$.  For example, the set of open balls
is not stably embedded in $\ACVF$, since the set of open balls containing a point $a \in K$
cannot in general be defined using a finite number of balls.

\<{lem} \lbl{eilem}  Let $D$ be a family of sorts of $L$; let $T|D$ be the theory induced on the sorts $D$; 
If $D$ is stably embedded and  $T|D$ admits elimination of imaginaries, then for any
definable $P$ and definable $S \subset P \times D^m$, viewed as a $P$-indexed family
of subsets $S(a) \subseteq D^m , a \in P$, we have a definable function 
$f: P \to D^n$, with $f(a)$
a canonical parameter for $S(a)$.   
  \>{lem}
  
 \prf   By stable embeddedness there exists a family $S' \subset P' \times D^m$ yielding
 the same family, i.e. $\{S(a): a \in P \} = \{S'(a'): a' \in P' \}$, and with $P' \subseteq D^n$;
 using 
   elimination of imaginaries we can take $S'$ to be a canonical family; now
   $a$ define $f(a)$ to be the unique $a' \in P'$ with $S(a)=S'(a')$.  \eprf 
  
  \<{cor} \lbl{eilem-c} If $D$ is stably embedded and admits elimination of imaginaries, then
  for any  
 substructure $A$,  \begin{enumerate}
  \item   $(T_A)|D  = (T|D)_{A \meet D}$
  \item  for $a \in A$, $tp(a / A \meet D) $ implies $tp(a/D)$.
 
\end{enumerate}
   
   \qed  \>{cor}
  
  Examples of definable sets of $\ACVF$ satisfying the hypotheses include 
  the residue field $\k$, or the value group $\G$, as well as $\RV \union \G$.   The stable embeddedness in this case  is an immediate consequence of quantifier elimination; cf. \lemref{Vplus1}.

If $M$ is saturated and $D$ is stably embedded in $T$, then we have an exact sequence
$$1 \to Aut(M/D(M)) \to Aut(M) \to Aut(D(M)) \to 1$$
 where   $Aut(M/D(M))$ is the group of automorphisms of $M$ fixing 
$D(M)$ pointwise, and $Aut(D(M))$ is the group of permutations of $D(M)$ preserving all definable relations.
Moreover $Aut(M/D(M))$ is has  a good Galois theory; in particular elements with a finite orbit are
almost definable over some finite subset of $D$.  This and some other characterizations    can be found in the appendix
to \cite{CH}.

\sssec{Generic types}  

Let $T$ be a complete theory with quantifier-elimination.  Let ${\mathcal C}$
be the category of substructures of models of $T$, with $L$-embeddings, and let  ${\mathcal S}$ 
be the category of pairs $(A,p)$ with $A \in \Ob {\mathcal C}$ and $p$ a type over $A$. 
We define $\Mor((A,p),(B,q)) = \{f \in \Mor _{\mathcal C} (A,B): f ^* (q)  = p\}$.

  By a 
 {\em generic type} we will mean a function $p$ on $\Ob  {\mathcal C} $, 
 denoted $A \mapsto (p|A)$, such that $A \mapsto (A,p|A)$ is a functor 
$ {\mathcal C}  \to  {\mathcal S}$.   One example, when $T$ is the theory of algebraically closed
fields, is provided by any absolutely irreducible variety $V$:  given a field $F$, 
let $p|F$ be the type of an $F$-generic point of $V$, i.e. the type of a point of $V(L)$
avoiding $U(L)$ for every proper $F$-subvariety $U$ of $V$,
where $L$ is some extension field of $F$.  Other examples will be given below, beginning with
\exref{O-gen}.  

\<{lem} \lbl{generic-acl} Let $p$ be a generic type of $T$, and let $M \models T$, $a,b \in M$.
Let $c \models p | M$.
  \begin{enumerate}
  \item If $a \notin \dcl(\emptyset)$, then $a \notin \dcl(c)$.
    \item If $a \notin \acl(\emptyset)$, then $a \notin \acl(c)$.
    \item  If $a \notin \acl(b)$, then $a \notin \acl(b,c)$.
\end{enumerate}
    
 \>{lem}
 
\prf  (1) 
Since $a \notin \dcl(\emptyset)$, there exists $a' \neq a$ with $tp(a) = tp(a')$.
Let $c' \models p | <\{a,a'\}>$.  Since $tp(a)=tp(a')$, there exists an isomorphism
$<a> \to <a'>$; by functoriality of $p$, $tp(a,c) = tp(a',c)$.  
 If $a \in \dcl(c)$, then $a$ is the unique realization of $tp(a/c)$, so $a=a'$; a contradiction.
 
(2) If $a \in \acl(c)$, then for some $n$ there are at most $n$ realizations of $tp(a/c)$.
Since  $a \notin \acl(\emptyset)$, there exist distinct realizations $a_0,\ldots,a_n$
of $tp(a)$.  Proceed as in (1) to get a contradiction.

(3)  By (2) for $T_{<b>}$. \eprf
\<{subsection}{Grothendieck rings}.

 \lbl{conv}  We define the Grothendieck group and associated objects of a theory $T$;
cf. \cite{denef-loeser01}. 
    $Def(T)$ is the category of definable sets and functions.  
   Let ${\fC}$ be a subcategory of $Def(T)$.   
  We assume
    $Mor(X,Y)$
 is a sheaf on $X$:  if $X_1=X_2 \union X_3$ are subobjects of $X$,
 and $f_i \in Mor(X_i,Y)$ with $f_1 | (X_2 \meet X_3) = f_2 | (X_2 \meet X_3)$, then there exists $f \in Mor(X_1,Y)$ with $f|X_i = f_i$.  Thus
 the disjoint union of two constructible sets in $Ob \fC$ is also the 
 category theoretic disjoint sum.

   If  only the objects are given, we will assume $\Mor \fC$ is the collection of all  definable bijections between them.

The {\em Grothendieck semigroup } $\SG ({\fC}) $ is defined to be the  semigroup 
 generated by the isomorphism classes $[X]$ of  elements  $X \in \Ob {\fC}$, subject to the 
 relation:
 $$[X]+[Y] = [X \union Y] + [X \meet Y]$$
In most cases, $\fC$ has disjoint unions; then the elements of $\SG(\fC)$
are precisely the isomorphism classes of $\fC$.  

If $\fC$ has Cartesian products, we have a semiring structure given by:
$$[X] [Y] = [X \times Y]$$
In all cases we will consider when products are present, the symmetry isomophism $X \times Y \to Y \times X$
will be in the category, as well as the associativity morphisms, so that $\SG(\fC)$
is a commutative semiring.  

(The assumption on Cartesian products is taken to include the presence
of an object   $\{p\} = X^0$ such that the bijections
$X \to \{p\} \times X$,
$x \mapsto (p,x)$, and $X \to X \times \{p\}$, $x \mapsto (x,p)$, are in $\Mor_\fC$ for all $X \in \Ob_\fC$.  All such $p$ gives the same
element $1=[\{p\}] \in \K({\fC})$, which serves as the identity element of 
the semiring.)

Let $\K(\fC)$ be the Grothendieck group, the formal groupification of $\SG(\fC)$.
When $\fC$ has products, $\K(\fC)$ is a commutative ring.

 We will often have dimension filtrations on out categories, and hence on the semi-ring.

By an {\em semi-ring ideal}   we mean a congruence relation, i.e. an equivalence
relation on the semi-ring $R$ that is a sub-semi-ring of $R \times R$.  
To show that an equivalence
relation
$E$ is a congruence on a commutative semi-ring $R$, it suffices to check that if $(a,b) \in E$ then
$(a+c,b+c) \in E$ and $(ac,bc) \in E$.  
 
{\bf Remark}  
When $T$ is incomplete, let $S$ be the (compact, totally disconnected) space of completions of
$T$.  Then $\{\K(t): t \in S \}$ are the fibers of a sheaf of rings over $S$.
$\K(T)$ can be identified with the ring of continuous sections of this sheaf.
In this sense, Grothendieck rings reduce to the case of complete theories.  

This last remark is significant even when $T$ is complete:    if one adds
a constant symbol $c$ to the language, $T$ becomes incomplete,
and so the Grothendieck ring of $T$ in $L(c)$ is the Boolean power
of $\K(T_a)$, where $T_a$ ranges over all $L(c)$-completions of $T$.
Say $c$ is a constant for an element of a sort $S$.  Then an $L(c)$-definable
subset of a sort $S'$ corresponds to an $L$-definable subset of $S \times S'$.   This allows an inductive analysis of
the Grothendieck ring of a structure, given good information about
definable sets in one variable  (cf. \lemref{collect}).  
  
 \sssec{Groups of functions into $\RR$} \lbl{functions} 
 Let $\fC(T)$ 
 be a subcategory of the category of definable sets and bijections, 
 defined systematically for $T$ and for expansions by constants   $T$.  
 Let $\RR(T) = \SG (\fC(T))$ be the Grothendieck semigroup of $\fC(T)$.  
 When $V$ is a definable set, we let $\fC_V$, $\RR_V$
denote the corresponding objects over $V$;  the objects of $\fC_V$
are definable sets $X \subseteq (V \times W)$ such that  for any $a \in V$, $X_a \in \fC_a$,
and similarly the morphisms.    
 In practice,
 $\RR$ will be the Grothendieck semgroup of all definable sets and definable isomorphisms
 satisfying some definable conditions, such as a boundedness condition on the 
 objects, or a ``measure preservation'' condition on the definable bijections.

To formalize the notion of  
 of ``definable function into $\RR$''  we will need to look at classes $X_a$ of parametrically definable sets.
The class of $X_a$ makes sense only in the 
Grothendieck groups associated with $\T_a$, not $\T$.    Moreover equality of such
classes, say of $X_a$ and of $X_b$, begins to make sense only in Grothendieck groups
of $\T_{(a,b)}$.  Expressions like
$$ [X] \oeq{a,b} [Y]$$
will therefore mean:   $X,Y$ are both definable in $\T_{a,b}$, $[X],[Y]$ denote
their classes in the Grothendieck group of $\T_{a,b}$, and these classes are equal.

  If $V$ is a definable set, we define 
the  {\em  semigroup of    definable functions $V \to \RR$}, denoted 
$Fn(V,\RR)$.  
An element of  $Fn(V,\RR)$ is represented
by a definable $X \in \fC_V$, viewed as the function:  $a \mapsto [X_a]$,
where $[X_a]$ is a class in $\RR_a$.    $X,X'$ represent the same function
if for all $a$, $[X_a],[X'_a]$ are the same element of $\RR_a$.  Note that despite
the name, the elements of $Fn(V,\RR)$   should actually be viewed as sections $V \to \Pi_{a \in V} \RR_a$.

Addition is given by disjoint union in 
the image (i.e. disjoint union over $X$.)

Usually $\RR$ has a natural grading by dimension; in this case $Fn(V,\RR)$ inherits the
grading.    

Assume 
   $V$ is a definable group and $\RR = \SG(T)$ is the Grothendieck semiring of
   all definable sets and functions of $T$,
   there is a natural  convolution product on 
$Fn(V,\Rr)$.  If $h_i(a) = [H_i(a)]$, $H_i \subset V \times B_i$,
  the convolution
$h_1 * h_2$ is   
represented by 
$$H = \{(a_1+a_2,(a_1,a_2,y_1,y_2)):  (a_i,y_i) \in H_i \} \subseteq V \times (V^2 \times B_1 \times B_2) $$
so that   $h_1 * h_2 (a) = H(a) = \{(a_1,a_2,y_1,y_2):   (a_i,y_i) \in H_i, a_1+a_2=a \}$.

\sssec{Grothendieck groups of orthogonal sets}

 \lemm{tensor0}  Let $T$ be a theory with two strongly orthogonal definable sets $D_1,D_2$,
 $D_{12} = D_1 \times D_2$.  
 Let $\SG D_i[n]$ be the Grothendieck semigroup of definable subsets of $D_i^n$.  Then 
 $\SG  D_{12}[n]  \iso  \SG D_1 [n] \tensor \SG D_2[n]$  \>{lem}
 
 \prf  This reduces to $n=1$.  
 Given  definable sets $X_i \subseteq D_i^n$, it is clear
 that the class of $X_1 \times X_2$ in $\SG D_{12}[n]$ depends only on the classes
 of $X_i$ on $D_i[n]$.  Define
 $[X_1] \tensor [X_2] = [X_1 \times X_2]$.  This is clearly $\Zz$-bilinear, and so extends
 to a homomorphism   $\eta:  \SG D_1[1] \times \SG D_2 [1] \to \SG D_{12}[1]$.  By 
 strong orthogonality, $\eta$ is surjective.    
 
 To prove injectivity, note that any element of $\SG D_1 [n] \tensor \SG D_2[n]$ can
 be written $\sum [X_1^i] \tensor [X_2^i]$, with $X_1^1,\ldots,X_1^k$ pairwise disjoint.
 To see this, begin with some expression $\sum [X_1^i] \tensor [X_2^i]$; use the
 relation $[X' \du X''] \tensor [Y] = [X'] \tensor [Y] + [X''] \tensor [Y]$ to replace
 the $X_1^i$ by the atoms of the Boolean algebra they generate, so that the new $X_i^i$ 
 are equal or disjoint; finally use the relation 
 $[X' \tensor Y'] + [X' \tensor Y''] = [X'] \tensor [Y' \du Y'']$ to amalgamate the terms with 
 equal first coordinate.  
 
 Hence it suffices to show that if $[\union_i X_1^i \times X_2^i]= [\union _j (Y_1^i \times Y_2^i]$,
 with the $X_1^i$ and the $Y_1^i$ pairwise disjoint, then
  $\sum [X_1^i] \times [X_2^i]= \sum [Y_1^i] \times [Y_2^i]$.  Let
  $F: \union_i X_1^i \times X_2^i \to \union _j  Y_1^i \times Y_2^i$ be a definable bijection.
 By strong orthogonality, the graph of $F$ is a disjoint union of rectangles.   Since $F$ is a bijection, 
 it is easy to see
 that each of these rectangles has the form $f_1^k \times f_2^k$,
 where for $\nu=1,2$, $f_\nu^k: X_\nu(k) \to Y_\nu(k)$ is a bijection from a subset of $\union_i X_\nu^i$ to a subset of $\union_j Y_\nu^j$.  The rest follows by an easy combinatorial 
 argument; we omit the details, since  a somewhat more complicated 
 case will be needed and proved later, see   \propref{tensor}. 
   \eprf

\sssec{Integration by parts}  The following will be used only in  \secref{exgam}, to study the Grothendieck semiring
of the valuation group.   
   
\<{defn}  \lbl{discrete}
Let us say that $Y \in \Ob \fC$ is {\em treated as discrete}  if
  for any $X \in \Ob \fC$  and any  definable $F \subset X \times Y$
such that $T \models F$ is the graph of a function,  the 
 projection map $F \to X$ is an invertible element of $\Mor_{\fC} (F, X)$.

\>{defn}

To explain the terminology, suppose each $X \in \Ob_{\fC}$
is endowed with a measure $\mu_X$, and $\fC$ is the category of 
measure preserving maps.   If $\mu_Y$ is 
the counting measure, and $\mu_{X \times Y}$   is the product measure, then for  any function
$f:X \to Y$, $x \mapsto (x,f(x))$ is measure preserving.

We will assume  $\fC$ is closed under products.

If $Y_1,Y_2$ are treated by $\fC$ as discrete, so is $Y_1 \times Y_2$:
  if $F \subset  X \times (Y_1 \times Y_2)$
is the graph of a function $X \to (Y_1 \times Y_2)$, 
then the projection to $F_1 \subset X \times Y_1$ is the graph of a function, hence the projection $F_1 \to X$ is  in $\fC$;
now $F \subset (F_1 \times Y_2)$ is the graph of a function,
and so $F \to F_1$ is invertibly represented too; thus so is the composition.)  In particular if $Y$ is discretely treated, any    bijection
$U \to U'$ between subsets of $Y^n$ is represented in $\fC$.

If $\RR$ is a Grothendieck group or semigroup, we  write $[X] \oeq{R} [Y]$ to mean that $X,Y$ have the same class in 
$\RR$.

\<{lem} \lbl{bc}  Let $f,f' \subset X \times L$ be objects of $\fC$ such that
  $[f(c)] \oeq{\K(\fC_c)} [f'(c)]$ 
 for any $c$ in $X$.   Then
$[f] \oeq{\K(\fC)} [f']$.   Similarly for $\SG$.   \>{lem}  

\proof  By assumption,  there exists 
$g(c)$ such that $f(c)+g(c),f'(c)+g(c)$ are   $\fC_c$-isomorphic.
By compactness  (cf. end of proof of \lemref{collect})  this must be uniform (piecewise in $L$, and
hence by glueing globally):  there exists a definable
$g \subset Z \times L$ and a definable isomorphism
$f+g \iso f'+g$, inducing the isomorphisms of each fiber.  By definition
of $\fC_c$, and since $\fC$ is closed under finite glueing,
$f+g,f'+g$ are in $\Ob \fC$  and the isomorphism between them is in $\Mor \fC$.  \qed

 Let $L$ be an object of $\fC$, treated as discrete in $\fC$, and assume given
 a definable partial ordering on $L$.

 \<{notation}  \lbl{sumnot} Let $f \subset X \times L$.   For $y \in L$, 
 let $f(y) = \{x: (x,y) \in f \}$.  Denote:
 $\sum_{\gamma < \beta} f(\gamma)  = [\{(x,y): x \in f(y), y < \gamma \}]$   \>{notation}
 
 \<{notation}  \lbl{sumnot2} Let $\phi: L \to \K(X)$ be a constructible
 function, represented by $f \subset X \times L$; so that
 $\phi(y) = [f(y)]$,    $f(y) = \{x: (x,y) \in f \}$.  Denote:
 $\sum_{\gamma < \beta} \phi(\gamma)  = [\{(x,y): x \in f(y), y < \gamma \}]$   \>{notation}

Note by \lemref{bc} that this is well-defined.

Below, we write $fg$ for the pointwise product of two functions
in $\K(\fC)$; $[fg(y)] = [f(y) \times g(y)]$.  

\<{lem}[Integration by parts]  \lbl{ibp}   Let $\G$ be an object of $\fC$, treated as discrete in $\fC$, and assume given
 a definable partial ordering of $\G$.    Let $f \subset X \times \Gamma$, 
$F(\beta) = \sum_{\gamma < \beta} f(\gamma)$,
$g \subset Y \times \Gamma$, 
$\bG(\beta) = \sum_{ \gamma \leq \beta} g(\gamma)$.

Then 
$$F\bG(\beta) = 
\sum_{\g < \beta} f\bG(\gamma)  + \sum_{\gamma \leq \beta} Fg(\gamma)$$
 
\>{lem}

\proof   Clearly
$$F\bG(\beta) = \Sum_{\g< \beta, \g' \leq \beta} f(\g)g(\g')$$
We split this into two sets: $\g < \g'$, and $\g' \leq \g$.  Now
$$\Sum_{\g < \g' \leq \beta} f(\g) g(\g') = \Sum_{\g' \leq \beta} F(\g')g(\g')$$
$$\Sum_{\g' \leq \g < \beta} f(\g) g(\g') = \Sum_{g < \beta} f(\g)\bG(\g')$$
   The lemma follows.  \qed

 This is particularly useful when $L$ is treated as discrete in $\fC$,
 since then, if the sets $f(\g)$ are disjoint, $[f] = [\union_\gamma f_\gamma]$.  Another version,
 with $G(\beta) =  \sum_{\gamma < \beta} g(\gamma)$:
 $$ FG(\beta) = \sum_{\gamma < \beta} (fG +gF + fg)(\g)$$

 \>{subsection}

\>{section}

\<{section}{Some $C$-minimal geometry.}
\lbl{minimality}

We will isolate the main properties of the theory $\ACVF$, and work with
an arbitrary theory $\T$ satisfying these properties.  This includes 
 the   rigid analytic expansions $\ACVFR$  of \cite{lipshitz}.  
 
The right general notion, $C$-minimality, has been introduced and studied in \cite{Cmin}.
They obtain many of the results of the present section.    
Largely for expository reasons, we will describe a  slightly less general version; it is essentially minimaility with respect to an ultrametric structure in the
sense of \cite{poizat}.  We will use notation suggestive of the case of valued fields;
thus denoting the main sort by $\VF$, and a binary function by $\val(x-y)$.  
Some additional assumptions will be made explicit later on.

Let $T$ be a theory in a language $L$, extending a theory $\mathsf{T}$ in a language $\mathsf{L}$.
$T$ is said to be {\em ${ \mathsf T}$-minimal} if   for any $M \models T$, any $L_M$-formula {\em in one variable} is $T_M$-equivalent to an  $\mathsf{L}_M$ formula.

More generally, if $D$ is a definable subset of $T$ (i.e. a formula of $L$), we say
that $D$ is {\em $\mathsf{T}$-minimal} if  if every for any $M \models T$, any $T_M$-definable
subset of $D$   is $T_M$-equivalent to one defined by an   $\mathsf{L}_M$ formula.

\vspace{4mm}
\paragraph{{\bf Strong minimality}} Let  $\mathsf{L}= \emptyset$.  The only atomic formulas of $\mathsf{L}$ are
  thus equalities $x =y$ of two variables.    $\mathsf{T}$ the theory of infinite sets.     $\mathsf{T}$-minimality is known as {\em strong minimality}; see 
  \cite{BL}, \cite{pillay}.   A theory $T$  is strongly minimal iff for any $M \models T$, 
  any $T_M$-definable subset of $M$ is finite or cofinite. 
  For us the primary example of a strongly minimal theory is $ACF$, the theory of algebraically closed fields.
  
  Let $M \models T$.  
 If $D$ is strongly minimal, and $X$ a definable subset of $D^*$, we define the $D$- {\em dimension} of $X$
 to be the least $n$  such that $X$ admits a $T_M$-definable map into $D^n$ with 
 finite fibers.  In the situation we will work in, there will be
  more than one  definable
 strongly minimal set up to isomorphism, and even up to definable isogeny; 
  in particular there will be the various sets of $\RES_M$.  However, between
 any of these, there exists an $M$-definable isogeny;
   so the $\k$-dimension agrees with the $D$-dimension for any of them. We will call it the $\RES$-dimension.  It   agrees with {\em Morley rank}, a notion defined in greater generality, that we will not otherwise need here.

  \vspace{4mm}
\paragraph{{\bf O-minimality}}  
 $\mathsf{L}= \{<\}$, $\mathsf{T}=DLO$ the theory of dense linear orders  without endpoints. (Cf. \cite{CK})  DLO minimality is known as $O$-minimality, and can also be stated thus:   any $T_M$-definable subset of $M$ is a finite union of points and intervals.  This 
    also forms the basis of an extensive theory; see \cite{vddries}.  
        
    Let $D$ be $O$-minimal.  Then the   $O$-minimal dimension of a definable set $X \subseteq D^*$ is   the 
    least $n$  such that $X$ admits a $T_M$-definable map into $D^n$ with bounded
 finite fibers.
 
The Steinitz exchange principle states that if $a \in \acl(B \union \{b\})$ but $a \notin \acl(B)$, then 
$b \in \acl(B \union \{a\})$.  

This holds for both strongly minimal and O-minimal structures, cf. \cite{vddries}.

For us the   relevant $O$-minimal theory is $\DOAG$   itself.   We will occasionally use
stronger facts valid for this theory.  
Quantifier elimination for
$\DOAG$ implies that

 \lemm{doag-fns}  
 \begin{enumerate}
  \item  any parameterically  definable function $f$ of one variable is piecewise affine; 
 there exists a finite partition of the universe into intervals and points,
such that on each interval $I$ in the partition, $f(x) = \alpha x + c$ for some rational
$\alpha$ and some definable $c$.  
  \item $\DOAG$ admits elimination of imaginaries.
 
\end{enumerate}
  \>{lem}

\prf  (1)  Follows from quantifier elimination for
$\DOAG$.

(2)  It   follows from (1) that  any   function definable with parameters in $\DOAG$  has a canonical code, consisting of the endpoints
of the intervals of the coarsest such partition, together with a specification of the rationals
$\alpha$ and the constants $c$.   But from this it follows on general grounds that
every definable set is coded (cf. \cite{acvf1} 3.2.2).  So $\DOAG$ admits elimination of imaginaries.  \eprf

\vspace{4mm}
\paragraph{{\bf $C$-minimality}}  

Let $\mathsf{T}=\mathsf{T}_{um}$ be the theory of ultrametric spaces, or equivalently chains of equivalence relations.    (cf. \cite{poizat}.)

  In more detail,   $\mathsf{L}$ has two sorts,   $\VF$ and $\G_\infty$.   
  The relations on $\G_\infty$ are a constant $\infty$ and a binary relation $<$.  In addition
  $\mathsf{L}$ has a function symbol $\VF^2 \to \G_{\infty}$, written $\val(x-y)$.  
  $\mathsf{T}$ states:
    \begin{enumerate}
  \item  $\G_\infty$ is a dense linear ordering with no least element, but
  with a greatest element $\infty$.
  \item $\val(x-y)=\infty$   iff $x=y$.
  \item  $\val(x-y) \geq \alpha$ is an equivalence relation; the classes are called {\em closed $\alpha$-balls.}   Hence so is the relation $\val(x-y)>\alpha$; whose classes are called
  {\em open $\alpha$-balls.}  
\item  Let $\G=\G_\infty \m \{\infty\}$.  For $\alpha \in \G$, every closed $\alpha$-ball 
contains infinitely many open $\alpha$-balls. \end{enumerate}

A $\mathsf{T}_{um}$-minimal theory will be said to be $C$-minimal.  The notion considered in \cite{Cmin} is a little  more general, but for  theories $\mathsf{T}_{um}$ they coincide.
Since we will be interested in fields, this level of generality will suffice.

 A theory $T$ extending $ACVF$  is $C$-minimal
iff for any $M \models T$, every $T_M$-definable subset of $\VF(M)$ is a Boolean combination 
of open balls, closed balls and points.  If $T$ is $C$-minimal, $A \leq M \models T$, and $b$ is an $A$-definable
ball, or an infinite intersection,  let $p^b_A$ be the collection of $A$-definable sets not contained in a finite union of proper sub-balls of $b$.
Then by $C$-minimality, $p^b_A$ is a complete type over $A$.

Let $T$ be $C$-minimal.  Then in $T$, $\G$ is O-minimal; and for any closed $\a$-ball $C$,
the set of open $\a$-sub-balls of $C$ is strongly minimal.   Denote it $C/(1 + \Mm)$.  
  (These facts are immediate
from the definition.)

Assume $T$ is $C$-minimal with a distinguished point $0$.  We define: $\val(x)=\val(x-0)$;
$\Mm = \{x: \val(x)>0 \}$.  Let
  $\cB$ be the family of all closed balls,
  including points.  Among them are
 ${\cB^c}_\alpha(0) = \{x: \val(x) \geq \alpha \}$.  Let
  $\RV = \dU _{\g \in \G} B^c_\g(0) / (1+\Mm)$, and let $\rv: \VF \m \{0\} \to \RV$
  and $\valr: \RV \to \G$ be  the natural map.
By  an $\rv$-ball we mean an open ball of the form $\rv \inv (c)$.  

The $T$-definable fibers of $\valr$ are referred to, collectively, as $\RES_T$. 
 Later we will fix a theory $\T$, and write $\RES$  for $\RES_\T$; we will also
 write $\RES_A$ for $\RES_{\T_A}$.  The unqualified notion ``definable'', as well as many derived notions,
 will implicity refer to $\T$.

A certain notion of genericity plays an essential role in these theories.

\<{example}  Let $T$ be a strongly minimal theory.  For any $A \leq M \models T$,
any $A$-definable set is  finite or has finite complement.  Therefore, the collection of co-finite 
sets forms a complete type.  A realization of this type is called a {\em generic} element of $M$, over $A$.
\>{example}

\<{example} \lbl{O-gen} Let $T$ be an O-minimal theory.  For any $A \leq M \models T$,
any $A$-definable set contains, or is disjoint from, 
an infinite interval $(b,\infty)$ for some $b \in M$.  The set of $A$-definable sets containing such an interval
is thus a complete type, the generic type of large elements of $\G$.  Similarly,
the set of $A$-definable sets containing an interval $(0,a)$ with $0<a$ is {\em the generic
type of small positive elements.}  More generally, given a subset $S \subseteq A$
$S' = \{b \in A: (\forall s \in S) (s < b )\}$; then the
definable sets $x > a (a \in S), x< b (b  \in S')$ generate a complete type over $A$, called the
type of elements {\em just bigger than $S$ }.  \>{example}

\<{defn}  Let $T$ be $C$-minimal.
Let  $b$ be a  $T_A$ -definable ball, or an infinite intersection of balls.
 The generic type $p_b$ of $b$ is defined by:  $p_b | {A' }= p^b_{A'}$, for any $A \leq A' \leq M \models T$.   
   \>{defn}

The completeness follows from $C$-minimality, since for any $A'$-definable subset  $S$ of $b$, either $S$
is contained in a finite union of proper sub-balls of $b$, or else the complement $b \m S$ is contained
in such a finite union.

A realization of $p_b  | A'$ is said to be a {\em generic point  of $b$ over $A'$}.  An $A'$- definable
set is said to be {\em $b$-generic} if it contains a generic point of $b$ over $A'$.
 
See \secref{generics} for some generalities about generic types.  For our purposes
it will suffice to consider generic types in one  $\VF$ variable.   
 For more information see   \cite{acvf1} \S 2.5.
 
\<{rem}\lbl{1types}  
If $A = \acl(A)$ then any type of a field element $tp(c/A)$ coincides with 
$p_b|A$, where $b$ is the intersection of all $A$-definable balls containing $c$. \>{rem}
This is intended to include the case of closed balls of valuative radius $\infty$, i.e. points;
these are the algebraic types $x=c$.    
Note also the degenerate case that $c$ is not in any $A$-definable  ball; then $b=\VF$ and $tp(c/A)$
is the generic type of $\VF$ over $A$.

Not
every generic 1-type is of the form $p_b$ for a ball $b$ as above.      For instance, let $b$ be an open ball, $c \in b$; then 
the generic type $p_b((x-c)^{-1})$ is not of this form.) 

For $\V$-minimal theories (defined below) it can be shown that
every generic 1-type  is of the form $p_b$ or $p_b((x-c)^{-1})$.

Let $\T$ be  a $C$-minimal theory. 
 Let $b$ be a definable ball, or an 
infinite intersection of definable balls.  
We say that $b$ is centered 
   if it contains a proper definable  finite union
of balls.  If $b$ is open, or a properly infinite intersection of balls, we have:

(*) If   $b$ contains a proper finite union of balls, then it contains a definable closed ball
(the smallest closed ball containing the finite set.) 

 For $C$-minimal fields of residue
characteristic $0$, 
(*) is true of closed balls: 
  the set of maximal open sub-balls of $b$ forms an affine
 space over the residue field $\k$, where the center of mass of a finite set is well-defined. 
 
%  group of affine automorphisms over $\k$  
 %fixing a finite set can only be a group of roots of unity, fixing a point too.  

Clearly $b$ is centered over $\acl(A)$ if and only if  it is centered over $A$.   The term ``centered''
will be justified to some extent by the assertion of \lemref{red1c1}, that when 
 $A$ is generated by elements of $\VF \union \RV \union \G$,  any $A$-definable closed ball  contains an $A$-definable point, and thus a centered ball has a definable ``center''.

\lemm{centeredtrans}  $b$ is centered over $A$ iff $b$ is not transitive over $A$. \>{lem}

This is immediate from the definition, and from $C$-minimality,
since any proper definable subset would have to be
a Boolean combination of balls.      
  
An often useful corollary of $C$-minimality:  
\<{lem} \lbl{vstar}  Let $T$ be $C$-minimal,  $X$   a definable subset of $\VF$,
and $Y$   a definable set of disjoint balls.  Then for all but finitely many $b \in Y$,
either $b \subseteq X$ or $b \meet X = \emptyset$.  \>{lem}

\proof $X$ is a finite Boolean combination of balls, so it suffices to prove this when $X$ is a ball;
then $X$ is contained in at most one ball $b \in Y$; for any other $b \in Y$, 
either  $b \subseteq X$ or $b \meet X = \emptyset$.  \qed

\<{lem} \lbl{transitive} Let $(b_t: t \in Q) $ be a  
definable family of pairwise disjoint balls.  Then for any nonalgebraic $t \in Q$, 
 $b_t$ is transitive over $<t>$     \>{lem}

\proof     
Consider a  definable $R' \subseteq Q \times \VF$ with 
$R'(t) \subseteq b_t$.    Let $Y = \union _{t \in Q} R'(t)$.  Then $Y$ is a definable subset of $\VF$, hence a finite   combination of a finite set $H$
of balls.   The $b_t$ are pairwise disjoint, so at most finitely many can contain 
an element of $H$, and thus no nonalgebraic $b_t$  contains an element of $H$.  
Thus each ball in $H$ is disjoint from,  or contains,   any given $b_t$.  It follows that
$Y$  is disjoint from,  or contains,   any given $b_t$. Thus $b_t \meet Y$
cannot be a non-empty proper subset of $b_t$.  \qed

\paragraph{{\bf Internalizing finite sets}}

The following lemma will be generalized later to finite sets of balls.
It is of such fundamental importance in this paper that we include it separately in its
simplest form.  The failure of this lemma in residue characteristic $p>0$ is  the  main
reason for the failure of the entire theory to generalize, in its present form.  Recall the definition of $\RV$,
\S \ref{language}.
 
\lemm{finite}  Let $\T$ be a  $C$-minimal theory of fields  of residue characteristic $0$ (possibly with additional structure), $A \leq M \models  \T$.   Let $F$ be a finite $\T_A$-definable subset of $\VF^n$.
  Then there exists $F' \subseteq \RV^m$,
and a $\T_A$-definable bijection $h: F \to F'$.
\>{lem}

\prf  First consider $F = \{c_1,\ldots,c_n\} \subseteq \VF$.  Let $c = (\sum_{i=1}^n c_i)/n$ be the average;
then $F$ is $\T_A$-definably isomorphic to $\{c_1-c,\ldots,c_n-c\}$.  Thus we may assume
the average is $0$.    If there is no nontrivial $A$-definable equivalence
relation on $F$, then $\val(x-y)=\alpha$ is constant on $x \neq y \in F$.
In this case $\rv$ is injective on $F$ and one can take $h=\rv$.  Otherwise,
let $E$ be a nontrivial $A$-definable equivalence
relation on $F$.  By an $E$-symmetric polynomial, we mean a polynomial $H(x_1,\ldots,x_n)$ with
coefficients in $A$, invariant under the symmetric group on each $E$-class.  
For any such $H$, $H(F)$ is a $\T_A$-definable set with $<n$ elements.
There exists $H$ such that $H(F)$ has more than one element.  By induction,
there exists an injective $A$-definable function $h_0: H(F) \to \RV^m$.
Let $h_1 = h_0 \circ H$. 
For $d \in h_0(H(F))$, and $d'=h_0 \inv d $,  let $F_d = H^{-1} h_0 \inv (d) = H \inv (d')$.   
   By induction again, there exists an $A(d)=A(d')$-definable 
injective function $g_d: F_d \to \RV^{m'}$.  (We can take the same $m'$ for all $d$.)
Define $h(x) = (h_1(x),g_{h_1(x)}(x))$.  Then clearly $h$ is $A$-definable and injective.  

The case $F \subseteq \VF^n$ follows using a similar induction, or by finding a linear projection
with $\Qq$-coefficients $\VF^n \to \VF$ which is injective on $F$.
\>{proof}

\ssec{Basic geography of $C$-minimal structures}

Let $\T$ be a $C$-minimal theory.   We begin with a rough study
of the   existence and non-existence of definable maps between various regions of the structure:  $\k,\G,\RV,\VF$ and $\VF / \Oo$.

 We will occasionally refer to   {\em stable} definable sets.
 
    A definable set $D$ of a theory $T$ is called
{\em stable} if there is no model $M \models T$ and $M$-definable relation $R \subseteq D^2$
and infinite subset $J \subseteq M(D)$ such that $R \meet J^2$ is a linear ordering.
This is a  model-theoretic  finiteness condition, greatly generalizing
finite Morley rank, and in turn strong minimality (cf. \cite{pillay}.)

 It is shown in 
 \cite{acvf1}, that
 a definable subset of $\ACVF_A^{eq}$   is stable if and only if it has finite Morley rank, if an only if it admits no parametrically definable map onto an interval of $\G$; and this is if and only it embeds, definably over $\acl(A)$, into a finite dimensional
$\k$-vector space.   These vector spaces have the general form 
$\Lambda \ /\Mm \Lambda$,  with $\Lambda \leq \VF^n$ a lattice.  Within the sorts we are using here, the  relevant ones   are the finite products of vector spaces of $\RES$.  
More generally in a $C$-minimal structure with sorts $\VF,\RV$, all stable sets are definably embeddable (with parameters) into $\RES$.   We will however make no use of these facts,
beyond justifying the terminology.  Thus ``$X$ is a stable  definable set'' can simply be read as `` there exists a  definable bijection between  $X$ and a subset of  $\RES^*$''. 
\vspace{2mm}

 The first fact is the unrelatedness of $\k$ and $\G$.

\<{lem} \lbl{orth1}  Let $Y$ be a stable definable set.  Then $Y,\G$ are strongly orthogonal.
In particular, any definable map from $Y$ to $\G$ has finite image.
\>{lem}

\<{proof} We prove the second statement first:  let $M \models \T$.  Let $f: Y \to \G$
be an $M$-definable map.  Then  $f(Y)$ is stable, and linearly ordered by $<_\G$;
hence by definition of stability it is finite.  
 
 Let $\g = (\g_1,\ldots,\g_m) \in \G$.  We have to show that for a $Y$-generated structure $A$, $tp(\g)$
implies $tp(\g/A)$.  It suffices to show that for any $a_,\ldots,a_n \in A$,
 $tp(\g_i/<\g_1,\ldots,\g_{i-1}>)$
implies $tp(\g_i / <\g_1,\ldots,\g_{i-1},a_1,\ldots,a_n>)$, for each $i$.  By
passing to $T_{<\g_1,\ldots,\g_{i-1}>}$ we may assume $m=1, \g \in \G$.
Similarly we may assume $n=1$; let $a=a_1 \in Y$.  
  To show that $tp(\g)$ implies $tp(\g/a)$, 
it suffices to show that any $T_a$-definable subset of $\G$ is definable.
By O-minimality, any $T_a$-definable subset of $\G$ is a finite union of intervals,
so (in view of the linear ordering) it suffices to show this for intervals $(c_1,c_2)$.  But if the interval is $T_a$-definable
then so are the endpoints, so $c_i = c_i(a)$ is a value of a definable map $Y \to \G$.
But such maps have finite images, so $c_i$ lies in a finite
definable set.  Using the linear ordering, we see that $c_i$ itself is definable, and hence so is the interval.  \>{proof}

\<{lem} \lbl{nosec} There are no definable
sections of $\valr: \RV \to \G$ over an infinite subset of $\G$.  In fact
if $Y \subset \RV^n$ is definable and $\valr$ is finite-to-one on $Y$,
then $Y$ is finite.
\>{lem}

\<{proof}  Looking at the fibers of the projection
of $Y$ to $\RV^{n-1}$, and using induction, we reduce the lemma to the case $n=1$.
     In this case, by \lemref{vstar},
every definable set is a Boolean combination of pullbacks
by $\valr$ of subsets of $\G$, and finite sets. 
                                                                                                     \>{proof}

\<{lem}  \lbl{fmrf}   Let $M \models \T$  and let $Y \subset \cB^n $ be an infinite definable set.  Then
   there exists
        a surjective $M$-definable map of $Y$  to a proper interval in $\G$.
  \>{lem}
 
\<{proof}    Since $\G$
is O-minimal, any infinite $M$-definable subset contains a proper interval.  Thus it suffices
to find an $M$-definable map of $Y$ onto an infinite subset of $\G$.

 If the projection of $Y$ to $\cB^{n-1}$ as well
as every fiber of this projection are finite, then $Y$ is finite.  Otherwise, replacing
$Y$ by one of the fibers or by the projection, we reduce inductively to the case $n=1$.

Let $v(y) \in \G$ be the   valuative radius of the ball $y$.  Then $v(Y)$ is an $M$-definable subset of $\G$.  If it is infinite, we are done; otherwise  we may assume   all elements of $Y$ have the same 
 valuative radius $\g$.
 
 Let  $W = \union Y$.   
 By $C$-minimality,
 $W$ is a Boolean combination of balls $b_i$ (open, of valuative radius $<\gamma$, or closed, of valuative radii $\delta_i \leq \g$.)  If 
 $W$ contains some $W'=b_i \setminus (b_{j_1} \union \cdots \union b_{j_l})$,
 where $b_{j_i}$ is a proper sub-ball of $b_i$, and  $\delta_i  <\gamma$, pick  a point $c$ in $W'$; then for any
 $\delta $ with $\gamma > \delta > \delta_i$ there exists $c' \in W'$
 with $\val(c-c') =\delta$.  It follows that the balls $b_\g(c),b_\g(c')$ of radius
 $\gamma$ around $c,c'$ are both in $Y$; but infinitely many such $\delta $ exist;
 fixing $c$, we obtain a map $\b_\g(c') \mapsto \val(c-c')$ into an infinite subset of $\G$.  
 
 Otherwise, $W$ can only be
a finite set of balls of valuative radius $\g$.  So $Y$ is finite.   
                                                                                                      \>{proof}

\<{cor}  \lbl{fmrf-r}  $\cB^n$ contains no  stable  definable set.   In particular $\VF$ contains no strongly minimal set.   \qed \>{cor}
By contrast,

\<{lem} \lbl{fmrf-rem}   
Any infinite definable subset of $\RV^n$ contains
a  strongly minimal $M$-definable subset.  \>{lem}
\prf  

By \lemref{nosec}, the
 inverse
image of some point in $\G^n$ must be infinite. \eprf

\<{lem}    \lbl{fmrf-rem2} Let $M \models \T$.   Let $Y \subseteq \cB$ be a definable set.
Let $\rad(y)$ be the  valuative radius of the ball $y$.  
Then either $\rad: Y \to \G$ is finite-to-one, or else 
there exists an $M$- definable map of an   $M$-definable $Y' \subseteq Y$  onto a strongly minimal set. \>{lem}

\proof If $\rad$ is not finite-to-one, then   $Y$ contains an infinite set $Y'$ of  balls
of the same  
radius $\alpha$.   Then $\union Y'$   contains a closed ball $b$    of valuative radius $\beta < \alpha$.  The set $S$ of open sub-balls $b'$ of $b$ of valuative radius $\beta$
forms a strongly minimal set; the map sending $y \in Y'$ to the unique $b' \in S$ containing $y$ is surjective.  \qed

The following lemma regarding  $\VF / \Oo$  will    be needed for integration with an additive character (\secref{additivecharS}).

\<{lem} \lbl{Om-st} Let $Y$ be a stable definable set,
  $Z \subset \VF \times Y$ a definable set such that for $y \in Y$, $Z(y) = \{x: (x,y) \in Z \}$ is
additively  $\Mm$ invariant.      Then 
for all but finitely many $\Oo$-cosets $C$,
 $Z \meet (C \times Y) $ is a rectangle  $C \times Y'$.  \>{lem}

\<{proof}   
For $y \in Y$, $Z(y)$ is a $\T_{y}$-definable subset
of $\VF$, hence a Boolean combination of a finite $<y>$-definable set of balls $b_1(y),\ldots,b_k(y)$.   
Let $B_i(y)$ be the smallest
closed ball containing $b_i(y)$.  According to \lemref{fmrf-r}, since the set of closed balls
occuring as $B_i(y)$ for some $y$ is stable, it is finite:

$$\{B_i(y): y \in Y \} = \{B_1,\ldots,B_l\}$$
All the $B_i$ are $\Oo$-invariant.  
Let $R$ be the set of $\Oo$-cosets $C$ that  are equal to  some $B_i$.  

 If $B_i(y)$ has valuative radius $< 0$ (i.e. it is bigger than an $\Oo$-coset),  
 then so is $b_i(y)$, so the characteristic function of 
such a $b_i(y)$ is constant on any closed $\Oo$-coset $C$.  
 If 
  $C \notin R$, then it is disjoint from any $B_i$ of valuative
 radius equal to (or greater than) $0$, so the characteristic functions of the
 corresponding $b_i(y)$ are also constant on it.  Thus with finitely many exceptional 
 $C$, any  such characteristic function is constant on $C$, and the claim follows.
                                                                                                     \>{proof}

\ssec{Generic types and orthogonality}  \lbl{generics}

Two generic types $p,q$ are said to be {\em orthogonal} if for any base $A'$, 
  if $c \models p |A'$, $d \models q |A'$,
then $p$ generates a complete type over $A(d)$, equivalently $q$ generates a complete
type over $A(c)$.   We will see that generics of different kinds are orthogonal 
(cf. \lemref{orth-gen}).     This orthogonality of types is weaker than the orthogonality of definable 
sets mentioned in the introduction, and   in the present case is only an indirect consequence
of the orthogonality between the residue field and value group; these types do 
not have  orthogonal definable neighborhoods.  

If $\g \in \G$ and $rk_{\Qq} ( \G(C(a)) / \G(C)) = n$, we say that $tp(\g/C)$ has $\G$-dimension
$n$.

\<{lem} \lbl{orthplus} Let $p_\G$ be a $\T_A$- type of elements of $\Gamma^n$
of $\G$-dimension $n$.  Let $P = \val^{\inv}(p_\G)$.  Then
\begin{enumerate}
  \item  $\val^{\inv} (p_\G)$ is a complete type over $A$.  In other words,
 for any $A$-definable set $X$, either $\val^{\inv} (p_\G) \subseteq X$ or 
 $\val^{\inv} (p_\G) \meet X  = \emptyset$.
  \item 
If $D$ is a stable 
 $A$-definable set and $d_1,\ldots,d_n \in D$, then $P$
implies a complete type over $A(d_1,\ldots,d_n)$.  
 \item  If $c \in P$ then $D(A(c)) = D(A)$.
   \item  $P$ is complete over $A$.  

\end{enumerate}

 \>{lem}

 \<{proof}   
 (1)  reduces inductively to the case $n=1$.
Since $\val^{-1}(p_\G)$ is
 a disjoint union of open balls,  (1) for $n=1$ follows from \lemref{vstar}: an $A$-definable set $X$ cannot
 intersect nontrivially each of an infinite family of open balls.  Therefore either $X$ is disjoint from almost all,or  $X$ contains almost all open
 balls $\val^{\inv} (c)$, $c \models p_\G$;
 in the former case 
  the complement of $X$ contains $\val^{\inv} (p_\G)$, 
  and in the latter $X$ contains $\val^{\inv} (p_\G)$ since $p_\G$ is complete.

 (2)    By strong orthogonality, $p_\G$ generates a complete type $q'$ over $A(d)$,
 of $\G$-dimension $n$.  

 By (1) over $A(d)$, $\val \inv (p_\G)$ is complete over $A(d)$.  But if $c \in P$ then $\val(c) \models p_\G$
 so $c \in \val \inv (p_\G)$.  Thus $P$ 
implies a complete type over $A(d)$.  

(3) follows from (2):  if $d \in D(A(c))$ then there exists a formula $\phi$ such that
$\models \phi(d,c)$ and such that $\phi(x,c)$ has a unique solution.  By (2)
$\phi$ is a consequence of $P(c) \union tp(d/A)$, and hence by compactness of 
a formula $\phi_1(x) \& \phi_2(c)$ where $\phi_2 \in tp(d/A)$.  So already $\phi_1(x)$
has the unique solution $d$, and thus $d \in D(A)$.  

(4) is immediate from (1).
                                         \>{proof}

\<{lem} \lbl{orthk}  Let $q$ be a $\T_A$-type of elements of $\RES_A^n$ of 
$\RES$- dimension $n$.  Let $Q = \rv^{\inv}(q)$.  Then $Q$ is complete over $A$.  
Moreover, if $\g_1,\ldots,\g_m \in \G$ then $Q$ implies a complete type over
$A(\g_1,\ldots,\g_m)$.
\>{lem}

\<{proof}  Again the lemma  reduces inductively to the case $n=1$, and for $n=1$
follows from  \lemref{vstar}, since $\val \inv (q)$ is a union of disjoint annuli;  the ``moreover''
also follows from orthogonality as in the  proof of \lemref{orthplus} (2). \>{proof}

\<{lem}[ \cite{acvf1} \S 2.5] \lbl{orth-gen}  

(1)   If $b$ is an open ball, or a properly infinite intersection of balls, and $b'$ a closed ball, 
  then $p_b, p_{b'}$ are orthogonal.  
  
  (2)  Any $b$-definable map to $\k$ is constant on 
  $b$ away from a proper sub-ball of $b$.
\end{lem}

\prf   We recall the proof from \cite{acvf1} \S 2.5:
The statement 
becomes stronger if the base set is enlarged.  Thus we may  assume that $b$ and $b'$ are centered;
by translating we may  assume both are centered at $0$, and by a multiplicative renormalization 
  that $b'$ is the unit closed ball.  So 
  
  (*) \ \  $c \models p_{b'}|A$ iff $c \in \Oo$
and $\res(c) \notin \acl(A)$.  

On the other hand let $p_\G $ be the type of elements
of $\G$ that are just bigger than the valuative radius of $b$ (cf. \exref{O-gen}).  Then 
  $d \models p_b | A$ iff $\val(d) \models p_\G$, i.e. $p_b$ is now the type $P$
  described in \lemref{orthplus}.  By \lemref{orthplus}, if $c' \in P$ then $\k(A(c'))=\k(A)$.
  It follows that if $c \models p_{b'} | A$ then $\res(c) \notin \acl(\k(A(c'))$.   By (*) $c \models p_{b'} | A(c')$.
  
  For the second statement, let $g$ be a definable map $b \to \k$; by   \lemref{orthplus} (3),
  $g$ is constant on the generic type of $b$; by compactness, $g$ is constant on $b$ 
  away from some proper sub-ball of $b$.
  \eprf

  \<{lem} \lbl{1gen} Let $a=(a_1,\ldots,a_n) \in \RV^n$,   and assume
  $a_i \notin \acl(A(a_1,\ldots,a_{i-1}))$ for $1 \leq i \leq n$.
    Then the formula
 $D(x) = \bigwedge_{i=1}^n \rv(x_i) = a_i$ generates a complete type over $A(a)$, and indeed over
 any $\RV \union \G $-generated structure  $A''$ over $A$.
 % can also add generators from any stable $A$-definable set. 
  
  In particular, if $q=tp(a/A)$, any $A$-definable function $f: \rv \inv (q) \to \RV \union \G$  factors through $\rv(x)=(\rv(x_1),\ldots,\rv(x_n))$.
  %can also add a stable set S to range

    \>{lem}

\<{proof}  This reduces inductively to the case $n=1$.  If we replace $A$ by a bigger set $M$
(such that $a_i \notin \acl(A(a_1,\ldots,a_{i-1}))$ for $1 \leq i \leq n$),  the assertion becomes stronger; so we may assume
$A=M \models \T$.  Let $\rv(c)=\rv(c')=a$.  Either $\val(c)=\val(c') \notin M$,
or else $\val(c)=\val(c')=\val(d)$ for some $d \in M$, and $\res(c/d)=\res(c'/d) \notin M$;
in either case, by \lemref{orthplus} or \lemref{orthk}, we have $tp(c/M)=tp(c'/M)$.
So $tp(c,\rv(c) / M) = tp(c',\rv(c')/M)$, i.e.$tp(c/M(a))=tp(c'/M(a))$.  This proves completeness
over $A(a)$.  

Let $A'$ be a structure generated over $A$ by finitely many elements of $\G$.
Then $A'(a) = A(\g_1,\ldots,\g_k,a)$, where  $\g_i \in \G$, and $\g_i \notin A(\g_1,\ldots,\g_{i-1},\val(a))$.  It follows that  $\rv(a) \notin A(\g_1,\ldots,\g_k)$, so $D(x)$ generates a complete type over $A(\g_1,\ldots,\g_k)(a)=A'(a)$.

Let $A''$ be generated over $A'(a)$ by elements of stable $A$-definable sets.  Since
$D(x)$ is the (unique, and therefore) generic type of an open ball over $A'(a)$, by
\lemref{orthplus},
it
 generates a complete type over $A''$.  

Now if $A'' = A(\g_1,\ldots,\g_k,r_1,\ldots,r_n, d)$ where $\g_j \in \G$, $r_i \in \RV$ and $d$ lies in a stable set over $A$,
let $A'=A(\g_1,\ldots,\g_k,\valr(r_1),\ldots,\valr(r_n))$; then $A'/A$ is $\G$-generated, and
$A''/A$ is generated by elements of stable sets (including $\valr^{-1}(r_i)$.)  Thus
the above applies.  

The last statement follows by applying the first part of the lemma over $A'' = A(f(c))$:  the formula
$f(x)=f(c)$ must follow from the formula $D(x)$, since $D(x)$ generates a complete type over $A''$.  \>{proof}

\ssec{Definable sets in group extensions} 
 \lbl{groupext}

 We will analyze the structure of $\RV$ in a slightly more abstract setting.   In the following lemmas we assume $R$ is a ring, and $0 \to A \to B \to C \to 0$  is a definable exact sequence of 
 $R$-modules in $T$.  This means that $A,B,C$ are definable sets, and that one is also given
 definable maps $+_A:A^2 \to A$, $f^r_A: A \to A$ for each $r \in R$, and similarly for $B,C$;
 and definable maps $\iota: A \to B, \rh: B \to C$, such that in any $M \models T$,
 $A(M),B(M),C(M)$ are $R$-modules under the corresponding functions, and 
 $0 \to A(M) \to_{\iota} B(M) \to _{\rh} C(M) \to 0$ is an exact sequence of homomorphisms of $R$-modules.  

\<{lem}  \lbl{seq} Consider a theory   with a sequence $0 \to A \to B \to_\rh C \to 0$ 
of definable $R$-modules and homomorphisms (carrying additional structure.)
Assume:\begin{enumerate}
  \item $A,C$ are stably embedded and orthogonal.
    \item Every almost definable subgroup of $A^n$ is defined by finitely many $R$-linear equations. 
   \item (``No definable quasi-sections''.)  If $P$ is a definable subset of $B^n$
  whose projection to $C^n$ is finite-to-one, then $P$ is finite.  
   \end{enumerate}
Then every  almost
definable subset $Z$ of $B^n$ is a finite union of sets of  the form
$$\{b: \rh(b) \in W, Nb \in Y \}$$
 where $N \in B_{n,k}(R)$ is an $n \times k$ matrix, 
 $Y$ is an almost definable subset of a single coset of  $A^k$,   $W$ is an almost
 definable subset of $C^n$.
\>{lem}

Note:   \begin{enumerate}
  \item  To verify (3), it suffices to check it for $n=1$ but for parametrically definable $P$.
  \item If $C$ is definably linearly ordered, and $Z$ is definable, then $Y,W$ may be taken definable.
\end{enumerate}
\<{proof}    Using a base change as in \secref{naming}, we may assume almost definable sets are definable.
Replacing $B$ by $B^n$ and $R$ by $M_n(R)$, we may assume $n=1$.
Let $Z$ be a definable
subset of $B$.  Given $X \subset A$, let $[X]$ denote the class of $X$ up to translation;
so $[X]=[X']$ if $X=X'+a$ for some $a \in A$.  Now a definable subset $U$ of a coset 
$b+A$ of $A$ has the form $b+X$, $X \subset A$;   the class $[X]$ is well-defined, 
and we will denote $[U]=[X]$.  
We obtain a map 
$$ c \mapsto [Z \meet \rh^{-1}(c)]$$
In more detail:  for any $b \in  (\rh \inv (c) \meet Z)$,
we have $(\rh \inv (c) \meet Z) - b \subseteq A$, and so by stable embeddedness
of $A$ we can write $(\rh \inv (c) \meet Z) - b = X(a)$ for some $a \in A^m$.
The tuple $a$ is not well-defined; but the class of $a$ in the definable equivalence
relation:

$$ x \sim x'  \iff  (\exists t \in A) (t+X(x))=X(x') $$

is obviously a function of $c$ alone.  
By the orthogonality assumption, this map is piecewise constant.  So we may assume
it is constant, and fix $C_0$ with $[Z \meet \rh^{-1}(c)] = [C_0]$.  Let $S$ be the stabilizer
$S= \{a \in A: a + C_0 = C_0 \}$.  Then for $a \in S$, $a+ (Z \meet \rh^{-1}(c)) =
 (Z \meet \rh^{-1}(c)) \}$ for any $c \in C$, so that also $S= \{a \in A: a + Z = Z \}$,
 and $S$ is definable.

Now  $ Z \meet \rh^{-1}(c)  =  C_0 + f(c)$ for   some $f(c) \in \rh^{-1}(c)$;
$f(c)+S$ is well-defined.
 
 By assumption (2), $S = Ker(r_1) \meet \ldots \meet Ker(r_m)$ for some $r_i \in R$.
 Let $I = \{r_1,\ldots,r_m  \}$.   For $r \in I$, 
 $f_r(c):=rf(c)$ is a well-defined element of $B$, and for all $c \in \rh(Z)$,
  $r(Z \meet \rh^{-1}(c)) = rC_0 + f_r(c)$.  
  
 We have $\rh f_r(c) = r c $.  If $d \in Ker(r: C \to C)$, then $f_r(d+c) = rc$ also,
 so $f_r(d+c) - f_r(c) \in A$.  By orthogonality, for fixed $r$,
 $f_r(d+c)-f_r(c)$ takes finitely many values as $c,d$ vary in $C$.  
 In other words,  $\{rf(c): c \in \rh(Z)$  is a quasi-section above $r \rh(Z)$. 
  By   (3), $r \rh(Z)$ is finite, for each $r \in I$.     Let $N=(r_1,\ldots,r_m)$,
  $Y'=NZ$.   Then $\rh(Y') $ is finite.  It follows that $Y'$ is contained in a finite union 
  of cosets of $A$, so $C,Y'$ are orthogonal.
  
Thus  $\{(\rh(z),Nz): z \in Z \}$ is a finite union of rectangles; upon dividing $Z$
further, we may assume this set is a rectangle $W \times Y$.  Now if $\rh(b) \in W$
and $Nb \in Y$ then for some $z \in Z$, $\rh(b)=\rh(z)$ and $Nb=Nz$;
it follows that $b-z \in A$ and   $b-z \in S$; so $b \in S+Z = Z$.  
     Thus $Z$ is of the required form.  
   \>{proof}

 \<{cor} \lbl{rvexpansions}
  $T $ be  a complete theory in a language $L$ satisfying   the assumptions of \lemref{seq}.
  Let  $L \subseteq L', T \subseteq T'$, and assume (1)-(3) persist to $T'$.  If $T,T'$
  induce the same structure on $A$ and on $C$, 
  up to constants they induce
  the same structure on $B$, i.e. every $T$-definable subset of $B^*$ is parameterically $T'$-definable.  \>{cor}
 \prf  Apply \lemref{seq} to $T'$, and note that every definable set in the normal form
 obtained there is already parametrically definable in $T$. \eprf

We will explicitly  use imaginaries in $\RV$ only rarely; but our ability to work with $\RV$, using $\G$ as an auxiliary, 
is partly explained by:

\<{cor} \lbl{ei-extns} 
Let  $0 \to A \to B \to_\rh C \to 0$ be as in   \lemref{seq}, and assume  $C$ carries a definable linear ordering.
{  
  Let ${\bar V}$ be the disjoint union of the definable cosets of $A$
in $B$, with structure induced from $T$.
Let $e$ be an   imaginary element of $B$.  Then $<e> = <(a',c')>$ for some pair $(a',c')$ consisting of an imaginary of ${\bar V}$ and an imaginary of $C$.  Thus if  ${\bar V} ,C$ eliminate imaginaries, so does $B \union C \union {\bar V}$.  }   \>{cor}

\proof {    Let $e$ be an imaginary element of $B$; let $E_0$
be the set of $A, {\bar V} $-imaginaries   that are algebraic
over $e$.   
 
 By \lemref{seq}, applied to a definable set with code $e$ in the theory $T_{E_0}$,
  there exist almost definable subsets of ${\bar V}, C^n$ from which $e$ can be defined.  These are coded by imaginaries permitted in the definition of $E_0$.
   Thus $e$ is $E_0$-definable.
     So $e=g(d)$ for some definable function $g$ and some tuple $d$ from $E_0$.
}  
Let us now show that     $e$ is equi-definable with a finite set, i.e. an imaginary
of the form $(f_1,\ldots,f_n) /Sym(n)$.  
   Let $W$ be the set of elements with the same type as $d$ over $e$;
   $W$ is finite by definition of $E_0$, and is $e$-definable.  But
     $e=g(w)$ for any element $w \in W$, so $e$ is definable from $\{W\}$.  
     
It remains to see that every finite set of elements of $E_0$ is coded by imaginaries of $A$ and $C$ and elements of
$B$.  Since  $C$ is linearly ordered, it suffices
to consider finite sets whose image in $C^m$ consists of one point.  These are subsets
of some definable coset of $A^m$, so again by elimination of imaginaries there they
are coded.  \qed 

\<{cor}\lbl{rvei} The structure induced on $\RV \union \G$ from $\ACVF$ eliminates imaginaries.   \>{cor}

\proof $\G_{E_0}$ eliminates imaginaries, and so does ACF (cf. \cite{poizat}).
Note that $\bar V$ is essentially a family of 1-dimensional $\k$-vector spaces,
closed under tensor products and roots and duals.  
Hence by \cite{EI},  ${\bar V}_{E_0}$ eliminates imaginaries too.
Our only application of this
lemma will be in a situation  when parameters can be freely added; in this case,
it suffices to quote elimination of imaginaries in ACF.
 \qed

\<{cor}\lbl{seq-c}  Let $T$ be a theory as in \lemref{seq},  with $R= \Zz$, and $C$ a linearly ordered group.  
  Then  every definable subset of $B^{n}$  is a disjoint union of $GL_n(\Zz)$-images of products $Y \times \rh^{-1}(Z)$, with   $|\rh Y|=1$.  In particular the Grothendieck semiring $\SG(B)$
(with respect to the category of all definable sets and functions of $B$) 
   is generated by the classes of elements $Y \subset B^n$ with $|\rh Y|=1$, and pullbacks $\rh^{-1}(Z)$, $Z \subset C^m$.
 \>{cor}
 
\proof  By   \lemref{seq}, the Grothendieck ring is generated by classes of sets
$X$ of the form $X = \{b \in B^n: \rh(b) \in W, Nb \in Y \}$.    
After performing row and column operations on the matrix $N$, we may assume it is the composition of a projection $p:R^n \to R^k$
with a diagonal $k \times k$ integer matrix with nonzero determinant.  The composition
$\rh p (X)$ is finite; since $C$ is ordered, each element of $\rh p(X)$ is definable,
and so we may assume $\rh p (X)$ has one element $e$.  Thus
$W = \{(e) \times W' \}$ for some $W'$, and $X = pX \times \rh^{-1}(W')$.   \qed

\<{lem} \lbl{seq-d} Let $T$ be a theory, and let
 $0 \to A \to B \to_\rh C \to 0$ 
 be an exact sequence o f
of definable Abelian groups and homomorphisms.  If $E \leq M \models T$,
we will write $E_A = A(E)$, etc.    
Assume: 
  \begin{enumerate}
  \item %1
  $A,C$ are orthogonal.  
  \item  %2
Any parametrically definable subset of $B$ is a Boolean combination 
of sets $Y$ with $\rh(Y)$ finite, and of full pullbacks $\rh^{-1}(Z)$.
  \item %3
   $C$ a uniquely divisible Abelian group, and for any $E \leq M \models T$, every divisible
  subgroup containing $E_C$ is algebraically closed in $C$ over $E$.
  \item   %4
   For any prime $p>0$,  $T \models (\exists x \in A) (px =0, x \neq 0)$
\end{enumerate}
     Let $Z \subset C^n$ and  $f: Z \to C$ 
be definable, and suppose there exists $E$ and   $E$-definable 
     $X \subset B^n$ and
$F: X \to B$  lifting $f$: $\rh X = Z, \rh F(x) = f(\rh x)$.     Then there exists
a partition of $Z$ to finitely many definable sets $Z_\nu$,  
such that for each $\nu$, for some $m \in \Zz^n$,
$f(x) - \sum_{i=1}^n m_ix_i$  is constant on $Z_\nu$.
  \>{lem}

The main point is the integrality of the coefficients $m_i$.

\<{proof}

It suffices to show that for any $M \models T$ and any
$c=(c_1,\ldots,c_n) \in Z(M)$, there exists  $m = (m_1,\ldots,m_n) \in \Zz$ such that
$f(c) - mc \in E^0$,
where $E^0 = \dcl(\emptyset)$
is the smallest   substructure of $M$.  
  For if so, there exists a formula of one variable of sort
$C$, such that $T \models (\exists ^{\leq 1} z) \psi(z)$,
$M \models \psi(f(c)-mc)$.  
By compactness 
there exists a finite set $F$ of such pairs $\nu=(m,\psi)$,
 such that for any $M \models T$ and $c \in Z(M)$, for some $(m,\psi) \in F$, 
$M \models \psi( f(c) - mc )$;   the required partition
is given by $X_{m,\psi} = \{z \in Z :  \psi( f(z) - mz ) \}$.

Fix $M$ and $c \in Z(M)$.    Let   $<c>$ be the
smallest divisible subgroup of $C(M)$ containing $E^0_C$
and 
$c_1,\ldots,c_n$.  By (3), $<c>$ is closed under $f$,
so $f(c) \in <c>$, i.e. $f(c)= \sum \alpha_i c_i + d$ for some $\alpha_i \in \Qq$ and some $d \in E^0_C$.  The only problem is to show that
we can take $\alpha_i \in \Zz$.

We will use induction on $n$.
Let $K = \{\b \in \Qq^n: \b \cdot c  \in E^0_C \}$.   $K$ is a $\Qq$-subspace of $\Qq^n$.  If  $K \neq (0)$, there exists a primitive
integral vector $\b_1 \in K$.  $\b_1$ may be completed to
a basis for a $\Zz$-lattice in $\Qq^n$.  Applying a $GL_n(\Zz)$
change of variables to $B^n$, we may assume $\b_1 = (1,0,\ldots,0)$,
i.e. $c_1 \in E^0_C$.  But then let   $f'(z_2,\ldots,z_n)=f(c_1,z_2,\ldots,z_n)$.
Then $f'$ lifts to a definable function on $B^n$ (with parameters,
of the form $F(b_1,y_2,\ldots,y_n)$) so by induction,
$f(c_1,\ldots,c_n) = f'(c_2,\ldots,c_n) = \sum_{i \geq 2} m_i z_i + d'$ for some  $m_2,\ldots,m_n \in \Zz$ and $d' \in  E^0_C$, as required.

Thus we can assume $K = (0)$.   
%
%Consider first the case:  $c=c_1$, $p f(c)=c$, $c \notin \acl(E)$.  Increasing $E$, we may
%assume $f$ lifts to an $E$-definable function $F$.  (We can still preserve $c \notin \acl(E)$,
%say by picking a different realization of the same type.)  So $p F(b) -b \in A$ for
%$b \in \rh \inv (c)$.  Now $b \mapsto pF(b) - b$ is an affine map from 

We can find $m, m_i \in \Zz, e \in \dcl(\emp)$ with 
$$m f(c) = \sum m_i c_i + e$$
   If $m | m_i$ we are done.   
We will now derive a   contradiction from the contrary assumption that 
$m$ does not divide each $m_i$ in such an equation, with $f$ a liftable  function.  We may assume
that the greatest common divisor of $m,m_1,\ldots,m_n$; so there exists a prime dividing $m$ but not (say) $m_1$.

Let $g(x) = f(x,c_2,\ldots,c_n)- e/m - \sum_{i=2}^n m_ic_i/m$; then $mg(c_1) = m_1 c_1$, $m$ does not divide $m_1$, $g$ is $E=\acl(c_2,\ldots,c_n)$-definable and liftable.   Since $K = (0)$, by assumption (3),
$c_1 \notin \acl(E)$.   Let $E' \supset E$ be such that 
$g$  lifts to an $E'$-definable function $G'$.  Enlarging the model if necessary, 
let $c_1'$ realize $tp(c_1/E)$, with $c_1' \notin E'$ (cf. \exref{alldepend}).  
  Therefore there exists $E''$
such that $E'',c_1$ and $E',c_1'$ have the same type.  In particular  
$g$     lifts to an $E''$-definable function $G$.   

  Consider any $b_1$ such that  $\rh(b_1)=c_1$.
Then   $m \rh G(b_1) - m_1 \rh(b_1) =0$.  
So $m G(b_1) - m_1 b_1 \in A$.   % in fact $m G(y)  - m_1y \in A$ for any $y$ with $\rh(y)=c_1$.  

Let $p$ be prime, 
$p | m$ but $p \not | m_1$.  
 Let $s,r \in \Zz$ be such that $sp - r m_1=1$, and let
$h(x)=sx - {\frac{rm}{p}} g(x)$.  Then
$ph(c_1) = ps c_1 -rm g(c_1) = ps c_1 -rm_1c_1 = c_1$.  
Also $h$ is liftable over $E''$:  indeed if $G$ is $E''$-definable and lifts $g$,
then  $H(x)=sx -{\frac{rm}{p}}G(x)$ lifts $h$.  

%Let $C_1 = \rh \inv (c_1)$,  a coset of $A$.  
%For $y \in C_1$, $pH(y)-y \in A$. 
% Pick $b_1 $ with $\rh (b_1)=c_1$. 
So $p H(b_1) = b_1 + d$, some   $d \in A$.   Let $b_2 = H(b_1)$;
then $b_1 = pb_2 -d$, or
$$ b_2 = H(pb_2-d)$$

Now let $c_2  =  h(c_1) = \rh(b_2)$.  Then $pc_2=c_1$, and so $c_2  \notin \acl(E'')$,
since by unique divisibility $c_1 \in \acl(h(c_1))$.  By (1), $c_2 \notin \acl(E''(d))$.
 Let $C_2 = \rh \inv c_2$.  By (2), any $E''(d)$-definable set either contains 
  $C_2$ or is disjoint from $C_2$.   Hence for any $y \in C_2$,  $y=H(py-d)$.  

By (4) there exists $0 \neq \omega_p \in A$ with $p \omega_p = 0$.
Let $b_2'=b_2 + \omega_p$.  Then $b_2 \in C_2$, so $b_2' = H(pb_2'-d)$.
But $pb_2' = pb_2$, so $b_2=b_2'$ and $\omega_p=0$, a contradiction.  
  
 \eprf

\<{remark} \begin{enumerate}
  \item It follows from \lemref{seq-d}  that a definable bijection between subsets of $C^n$ that lifts to 
  subsets of $B^n$ is piecewise given by an element of $GL_n(\Zz) \ltimes C^n$.
  (cf. \lemref{rvfn}.)
  \item The assumption (4) on torsion does not hold in characteristic $p>0$
for the sequence $\k^* \to \RV \to \G$.  In this case there is $l$-torsion for $l \neq p$,
but no $p$-torsion, and the corresponding group is $GL_n(\Zz[1/p]) \ltimes C^n$.
\end{enumerate}
   \>{remark}

Note as a corollary that there can be no definable sections 
of $B \to C$ over an infinite definable subset of $C$.

\<{lem} \lbl{rvfn} Let $0 \to A \to B \to C \to 0$ be as in \lemref{seq-d}.
Let $X \subset B^n$ be definable, and let
 $f: X \to B^l$ be a  definable function.   $X$ may be partitioned
 into finitely many pieces $X'$, such that on each $X'$
 \begin{enumerate}
  \item  $f(x)= Mx + b(x)$, where $M$ is a $l \times n$-integer matrix
  and $\rh b(x)$ is constant.
  \item    There exists $g \in GL_n(\Zz)$ such that 
  $b \circ g$ factors through a   projection
  $B^n \to_\pi B^k$, where $\rh \pi(X') $ is one point of $C^k$.
\end{enumerate}
\>{lem}

\proof  We first prove (1,2) for complete types.  

 (1)  This reduces to $l=1$.  Let $P$ be a complete type
of elements of $X$.  
Then on $P$ we have $\rh \circ f(x) =\sum m_i \rh(x_i) + d$ for some constant
$d$.  (\lemref{seq-d}).

Thus $f(x)= \sum m_ix_i + b(x)$, where 
$b(x)= f(x) - \sum m_ix_i$, and $\rh b (x) = d$ is constant. 

(2).  %Suppose now that $\rh(b(x))$ is constant on $P$.
  Let 
$\pi: B^n \to B^k$ be a projection such that $\rh \pi(X) $ is one point of $C^k$, and with $k$ maximal.  So $P \subset P'  \times P''$,
$P' \subset B^{n-k}, P'' \subset B^k$, 
and $\rh(P'')$ is a single point of $C^k$, while $\rh(P')$
is not contained in any proper hypersurface $\sum n_i x_i ={\rm constant}$ with $n_i \in \Zz$.     Pick $b'' \in P''$.
Let $\g=(\g_1,\ldots,\g_k) \in \rh(P')$,
$\g$ not in any such hypersurface.   Let $a=(a_1,\ldots,a_k)$,
$\rh(a_i)=\g_i$, and let $a'$ be another point with $\rh(a')=\g$.
Let $e = f(a,b)$. 
Then $tp(a/b,e)=tp(a'/b,e)$, so $f(a',b)=e$.  Thus $f(a,b)$
depends only on $b \in P''$ and not on $a$ (with $(a,b) \in P$.)

 Since (1),(2) hold on each complete type, there exists
a definable partition such that they hold on each piece.   \qed

\ssec{$\V$-minimality} \lbl{V+}

We assume from now on that $\T$ is a theory of $C$-minimal valued fields, of residue
characteristic $0$.     When using the many-sorted language, we will still say that $\T$ is {\em a theory of valued fields}
when $\T = Th(F,\RV(F))$ for some valued field $F$, possibly  with additional structure.
 A $C$-minimal $\T$ satisfying assumption \eqref{V+c} below will be said to have  {\em centered closed balls}.  If in addition \eqref{V+RV}, \eqref{V+comp}  hold, we will say $\T$ is {\em  $\V$-minimal.}
%Many lemmas use only a subset of the assumptions, but for readability's sake we will
%not usually indicate these. 
%for instance the comparison of $\VF$ and $\RV$ derivatives \propref{rv-compat}  works 
%for $C$-minimal fields satisfying   \eqref{V+RV}.   
Expansions by definition of the language,
i.e. the addition of a relation symbol $R(x)$ to the language along with a definition
$(\forall x)(R(x) \iff \phi(x))$ to the theory, do not change any of our assumptions. 
Thus  we can    assume that $\T$ eliminates quantifiers.

\<{enumerate}

\item{\bf Induced structure on $\RV$ }  \lbl{V+RV}
$\T$ contains $\ACVF(0,0)$,  and every parametically $\T$-definable relation on $\RV^*$ is parametrically definable in $\ACVF(0,0)$.

\item{\bf Definable completeness}  \lbl{V+comp}  Let $A \leq M \models T$, and let $W \subset \fB$ be a 
$\T_A$-
definable family of closed balls linearly ordered by inclusion.
Then $\meet W  \neq \emptyset$.

\item{\bf Choosing points in closed balls}    Let $M \models \T$, $A
 \subseteq \VF(M)$, and
let $b$ be an almost $A$-definable closed ball. 
 Then 
 $b$ contains an  almost $A$-definable point.
 \lbl{V+c}

 \>{enumerate}

   $\T$ will be called {\em effective}  if every definable finite disjoint union of balls contains a definable set, with exactly one point  in each.   
   A substructure $A$ of a model of $\T$ will be called effective if $\T_A$ is effective. 

%By an $\rv$-ball, we mean one of the form $\rv \inv (a)$, $a \in \RV$. 
 If every 
    definable finite disjoint union of $\rv$-balls contains a definable set, with exactly one point  in each,
    we can call $\T$ $\rv$-effective.   However:
    \<{lem}\lbl{rveffective}  Let $\T$ be $\V$-minimal.  Then $\T$ is effective iff it is $\rv$-effective.  \>{lem}

\prf  Assume $\T$ is $\rv$-effective.  Let $b$ be an algebraic ball.  If $b$ is closed, it has
an algebraic point by \aref{V+} \ref{V+c}.  If $b$ is open, let $\bar{b}$ be the closed ball
surrounding it.  Then $\bar{b}$ has an algebraic point $a$.  Let $f(x)=x-a$.  Then $f(b)$ is an $\rv$-ball, so by $\rv$-effectivity it has an algebraic point $a'$.
Hence $a'+a$ is an algebraic point of $b$.
\eprf

%    
 
% surjectivity of the map $\VF^* \to \RV$ on definable points; 
% cf. \propref{eff}. 
  In general, 
effectivity is needed for lifting morphisms from $\RV$ to $\VF$,   not for the 
  ``integration'' direction.

   If $\T$ is $\V$-minimal
 and $A$ is a $\VF \union \RV \union \G$ - generated structure, we will see that
 $\T_A$ is $\V$-minimal too.  The analog for points in open balls is 
 true  but only
 for $\VF \union \G$-generated substructures; for thin annuli it is true only for $\VF$-generated
 structures.    For this reason the condition on closed balls   is
 more flexible; luckily we will be able to avoid the others.
 
\lemm{Vplus1}  Let $\T$ be a C-minimal theory of valued fields.

 Then $(1) \implies (2) \implies (3) \implies (4)$

  \begin{enumerate}
  \item  $\T$ admits quantifier-elimination in a three-sorted language $(\VF,\k,\G)$, such that
  for any basic function symbol $F$ with range $\VF$, the domain is a power of $\VF$; and 
  no relations on $\k,\G$ beyond the  field structure on $\k$ and the ordered
  Abelian group structure on $\G$.
  \item  Every parametrically definable relation on $\k$ is parameterically definable
  in ACF(0); and every parametrically definable relation on $\G$ is parameterically definable
  in DOAG.   
  \item   \aref{V+}  \eqref{V+RV}. 
  \item $\k,\G$ and $\RV$ are stably embedded.  
\end{enumerate}

\>{lem}

\prf $(1) \implies (2):$
Let $\phi(a,x)$ be an 
 atomic formula  with paramaters $a=(a_1,\ldots,a_n)$ from $\VF$ and $x=(x_1,\ldots,x_m)$ variables for the $\k,\G$ sorts.  Then $\phi$  must have the form $\psi(t(a),x)$, 
where $t$ is a term (composition of function symbols) $\VF^* \to (\k \union \G)$. So
$\phi(a,x)$ defines the same set as $\psi(b,x)$ where $b = t(a)$.  Since every formula
is a Boolean combination of atomic ones,  (2) follows.      

$(2) \implies (3):$  
This follows from \corref{rvexpansions}.   The assumptions of \lemref{seq-d} are satisfied: (1) is automatic
since by C-minimality $\k$ is strongly minimal and $\G$ is $O$-minimal; (2) follows from C-minimality; (3),(4) follow from the assumptions on $\k,\G$.  

$(3)$ immediately  implies (4).

 \eprf

\lemm{average}  Let $\T$ be a theory of valued fields satisfying \aref{V+} \eqref{V+RV},
such that $\res$ induces a surjective map on algebraic points.
  Then $(1) \implies (2) \implies (3) \implies (4)$:\begin{enumerate}
  \item   for any $\VF$-generated substructure $A$ of
a model $M$ of $\T$, if $\G(A) \neq (0)$, then $\acl(A) \models \T$.  
  
  \item for any $\VF$-generated substructure $A$ of
a model of $\T$, any $\T_A$-definable nonempty finite union of balls  contains a nonempty  $\T_A$-definable finite set.

  \item \aref{V+} \eqref{V+c} holds.
  
  \item Let $A$ be $\VF$-generated, and $Y$  a finite $A$-definable set of disjoint closed balls.  Then there exists an $A$-definable
finite set $Z$ such that $|b \meet Z|=1$ for each $b \in Y$. 
\end{enumerate}
\>{lem}

\prf  We first show:

\Claim{}      For any $\VF$-generated $A$ with $\G(A)=(0)$,  
$\res: \VF(\acl(A)) \to \k(\acl(A))$  is surjective.

\prf    It suffices to prove the claim
for finitely generated $A$.
For $A = \emptyset$ this is true by assumption.  Using induction on the number of generators,
it suffices to show that if the Claim holds for $A_0$ and $c \in \VF$  then it holds for
 $A=A_0(c))$.
%Let $A'=\acl(A_0)$.  
 
 Since $\G(A)=(0)$,
   $\res$ is defined and injective on $\VF(A)$.     If $c \in \acl(A_0)$ there is nothing to prove.  Otherwise,
 by injectivity, $\res(c) \notin \acl(A_0)$.   As a consequence of  \aref{V+} \eqref{V+RV},
 both $\dcl$ and $\acl$ agree with the corresponding field-theoretic notions
  on $\RV$, and in particular 
 on the residue field.
  
 By \lemref{1gen},
  $$\k(A_0(c)) \subseteq \dcl(\RV(A_0),\rv(c)) = \dcl(\k(A_0),\res(c))  = \k(A_0)(\res(c))$$
 
Now if $d \in \k(\acl(A))$ then $d \in \k$ and $d \in \acl(A)$, so by stable embeddedness of $\k$ we 
have $d \in \acl(\k(A))$; but $\acl(\k(A)) = \k(A)^{alg}$ by Assumptions \ref{V+} and \ref{V+RV};
so $d \in \k(A_0)(\res(c))^{alg} \subseteq \res(A_0(c)^{alg})$.    \eprf

Assume (1).  If $\G(\acl(A)) \neq (0)$ then by (1)
$\acl(A) \models T$, and in particular every $\acl(A)$-definable 
ball has a point in $\acl(A)$, so (2) holds.  Assume therefore that $\G(\acl(A))=0$.      
Let $b$ be an $\acl(A)$-definable ball.  Then $b$ must have valuative radius $0$.
If some element of $b$ has valuation $\g<0$ then all do, and $\g \in A$, a contradiction.
So $b$ is the (open or closed) ball of radius $0$ around some $c \in \Oo$.  If $b$ is closed,
then $b = \Oo$ and $0 \in b$.  If $b$ is open, then $b= \res \inv (b')$ for some element
$b'$ of the residue field $\k$; in this case $b$ has an $\acl(A)$-definable   point by the Claim.  

  (3) is included in  (2), being the case of closed balls.
  
 Assume (3).  In expansions of ACVF(0,0), the  {\em average} of a finite subset of a ball 
remains within the ball.  Thus if  $Y$ is a finite   $A$-definable set of disjoint balls, by (3), 
  there exists a finite $A$-definable set $Z_0$ including
a representative of each ball in $Y$.  Let $Z = \{\av(b \meet Z_0): b \in Y\}$, where   $\av(u)$ denotes the average of a finite set $u$. 
 
\eprf

\lemm{complete} 
 When $\T$ is a complete theory, definable completeness  is true as soon as  $T$ has a single {\em spherically complete} model $M$ in the sense of Ribenboim and Kaplansky: every intersection of nested closed balls is nonempty.  \>{lem}
 
 \prf Clear. \eprf

Let $ACVF^{an}$ denote any of the rigid analytic theories of \cite{lipshitz}.    For definiteness, let us
assume the power series have coefficients in $\Cc((X))$.  See \cite{clr} for variants living
over $\Zz_p$.

\<{lem}  \lbl{3.32}

 $ACVF(0,0)$ is $\V$-minimal and effective.  So is  $ACVF^{an}$. \>{lem}

\prf   $C$-minimality is proved in \cite{lipshitz-robinson}.   \lemref{Vplus1}(1)  for ACVF is a version
of Robinson's quantifier elimination; cf. \cite{acvf1}.  

$ACVF^{an}$ admits quantifier elimination in the sorts $(\VF,\G)$ by \cite{lipshitz} Theorem 3.8.2.  The  residue field sort is not explicit in this language, but 
one can argue
as follows.  %Assume for example the power series have coefficients in $\Cc((X))$. 
Let $\k_1$ be a large algebraically closed field containing $\Cc$, and let $K = \union_{n \geq 1} \k_1((X^{1/n}))$ be the Puiseux series ring. 
Then $K$ admits a natural expansion to a model of the theory.  $K$ is not saturated,
but by C-minimality the induced structure on the residue field is strongly minimal, so $\k_1$ is 
saturated.  Now any automorphism of $\k_1$ as a field extends to an automorphism
of $K$ as a rigid analytic structure.  Thus every $K$ - definable relation on $\k_1$ is algebraic.
(This could be repeated over a larger value group if necessary.)  \lemref{Vplus1} (2)
thus holds in both cases, hence \aref{V+} \eqref{V+RV}.

Condition   \lemref{average} (1) is obviously  true for $\ACVF$.  For   $ACVF^{an}$ it is proved   in \cite{lipshitz-robinson}.)   It is also evident that these theories have a spherically complete model.  
Thus by \lemref{average} and \lemref{complete}, \aref{V+} \eqref{V+c} and \eqref{V+comp}
hold too.

\eprf

{\bf Remarks} \begin{enumerate}
  \item  \lemref{average} (1,2,3) remain true  for $\ACVF$ in positive residue characteristic, but (4) fails.

  \item 
 ACVF(0,0) also admits quantifier-elimination in the two sorted language
with sorts $\VF,\RV$; so  \aref{V+}  \eqref{V+RV} can also be proved directly,
without going through $\k,\G$ as in \lemref{Vplus1}.
 
  \item
\eqref{V+RV} is  needed  for lifting definable bijections of $\RV$ to $\VF$,
\propref{rvlift}, \lemref{rvlift-m}.  Specifically it implies the truth of assumptions (2) of    \lemref{seq} and (4) of \lemref{seq-d}.  
These lemmas are only needed for the injectiveness
of the Euler charactersitic and integration maps, not for their construction and main properties.
It is also needed for the theory of differentiation and for 
comparing derivations in $\VF$ and $\RV$; indeed even for
posing the question, since in general there is no notion of differentiation on $\RV$.
 The theory of differentiation itself  is   needed neither for the Euler characteristic nor for integration of  definable sets with a $\G$-volume form.  
 They are required only for the finer theory introduced here of 
  integration of $\RV$-volume forms. 
  
  \item  
   We know no examples of $C$-minimal fields where \eqref{V+comp} fails.   

\item    Beyond effectivity of $\dcl(\emptyset)$, \eqref{V+c},
imposes a condition
on liftability of definable functions from $\VF$ to $\fB^{cl}$.   
Let $\T_1$ be the theory, intermediate between
$\ACVF(0,0)$ and a Lipshitz rigid analytic expansion,
generated over $\ACVF(0,0)$ by the relation 
$$\val(f(t_0x)-y) \geq \val(t_1)$$ on $\Oo^2$ 
where $t_0,t_1$ are constants with $\val(t_1) >> \val(t_0)>0$ and $f$ is an analytic function.
   It appears that balls do not necessarily remain pointed upon adding $\VF$-points to $\T_1$;
   so \eqref{V+c} is not redundant.  
\end{enumerate}

\ssec{Definable completeness and functions on the value group}

We assume  $\T$ is $C$-minimal and definably complete.  We show that the property of having
   centered closed balls  is preserved under passage to $\T_A$ if $A$ is $\RV,\G,\VF$-generated;
similarly for open balls if $A$ is $\G,\VF$-generated.     Also included is a lemma
stating that every image of an $\RV$-set in $\VF$ must be finite; 
 from the point of view of content this belongs to the   description of the "basic geography",
 but we need the lemmas on functions from $\Gamma$ first.

\<{prop}  \lbl{V+11}   Let $M \models \T$, 
    $\g=(\g_1,\ldots,\g_m)$
  a tuple of elements of $\G(M)$.  Any almost $A(\g)$-definable
ball $b$ contains an almost $A$-definable ball $b'$. \>{prop}
\prf
      See  \cite{acvf1} Prop. 2.4.4.  While the Proposition is stated
     for ACVF there, the proof uses only $C$-minimality and definable completeness.   We review the proof 
 in the case that $b \in A(\g)$, i.e. 
 $b=f(\g)$ for some definable function $f$ with domain $D \subseteq \G^M$.

 Let $P=tp(\g /A)$.  
 Let $r(\g)$ be the valuative radius of $f(\g)$. 
  By O-minimality, $r$ is piecewise monotone;
 since $P$ is a complete type, $r$ is monotone, say decreasing.   
 For $a\in P$   let   $P_a = \{b \in P: b<a \}$, and for $b \in P_a$ let $f_a(b)$ be the open ball of size $r(a)$ containing  $f(b)$.  
    By \lemref{fmrf-rem2}, the valuative radius map $\rad$ is finite-to-one on  $f_a(P_a)$; but by definition it is constant,
 so $f_a(P_a)$ is finite.  Using the linear ordering, $f_a(P_a)$ is constant on 
 each complete type over $a$.    Pick $b_1 \in P$, $\e \in \G$ with $\e >0$ but
 very small (over $A(b_1)$), and $\e' \in \G$ with $\e'>0$ but $\e'$ very small
 (over $A(b_1,\e)$).  Let $b_2 = b_1 + \e$, $a=b_2 + \e'$.  Then $tp(b_1,a/A)=tp(b_2,a/A)$,
 so $f_a(b_1)=f_a(b_2)$.  Now if $f(b_1),f(b_2)$ are disjoint, let $\delta = \val(x_1-x_2)$
 for (some or any) $x_i \in f(b_i)$.  Then $r(b_2) > \delta$.  Since $\e'$ is very small, 
 $r(a) > \delta$ also.  So $f_a(b_1),f_a(b_2)$ are distinct, a contradiction.  Thus
 $f(b_1) \subset f(b_2)$.    Since $tp(a/A)=tp(b_2/A)$, we have $f(y) \subset f(a)$ for some $y \in P_a$.
   If $f(y) \subset f(a)$ for all $y \in P_{a}$, we are done;
 otherwise let $c(a)$ be the unique smallest element such that $f$ is monotone on 
 $(c(a),a)$.  We saw however that $f$ is monotone on $(d,c(a))$ for some $d< c(a)$,
 hence also on $(d,a)$, a contradiction.  So $f$ is monotone with respect to inclusion.  By compactness,
 this is true on some $A$-definable interval 
 hence on some interval
 $I$ containing $P$.  
 
 Let $U = \meet _{a \in I} f(a)$.  By definable completeness  \eqref{V+comp}
 $U \neq \emptyset$.   
Clearly $U$ is a ball, and $U \subseteq b$.  
   \eprf

 \<{lem}  \lbl{fnsg}    Let $M \models \T$, 
    $\g=(\g_1,\ldots,\g_m)$
  a tuple of elements of $\G(M)$.  Then any $A(\g)$-definable
ball contains an $A$-definable ball.   If $Y$ is a finite $A(\g)$-definable
set of disjoint balls, then there exists a finite $A$-definable set $Y'$
of balls, such that each ball of $Y$ contains a unique ball of $Y'$. \>{lem}

\<{proof}  This reduces immediately to $m=1$.  For $m=1$, by \propref{V+11},
any almost $A(\g)$-definable
ball $b$ contains an almost $A$-definable ball $b'$. 
Thus given a finite $A(\g)$-definable set $Y$ of disjoint balls,
there exists a finite $A$-definable set $Z$ of balls, such that
any ball of $Y$ contains a ball of $Z$.  Given $b \in Y$, let
$b'$  be the smallest ball containing every sub-ball $c$ of 
$b$ with $c \in Z$.  Then  $Y' = \{b': b \in Y\}$ is $A(\g)$-definable,
finite, almost $A$-definable, and (since $b_1'$ is disjoint from $b_2$
if $b_1 \neq b_2 \in Y$) each ball  of $Y$ contains a unique ball of $Y'$.
Using elimination of imaginaries in $\G$,  by \exref{def-lin},
being $A(\g)$-definable
and almost $A$-definable, $Y'$ is $A$-definable.   
                                                                                                     \>{proof}

The following corollary of 
 \lemref{fnsg} concerning definable functions from $\Gamma$ will be
    important for the theory of integration with an additive character in \secref{additivecharS}.

 \<{cor} \lbl{fns-from-gamma} Let $Y$ be a definable set admitting a finite-to-one map
 into $\G^n$, and let  
  into $h$ be a definable map on $Y$ into $\VF$ or 
 $\VF/\Oo$ or $\VF/\Mm$.  Then $h$ has finite image. \>{cor}
 
 \proof One can view $h$ as a function from a subset of $\G^n$
 into finite sets of balls.      Since a ball whose radius is  definable containing  a definable ball is itself definable, 
 \lemref{fnsg} implies that  $h(\g) \in \acl(\emptyset)$ for any $\g \in \G^n$.
  By \lemref{tp-f},   the corollary follows.  \qed

\<{cor} \lbl{Om-rv} Let  $Y \subseteq (\RV \union \G)^n$ and 
$Z \subseteq \VF \times Y$ be definable sets, with $Z$  invariant for the action
of $\Mm$ on $\VF$.  Then for all but finitely many $\Oo$-cosets $C$, 
$Z \meet (C \times Y) $ is a rectangle  $C \times Y'$.  \>{cor}

\prf    Let $p:  (\RV \union \G)^n \to \G^n$
be the natural projection, and for $\g \in \G^n$ let $Z_\g$ be the fiber.  For each $\g$,
by \lemref{Om-st}, there exists a finite $F(\g) \subseteq \VF/\Oo$ such that for any
$\Oo$-coset $C \notin F(\g)$, $Z_\g \meet (C \times Y)$ is $\Oo$-invariant.
Now $\{(u,\g): u \in F(\g)\}$ projects finite-to-one to $\G^n$, so by 
  \lemref{fns-from-gamma}, this set projects to a finite subset of $\VF/\Oo$.  Thus there
  exists a finite $E \subset \VF/\Oo$ such that for any $\g$, and 
  any  $\Oo$-coset $C \notin E$,  $Z_\g \meet (C \times Y)$ is $\Oo$-invariant.
In other words, for any $C \notin E$, $Z \meet (C \times Y)$ is $\Oo$-invariant.
\eprf

  \<{lem}  \lbl{fnsg2}  Let $M \models \T$, $A$ a substructure of 
 $M$ (all imaginary elements allowed), and let $r=(r_1,\ldots,r_m)$
be a tuple of elements of $\RV(M) \union \G(M)$.  Then any closed ball almost defined
 over 
$A(r)$ contains a ball almost defined over $A$. \>{lem}
  
\proof This reduces to $m=1$, $r=r_1$; moreover using \lemref{fnsg},
to the case $r \in \RV(M), \valr(r) = \g \in A$.  Let $E = \{y \in \RV: \valr(y) = \g \}$.
Then $E$ is a $k^*$-torsor, and so is strongly minimal within
$M$.  If $c$ is almost  defined over $A(r)$, there exists
an $A$-definable set $W \subset  E   \times \cB$, with 
$W(e) = \{y: (e,y) \in W\}$ finite, and $c \in W(r)$.  But then $W$ is a finite union of 
strongly minimals,  hence so is the projection $P$ of $W$ to $\cB$.
But any strongly minimal subset of $\cB$ is finite (otherwise it admits a definable
map onto a segment in $\G$; but $\G$ is linearly ordered and cannot have a strongly minimal segment.)   So $c \in P$
is almost defined over $A$.  \qed

 \<{lem} \lbl{red1c1} Let $M \models \T$, $\T$  $C$-minimal with centered closed balls.
  Let $B$ be substructure
 of $\VF(M) \union \RV(M) \union \G(M)$.   Then every $B$-definable 
 closed   ball  has a $B$-definable point.   If 
 $Y$ is a finite $B$-definable  set of disjoint closed balls, there exists
 a finite $B$-definable set $Z \subset M$, meeting each ball of 
 $Y$ in a unique point.
 \>{lem}
 
\proof  We may take $B$
to contain a subfield $K$ and be generated over $K$ by finitely
many points $r_1,\ldots,r_k \in RV$.  Let $Y$ be a a finite $B$-definable  set of disjoint closed balls, and let $b \in Y$.  We may assume
all elements of $Y$ have the same type over $B$.  
By \lemref{fnsg2},  there exists a closed ball $b'$ defined almost
over $K$ and contained in $b$.  By \aref{V+} \ref{V+c} 
there exists a finite $K$-definable set $Z'$ meeting $b'$ in a unique
point.  Let $Y' = \{b'' \in Y: b'' \meet Z' \neq \emptyset \}$, and
$Z= \{\av(Z' \meet b''): b'' \in Y' \}$.  Then $Z$ meets
each ball  of $Y'$  in a unique point, and $Z,Y'$ are $B$-definable.
As for $Y \setminus Y'$, it may be treated inductively.  \qed

 \<{cor} \lbl{red1c0} Let $M \models \T$,  $\T$  $C$-minimal with 
 centered closed balls, and effective.   Let $B$ be an almost $\G$-generated  substructure.   Then $\T$ is effective. 
 \>{cor}

\prf  Same as the proof of \lemref{red1c1}, using \lemref{V+11} in place of \lemref{fnsg2}.
\eprf

\<{lem}\lbl{fmrfr}   Let $Y$ be a $\T$ - definable set   admitting
 a finite-to-one map  into $\RV^n$.
Let $g: Y \to \VF^m$ be another definable
map.  Then $g(Y)$ is finite.
\>{lem}  
 
\<{proof}    It suffices to prove this for $\T_A$, where   $A \models \T$.   We may also assume $m=1$.     We will use the equivalence $(3) \iff (4)$ of
\lemref{tp-f}. If $g(Y)$ is infinite, then by compactness
there exists $a \in g(Y)$ $a \notin \acl(A)$.  But for some $b$ we have
$a = g(b)$, so if $c=f(b)$, we have $c \in \RV^n, a \in \acl(c)$.  Thus
it suffices to show:  

(*) If $a \in \VF$, $c \in \RV^n$ and $a \in \acl(A(c)))$, then $a \in \acl(A)$. 

This clearly reduces to the case $n=1$, $c \in \RV$.  Let $d = \valr(c)$, $A' = \acl(A(d))$.
Then $c$ lies in an $A'$-definable strongly minimal set $S$
 (namely $S=\valr \inv (d)$).  
 Using \lemref{tp-f}
in the opposite direction, since $a \in \acl(A'(c)))$
there exists a finite-to-one map $f:S' \to S$ and a definable
map $g': S' \to \VF$ with $a \in g'(f \inv (S'))$.   By  \corref{fmrf-r}, $g'(f \inv (S'))$ is finite.
Hence $a \in \acl(A(d))$.  But then by  \lemref{fns-from-gamma}, $a \in \acl(A)$.     \>{proof}

     In particular, there can be no definable isomorphism between an infinite subset of
$\RV^n$ and one of $\VF^m$.

 \<{lem}\lbl{fmrfc}  Let $M \models \T$, $\T$ $C$-minimal with centered closed balls,
  and let $A$ be a substructure
 of $M$.  Write $A_\VF$ for the field elements of $A$, $A_{\RV}$
 for the $\RV$-elements of $A$.

 Let
 $c \in \RV(M)$, and let $A(c) = dcl(A \union \{c\})$.
   Then $A(c)_\VF \subset (A_\VF)^{alg}$, and 
   $\rv(A(c)_\VF) \meet A_{\RV} = \rv(A_\VF)$.  
  \>{lem}

\proof  Let $e \in A(c)_\VF$.   Then $e = f(c)$
for some $A$-definable function $f: W \to \VF$, $W \subseteq \RV$.
By   \lemref{fmrfr}, the image of $f$ is finite, 
$e \in acl(A)$.    This proves the first point.    Now if $d \in \RV_A$ and $\rv^{-1}(d)$ has a point
in $A(c)$, then it has a point in $(A_\VF)^{alg}$, by
 \aref{V+} \eqref{V+c}  \qed

\ssec{Transitive sets in dimension one}

 Let $b$ be a closed ball in a valued field.  Then the set $\Aff(b)$ of maximal open sub-balls of $b$
 has the structure of an   affine space over the residue  field.  We will now begin using
 this structure.  
Without it,  more general transitive annuli (missing more than one ball) could exist.

\<{lem} \lbl{ballrec} Let $X \subseteq \VF$ be a transitive $\T_B$-definable set, where $B$ is some 
set of imaginaries.   
Then $X$ is a finite union of open balls  of equal size, or a finite union of closed
balls of equal size, or a finite union of  thin annuli.  
  \>{lem}
 
\proof   By $C$-minimality, $X$ is a finite Boolean combination of balls.  There
are finitely many distinct balls $b_1,\ldots,b_n$ that are almost contained in $X$
(i.e. $b_i \setminus X$ is contained in a finite union of proper sub-balls of $b_i$)
but such that no ball larger than $b_i$ is almost contained in $X$.  These $b_i$ must be disjoint.   If some of the $b_i$
have different type than the others, their union (intersected with $X$) will be a proper
$B$-definable subset of $X$.  Thus they all have the same type over $B$; in particular
they have the same radius $\beta$. 

Consider first the case where  the balls $b_i$ are open, then $b_i \subseteq X$:
otherwise $b_i \m X$ is contained in a unique smallest ball $c_i$; say $c_i$
has radius $\alpha$; then $\alpha > \beta$.  Let $b_i'$ be the open ball
of radius $(1/2)(\alpha+\beta)$ around $c_i$; then $\union_i b_i'$ is a $B$-definable
proper subset of $X$, a contradiction.  Thus in the   case of open balls,
$X \supseteq \union _i b_i$ and therefore $X = \union _i b_i$.

If the balls $b_i$ are closed, let $c_{ij}$ be a minimal 
finite set of sub-balls of $b_i$ needed to cover 
$b_i \m X$.  The same argument shows that no $c_{ij}$ has radius $< \beta$.
Thus all $c_{ij}$ are elements of the set $V_i$ of open sub-balls of $b_i$ of radius $\beta$.  
Now $V_i$ is a $\k$-affine space, and if there is more than one $c_{ij}$ then
over $\acl(B)$, $V_i$ admits a bijection
with $\k$; so there is a finite $B$- definable set of bijections $V_i \to \k$; 
since any finite definable subset of $\k$ is contained in a strictly bigger one, the union of the
pullbacks gives a   $B$-definable subset of $V_i$ properly containing
the ${c_{ij}}$, leading to a proper $B$-definable subset of $X$.  Thus  
either $b_i \subseteq X$ (and then $X =  \union _i b_i$), or else $b_i \m c_i \subseteq X$
for a unique maximal open sub-ball $c_i$.  Now $\union c_i$ intersects $X$ in a proper
subset, which must be empty.  So in this  case, $X = \union_i (b_i \m c_i)$.  \qed

Let $X$ be a transitive $B$-definable set.  Call $Y \subseteq X$ {\em potentially transitive} if
 there exists $B' \supset B$ such that  $Y$ is $B'$-definable  and $B'$-transitive.   
Let  ${\mathcal F}(X)$ be the collection of all proper potentially transitive subsets $Y$ of $X$.   Let $ {\mathcal F}_{max}(X)$ be the set of maximal elements of ${\mathcal F}(X)$.    
  
\lemm{maxtrans}  \begin{enumerate}
  \item If $X$ is an open ball, $ {\mathcal F}_{max}(X) = \emptyset$
  \item If $X$ is a  closed ball, $ {\mathcal F}_{max}(X) =\{X \m Y: Y \in \Aff(X) \}$
  \item If $X$ is a thin annulus $X' \m Y$ with $X'$ closed, then $ {\mathcal F}_{max}(X) = \Aff(X) \m \{Y\}$
\end{enumerate}
\>{lem}
 
 \prf  Any element of ${\mathcal F}(X)$ must be a ball or a thin annulus, so the lemma follows by inspection \eprf

\<{lem} \lbl{image.}   Let $b$ be a transitive  closed ball (respectively thin annulus).  Let $Y=\Aff(b)$ be the set of maximal open sub-balls of $b$.  Then the group of automorphisms
of $Y$ over $\k$ is definable, acts transitively on $Y$, and in fact contains 
$G_a(\k)$ (respectively $G_m(\k)$). 

If $b,b'$ are transitive definable closed balls, and 
$F: b \to b'$ a definable bijection, let $F_*:Y(b) \to Y(b')$ be 
the induced map. Then $F_*$ is a homomorphism of affine spaces, i.e. there exists
a vector space isomorphism $F_{**}: V(b) \to V(b')$ between the corresponding vector
spaces, and $F_*(a+v)=F_*(a)+F_{**}(v)$.     If $b=b'$ then $F_{**}=Id$.   \>{lem}

\<{proof}  $Y=\Aff(b)$ is transitive,
and  there is a $\k$- affine space structure on $Y$ (respectively a $\k$-vector space structure on $V = Y' \du \{0\}$.)    Let $G=Aut(Y/\k)$ be the subgroup of the group $Aff = (G_m \sd G_a) (\k)$ of 
affine transformations of $Y$ that preserve all definable relations.
By definition, this is an intersection of definable subgroups of $Aff$.  However there is no infinite 
descending chain of definable subgroups of $Aff$, so $G$ is definable.  

   If $G$ is finite, then $Y \subseteq \acl(\k)$, and it follows (cf. \secref{stab-emb} that 
there are infinitely many algebraic points of $Y$, contradicting transitivity.  
Thus $G$ is an infinite subgroup of $(G_m \sd G_a) (\k)$ such that
 the set of fixed points $Y^G$ is empty.
 So $G$ must contain a translation, and by strong minimality it must
contain $G_a(\k)$.  
Similarly in the case of the annulus $G$ is an infinite definable subgroup of $G_m(\k)$, 
so it must equal $G_m(\k)$.  

As for the second statement, $F$ induces a group isomorphism $Aut(Y(b)/k) \to Aut(Y(b')/k)$,
and hence an isomorphism $G_a(k) \to G_a(k)$, which must be multiplication by some
$\g \in \k^*$.  Since $G_a(k)$ acts by automorphisms on $(Y(b),Y(b'))$, 
any definable function $Y(b) \to Y(b')$ commutes with this action and hence has the
specified form.  If $b=b'$ then $Y(b)=Y(b')$; now if $F_{**} \neq Id$ then 
$F_*$ would have a fixed point, contradicting transitivity.   \>{proof}

\<{lem}  \lbl{image1}  Let $b$ be a transitive $\T_B$-definable closed (open) ball.
  Let $F$ be a $B$-definable   function, injective on $b$.  Then 
$F(b)$ is a closed (open) ball. \>{lem}

\prf   %Note that $b$ is actually finitely primitive over $B$ (\exref{exprim}).  
 
By \lemref{ballrec}, since $F(b)$ is also transitive, it is either a closed ball, or an open ball, or a thin annulus.  We must
rule out the possibility of a bijection between such sets of distinct types.

Consider the collection ${\mathcal F}_{max}(b)$ defined above.   Any definable bijection between $b$ and $b'$
clearly induces a bijection ${\mathcal F}_{max}(b) \to {\mathcal F}_{max}(b')$.   By \lemref{maxtrans}, the bijective image of an open ball is an open ball.
 
Let $b$ be a closed ball, $b'=b'' \m b'''$ a closed ball minus an open ball, $A=  {\mathcal F}_{max}(b) \iso \Aff(b)$
$A' =  {\mathcal F}(b') \iso \Aff(b'') \m \{b'''\}$, $G=Aut(A/\k),
G' = Aut(A'/\k)$.   Then a definable bijection $A \to A'$ would give a definable group isomorphism $G \to G'$.  But
by \lemref{image.}, $G'=G_m(\k)$ while $G$ contains $G_a(\k)$, so no such isomorphism
is possible (say because $G_m(\k)$ has  torsion points.)

Thus the three types are distinct. \eprf

We will see later that there can be no definable bijection between an open and a closed ball, whether
transitive or not.

  \<{lem}\lbl{rv-transitive}  Let $b$ be a transitive ball.  Then every definable function
 on $b$ into $\RV \union \G$ is constant.  If $b$ is a transitive thin annulus,
 every  definable function
 on $b$ into $\k \union \G$ is constant.   More generally, this is true for definable
 functions into definable  cosets $C$ of $\k^*$ in $\RV$ that contain
algebraic points.

 \>{lem}
 
 \prf    When a ball $b$ is transitive, it is actually   finitely primitive.  
  For if $E$ is 
a $B$-definable equivalence relation with finitely many classes, then exactly one of these
classes is generic (i.e. is not contained in a finite union of proper sub-balls of $b$.)  
This class is $B$-definable, hence must equal $b$.
 
   Thus a definable function on $b$ with finite image is 
 constant.
 
 Let $F$ be a definable function on $b$ into $\G$.  If $F$ is not constant, then for some 
 $\g \in \G$, $F \inv (\g)$ is a proper subset of $b$; it follows that some finite union of 
 proper sub-balls of $b$ is $\g$-definable.  By \lemref{fnsg}, it follows that some such finite union is already definable, a contradiction.
 
 Thus it suffices to show that functions into a single coset $C=\valr \inv (\g)$ of $\k^*$ are constant on $b$. 
 
 Assume first that  $b$ is open, or a properly infinite intersection of balls. By \lemref{orth-gen}
 definable 
 functions on $b$ into $C$ are generically constant; but then by transitivit  they
 are constant.   
 
 Now suppose $b$ is closed, or a thin annulus.  Let $Y$ be the set of 
 maximal open sub-balls $b'$ of $b$.  Each $b' \in Y$ is transitive over $\T_{b'}$, so $F | b'$
 is constant.  Thus $F$ factors through $Y$.  
  
   In the case of the annulus, by \lemref{image.}, $G_m(\k)$  acts transitively
 on $Y$ by automorphisms over $\k$.  This
   suffices to rule out nonconstant functions into $\k$.  More generally, 
 if a coset $C$ of $\k^*$ has algebraic points, then $Aut(C/\k)$ is finite.  Since 
 $Aut(Y/\k)$ is transitive, it follows that if $f: Y \to C$ is definable then 
 $f(Y)$ is finite.  But $Y$ is finitely primitive, so $f(Y)$ is a point.
 
 Assume finally
 that $b$ is a closed ball.   Using \lemref{image.}, we can view $G_a(\k)$ as a subgroup
 of $Aut(Y/\k)$. 
  $Aut(C/\k)$ is contained in $G_m(\k)$. 
 Let $S=Aut(Y \times C / \k) \meet  (G_a(\k) \times G_m(\k))$.  Then $S$ projects onto $G_a(\k)$.  
By strong minimality, $S \meet (G_a(\k) \times (0))$ is either $G_a(\k)$ or a finite group.
%hence zero since the residue field has characterstic $0$ by \aref{V+} \eqref{V+RV} .
  In the first case,  
$S=G_a \times T$ for some $T \leq G_m$. 
 In the latter,   
$S$ is the graph of an finite-to-one homomorphism $G_a \to T$; 
but this is impossible.
% since
%  $G_a$ has too few or too many torsion points.  
  Thus  $G_a \times (0) \leq S$ and $G_a$ acts transitively on $Y$ by automorphisms 
  fixing $C$; it follows that 
 $F$ is constant.    
   \eprf

\ssec{Resolution and finite generation}

\<{lem}\lbl{fg}  Let $A \leq B$ be   substructures of  a model of  $\T$.     
 Assume   $B$ is  finitely generated over $A$.  Then $\RV(B)$ is
finitely generated over $\RV(A)$.  Also, if $\RV(A) \leq C \leq \RV(B)$
then $C$ is finitely generated over $\RV(A)$. \>{lem}

\<{proof}    Suppose $\G(B)$ has infinitely many elements, $\Qq$-linearly independent
modulo $\G(A)$.  By \lemref{doag-fns}, they are algebraically independent.  
 By \lemref{1gen}, they lift to algebraically independent elements of $B$ over $A$, 
 contradicting the assumption of finite generation.  Thus 
    $rk_{\G} \G(B)/\G(A) < \infty$. 
 It is thus clear that any substructure of $\G(B)$ containing $\G(A)$ is 
finitely generated over $\G(A)$.   So it suffices to show that $\RV(B)$ is
finitely generated over $A \union \G(B)$; replacing $A$ by $A \union \G(B)$, 
 we may assume $\G(B) = \G(A)$.
In this case   $\RV(B) \subset  \RES$.  See   \cite{acvf2} Proposition 7.3 for a proof
  stated for $\ACVF_A$, but valid in the present generality.  Here is a sketch.   One looks at 
  $B=A(c)$ with $c \in \VF$.  If $c \in \acl(A)$ then the Galois group 
  $Aut(\acl(A) / A(c))$ has finite index in $Aut(\acl(A)/ A)$.  Hence the same is true of their images
  in $Aut(\acl(A) \meet \RV)$, and since $\RV$ is stably embedded (by clause (1) of the definition of
  $V$-minimality) it follows that there exists a finite subset
  $C'$ of $A(c) \meet \RV$ such that any automorphism of $\acl(A)$ fixing $A(C')$
  fixes $A(c) \meet \RV$.   By   Galois theory for saturated structures (\secref{Types}) $C'$ generates $A(c) \meet \RV$ over $A$. 
  
   On the
  other hand if $c \notin \acl(A)$, then $tp(c/\acl(A))$ agrees with the generic type over $A$
  of either a closed ball, an open ball, or an infinite intersection of balls.  In the latter two cases, 
    $\RES (A) = \RES(B)$
  using \lemref{orth-gen}.  In  the case of a closed ball $b$, let $b'$ be the unique maximal open sub-ball of $b$ containing $c$.  Then $b' \in A(c)$, and $tp(c/A(b'))$ is generic in the open ball $b'$.  Thus by \lemref{orthplus}, $\RES(B) = \RES(A(b'))$ so it is 1-generated. 
  % (Cf. \lemref{1gen} in case $b$ is a pullback from $\RV$.)
     \>{proof}

Recall $\fB = \fB^o \union \fB^{cl} $ is the sort of closed and open balls.

 We require a variant of a result  from  \cite{acvf2} on  canonical resolutions.  We state it
 for $\fB$-generated structures, but it can be generalized to 
 arbitrary $\ACVF$-imaginaries (\cite{acvf1}).

 The proposition and corollaries  will have the effect
 of allowing free use of the technology constructed in this paper over   arbitrary base.
 (cf. \propref{eff}.)

For this proposition, we   allow $\fB$ (and $\G$) as  sorts, in addition to $\VF$ and $\RV$;
 so that a structure is a subset of $\fB , \G$ of a model of $\T$, closed under definable functions.

Assume for simplicity that $\T$ has quantifier elimination.    (cf. 3.4)

Let us call a structure $A$ {\em resolved} if any ball and any thin annulus defined
over $\acl(A)$ has a point over $\acl(A)$.  

\lemm{curve-selection} Let $\T$ be $\V$-minimal.  Let $M \models \T$,
and $A$ be a substructure of $M$.
Then  (1),(2) are equivalent; if $\G(A) \neq (0)$ then (3) is equivalent to both.
\begin{enumerate}  
\item       $A$ is  effective and $\VF(\acl(A)) \to \G(A)$ is surjective.
  \item  $A$ is resolved.
 \item      $\acl(A)$ is an elementary submodel of $M$.  
\end{enumerate}

\>{lem}  

\prf  Clearly (3) implies (1) and (2) implies (1).
 To prove that (1) implies (3) it suffices to show that every definable $\phi(x)$ of $\T_A$ in one variable, with a solution in $M$, has a solution in $A$.  If $x$ is an $\RV$-variable it
 suffices to show that $\phi(\rv(y))$ has a solution; so  we may assume
$x$ is a $\VF$- variable, so $\phi$ defines $D \subseteq \VF$.  By $C$-minimality $D$ is a finite Boolean combination of balls.  $D$ can be written as a finite union of definable sets
of the form $\union_{j=1}^m D_j \m E_{j}$, where for each $j$,  $D_j$ is a closed ball, and $E_j$ a finite union 
of maximal open sub-balls of $D_j$, {\em or} $D_j$ is an open ball and $E_j$ is a proper
sub-ball of $D_j$, {\em or} $E_j = \emptyset$, {\em or} $D_j=K$.     In the third case, 
by effectivity there exists a finite set meeting each $D_j$ in a point; since $A=\acl(A)$, this finite set is contained in $A$; so $D(A) \neq \emptyset$, as required.  In the first and second cases,
there exists similarly a finite set $Y$ meeting each  $E_j$.  Since $A=\acl(A)$,
$Y \subseteq A$.  By picking a point
and translating by it, we may assume $0 \in E_j$ for some $j$.  Say $E_j$ has valuative radius $\alpha$; picking a point $d \in A$ with $\val(d) = \alpha$ and dividing, we may assume $\alpha=0$.  Now in the open case any element of valuation $0$ will be in $D_j$.  In the closed case,
the image of   $E_j$ under $\res$ is a finite subset of the residue field; pick some
element $\bar{a}$ of $\k(A)$ outside this finite set; by effectivity, pick $a \in A$ with 
$\res(a) = \bar{a}$; then $a \in D$.     In the fourth case, we use the assumption
that $\G(A) \neq (0)$.  This proves (3).

   It remains to show that (1) implies (2).  Let $b$ be a thin annulus defined
over $\acl(\emptyset)$; so $b = b' \m b''$ for a unique closed ball $b'$ and maximal open sub-ball $b''$.
By effectivity, $b''$ has an algebraic point, so translation we may assume $0 \in b''$.
In this case the assumption that $\VF(\acl(A)) \to \G(A)$ is surjective gives a point of 
$b' \m b''$.
  \eprf

If $\T_0$ is $\V$-minimal, $A$ is a finitely generated structure (allowing $\fB$, 
or even $ACVF$-imaginareis), and $T= (\T_0)_A$, we will call $T$ a finitely generated  extension
of a $\V$-minimal theory.

\<{rem}  If $A$ is effective, then $A$ is $\VF \union \G$-generated.
If $A$ is resolved, then $A$ is $\VF$-generated.  \>{rem}

\<{prop} \lbl{resolve}   Let $\T$ be $\V$-minimal.  
 \begin{enumerate}
\item There exists an effective structure $E_{eff}$ admitting an embedding into any effective
structure $E$.  We have $\RV(E_{eff}),\G(E_{eff}) \subseteq \dcl(\emptyset)$.  

\item There exists a  resolved $E_{rslv}$ embedding into any resolved
structure $E$.  We have $\k(E_{rslv}),\G(E_{rslv}) \subseteq \dcl(\emptyset)$.  
In fact $C(E_{rslv}) \subseteq \dcl(\emptyset)$ for any  cosets $C$ of $\k^*$ in $\RV$ that contain
algebraic points.

\item  Let    $A$ be a  finitely  generated 
substructure  of a model of $\T$, in the sorts $\VF \union \fB$.   Then (1),(2) hold
for $\T_A$.
   \end{enumerate}
 \>{prop}

\prf (1)  Let $(b_i)_{i < \lambda}$ enumerate the definable balls.
Define a tower of $\VF$-generated 
structures $A_i$, and a sequence of balls $b_i$, as follows. 
 Let $A_0 = \dcl(\emptyset)$; 
if $\kappa$ is a limit ordinal, let $A_\kappa = \union_{i < \kappa} A_i$.  Assume $A_i$ has been defined.  If possible, let $b_i$ be an $A_i$-definable, $A_i$-transitive ball, not a point; and  let $c_i$ be
any point of $b_i$.  If no such ball $b_i$ exists, 
 the construction ends, and we let $E_{eff}=A_i$ for this $i$.

Suppose $E$ is any effective substructure of a model of $\T$.  We can inductivelly
define a tower of embeddings $f_i: A_i \to E$.  At limit stages $\kappa$ let
$f_\kappa= \union_{i<\kappa} f_i$.  Given $f_i$ with $A_i \neq E$,
let $b_i' $ be the image under $f$ of $b_i$.  By effectivity, $b_i'$ has a point $c'_i \in E$.
Since $b_i$ is transitive over $A_i$, the formula $x \in b_i$ generates a complete
type; so $tp(c_i/A_i)$ is carried by $f$ to $tp(c_i' / A_i')$.  Thus there exists an embedding
$f_{i+1} : A_{i+1} \to E$ extending $f_i$, and with $c_i \mapsto c_i'$.

Each $A_i$ is $\VF$-generated; by \lemref{average} (3) $\implies$ (4), the process 
can only stop when $A_i=E_{eff}$.  This shows that $E_{eff}$ embeds into $E$,
and at the same time that the construction of $E_{eff}$ itself must halt at some stage
(of cardinality $\leq |\T|$.)  

By construction, $E_{eff}$ is $\VF$-generated; and hence $\T_{E_{eff}}$ is 
$\V$-minimal. 
 Moreover there   are no $E_{eff}$-definable
$E_{eff}$-transitive balls (except points).    In other words all $E_{eff}$-definable
balls are centered.   By $\V$-minimality \eqref{V+c}
every closed ball has a definable point, so 
every centered ball has  one.  So $E_{eff}$ is effective.  

It remains only to show that  $\RV(E_{eff}),\G(E_{eff}) \subseteq \dcl(\emptyset)$.  We show
inductively that   $\RV(A_i) , \G(A_i) \subseteq \dcl(\emptyset)$.  At limit stages
this is trivial, and at successor stages it    follows
from \lemref{rv-transitive}.      
  
(2) Identical to (1), but using thin annuli as well as balls.  If a thin annulus is not transitive, it contains a proper nonempty  finite union of balls, so by $\V$-minimality it contains 
a proper nonempty finite set.  Hence the construction of the $A_i$  stops only when $A_i$ is resolved.

(3)      Let $A_0 = (A \meet (\VF \union \G))$.    $A$ is generated over $A_0$ by some
$b_1,\ldots,b_n \in \fB$ with $b_i$ of  valuative radius $\g_i \in A_0$.  
Since $\T_{A_0}$  is $\V$-minimal, we may assume $\T=\T_{A_0}$ and $A$ is
generated
  by  $b_1,\ldots,b_n$, with $\g_i$ definable.
  
    Let $J$ be a subset of $\{1,\ldots,n\}$ 
of smallest size such that $\acl(\{b_j: j \in J\}) = \acl(\{b_1,\ldots,b_n\})$.
By minimality, no $b_j$ is algebraic over 
$\{b_{j'}: j' \in J, j' \neq j\}$.  Let $j \in J$, and let $Y_j$ be the set of balls of radius $\g_j$;
then $Y_j$ is a definable family of disjoint balls.  By \lemref{transitive} for
$\T'=\T_{<\{b_{j'}: j' \in J, j' \neq j\}>}$,   $b_j$ is transitive in $\T'_{b_j}$, i.e. in $\T_{<b_{j'}: j' \in J>}$;
 hence $b_j$ is transitive over $\acl(b_1,\ldots,b_n) = \acl(A)$.  
Let us now show, using induction on $|J|$, that $\Pi_{j \in J} b_j$ is transitive over $A$.
 Let $c_j \in b_j$.     By \lemref{generic-acl}
the $ (\{b_{j'}: j' \in J, j' \neq j\})$ remain algebraically independent over $<c_j>$.
Thus by induction, $\Pi_{j \neq j'} b_{j'}$ is transitive over $A(c_j)$; since
$b_j$ is transitive over $A$, $\Pi_{j \in J} b_j$ is too.  Let $A' = A(c_j: j \in J)$.

\Claim{} If $B$ is a $\VF \union \G$-generated structure containing $A$, then $A'$ embeds into $B$ over
$A$.
 
\proof   Since $B$ is $\VF \union \G$-generated, every ball of $\T_B$ is centered, in particular $b_j$
has a point $c_j'$ defined over $\T_B$.  Let $c'  = (c_j': j \in J)$.  By transitivity of $\Pi_{j \in J} b_j$ we have
 $tp(c/A)=tp(c'/A)$.
So $A'$ embeds into $B$.

Note that  $A'$ is almost $\VF \union \G$-generated; indeed since $\g_i$ is definable,
$b_i \in \dcl(c_i)$ so $A' \subseteq \acl((c_j)_{j \in J})$.  Thus $\T_{A'}$ is $\V$-minimal.
So (1,2) applies and prove  (3).  
 
  \eprf

See \lemref{rslv-uniqueness}  for a uniqueness statement. 

\<{cor} \lbl{resolve-cor}  Let $f: \VF \to (\RV \union \G)^*$ be a definable map.  
\begin{enumerate}
  \item There
exists a definable $\tf: \RV \to (\RV \union \G)^*$ such that for any $\bx \in \RV$, for
some $x \in \VF$ with $\rv(x)=\bx$, $\tf(\bx) = f(x)$.
  \item Let $\Omega = \VF/\Mm$.
  There
exists a definable  map $\tf: \Om \to (\RV \union \G)^*$ such that for any $\bx \in \Om$, for
some $x \in \VF$ with $x+\Mm=\bx$, $\tf(\bx) = f(x)$.
\end{enumerate}
     \>{cor}

\prf   (1) In view of   \lemref{collect}, it suffices to show that 
 for a given complete type $P \subseteq \RV$, there exists such a function $\tf$ on $P$.  
We fix $\ba \in P$, and show the existence of $\bc \in \dcl(\ba)$
such that for some $a$ with $\rv(a) = \ba$, $f(a) = \bc$. 

   By \propref{resolve},
there exists an effective substructure $A$ with 
$\ba \in A$ and $ (\RV \union \G)(A) = (\RV \union \G)(<\ba>)$.  Thus the open ball $\rv \inv (\ba)$ has an $A$-definable point $a$.
Set $\bc = f(a)$; since $f(a) \in \RV(A) = \RV(<\ba>)$ we have
$\bc = \tf_P(\ba)$ for some definable function $\tf_P$.   Clearly 
$\tf_P$ satisfies the lemma for the input $\ba$, hence for any input from $P$.

  The proof of (2) is identical, using \lemref{resolve} (3). \eprf

\<{cor}     \lbl{resolve-total}   Let $\T$ be $\V$-minimal.   Assume
every definable point of $\G$ lifts to an algebraic point of $\RV$.  
 Then  there exists a resolved structure $E_{rslv}$ such that $E_{rslv}$ can be embedded
  into any resolved
structure $E$, and  $\RV(E_{rslv}),\G(E_{rslv}) \subseteq \dcl(\emptyset)$.  
If $A$ is a   finitely  generated 
substructure  of a model of $\T$, in the sorts $\VF \union \fB$, the same is true
for $\T_A$.  \>{cor}

\prf 
Under the assumption of the corollary, 
 the conclusion of \propref{resolve} implies  $\RV(E_{rslv}) \subseteq  \dcl(\emptyset)$.
 \eprf
 
 \<{rem} It is easy to see using the description of imaginaries in  \cite{acvf1} that 
 in a  resolved structure,
  any definable $\ACVF-$ imaginary   
 is resolved.  In other words, 
 if $A$ is a resolved, and $\sim$ is a definable equivalence relation on a definable set $D$, 
 then $D(A) \to (D/\sim)(A)$ is surjective. 
 
If $A$ is only effective, then there exists $\g \in \G(A)^n$ such that for any
$t$ with $\val(t) = \g$, 
 $ (D/\sim)(A) \subseteq \dcl(D(A)/\sim, t)$; this can  be
seen by embedding $D/\sim$ into $B_n(K)/H$ for an appropriate $H \leq B_n(\Oo)$,
and splitting $B_n=T_nU_n$.  
 
 \>{rem}

 \ssec{Dimensions} 
 
We define the {\em $\VF$-dimension} of a $T_M$-definable set $X $ to be the smallest
$n$ such that for some $n$, $X$ admits a $T_M$-definable map with finite fibers into $\VF^n \times (\RV \union \G)^*$. 

By {\em essential bijection $Y \to Z$} we mean
a bijection $Y_0 \to Z_0$, where $\dim_{\VF}(Y \m Y_0), \dim_{\VF}(Z \m Z_0) < \dim_{\VF}(Y) = \dim_{\VF}(Z)$;
and where two such maps are identified if they agree away from a set of dimension $<\dim_{\VF}(Y)$.

We say that a map $f: X \to \VF^n$ has $\RV$-fibers if there exists 
  $g: X \to (\RV \union \G)^*$ with $(f,g)$ injective.

\<{lem} \lbl{vfdim}  Let $X \subseteq \VF^n \times (\RV \union \G)^* $ be  a definable set.  Then 
\begin{enumerate}
  \item 
$X$ has  $\VF$-dimension $\leq n$ iff
there exists a definable map $f: X \to \VF^n$ with $\RV$-fibers.
  \item If it exists, 
the map $f$ is   ``unique up to isogeny'':  
 if $f_1,f_2: X \to \VF^n$ have $\RV$-fibers, then there exists a definable
 $h: X \to Z \subseteq \VF^n \times (\RV \union \G)^*$ and $g_1,g_2: Z \to \VF^n$
 with finite fibers, 
 such that $f_i = g_ih$. 
 
\end{enumerate}

 \>{lem}
 
\prf  (1) If $f: X \to \VF^n$ has $\RV$-fibers, let $g$ be as in the definition of $\RV$-fibers;
then $(f,g): X \to \VF^n \times (\RV \union \G)^*$ is injective, so certainly finite-to-one.
If $\phi: X \to \VF^n \times \RV^*$ is  finite-to-one, by \lemref{finite},
each fiber $\phi \inv(c)$ admits a $c$-definable injective map $\psi_c: \inv(c) \to \RV^*$.
By \lemref{collect} we can find $\theta: X \to \VF^n \to \RV^*$ that is injective on each
$\phi$-fiber.  Let $f(x)=(\phi,\theta)$.    This proves the equivalence.

(2) Now suppose $f_1,f_2: X \to \VF^n$ both have $\RV$- fibers.

Let  $h(x)=(f_1(x),f_2(x))$, $Z'=h(X)$, 
and define  $g_i: Z' \to \VF^n$ by
 $g_1(x,y) = x$, $g_2(x,y)=y$.  
Then $g_i$ has finite fibers.  Otherwise  we can find $a \in X$
such that $f_1(a) \notin \acl(f_2(a))$ (or vice versa.)  But 
for any $a \in X$, we have  $f_1(a) \in \acl(f_2(a),c)$
for some $c \in (\RV \union \G)^*$.  By \lemref{fmrfr}, $f_1(a) \in \acl(f_2(a))$, a contradiction.
By \lemref{finite} (cf. \lemref{collect}), there exists a definable bijection between $Z'$ and a 
subset $Z$ of $\VF^n \times \RV^*$.  Replacing $Z$ by $Z'$ finishes the proof of the lemma. \eprf

\<{cor}\lbl{vfrvdim} Let $f: X \to \RV \union \G$, $X_a = f \inv (a)$.
  Then $\dim(X) = \max_{a } \dim X_a$. \>{cor}
  
\prf  Let $n= \max_{a } \dim X_a$.  
For each $a$ there exist  definable functions $g_a: X_a \to \VF^{n}$ and $h_a: X_a \to (\RV \union \G)^*$ with $(g_a,f_a)$
injective on $X_a$.  Thus by the compactness argument of \lemref{collect},
 there exists  definable
functions $g: X \to \VF^n$ and $h: X_a \to (\RV \union \G)^*$ such that $(g,h)$
is injective when restricted to each $X_a$.  But then clearly $(g,h,f)$ is injective,
so $\dim(X) \leq n$.  The other inequality is obvious.  \eprf

We continue to assume $\T$ is $\V$-minimal.

\lemm{vf-exchange}   Let $a,b \in \VF$.  If $a \in \acl(b) \m \acl(\emptyset)$ then
$b \in \acl(a)$.
\>{lem}

\prf Suppose $b \notin \acl(a)$.  Let $A_0= \G(\acl(a,b))$.  Then by 
\lemref{fns-from-gamma}, $b \notin \acl(A_0(a))$.  

Let $C$ be the intersection of all $\acl(A_0)$-definable balls such that $b \in C$,
and let $C'$ be the union of all $\acl(A_0)$-definable proper sub-balls of $C$.  
Let $B= \meet_i \{B_i\}$ be the set of all balls defined over $\acl(A_0(a))$ with $b \in B_i$,
and let $B' = \union_j \{B'_j\}$ be the union of all  $\acl(A_0(a))$-definable proper sub-balls of $B$.  

Since $a \in \acl(b)$, we have  $a \in \acl(b')$ for all $b' \in B \m B'$, outside some proper sub-ball.
It follows by compactness that for some $i,j$,  $a \in \acl(b')$ for all $b' \in B_i \m B'_j$. 
Say $i=j=1$, $B'_1 \subset B_1$.   
By \exref{vf-exchange},  $a \in \acl(A_0(f_1))$ where $f_1 \in \fB$ codes the ball $B_1$.  

If $B_1$ is a point, we are done.  Othewise $B_1$ 
has valuative radius $\alpha_1 < \infty$  defined over $A_0$.  It follows that if 
 $B_1 \supseteq C$ then 
$B_1$ is $\acl(A_0)$-definable; but then 
  $a \in \acl(A_0)$, contradicting the assumption.  Since $B_1$ meets $P$ nontrivially,
we have therefore $B_1 \subset C$.  Similarly $B_1$ cannot contain any ball in $C'$
since it is not $\acl(A_0)$-definable, 
but it cannot be contained in $C'$ since $B_1 \meet P \neq \emptyset$.   so $B_1 \meet C' = \emptyset$.  Thus $B_1 \subset  P$.
 
 Let $\bar{B_1}$ be the closed ball of radius $\alpha_1$ containing $B_1$, and let
 $e_1$ be the corresponding element of $\cB$.     Since $\bar{B_1}$ is almost definable
 over $A_0(a)$,   it follows from  $\V$-minimality
  
 that
there exists an almost $A_0(a)$-definable point $c(a)$ in $\bar{B_1}$.
Now if  $a \in \acl(A_0(e_1))$, then $\bar{B_1}$ contains an $A_0(e_1)$-definable finite set $F_1=F_1(e_1)$.  But since $B_1$ is a proper subset of $P$, $e_1 \notin \acl(A_0)$.  This 
contradicts \lemref{transitive}.   So $a \notin \acl(A_0(e_1))$.

 Nevertheless we have seen that $a \in \acl(A_0(f_1))$.  Thus $B_1 \neq \bar{B_1}$, 
 so $B_1$ is a maximal open sub-ball of $\bar{B_1}$.
Let $b_1$ be the point of $\Aff(\bar{B_1})$ representing $B_1$.  Then $a \in \acl(b_1)$.
It follows that $tp(a / acl(A_0(e_1)))$ is strongly minimal, contradicting \lemref{fmrf-r}.    
We have obtained a contradction in all cases; so $b \in \acl(a)$.
\eprf

Since the lemma continues to apply over any $\VF$-generated structure, 
algebraic closure is a dependence relation in the sense of 
Steinitz (also called a pre-matroid or combinatorial geometry, cf.   \cite{rota}.)  
Define the $\VF$-transcendence degree of a finitely generated structure $B$
 to be the maximal number of elements of $\VF(B)$ that are algebraically independent
over $\VF(A)$.  This is the size of any maximal independent set, and also the minimal
size of a subset whose algebraic closure includes all $\VF$-points.    Hence:

\<{cor}  \lbl{steinitz}
 The $\VF$ dimension of a definable set $D$ is the maximal transcendence degree
of $<b>$.  \qed \>{cor}

  We can now obtain a strengthening of \lemref{fmrfr}, and a uniqueness
statement in \propref{resolve}:

\<{cor}\lbl{fmrfr+}   Let $Y$ be a $\T$ - definable set   admitting
 a finite-to-one map $f$ into $\fB^n$.
Let $g: Y \to \VF^m$ be a  definable
map.  Then $g(Y)$ is finite.
\>{cor}  
 
\<{proof}   We may   assume $m=1$.     We will use the equivalence $(3) \iff (4)$ of
\lemref{tp-f}. If $g(Y)$ is infinite, then by compactness
there exists $a \in g(Y)$ $a \notin \acl(A)$.  But for some $b$ we have
$a = g(b)$, so if $c=f(b)$, we have $c \in \fB^n, a \in \acl(c)$.  Thus
it suffices to show:  

(*) If $a \in \VF$, $c \in \fB^n$ and $a \in \acl(A(c)))$, then $a \in \acl(A)$. 

This clearly reduces to the case $n=1$, $c \in \fB$.   Let $\g$ be the valuative
radius of $c$.  As follows from 
\corref{fns-from-gamma}, it suffices to show that  $a \in \acl(A(\g))$.  So in (*)
we may assume $\g \in A$.

Finally, to prove (*) (using again the equivalence of \lemref{tp-f}), we may 
enlarge $A$, so we may assume $A \models \T$.  

Since $\g \in A$, $c \in \dcl(A(e))$ for any element $e$ of  the ball $c$.  Thus 
  $a \in \acl(A(e))$.   Suppose $a \notin \acl(A)$;
then by  exchange for algebraic closure in $\VF$,
$e \in \acl(A(a))$.  Thus any two elements of the ball $c$ are algebraic over each other.
By \exref{alldepend},  $c$ has finitely many points; which is absurd.  This contradiction 
shows that $a \in \acl(A)$.
   \>{proof}

 \<{lem}\lbl{rslv-uniqueness} (cf. \propref{resolve}.) Let $\T$ be a 
  finitely generated extension of an effective $\V$-minimal theory. 
 Then if $E_1,E_2$ are effective
 and both embed into any effective $E$, then they are finitely generated, and  $E_1 \iso E_2$.    
  \>{lem}
 
 \prf The finite generation is clear.  Since $E_1,E_2$ embed into each other, they
 have the same $\VF$-transcendence degree  
 We may assume $E_1 \leq E_2$.  But then by \lemref{steinitz}, 
 $E_2 \subseteq \acl(E_1)$.  By \lemref{finite}, $E_2 \subseteq \dcl(E_1,F)$
 for some finite $F \subseteq \RV^* \meet \dcl(E_2)$.  But $\RV(E_1)=\RV(E_2)$,
 so $F \subseteq \dcl(E_1)$, and thus $E_2 = E_1$.  
 
 \<{remark}  The analogous statement is true for  resolved
 structures.  Note that if $F$ is a finite definable
 subset of $\RV^n$ then automatically the coordinates of the points of $F$ lie in cosets of 
 $\k^*$ that have algebraic points. \>{remark}
 
  \eprf

{\bf Remark} The hypothesis of \lemref{rslv-uniqueness} can be slightly weakened to:  $\T$ is  
 finitely generated over a $\V$-minimal theory, and   there exists a finitely generated effective $E$.

\<{example}  In $\ACVF$, when $X \subseteq \VF^n$, the $\VF$-dimension equals the dimension of the Zariski closure of $X$.  \rm
 This is proved in \cite{vddries-alg-bdd}.  The idea of the proof: the $\VF$-dimension is
 clearly bounded by the Zariski dimension.   For the opposite inequality, in the case of dimension $0$, if $X$ is a finite $A$-definable
subset of $\VF$, then using quantifier elimination there exists a nonzero polynomial $f$ with coefficients in $A$, such that $f$ vanishes on $X$.   In general, if 
 a definable $X \subseteq \VF^n$ has $\VF$ dimension $<n$, one can reduce to the case where all   fibers of the projection
$\pr: X \to \pr X \subset \VF^{n-1}$ are finite, then $X$ is not Zariski dense
in $\VF^n$, using the $0$-dimensional case. \>{example}

The  {\em $\RV$-dimension} of a definable set $X \subseteq \RV^*$ is the smallest integer $n$
(if any) such that $X$ admits a parametrically definable finite-to-one map 
into $\RV ^n$.   More generally for $X \subseteq (\RV \union \G)^*$, $\dim_{\RV}(X)$ is
the smallest integer $n$
(if any) such that $X$ admits a parametrically definable finite-to-one map 
into $(\RV \union \G)^n$

  Note that $\RV$ is one-dimensional, but $\G$ and every fiber of $\valr$ are also
  one-dimensional.  In this sense $\RV \union \G$-dimension is not additive; model-theoretically it is closer to weight than to rank.   We do have $\dim(X \times Y) = \dim(X) + \dim(Y)$.

Dually, if a structure $B$ is $\RV$-generated over a substructure $A$, we can define the {\em weight} of $B/A$ to be the least $n$ such that $B \subseteq \acl(A,a_1,\ldots,a_n)$, 
with $a_i \in \RV$.

For subsets of $\RV$, $\RV$-dimension can be viewed as the size of a Steinitz basis with respect to algebraic closure.  One needs to note that the exchange principle holds:

\<{lem}[exchange]  \lbl{exchange} Let $a,b_1,\ldots,b_n \in \RV$; assume $a \in \acl(A,b_1,\ldots,b_n)
\setminus \acl(A,b_1,\ldots,b_{n-1})$.  Then $b_n \in \acl(A,b_1,\ldots,b_{n-1},a)$.
\>{lem}

\<{proof}    We may take $n=1, b_n=b$, and $A = \acl(A)$.  Let $\a = \valr(a) \in \G$, $\b = \valr(b)$.  
If $\b \in A$ then $\G(A(a,b)) = \G(A(b)) = \G(A)$.  The first equality is true since $a \in \acl(A(b))$
so $A(a,b) \subset \acl(A(b))$, and using the stable embeddedness of $\G$ (\secref{stab-emb}) and the linear ordering on $\G$.   The second equality follows
from \lemref{orth1}.   So if $\b \in A$ then 
 $a,b$ lie in $A$-definable strongly minimal sets, cosets of $\k^*$, and the lemma is clear.  
  
  Assume $\b \notin A$.  
If $\a \in A$,
then $tp(a/A)$ is strongly minimal, and $tp(a/A)$ implies $tp(a/ A(b))$ by \lemref{orth1};
; but then $a \in \acl(A)$,
contradicting the assumption.  So $\a,\b \notin A$; from the exchange principle in $\G$
it follows that $A' := \acl(A,\a)= \acl(A,\b)$.  
Moreover $a \notin \acl(\a) $ by \lemref{nosec} and \lemref{tp-f}. 
By the previous case, $b \in \acl(A',a)$, so $b \in \acl(A,a)$.         \eprf

\<{lem}\lbl{rvd1}  A definable $X \subseteq \RV^n$ has $\RV$-dimension $n$ iff it contains an $n$-dimensional 
definable subset of some coset of ${\k^*}^n$. 
\>{lem}

\<{proof}  Assume $X$ has $\RV$-dimension $n$.  Then there exists $(a_1,\ldots,a_n) \in X$ with
$a_1,\ldots,a_n$ algebraically independent.  Let $c \in \G$; then since $a_n \notin \acl(a_1,\ldots,a_{n-1})$,
it follows as in the proof of \lemref{exchange} that $a_n \notin \acl(a_1,\ldots,a_{n-1},c)$.  This applies
to any index, so $a_1,\ldots,a_n$ remain algebraically independent over $c$; and inductively we may add
to the base any finite number of elements of $\G$.  Let $c_i = \valr(a_i)$, and let $A'=A(c_1,\ldots,c_n)$.
Then $a_1,\ldots,a_n$ are algebraically independent over $A'$, and the lie in $X'=X \meet \Pi_{i=1}^n \rv \inv (c_i)$;
thus $X'$ is an $n$-dimensional definable subset of a coset of ${\k^*}^n$.   
\>{proof}

\<{defn} \lbl{VFcat}  $\VFni$ be the category of definable subsets of $\VF^* \times \RV^*$ of dimension
$\leq n$.  Morphisms are definable maps.  \>{defn}

Let $X \in \Ob \VFni$.  By \lemref{vfdim}, there exists a definable 
   $f: X \to \VF^n$  
with $\RV$-fibers; and the   maximal $\RV$-dimension of  a fiber is a well-defined
quantity, depending only on the isomorphism type of $X$ (but not on the choice of $f$).  In particular,
the   subcategory of definable sets of maximal fiber dimension $0$ will be denoted $\VF[n]$.

\<{defn} \lbl{RVcat} \rm                                                          
We define 
  $\RVni$ to be the category of definable pairs $(U,f)$, with   $U \subseteq  \RV^*$,
  $f: U \to \RV^n$.   If $U,U' \in \Ob \RVni$, a morphism $h: U \to U'$
 is a definable map, such that $U''=\{(f(u),f'(h(u))): u \in U \}$ has finite-to-one first projection
to $\RV^n$.     $\RV[n]$ is the full subcategory of pairs $(U,f)$
 with $f: U \to \RV^n$ finite-to-one.
 
 $\RES[n]$ is the full subcategory of $RV[n]$ whose objects are pairs $(U,f) \in \Ob \RV[n]$
 such that  $\valr (U)$ is finite, i.e. $U \subseteq \RES^*$.
  \>{defn}
 
\<{remark}\lbl{rvcat-rem}  \begin{enumerate}
  \item For $X,Y \on \Ob \RV[n]$, any definable bijection $X \to Y$ is in 
 $\Mor_{\RV[n]}(X,Y)$.
  \item The forgetful map $(X,f) \mapsto X$ is an equivalence of categories between
  $\RV[n]$ and the category of all definable  subsets of $\RV^*$ of $\RV$-dimension $\leq m$, with all maps between them.  The presentation with $f$ is nonetheless useful for 
  defining $\L$.
\end{enumerate}
\end{remark}

By \remref{rvcat-rem}  $\SG(\RV[m])$ is isomorphic to the Grothendieck semigroup of 
 definable subsets of $\RV^*$ of $\RV$-dimension $\leq m$.  If $\dim(X)\leq m$, 
 let $[X]_m$ denote the class $[X]_m = [(X,f)]_m \in \RV[m]$, where $f : X \to \RV^*$ is any finite-to-one
 definable map.
  
Unlike the case of $\VF[n, \idot]$ or $\RV[n]$, for   $(U,f) \in \Ob \RVni$ the map $f$ cannot be reconstructed from $U$ alone, even up to isogeny,
 so it must be given as part of the data.  We view $(U,f)$ as a cover of $f(U)$ with 
  ``discrete'' fibers.
 
 We denote 
\[\RV[\leq N, \idot] := \oplus_{0 \leq n \leq N} \RVni, \ \ 
\RV[\leq N] = \oplus_{0 \leq n \leq N} \RV[n]  \]
\[\RV[*,\idot]  := \oplus_{0 \leq n} \RVni, \ \  \ \
 \RV[*]  := \oplus_{0 \leq n} \RV[n]   \ \ \ \ 
 \RES[*] := \oplus_{0 \leq n} \RES[n]
   \]
 We have natural multiplication maps $\SG \RV[k,\idot] \times \SG \RV[l,\idot] \to \SG[k+l,\idot]$,
 $([(X,f)],[(Y,g)]) \mapsto [(X \times Y, f \times g)]$.  This gives a semiring structure
 to  $\SG(\RV[*])$.   This differs from the Grothendieck ring $\SG(\RV)$.

 \sssec{Alternative description of $\RV[\leq N,\idot]$. }  \lbl{altRV}
An object of  $\RV[\leq N, \idot]$ thus consists of a formal sum $\sum_{n=0}^N \bX_n$ 
of objects $\bX_n = (X_n,f_n)$ of $\RVni$.  This can be explained from another 
angle if one adds a formal element  $\infty$ to $\RV$, and extends $\rv$ to $\VF$
by $\rv(0)=\infty$.  Define a function $f[k]$ by $f[k](x)=(f_n(x),\infty,\ldots,\infty)$ ($N-k$ times).  
If $\bX = (X,f)$, let   $\bX[k] = (X,f[k])$.  
Then
$\sum_{n=0}^N \bX_n$  can be viewed as the disjoint union 
$\union_{i=0}^N X_i \times \{\infty\}[N-i]$.   The $\rv$ pullback is then a set of $\VF$-dimension
$N$, invariant under multiplication by $1+\Mm$; the sum over dimensions $\leq N$
is necessary to ensure that any such invariant set is obtained.  (Cf. \lemref{tr2-surjectivity}.)
  From this point of view,
an isomorphism is a definable bijection preserving the function ``number of finite coordinates''.
We will use $\RV[\leq N, \idot]$ or $\RVi[N,\idot]$ interchangeably.

\lemm{RVVF0}   Let $X,X' \in \Ob \RVni$, and assume a bijection $g: X' \to X$ lifts
to $G:   \L X' \to \L X$.  Then $g \in \Mor_{\RVni} (X',X)$.  \>{lem}

\prf   We only have to check the isogeny condition,
i.e. that $f(g(a)) \in \acl(f'(a))$ for $a \in X'$ (and dually).    
By \lemref{fmrfc},  for $x \in \rho_{X'} \inv (a)$, $G(x)_{\VF} \in \acl(x_{\VF})$, 
i.e. the $\VF$-coordinates of $G(x)$ are algebraic over those of $x$.   Thus
$f(g(a)) \in \acl(x_{\VF})$.  This is true for any $x \in \rho_{X'} \inv (a)$, so
$f(g(a)) \in \acl(a)$.  \eprf

\>{section}

\<{section}{Descent to $\RV$:  objects}  \lbl{descent-objects}
 
 We assume $\T$ is $C$-minimal 
 with centered closed balls.   We will find a very restricted set of maps that transform any definable set to a pullback from $\RV$. 
 This is   related to Denef's cell decomposition theorem;  since we work in $C$-minimal theories it takes a simpler form.
 Recall that this assumption is preserved
 under passage to $\T_A$, when $A$ is a $(\VF,\RV,\G)$-generated substructure of a model of
 $\T$ (\lemref{red1c1}).
  
  Recall that $\RV = \VF^\times/ (1+\Mm)$, $\rv: \VF^\times \to \RV$ the quotient map.
Let $\RVi = \RV \union \{\infty\}$,
and define $\rv(0)=\infty$.  
We will also write $\rv$ for the induced map ${\rv}^n: (\VF^\times)^n \to (\RV)^n$.

 \<{defn}  \lbl{fc0}  Fix $n$.  
 Let $\fC^0$ be the category whose objects are the definable subsets
 of $\VF^n \times \RVi^*$, and whose morphisms are generated by
 the inclusion maps together with   functions
 of one of the following types:
 \begin{enumerate}
  \item Maps 
  $$(x_1,\ldots,x_n,y_1,\ldots,y_l) \mapsto (x_1,\ldots,x_{i-1},x_i+a,x_{i+1},\ldots,x_n,y_1,\ldots,y_l)$$
with $a=a(x_1,\ldots,x_{i-1},%x_{i+1},\ldots,x_n,
y_1,\ldots,y_l): \VF^{i-1} \times \RVi^l \to \VF$
an $A$- definable function of the coordinates $y,x_1,\ldots,x_{i-1}$. 

%We can also demand that $a$ is invariant under multiplication
%by $1+\Mm$ in the first $i-1$ coordinates, i.e. $a(m_1x_1,\ldots,m_{i-1}x_{i-1},x_{i+1},\ldots,x_n,y) = a(x_1,\ldots,x_{i-1},x_{i+1},\ldots,x_n,y)$.
%  \item Maps $(x_1,\ldots,x_n,y_1,\ldots,y_l) \mapsto (x_1,\ldots,x_n,y_1,\ldots,y_l,\rv(x_i))$
  \item Maps $(x_1,\ldots,x_n,y_1,\ldots,y_l) \mapsto  (x_1,\ldots,x_n,y_1,\ldots,y_l,\rv(x_i))$
\end{enumerate}    \end{defn}
 
 The above functions are called {\em elementary admissible transformations over $A$};
 a morphism in $\fC_A^0$  generated by elementary  admissible transformations
 over $A$ will be called an {\em admissible  transformation} over $A$.   Taking $l=0$, we see that     all
 $A$-definable additive translations of $\VF^n$ are admissible.

 Analogously, if $Y$ is a given definable set, one defines the notion of a $Y$-family of admissible
 transformations.

If $e \in \RV$ and $T_e$ is an $A(e)$-admissible transformation, then there
exists an $A$-admissible $T$ such that
$\iota_e T_e  = T \iota_e$, 
where $\iota_e(x_1,\ldots,x_n,y_1,\ldots,y_l) = (x_1, \ldots,x_n,e,y_1,\ldots,y_l)$.
This is easy to see for each generator and follows inductively.

 Informally, note  that admissible maps preserve volume for any 
product  
 satisfying Fubini's theorem of   %(a number of copies of)
 translation invariant measures on $\VF$  and counting
 measures on $RV$.

We will now see  that any  $X \subset \VF^n$ is a finite disjoint union of admissible transforms
of pullbacks from $\RV$.   We begin with $n=1$.

 \<{lem}  \lbl{tr1}  Let $\T$ be $C$-minimal with centered closed balls.   Let $X$ be a   definable subset of $\VF$.
                 %%   $\VF \times \RV^n$.  
  Then $X$ is the disjoint
  union of finitely many     definable 
 sets $Z_i$, such that for some admissible transformations $T_i$, and   
 definable subsets $H_i$ of $\RVi^{l_i}$,
 $T_iZ_i =   \{(x,y): y \in H_i, \rv(x)=y_{l_i}\}$

If $X$ is bounded,   $H_i$ is bounded below;
in fact for any $h \in H_i$, $\valr(h) \geq \val(x)$ for some $x \in X$.
 \>{lem}
 
 Here $\VF$ will be considered a ball of valuative radius $-\infty$, points
 as balls of valuative radius $\infty$.
 
 \<{proof} % We first prove the lemma for $X \subseteq \VF$.  In this case  
 We may assume $X$ is a finite union of disjoint balls of the same valuative radius
 $\alpha \in \Gamma \union \{\pm \infty \}$, each minus a finite union of proper sub-balls, 
 since any definable set
 is a finite union of definable sets of that form.   
 
 {\bf Case 1} $X$ is a closed ball.

In this case, by the assumption of centered closed balls,
 $X$ has a definable point $a$.  Let $T(x)=x-a$.  Then 
 $TX \m \{0\}$ is the pullback of a subset of $\G$, 
the  semi-infinite interval $[\alpha, \infty)$ (where $\alpha$ is the valuative radius of $X$.)
  Thus $TX = \rv^{-1}(H)$ where $H = \valrv^{-1}([\alpha,\infty)) \union \{\infty\}$.

{\bf Case 2} $X$ is an open ball.

Let $\bX$ be the surrounding closed ball of the same radius $\alpha$, and
as in Case 1 let $a \in \bX$ be an definable point, $T(x)=x-a$.  If $0 \in TX$
then $TX = \rv^{-1}(H)$ where $H = \valrv^{-1}((\alpha,\infty)) \union \{\infty\}$.  If $0 \notin TX$,
then $TX = \rv^{-1}(H)$ where $H = \rv(TX)$ is a singleton of $\RV$.

{\bf Case 3} $X = C \m F$ is a ball with a single hole, the closed ball $F$.

Let $\beta$ be the valuative radius of $F$.
Let $a \in F$ be a definable point, $T(x)  = x-a$.  Then $TX = \rv^{-1}(H)$,
$H = \valrv^{-1}(I)$, where $I$ is the open interval $(\alpha,\beta))$ of $\G$ in case $C$
is closed, the half-open interval $[\alpha,\beta)$ when $C$ is open.  

{\bf Case 4}  $X = C \setminus \union_{j \in J} F_j$ is a closed ball, minus a finite  union
of maximal open sub-balls.  

As in Case 1, find $T_1$ such that $0 \in T_1X$.  Then $T_1X$ is the union of the
maximal open sub-ball $S$ of radius $\alpha$, with $\rv^{-1}(H)$, where $H = \rv(X \m S)$.
$S$ can be treated as in Case 2.  Here $H$ is a subset of $\valr^{-1}(\alpha)$, consisting
of $\valr^{-1}(\alpha)$ minus finitely many points.

{\bf Cases 3a,4a}: $X$ is a union of $m$ balls (perhaps with holes) of types 1-4 above.
%$C_j$ of size $\alpha$ (with no holes.)

Here we use induction on $m$; we have $m$ balls $C_j$   covering $X$.
Let $E$ be the smallest ball containing all $C_j$.  As we may assume $m>1$, 
 $E$ must be a closed ball; and each $C_j$ is contained
in some maximal open sub-ball $M_j$ of $E$.  By the choice of $E$, not all $C_j$ can be contained in  the same maximal open ball of $E$.   Let $a \in E$ be a definable point,
$T_1 (x) = x-a$.  If $0 \in T_1C_j$ for some $j$, the lemma is true by induction for this $C_j$
and for the union of the others, hence also for $X$.  Otherwise, $F =  \rv(T_1(X))$ is a finite
set, with more than one element.  For $b \in F$, let $Y_b = T_1 X \meet \rv^{-1}(b)$.  
By \lemref{collect}, we can in fact find a definable $Y$ whose fiber at $b$ is $Y_b$.
By induction again, there exists
an admissible transformation $T_b$ such that $T_b(Y)$ is a pullback of the required form.  
Let $T_2(x)= (x,\rv(x))$, $T_3((x,b)) = ((T_b(x),b))$.  Then $T_3T_2T_1$ solves the problem.

 {\bf General subsets of $\VF$}
 
Let $\beta \geq \alpha$ be the least size (i..e greatest element of $\G$)
 such that some ball of radius $\beta$ contains more than one hole of $X$.  Let 
 $\{C_j: j \in J \}$ be the balls of radius $\beta$ around  the holes $W$ of $X$,
 and let  $C = \union_{j \in J} C_j$.
  Then $X = (X \setminus C) \du (C \setminus W)$.  Now $X \setminus C$
 has fewer holes than $X$, so it can be dealt with inductively.  Thus
 we may assume $X  = C \setminus W$;  and any proper sub-ball of 
 $C$ of less than maximal size contains at most one hole of $X$.   We may assume
 the $\{C_j\}$ form a single Galois orbit; so they each contain 
 two or more holes of $X$.   Since these holes are not contained in a proper
 sub-ball of $C_j$,  each  $C_j$ must be closed, 
 and the maximal open sub-balls of $C_j$ separate holes.  Let $D_{j,k}$
 be the maximal open sub-balls of $C_j$ containing a  hole $F_{j,k}$.  Let $\bar{F}_{j,k}$
 be the smallest closed ball containing $F_{j,k}$.  Then 
 $X = (C \setminus \union_{j,k} D_{j,k}) \du \union_{j,k} (D_{j,k} \setminus \bar{F}_{j,k}) \du \union_{j,k} \bar{F}_{j,k}\setminus {F}_{j,k}) $.
 The second summand in this union falls into Case 3a, the first and third (when non-empty) into
 case 4a.  
 
% 
%The case of subsets of $\VF$ is thus proved.  Now in general,
%let  $b \in \RV^n$; then  $X(b) = \{x: (x,b) \in X \}$ satisfies the conclusion of the lemma
% (applied for $\T_{A(b)}$.)  
%So for each $b$, there exists a finite partition   $X(b) = \union_{i} Z_i(b)$, and 
%$\T_A$-
%admissible transformations
%$T_i(b)$, and $H_i(b)$, such that $T_i(b) Z_i(b) = \rv \inv H_i(b)$.  
%Note that $T_i(x,b) := (T_i(b)(x),b)$ is also an admissible transformation.
%By compactness,
%we can partition $X = \union_j Z_j$, and find admissible transformations $T_i$, such
%that for each $b$, $T_i X (b) = \rv \inv  H_i(b)$.  Letting $H_i = \{(h,b): h \in H_i(b) \}$,
%we see that the lemma holds for $X$.  

                                                                                                  \>{proof}
 
  {\bf Remark}  If we allow arbitrary   Boolean combinations (rather than disjoint unions
 only), we can demand in \lemref{tr1} that the sets $H_i$ be finite.  More precisely,     let $X$ be a   definable subset of $\VF$.  Then there
  exist definable
 sets $Z_i$, admissible transformations $T_i$, and finite 
 definable subsets $H_i$ of $\RVi^{l_i}$ such that:
 
   $X$ is a  % good
    Boolean combination of the sets $Z_i$, 
  and $T_iZ_i$ is one of the following:
   
  \begin{enumerate}
\item  $\VF $
\item  $(0) \times H_i$
\item  $b_i \times H_i$, with   $b_i$ a  definable   ball containing $0$
\item     $\{(x,y): y \in H_i, \rv(x)=f_i(y)\}$, for  
some definable function  $f_i: H_i \to \RVi$
%(such that $\val \circ f $ is constant.)
\end{enumerate}

\<{cor} \lbl{cell1}  Let $X \subseteq \VF \times \RV^*$ be definable.  Then there 
exists a definable $\rho: X \to \RV^*$ and $c: \RV^* \to \VF$, $c': \RV^* \to \RVi$, $c'': \RV^* \to \RV^*$
such that every fiber $\rho \inv (\alpha)$
has the form $(c(\a)+ \rv \inv(c'(\a))) \times \{c''(a)\}$.  Moreover
$c$ has finite image.
\>{cor} 

\<{proof}  The finiteness of the image of $c$ is automatic, by \lemref{fmrfr}.  
The corollary is obviously true for sets of the form $\L (H,h) = \{(x,u) \in \VF \times H: \rv(x)=h(u)\}$;
take $\rho(x,u) = (\rv(x),u)$.   If the statement holds for $T X$ where $T$ is an admissible
transformation, then it holds for $X$.  
 If true for two disjoint sets, it is also true for their union (add to $\rho$ a map
 to $\{1,-1\} \subseteq \k^*$ whose fibers are the two sets.)  Hence by \lemref{tr1} is
 is true for all definable sets.  \>{proof}

  \<{cor} \lbl{cell1.5}  Let $\T$ be $\V$-minimal, $X \subseteq \VF$ and let $f: X \to \RV \union \G$ be a definable function.
 Then there exists a definable finite partition of  $X= \union_{i=1}^m X_i$ such that either
 $f$ is constant on $X_i$, or else $X_i$ is   a finite union of balls of equal radius (possibly missing some sub-balls), there is a definable set $F_i$ meeting each of the balls $b$ in a single point, and for $x \in X_i$, letting $n(x)$ be the point of $F_i$ nearest $x$, for some function 
 $H$,
 $f(x) =H( \rv(x-n(x)))$. \>{cor}
 
 \prf  The conclusion is so stated that it suffices to prove it over $\acl(\emptyset)$, i.e. we may
 assume every almost definable set is definable; cf. \secref{naming}.  By compactness
 it suffices to show that for each complete type $p$, $f| p$ has the stated form.  Let $b$ be the intersection of all balls containing $p$.  If $b$ is transitive then by  
  \lemref{rv-transitive} $f|p$
 is constant.  Otherwise by $\V$-minimality $b$ contains a definable point, and so we may assume
 $0 \in b$.  It follows that $\rv(p)$ is infinite.  Thus by \lemref{1gen}, $f$ factors through $\rv$.  \eprf

 \<{prop}  \lbl{tr2}  Let $\T$ be $C$-minimal with centered closed balls,
  and let $X$ be a definable subset of $\VF^n  \times \RV^l$.   
  Then 
 $X$ can be expressed as a finite disjoint union of $A$-definable sets $Z$,  with each $Z$  
  of the following form.  For some 
 $A$-  admissible transformation  $T $,    
 $A$-definable subset  $H$ of $\RVi^{l^*}$, and map   of
 indices
 $\nu \in \{1,\ldots,n\} \mapsto \nu' \in \{1,\ldots,l^*  \}$, 
     $$TZ = \{(a,b):   b \in H, \rv(a_\nu) = b_{\nu'} \, (\nu=1,\ldots,n)  \}$$

  If $X$ projects finite-to-one to  $\VF^n$, then the projection of $H$ to
  the primed coordinates $1',\ldots,n'$ is finite to one.
 
  If $X$ is bounded, then  $H$ is bounded  below  in $\RVi$.

   \>{prop}

\proof By induction on $n$; the case $n=0$ is trivial.  
Let $\pr: X \to \pr X$ be the projection of $X$ to $\VF^{n-1} \times \RV^{l}$,
so that $X \subset \VF \times \pr X$.  

 {     Let $\pr^*(Y) = \{v: (\exists y \in Y) (x,y) \in Y \}$. }
For any $c  \in \pr X$, according to   \lemref{tr1},
we can write  {   $\pr^{*}(c)  =  \dU_{i=1}^k Z_i(c)$, }
where  
 $$T_i(c)Z_i(c)  = \{(a,b): b \in H_i(c), \rv (a) = b_{1'}  \}$$
  for some $A(c)$-admissible $T_i(c)$, $A(c)$-definable $Z_i(c)$,
  and $H_i(c) \subseteq \RV = \RV^{{1'}}$.  We can write
  {     $Z_i(c) = \{x: (x,c) \in Z_i \}$, 
  $H_i(c) = \{x: (x,c) \in H_i \}$   }
  for some  definable $Z_i$ and 
  $H_i \subset \VF^{n-1} \times \RV^{{1'}}$.
 By compactness, as in \lemref{collect}, one can asssume that the $Z_i(c),H_i(c),T_i(c)$ are
 uniformly definable:   there exists a partition
of $\pr X$ into finitely many definable sets $Y$, 
and for each $Y$ families $Z_i,H_i,T_i$ over $Y$ of definable sets and admissible transformations over $Y$,
such that the integer $k$
is the same for all $c \in Y$, and the $Z_i(c),H_i(c),T_i(c)$ are fibers over $c$ of $Z_i,H_i,T_i$.   
  In this case, {   $\pr^{*}(Y) = \dU_{i=1}^k Z_i$.  }
We can express $X$ as a disjoint union of the various $\pr^{*}(Y)$;  so we may as
well assume $\pr X = Y$ and  $X=Z_1$.   Let $T_1$ 
be such that $\iota_c T_1(c) = T_1 \iota_c$.   Then
$$T_1X = \{(a,c,b):  (c,b) \in H_1,  \rv(a) = b_{1'} \}$$

Any   admissible transformation 
   is injective and so commutes with disjoint unions.

 Now by induction, $H_1$ itself is a disjoint union
  $H_1 = \dU_{j=1}^{k'} Z_{j}$,   with     
        $$T_i'Z_i'
         = \{(d,b):   b \in H_i', \rv(d_\nu) = d_{\nu'} \, (\nu=2,\ldots,n)  \}$$ 

Notational remarks:  Here $d=(d_2,\ldots,d_n)$ are the $\VF$-coordinates of $c$
above.      The $'$ depends on $i$ but we will
not represent this notationally

Let $T_i^*(a,d,b) = (a,T_i'(d,b))$, i.e $T_i^*$ does not touch the
first coordinate.  Note that $T_i^*$ also does not move the
$1'$ coordinate, since in general admissible transformations can
only add $RV$ coordinates but not change existing ones.  
Let $$Z_i = \{x  : T_1(x) = (a,d,b),  (d,b) \in Z_i' , \rv(a)=b_{1'} \}$$  
Then (as one sees by applying $T_1$)
$X = \dU_{i=1}^k Z_i$,   and if $T_i = T_i^*T_1$, we have:
$$ T_i Z_i = \{(a,d',b'): (d',b') \in T_i'Z_i', \rv(a)=b'_{1'} \} 
 = \{((a,d',b'): b \in H_i', \rv(a)=b_{1'}, \rv(d_\nu) = b_{\nu'} \}$$

 As for the finiteness of the projection, if $X$ admits a finite-to-one projection to $\VF^n$,
 so does  each $Z$ in the statement of the Proposition, and hence the isomorphic set $TZ$. 
 We have $H \subset \RV^{n+l}$, $\pi: \RV^{n+l} \to \RV^n$,
so $ TZ = \{(a,b,b'): (b,b') \in H,  \rv(a)= b' \}$. 
%.
 For fixed $a$, this yields   an   $a$-definable
finite-to-one map $TZ'(a) = \{b': (a,b,b') \in TZ \} \to \VF^n$.  By  \lemref{fmrfr}),    $TZ'(a)$ is finite.  Now fix $b$ and
suppose $(b,b') \in H$ with $b'$ not algebraic over $b$.  Then
for generic $a \in \rv^{-1}(b)$, $b'$ is not algebraic over $b,a$.
Yet $(a,b,b') \in TZ$ and so $b' \in TZ'(a)$, a contradiction.

The statement on boundedness is obvious from the proof; if 
$X \subseteq \{x: \val(x) \geq -\gamma\}^n \times \RV^m$,
then $H$ is bounded below by $-\gamma$ in each coordinate. \qed

\sssec{A remark on more general base structures}

  \<{lem}  \lbl{tr2p}  Let   
   $\T$ be $\V$-minimal,
    $A$  an $\fB$-generated substructure  of a model of $\T$.  
Let $X$ be a $\T_A$- definable subset of $\VF^n  \times \RV^l$.
Then there exist  $\T_A$-definable subsets $Y_i \subset \RV^{m_i}$ and (projection) maps $f_i: Y_i \to \RV^n$,   a disjoint union  $Z$ of
$$Z_i = Y_i \times _{f_i,\rv} \VF^n $$
and a nonempty $A$-definable family ${\mathcal F}$ of admissible transformations
$X \to Z$.  ${\mathcal F}$ will have an $A'$-point for any $\VF \union \RV \union \G $-generated structure
 containing $A$. \>{lem}.
 
\prf  We may assume $A$ is finitely generated.   By \propref{resolve}  there exists an almost $\VF \union \G$-generated 
 $A' \supset A$  embeddable over $A$ into any  $\VF \union \G$-generated structure containing $A$,
  and with $\RV(A')=\RV(A)$.  By \propref{tr2}, the required objects $Y_i,f_i$ exist
  over $A'$.  But since $\RV$ is stably embedded, this data is defined over $\RV(A') \subseteq A$.  The admissible transformations $X \to Z  = \du (Y_i \times _{f_i,\rv} \VF^n)$
  exist over $A'$; so one can find a definable set $D$ %parameterizing admissible transformations  $X \to Z$
  with an $A'$-point,
  and such that any element of $D$ codes an admissible transformation  $X \to Z$.  \eprf
  
{\bf Remark} In fact, arbitrary ACVF-imaginaries may be allowed here.
%  
% We may therefore  assume that the data $Y_i, \ldots$ exists over $\dcl(A \union \{a\})$,
%where $a$ is a point of an open ball $c$ defined over $A$ and
%containing no $A$-definable sub-balls, and show it can be found already over $A$.
%Let $p_c$ be the type of elements of $c$.  Then $A(d) \meet (\RV \union \G) = A \meet (\RV \union \G)$.  Thus indeed the data $Y_i$ is already defined over $A$.
% % In effect the projections
%of this data to $\G$ may be taken to be $A$-definable, using the Skolemization of $\G$.  But then the choice of the $Y_i,f_i$ themselves 
%is a choice of an element in a stable definable set of parameters; 
%%since the unique type of elements of $c$ over $A$ is $\G$-like,
%any function from elements of $c$ into $P$ must be constant.   This
%shows that the $Y_i,f_i$ can be chosen over $A$; the rest follows from the construction.  \qed

\<{example}   \lbl{tr2m}
${\mathcal F}$ need not have an $A$-rational point. \rm  
For instance if $A$ consists
of an element of $\VF / \Mm$, i.e. an open ball $c$, then we can take
$Y=Y_1$ to be the point $0 \in \RV$ (since $c$ can be transformed
to $\Mm$); but there is no $A$-definable bijection of $c$ with $\Mm$.
\>{example}

\sssec{A statement  in terms of  Grothendieck groups}

 Recall definitions \ref{VFcat}, \ref{RVcat}.  
\<{defn} \lbl{LL}  Define $\L: \Ob \RVni \to \Ob \VFni$ by: 
 
$$\L (X,f) = (\VF^\times)^n \times_{\rv^n,f} X \subset  \VF^n \times \RV^m$$
where $\VF^\times = \VF \m \{0\}$.

For $\bX = \sum_i \bX_i \in \RV[*]$, we let $\L(\bX)$ be the disjoint sum   $\sum_i \L(\bX_i)$
over the various components in $\RV[i]$.

Let $\rho$ denote the natural map $\L  (X,f) \to X$.
\>{defn}

\lemm{tr2-surjectivity} The image of  $\L: Ob \RV[\leq n,\idot] \to \Ob \VFni$ meets every isomorphism 
class of $\VFni$.  \>{lem}
 
 \prf  %This is just the assertion of \propref{tr2}, but stated in terms of $\RVi$.
For $X \subseteq \RV^*$ and $f: X \to \RVi$, define $\rv(0)=\infty$ and 
$$\L (X,f) = \VF ^n \times_{\rv^n,f} X \subset  \VF^n \times \RV^m$$
 Then in the statement of \propref{tr2}, we have 
 $TZ= \L (H,h)$ where    $h$ is the projection to the primed coordinates.
 For $x \in H$, let $s(x)=\{i: h_i(x)=\infty\}$.  For $w \subseteq \{1,\ldots,n\}$,
 let $H_w = \{x \in H: s(x)=w \}$.   Let $\bar{H}_w = (H_w,h'_w)$
 where $h'_w = (h_i)_{i \notin w}$.  Then $\bar{H}_w \in \RV[|w|,\idot]$, and
    $\L(H_w,h | H_w) \iso \L(\bar{H}_w)$.  So
    $\L(H,h) \iso  \L (\sum_w \bar{H}_w)$.  \eprf

\sssec{A restatement in terms of $\VF$ alone}

This restatement will not be used later in the paper. 

\<{defn}  \lbl{fc1}  Let $A$ be a subfield of $\VF$.  
 Let $\fC^1_A(n,l) $ be the category  of definable subsets of
 $\VF^n \times (\VF^\times)^l$, generated by composition and restriction
 to subsets by maps of  one of the following types:
 \begin{enumerate}
  \item Maps 
  $$(x_1,\ldots,x_n,y_1,\ldots,y_l) \mapsto (x_1,\ldots,x_{i-1},x_i+a,x_{i+1},\ldots,x_n,y_1,\ldots,y_l)$$
with $a=a(x_1,\ldots,x_{i-1}, 
y_1,\ldots,y_l): \VF^{i+l-1} \to \VF$
an $A$- definable function of the coordinates $y,x_1,\ldots,x_{i-1}$. 
  
  \item Maps $(x_1,\ldots,x_n,y_1,\ldots,y_l) \mapsto  (x_1,\ldots,x_n,y_1,\ldots,y_{i-1},x_iy_i,y_{i+1},\ldots,y_l): X \to Y$
 assuming $x_i \neq 0$ on $X$, and that this function takes $X$ into $Y$.
\end{enumerate}    \end{defn}
 
\<{remark}  The morphisms in this category are measure preserving with respect to Fubini products of invariant measures (additively
for $\VF$,  multiplicatively for $\VF^\times$), viz.   
  $dx_1 \wedge \cdots \wedge dx_n \wedge {\frac{dy_1}{y_1}}
\cdots \wedge  {\frac{dy_l}{y_l}}$ \>{remark}

 \<{lem}  \lbl{tr2v2}   
 Let $\T$ be $C$-minimal with centered closed balls,
    $X$   a definable subset of $\VF^n$.   
  Then 
 $X$ can be expressed as a disjoint union of  $A$-definable sets $Z$   with   the following property.  For some $l \in \Nn$, 
 there exists an  $\fC^1_A(n,l)$-transformation    $T$
 and a definable subset  $H$ of $\RVi^{n} \times \RV^l$, such that 
     $$T(Z \times (\Mmm)^{l}) =  \rv^{-1}(H)$$

Moreover, the projection of $H$ to $\RVi^n$ is finite-to-one.

If $\val(x)$ is bounded below, then $\val(H)$ may be taken to be bounded below in the
$\RV$-coordinates, and bounded in the $\RVi$-coordinates.
 
   \>{lem}

\<{proof}

\propref{tr2}  \>{proof}

\>{section}

\<{section}{$\V$-minimal geometry:  Continuity and differentiation}

  We work with a $\V$-minimal theory.

 \ssec{Images of balls under definable functions.}

\<{prop} \lbl{image2} Let $X,Y$ be definable subsets of $\VF$, and let $F: X \to Y$ be a definable
bijection.  Then there exists a partition of $X$ to finitely many definable equivalence classes,
such that for any open ball $b$ contained in one of the classes, $F(b)$ is an open ball; 
and dually, if $F(b)$ is an open ball, so is $b$. 
\>{prop}

\proof It suffices to show that such a partition exists over $\acl(\emptyset)$; for any finite
almost definable   partition has a finite definable refinement (cf. discussion of Galois theory in
\secref{types}).  Thus as in \secref{naming} we may assume every almost definable set is named. 

We will show that if $p$ is a complete type, and $b$ is 
an open sub-ball
of $p$, then $F(b)$ is an open ball; and that if $b' $ is an open sub-ball of $F(p)$,
then $b$   is an open ball.  From this it follows by compactness that there exists a definable $D_p$
containing $p$ with the same property; by another use of compactness, finitely many $D_p$
cover $X$; it then suffices to choose any partition, such that any class is contained in some $D_p$.

 When $p$ has a unique solution, the assertion is trivial.  When $p$ is the generic type of
 a closed ball, or of $\VF$, or of a transitive open or $\infty$-definable ball,
 for any $\alpha \in \G$, $p$ remains complete over $<\alpha>$.  
  In the transitive cases, this follows from  \lemref{rv-transitive}, while in the centered closed case it follows
  from   \lemref{orthk}.
 
 Thus
 all open sub-balls $b_t$ of $p$ of any radius $\alpha$ have the same type over $<\alpha>$;
 hence  they are all transitive over $<t>$, where $t \in K /\Mm _\alpha$, where $\Mm_\alpha = \{x: \val(x)>\alpha\}$   (\lemref{transitive}, with $Q=p$.)  
 Thus by \lemref{image1},  $F(b_t)$  is an open ball. 
 
  The remaining case is
 that $p$ is the generic type of a centered open or $\infty$-definable ball $b_1$. 
So $b_1$ contains a definable proper sub-ball $b_0$.
  If $b$ is
 an open sub-ball of $p$, of radius $\alpha$, then $b \meet b_0 = \emptyset$;
  let $\bar{b}$ be the smallest closed ball
 of  containing $b$ and $b_0$.   Then $b$ is contained in the generic type
 of $\bar{b}$, and so by the case of closed balls, $F(b)$ is an open ball.   \qed

\<{remark} \lbl{balldef} When $X \subseteq \VF \times \RV^n$,  by a {\em ball contained in $X$}
we will mean a subset of $X$ of the form $b \times \{e\}$, where $b \in \fB$ and $e \in \RV^n$.
With this understanding, the proposition extends immediately to such sets $X$.  \rm
Indeed for each $e \in \RV^n$, according to the proposition there is a finite partition 
of $X(e)$ with the required property; as in \lemref{collect} these can be patched
to form a single partition of $X$.
\>{remark}
 
%note:  in the proof of image2, do not need that  all open sub-balls $b_t$ of $p$ of radius $\alpha$ have the same type over $A(\alpha)$;
%merely need:  none are algebraic.   
\<{remark}  When $X \subseteq \VF$   there exists a finite  set of
points $F$ (not necessarily $A$-definable) such that $F(b)$ is an open ball whenever $b$ is an open ball disjoint from 
$F$.  (This does not extend to $X \subseteq \VF \times \RV^*$.)  \>{remark}
  Indeed by \propref{image2} there is a finite number of closed and open balls $b_i$ and points,
such that $F(b)$ is an open ball for any open  ball $b$ that is either contained in or is disjoint from each $b_i$.  Now let $c_i$ be a point of $b_i$.  If $b$ is an open ball and no $c_i \in b$,
then $b$ must be disjoint from, or contained in, each $b_i$; otherwise $b$ contains $b_i$,
hence $c_i$.

\ssec{Images of balls, II}

\<{lem}   \lbl{image2l}  Let $X,Y$ be balls, and $F: X \to Y$ a definable bijection taking open balls
to open balls.  Then  
 for all $x,x' \in X$ $$\val(F(x)-F(x')) = \val(x-x')+v_0$$
 where $v_0$ is 
the difference of the valuative radii of $X,Y$.
 \>{lem}  

\<{proof}  Translating by some $a \in X$ and by $F(a) \in Y$,  we may assume $0 \in X, 0 \in Y$, 
$F(0)=0$; and by multiplying we may assume and both $X,Y$ have valuative radius $0$,
i.e. $X=Y = \Oo$.   Let $M(\alpha) = \{x: \val(x)<\alpha\}$.  Then $F(M(\alpha))= M(\beta)$
for some $\beta = \beta(\alpha)$.  $\beta$ is an increasing definable surjection 
from $\{\alpha \in \Gamma: \alpha > 0 \}$ to itself; it must have the form
$\beta(\alpha) = m \alpha$ for some rational $m>0$.  
 By \lemref{seq-d}, we have $m \in \Zz$.  
 Now reversing the roles of $X,Y$ and using $F^{-1}$ will transform $m$ to $m^{-1}$,
so $m ^{-1} \in \Zz$ also, i.e. $m = \pm 1$.  Since $m>0$, we have $m=1$.   \>{proof}

\<{lem} \lbl{image2rd}   Let $X$ be a transitive open or closed ball (or infinite intersection of balls), and $F: X \to Y$ a  definable bijection.  Then there exists definable $e_0 \in \RV$ such that
for  $x \neq x' \in X$,
$\rv(F(x)-F(x')) = e_0  \rv(x-x')$.  \>{lem}

\prf  We first show a weaker statement:

\Claim{}   For some definable $e_0: \G \to \RV$ %with finitely many values,
$\rv(F(x)-F(x')) = e_0(\val(x-x')) \rv(x-x')$ for all    $x \neq x' \in X$.

   \<{proof}   Fix $a \in X$.  For $\delta \in \G$, let $b_\delta=b_\delta(a)$, the closed ball around $a$ of valuative
   radius $\delta$.  Consider those $b_\delta$ with $b_\delta \subseteq X$. 
    As we saw in the proof of \lemref{image2}, as any $a \in X$ is generic, $b_\delta$ is transitive  in $\T_{b_\delta}$.  
  By \lemref{image.},% when $\val(x-a)=\delta$, 
  $\rv(F(x)-F(a)) = f_a(\delta) \rv(x-a)$,  where $\val(x-a)=\delta$, and
    $f_a(\delta)$ is a function of $a$ and $\delta$.  But then $f_a$
  is a function $\G \to \RV$,  so by \lemref{nosec} it takes finitely many values $v_1,\ldots,v_n$.
    Let $Y_i = f_a \inv (v_i)$.     $Y_i$ has a canonical code $e_i \in \G^*$, consisting of the endpoints
    of the intervals making up $Y_i$.
    Using
  the linear ordering on $\G$, each individual $e_i$ is definable from the set $\{e_i\}_i$, and hence
  from $a$; thus $v_i = f_a(Y_i)$ is also definable from $a$.  So $f_a$ 
   $f_a$ is 
  definable from $(e_i,v_i)_i$.  (This last argument could have been avoided by quoting elimination of imaginaries
  in $\RV \union \G$.)  
    However as $X$ is transitive, 
    every definable function $X \to (\RV \union \G)$ is constant, and so  $f_a=f_b$
for any $a,b \in X$.)  Let $e_0(\delta) = f_a(\delta)$.  \>{proof}

 We now have to show that  the function $e_0$ of the Claim is constant.  Using
the O-minimality of $\G$, it suffices to show for any definable $\delta \in \dom(e_0)$

1)  If $e_0(\delta)=e$ then    $e_0(\gamma)=e$ for sufficiently small $\gamma > \delta$.

and if  $\delta$ is not a minimal element of $ \dom(e_0))$, then also:

2)  if $e_0(\gamma) = e$ for sufficiently large $\gamma < \delta$, then 
$e_0(\delta)=e$.

To determine $e_0(\delta)$ it suffices to know $\rv(F(x)-F(x'))$ and $\rv(x-x')$
for one pair $x,x'$ with $\val(x-x')=\delta$.  Thus in 1)  we may replace $X$ by
 a closed sub-ball $Y$ of valuative
radius $\delta$, and in 2) by any     closed sub-ball $Y$ of $X$ of valuative radius $< \delta$.
Since such closed balls $Y$ are transitive (over their code),  we may assume $X$ is a closed ball.

Fix $a \in X$.   Pick a generic $c$ (over $a$)  with $\rv(c) = e$.

To prove 1), note that  type of such $c$ is generic in an open ball, whereas the elements of $X$ are generic in a closed ball; these generic types
 are orthogonal by \lemref{orth-gen}; so $X$ remains transitive in $\T_c$.  Thus we may assume (by passing to $\T_c$) that $c$ is definable.

Let $q_a$ be the generic type of the closed ball $\{x: \val(a-x) \geq \delta\}$.   For   $x \models q_a$, let $v_0 =  \val(F(a)-F(x) - c(a-x)) - \val(c)$.

By the definition of $e$, $\val(F(a)-F(x) - c(a-x)) > \val(F(a)-F(x))$, so 
we have  
\beq{im1}  v_0 + \val(c) = \val(F(a)-F(x) - c(a-x))  > \val(F(a)-F(x))=\val(c(a-x)) = \delta + \val(c)  \eeq
 If $\delta < \gamma < v_0$,
find $x,x' \models q_a$ with $\val(x-x') = \gamma$.  Then $\val(F(x)-F(x') - c(x-x')) \geq v_0 +\valr(e)
   > \gamma + \valr(e) = \val(c(x-x'))$, so $\rv(F(x)-F(x')) = \rv(c(x-x'))$ showing that $e_0(\gamma)=\rv(c)=e$. 
   This proves 1).

For 2),  
let $Q_0 = \{\g: \g< \delta\}$, $Q_0^{def}$ the set of definable
elements of $Q_0$,  and $Q= \{\g \in Q_0:  (\forall y \in Q_0^{def}) ( \g > y) \}$.  
So $Q$ is a complete type of elements of $\G$.   For $\g \in Q$,
according to  \lemref{orthplus}, the formula  $\val(x-a) = \g$ generates a complete type 
  $q_{\g;a}(x)$.  By \lemref{rv-transitive}) $X$ is transitive over $\g$, so the formula $x' \in X$
generates a complete $\T_\g$-type.   Thus by transitivity
a complete $\T_\g$-type $q_\g(x,x')$    is generated by:  $x,x' \in X, \val(x-x')=\g$;  namely 
  $(a,b) \models q_\g$ iff $b \models q_{\g;a}$.  %$x,x' \in X$ and $\val(x-x')=\g$.  

For some definable
$v_0$, for $(a,x) \models q_\g$ we have, as in 1): 
\beq{im2}  \val(F(a)-F(x) - c(a-x)) =v_0(\gamma)+ \val(c) > \gamma +\val(c) \eeq
If we show that  $v_0(\g)   > \delta$
%\beq{im-g}  v_0(\g) + \val(c) > \delta + \val(c) \eeq
 we can finish as in 1).

Now 
$v_0(\g) = m \g + \g_0$ for some definable $\g_0 \in \G$, and some rational $m$.
Letting $\gamma \to \delta$ in \eqref{im2} gives
$ m \delta + \g_0  \geq \delta $.  If $m < 0$ then 
$v_0(\g) = m \gamma + \g_0 > m \delta + \g_0 \geq \delta$ so we are done; hence we may  take $m \geq 0$.

By \lemref{rv-transitive} $\RV(<\emptyset>) = \RV(<a>) $; by \lemref{1gen}, when $x \models q_{\g;a}$, 
$\RV(<a,x>)$ is generated over $\RV(<a>)$ by $\rv(a-x)$. 

 In particular on $q_{\g,a}$, $x \mapsto \rv(F(a)-F(x)-c(a-x))$ 
is a  function of $\rv(a-x)$.  This function lifts $v_0$ to a function on $\RV$; hence
by   \lemref{seq-d}, $m \in \Zz$.   (This, and $m \geq 0$, are simplifications rather than essential points.)
We have:
$$ \val ((F(a)-F(x) - c(a-x) ) (a-x)^{-m}) = \g_0$$

By \lemref{rv-transitive},
 $(\RV \union \G) (<a>) = (\RV \union \G) (<\emptyset>) $.  By \lemref{orthplus},   then $\valr \inv (\g_0) \meet \dcl(a,x) = \valr \inv (\g_0) \meet \dcl(a)$.
Thus $\valr \inv (\g_0) \meet \dcl(a,x) = \valr \inv (\g_0) \meet \dcl(  \emptyset)$.  
So $\rv((F(a)-F(x) - c(a-x) ) (a-x)^{-m}) \in \dcl(\emptyset)$; i.e.
$$\rv((F(a)-F(x) - c(a-x) ) (a-x)^{-m})  = e_1$$
for some definable $e_1$.  As in 1) we may assume there exists a definable $c_1$
with $\rv(c_1)=e_1$.  
%So $\rv((F(a)-F(x) - c(a-x) ) (a-x)^{-m})  = \rv(c_1)$
Thus for $(a,x) \models q_\g$, 
\beq{im3}   \val((F(a)-F(x) - c(a-x)  - c_1 (a-x)^m)) >  \val(F(a)-F(x) - c(a-x)) = v_0(\g)+\val(c) \eeq

Let $x' \models q_{\g,a}$ be generic over $ \{\g,a,x\}$, so in particular 
$\val(x-x')=\val(x-a)=\val(a-x')=\g$.  We have
$$ \val((F(a)-F(x') - c(a-x')  - c_1 (a-x')^m)) >  \val(F(a)-F(x) - c(a-x)) = v_0(\g)+\val(c)$$
Subtracting from \eqref{im3} we obtain:
\beq{im4} \val((F(x')-F(x) - c(x'-x)  - c_1 [(a-x)^m)-(a-x')^m])    > v_0(\g) + \val(c) =  \val(c_1(a-x')^m)   \eeq
But since $(x,x') \models q_\g$, by \eqref{im3}  we have
\beq{im5} \val((F(x)-F(x') - c(x-x')  - c_1  (x-x')^m))    > v_0(\g) + \val(c) =  \val(c_1(x-x')^m)   \eeq
Comparing \eqref{im4}, \eqref{im5} (and subtracting $\val(c_1)$) we see that
$$  \val(  (a-x)^m -(a-x')^m - (x'-x)^m )>   \val((x-x')^m) = \val((a-x')^m)= \val((a-x)^m)$$
Let $u = (a-x')/(x'-x)$; then $(a-x)/(x'-x) = u+1$, $\val(u)=0=\val(u+1)$, and
$\val( (u+1)^m - u^m-1) >0$.  If $U=\res(u)$ we get $(U+1)^m =U^m+1$.  Since the residue
characteristic is $0$ this forces $m=1$.  (Note that $U$ is generic.)  
So $v_0(\gamma)=\gamma+\gamma_0$.

From \eqref{im2}, $\gamma+\gamma_0+\val(c) > \gamma+ \val(c)$, or 
$\gamma_0 > 0$.  But $\delta - \gamma_0 \in Q_0^{def}$, so 
 since $\gamma \in Q$ we have $\g > \delta - \gamma_0$, or
 $v_0(\gamma) = \g + \gamma_0 > \delta$.   As noted below \eqref{im2} this proves the lemma.\eprf
 
 \<{remark}  In ACVF(p,p), the Claim of \lemref{image2rd} remains true, but 
it  is possible for $e_0$ to take more than one value;
consider $x-cx^p$ on a closed ball of valuative radius $0$, where $\val(c)<0$. \>{remark}

\<{lem}  \lbl{image2rr} Let $X$ be a transitive open ball, and let $F: X \to X$ be a  definable bijection.  Then $\rv(F(x)-F(y)) = \rv(x-y)$ for all  $x \neq y \in X$. \>{lem}

\prf This follows from the second assertion in \lemref{image.}, and from \lemref{image2rd}.  \eprf

 At this point, \lemref{image2} may be improved.  
 
 \<{defn} \lbl{nice} Call a function $G$ on an open ball 
{\em nice} if for some $e_0$, for all $x \neq x' \in \pr X$, $\rv(G(x)-G(x')) = e_0 \rv(x-x')$. \>{defn}
 
\<{prop} \lbl{image2+} Let $X,Y$ be definable subsets of $\VF$, and let $F: X \to Y$ be a definable
bijection.  Then there exists a partition of $X$ to finitely many definable classes,
such that on any open ball $b$ contained in one of the classes, $F(b)$ is an open ball,
and $F| b$ is nice.   
\>{prop}

\prf The proof of \propref{image2} goes through verbatim, only quoting
\lemref{image2rd} along with \lemref{image1}.  \eprf

A definable translate of a ball $\rv \inv (\a)$ will be called a {\em basic 1-cell}.  Thus
\corref{cell1} states that every fiber of $\rho$ is a basic 1-cell.  By
 a {\em basic 2-cell} we mean a set of the form: 
 $$X=  \{(x,y): x \in pr X, \rv(y-G(x)) = \a$$
 where $pr X $ is a basic 1-cell, and $G$ is  nice.

\<{cor} \lbl{cell2} 
Let $X \subseteq \VF ^2$ be definable.  Then there 
exists a definable $\rho: X \to \RV^*$ such that every fiber is a basic 2-cell.\>{cor}

\prf  Let $X(a) = \{y: (a,y) \in X \}$.  By \corref{cell1} there exist  an $a$-definable
$\rho_a: X(a) \to \RV^*$ and functions $c,c'$ such that every fiber $\rho_a \inv(\a)$
 is a basic 1-cell $\rv \inv(c'(a,\alpha))+c(a,\a)$.  By \lemref{collect} we can
 glue these together to a function $\rho_1: X \to \RV^*$ with 
 $\rho_a(y)= \rho_1(a,y)$.  Let $\rho_2(x,y) = (\rho_1(x,y), c'(x,\rho(x,y)))$.
Then any fiber $D$ of $\rho_2$ has the form 
$\{(x,y): x \in pr_1 D, \rv(y-G_D(x)) = \a$ where $G_D(x)=c(x,\a)$, $\a$ depending on the fiber $D$.
Combining $\rho_2$ with a function whose fibers yield a partition as in \propref{image2+},
we may assume $G$ takes open balls to open balls (cf. \remref{balldef}).  Now apply \corref{cell1}
to $\pr X$ to obtain a map $\rho': pr X \to \RV^*$ with nice fibers. \eprf

\ssec{Limits and continuity}

We now assume $\T$ is a $C$-minimal theory of valued fields, satisfying 
\aref{V+} \eqref{V+RV}.

Let $V$ be a $\VF$-variety.   By  ``almost all $a$'' we will mean:  all $a$ away from a set of smaller $\VF$-dimension.
%
%By the  {\em $n$- germ}   of a definable function on  a definable set $X$ of     dimension $n$, we mean here:  a definable function  on a definable subset $X'$, where $\dim(X \m X')< n$;
%where two such have the same germ if they agree on some $X''$ with $\dim(X \m X'') < n$.  
%A germ $f: X \to Y$ is {\em invertible} if it is represented by a bijection (from $X'$ to some $Y'$.)  Such germs can be composed.  
 
\<{lem}\lbl{limitexists} Let $g$ be a definable function on a  ball around $0$.  Then
either $\val g(x) \to -\infty$ as $\val(x) \to \infty$, or 
 there exists
a unique $b \in \VF $  such that $b=\lim_{x \to 0, x \neq 0} g(x)$, i.e.   
$$ (\forall \e \in \G)(\exists  \delta \in \G)  (  0 \neq x \ \& \  \val(x)>\delta \implies \val(g(x)-b) > \e$$
   \>{lem}

\<{proof}  
  
  Let $p$ be the generic  type of an element of large valuation; so $c \models p|A$
  iff $\val(x) >  \G(A)$.  
and let $q = tp(g(c)/A)$ where $c \models p | A$.    By \remref{1types},  $q$ coincides  with the generic type of $P$ over $A$
where $P$ is a closed ball, an open ball, or an infinite intersection of balls, or $P=\VF$.
The last case means that $\val g(x) \to -\infty$.  
   The existence of $g$ shows that $p,q$ are non-orthogonal, so 
it follows from  \lemref{orth-gen} that the first case is impossible.

We begin by reducing to the case $P$ is centered.  Assume therefore that $P$ is transitive.  For $b \in P$,  let $q_b = tp(g(c') / A(b))$ where
$c' \models p| A(b)$.  If $q_b$ includes a proper $b$-definable
sub-ball $P_b$ of $P$, or a finite union of such balls, we may take them all to have
the same radius $\alpha(b)$; so $\alpha(b)$ is $b$-definable.  
By \lemref{rv-transitive}, $\alpha$ is constant.   If as $b$ varies there
are only finitely many ball $P_b$,
then $P$ is after all centered.  If not, then there are two disjoint $P_b,P_{b'}$;
but this is absurd since if $c'' \models p| A(b,b')$ then $g(c'') \in P_b \meet P_{b'}$.
Thus $q_b$ cannot include a proper sub-ball $P_b$ of $P$; so $q_b$ is just the generic type of 
$P$, over $A(b)$.  Moving from $A$ to $A(b)$ we may thus assume that $P$ is centered. 

So $P$ is a centered open or infinitely-definable ball; therefore it has a proper definable sub-ball $b$.
If $y \notin b$, write  $\val(b-y)$ for the constant value of $\val(c-y)$, $c \in b$.  
By the definition of a generic type of $P$, $\val(b-g(c)) \notin \G(A)$.  Now
$\val(b-g(c)) \in \G(A(c)) = \G(A) \oplus \Qq \val(c)$ (by \eqref{V+comp} of the definition of $\V$-minimality), and $\val(c)  > \G(A)$;
it follows that $\val(b-g(c))  < \G(A)$ or $\val(b-g(c))  > \G(A)$.  
The first case is 
again the case of $P=\VF$, while the second implies that   $P$ is an infinite intersection
of balls $P_i$, whose radius is not bounded by any element of $\G(A)$.  In other words,
$P = \{b\}$.  Unwinding the definitions shows that $b=\lim_{x \to 0, x \neq 0} g(x)$.  \>{proof}

{\bf Remark}  In reality, the transitive case considered in the proof above cannot occur.

 By an (open, closed) polydisc, we mean a product of (open, closed)  balls.  Let $B$ be a
 closed polydisc.  Let $M \models T$.  Let $b \in B(M), a \in B(\acl(\emptyset))$.
 Write $b \to a$ if for any definable $\g \in \G$, and each coordinate $i$, 
 $\val(b_i-a_i) > \g$.  Let $p_0$ be the   type of elements of $\G$
 greater than any given definable element.  % and $P_0(x) = p_0(\val(x))$.   
  Then \lemref{limitexists} can also be stated thus: given a definable $g$ on a ball $B_0$
  around $0$ into $B$, there exists
 $b \in \dcl(\emptyset)$ such that 
  if $\val(t) \models p$, then $(t,g(t)) \to (0,b)$.
 
 Stated this way, the lemma generalizes to functions defined on a finite cover of $B_0$:
 
 \lemm{limitexists2}  Let $B_0$ be a ball around $0$, and $B$ a closed polydisc, 
 both $0$-definable.  Let $t \in B_0$ have $\val(t) \models p_0$,
 and let $a \in \acl(t)$, $a \in B$.  Then there exists $b \in B$, $b \in \acl(\emptyset)$
 with $(t,a) \mapsto (0,b)$.  \>{lem}
 
 \prf The proof of \lemref{limitexists} goes through.  \eprf

 The following is an analogue
 of a result of Macintyre's for the $p$-adics.  By the {\em boundary} of a set $X$, we mean the closure
 minus the interior of  $X$.
 
\<{lem} \lbl{opensets} \begin{enumerate}
  \item Any definable $X \subseteq \VF^n$ of dimension $n$
 contains an open polydisc. 
  \item Any definable function $\VF^n \to \RV \union \G$ is constant on some open polydisc.
  
 \item The boundary of any definable $X \subseteq \VF^n$ has dimension $<n$.
\end{enumerate} 
\>{lem}
  \<{proof}   
  Given (1), (3) follows since the boundary is definable; so it suffices to prove (1,2).  For a given $n$, (2)
  follows from (1):   by \lemref{vfrvdim}, the fibers of the function cannot all have dimension $<n$.
   
     For $n=1$, (1) is immediate from $C$-minimality.  Assume (1-2) true for $n$ and let $X \subseteq\VF \times \VF^{n}$ be a definable set of dimension $n+1$.  For any $a \in \VF^{n}$ such that
 $X_a = \{b: (a,b) \in X \}$ contains an open ball, let $\g(a)$ be the infimum of all $\g$
 such that $X_a$ contains an open $\g$-ball.   By (2) for $n$, $\g$ takes a constant value $\g_0$ on some polydisc $U$; pick $\g_1 > \g_0$.  Let 
 $X' = \{(u,z) \in X: u \in U, \& (\forall z')(\val(z-z') >  \g_0 \implies (u,z') \in X \}$.
 Then $\dim(X')=n+1$.  
 Now consider the projection $(u,z) \mapsto z$.  For some $c \in \VF$, the
 fiber $X'_c = \{u: (u,c) \in X' \}$ must have dimension $n$.  By induction, $X'_c$
 contains a polydisc $V$.  Now clearly $V \times B_{\g_1}^o(c) \subseteq X$.     \eprf
    
If $x = (x_1,\ldots,x_n), x'=(x_1',\ldots,x_n')$,  write $\val(x-x') $ for $\min \val(x_i-x_i')$.
Say a function $F$ is $\d$-Lipschitz
  at $x$ if whenever $\val(x-x')$ is sufficiently large,
$\val(F(x)-F(x')) > \d + \val(x-x')$.   Say $F$ is locally Lipschitz on $X$ if
for any $x \in X$, for some $\d \in \G$, $F$ is $\d$-Lipschitz at $x$.  

\<{lem}\lbl{continuity}  Let $F: X \subseteq \VF^n \to \VF$ be a definable function.  Then
 $F$ is continuous 
away from a subset $X'$ of dimension $<n$.  Moreover $F$ is 
locally  Lipschitz on $X \m X'$. 
  \>{lem}
\<{proof}  Let $X'$ be the (definable) set of points $x$ where $F$ is not Lipschitz.
 We must show
that $X'$ has dimension $<n$ (in this case, by \lemref{opensets}, the closure of $X'$
has dimension $<n$ too.)    Suppose otherwise.   For $n=1$ the lemma follows from 
Lemmas \ref{image2} and \ref{image2l}.   Let $\pi_i: X' \subseteq \VF^n \to \VF^{n-1}$ be the projection
along the $i$'th coordinate axis.    Let $Y$ be
the set of $b \in \VF^{n-1}$ such that $\pi_i \inv (b)$ is infinite, or equivalently contains a ball;
it is a definable set.  For  $b \in Y $, let 
$$  D_i(b) = \{ x \in \pi_i \inv (b): (\exists \d \in \G)(F |  \pi_i \inv (b)  \text{ is }   \d  \text{-Lipschitz near }  x \}$$
By the case $n=1$, $\pi_i \inv(b) \setminus D_i(b)$ is finite.  Thus if $D_i= \union_{b \in Y} D_i(b)$,
then $\pi_i$ has finite fibers on $X \m D_i$, so $\dim(X \m D_i) < n$.  Let $X^* = \meet_i D_i$, and for $x \in X^*$ let  $\d(x)$ be the infinimum of all such Lipschitz constants $\d$
(for all $n$ projections.)
By \lemref{opensets}, $\d$ is constant on some open polydisc $U \subseteq X^*$.  Let $\d'$
be greater than this constant value.  Then at any $x \in U$, the restriction of 
$F$ to a line parallel to an axis is $\d'$-Lipschitz.  It follows immediately (using the ultrametric inequality) that $F$ is $\d'$-Lipschitz on $U$; but this contradicts the definition of $X'$.  \eprf

\<{remark}  Via  \aref{V+} \eqref{V+RV}, we used the existence
of $p$-torsion points in the kernel of $\RV \to \G$  for each $p$.  In $ACVF(p,p)$ this fails; one can still show that
$F$ is locally logarithmically Lipschitz, i.e. for some rational $\a >0$, for any
$x \in X \m X'$, for sufficiently close $x'$, $\val(F(x)-F(x')) > \d   \val(x-x')$. \>{remark}

 \ssec{Differentiation in $\VF$}
  
Let $F: \VF^n \to \VF$ be a definable function, defined on a neighborhood of $a \in \VF^n$.
We say that $F$ is differentiable at $a$ if there exists a  linear map $L: \VF^n \to \VF$
 such that for any $\gamma \in \G$, for large enough $\delta \in \G$, 
if $\val(x_i) > \delta$ for each $i$, $x=(x_1,\ldots,x_n)$, then 
$\val(F(a+x)-F(a)-Lx) > \delta + \gamma$.   If such an $L$ exists it is unique, and we denote it 
   $dF_a$.

\<{lem} \lbl{diffn}  Let $F: X \subseteq \VF^n \to \VF^m$ be a definable function.  Then each partial derivative is
defined at almost every $a \in X$. \>{lem}

\<{proof}   We may assume $n=m=1$. 
%% Moreover it suffices to show, by the usual proof, that  the partial derivatives
%exist and are continuous at almost every point.    As every definable function is continuous
%at almost every point, we reduce to $n=1$.  In this case let 
Let
$g(x) = (F(a+x) - F(a)) / x$.  By \lemref{image2l}, for almost every $a$, for some $\delta \in \G$,
for all $x$ with $\val(x)$ sufficiently large, $\val (F(a+x)-F(a)) = \delta+\val(x)$;
so $\val g(x)$ is bounded.
 By  \lemref{limitexists}, and \propref{image2},  $g(x)$ approaches a limit $b \in \VF$
as $x \to 0$ (with $x \neq 0$); the lemma follows.
  \>{proof}

\<{cor}    \lbl{diffble}  Let $F: \VF^n \to \VF$ be a definable function.  Then $F$ is continuously differentiable  
away from a subset of dimension $<n$ \>{cor}
\<{proof} $F$ has partial derivatives almost everywhere, and these are continuous almost
everywhere, so the usual proof works.  \>{proof}

%		We will later define a category $\VF [ n]$, whose objects are  definable sets   $X \subseteq \VF^n \times \RV \st$, such that the projection $X \to \VF^n$ has (bounded) finite fibers.  
%		We will see that $\VF[ n]$ is equivalent to the category of definable subsets of
%		$n$-dimensional varieties.
%		 We wish already to discuss Jacobians and more generally partial derivatives
%		of morphisms of $\VF [ \leq n]$.  

\<{lem}  \lbl{vf-local} Let 
$X \subseteq \VF^n \times \RV ^m$    be  definable, $pr: X \to \VF^n$ the projection.
Then 
  for almost every $p \in \VF^n$, there exists an open neighborhood $U$ of $p$ and
   $H \subseteq \RV^m$ such that $pr \inv (U) = U \times H$.  
   If $h: X \to \VF$, then for almost all $x \in X$, $h$ is differentiable with respect to each $\VF$-coordinate.
    \>{lem}

 \<{proof} For $x \in \VF^n$, let $H(x) = \{h \in \RV^m: (x,h) \in X \}$.  
 By \corref{rvei},  \lemref{eilem},
 there exists $H' \subseteq \RV^m \times \RV^l \times \G^k$ such that for any
 $x \in \VF^n$, there exists a unique $y=f(x) \in \RV^l \times \G^k$ with $H(x)=H'(y)$.
 By \lemref{opensets}, $f$ is locally constant almost everywhere.  Thus for almost all $x$, for some neighborhood $U$ of $x$, for all $x' \in U$,  $H(x)=H(x')$; so $pr \inv U = U \times H(x)$.
The last assertion is immediate. 
  \>{proof}

This
allows the definition of partial derivatives of a definable map $F: X \to \VF$ (almost everywhere);
we just take them with respect to the $\VF$-coordinates, ignoring the $\RV$-coordinates.

Given  
$h: X \to \VF^n, h': X'  \to \VF^n$ with $\RV$-fibers,  and a definable map $F: X \to X'$,
we define the partials of $F$ to be those of $h' \circ F$.  
 Then   the differential $dF_x$ exists at almost every point $x \in X$ by \corref{diffble}, and we denote the determinant by $\JVF$, and refer
to it as usual as the Jacobian.

\<{defn} \lbl{vfvol} Let $X,X' \in \VFni$ and let $F: X \to X'$ be a definable bijection.  $F$ is {\em  measure preserving} if $\rv \JVF (x)  = 1$ for almost all $x \in X$.  $\VFvol[n,\idot]$ is the subcategory of $\VFni$     with the same objects,
and whose morphisms are the  measure preserving morphisms of $\VFni$.

 Let $\VFvol$ be the category whose objects are those of $\VFni$, and whose 
 morphisms $X \to Y$ are  the essential bijections $f: X \to Y$ that are measure preserving.

 \>{defn}

\ssec{Differentiation and Jacobians in $\RV$.}

%\ssec{Jacobians of definable maps on $\RV$}
\lbl{RVjacobian}

Let $X, Y$ be definable sets, together with finite-to-one
definable maps $f_X: X \to \RV^n$, $f_Y: Y \to \RV^n$.  Here $X,Y$ can be subsets
of $\RV^{*}$ or of $\RV^* \times \VF^*$, etc.; the notion of Jacobian will not depend on the particular realization 
of $X,Y$.  

Let $h: X \to Y$ be a definable map.  

The notion of Jacobian will depend not only on $h,X,Y$ but also on $f_X,f_Y$;
to emphasize this we will write  $h: (X,f_X) \to (Y,f_Y)$.

We first define  smoothness.  When  $A= f_X(X),B=f_Y(Y) $ are definable subsets of  $\k^n$, we say $h,X,Y$ are smooth if 
$A,B$ are Zariski   open,   
$\{(f_X(x),f_Y(h(x))): x \in X \}  \meet (A \times B) = Z$ for some 
  nonsingular Zariski closed set $Z \subset A \times B$, and
  the differentials of the projections to $A$ and to $B$ are isomorphisms
  at any point $z \in Z$.  In this case, composing the inverse of one of these
  differentials with the other, we obtain a linear isomorphism $T_a(A) \to T_b(B)$
  for any $a= f_X(x), b=f_Y(h(x))$; since $T_a(A) = k^n = T_b(B)$, this linear
  isomorphism is given by an invertible matrix, whose determinant is the Jacobian $J$.

  In general, to define smoothness of $X,Y$ at $(x,y=h(x))$, we restrict to the cosets
  of $(\k^*)^n$ containing $x$ and $y$, translate multiplicatively by $x$ and $y$
  respectively,  and pose the same condition.  

Any $X,Y, h$ are smooth  outside of a set $E$ where
$E \meet C$ has dimension $<n$ for each coset $C$ of $(\k^*)^n$.
Equivalently (by \lemref{rvd1}) $E$ has $\RV$-dimension $<n$.  

Assume now that $X,Y,h$ are smooth.   Define
$$ \JRV(h)(q) = \Pi(f_X(q)) \inv \Pi(f_Y(q')) J(1, 1) \in \RV$$  
where $\Pi(c_1,\ldots,c_n) = c_1 \cdot \cdots \cdot c_n$.

%\<{remark}\lbl{rv-jacobian}  The valuation $\val \JRV$ of $h$ depends directly on the values of $h$, rather than on derivatives.  We can thus define $\val \JRV(h)(q) = \sum \val f_Y(q') - \sum \val f_X(q) $, {\em even if $h$ is defined only at the point $q$.}  \>{remark}

At times it is preferable not to use a different translation at each point of a coset of $(\k^*)^n$.
The Jacobian $\JRV(h)$ 
of $h$ at $q \in X$ can also be defined as follows.   Let $q'=h(q)$, $\g = \valr(q), \g'=\valr(q') \in \G^n$.
  Pick any $c,d \in \RV^n$ with $\valr(c) = \g, \valr(d)=\g'$ (one can take
  $c=f_X(q), d=f_Y(q')$. )  Let 
  $$W(\g,\g')= \{a: f_X(a) \in \valr^{-1} (\gamma), f_Y(h(a)) \in \valr \inv (\g ')) \}$$
$$H' = \{ (c \inv f_X(a), d\inv f_Y(h(a))): a \in W \}$$  
Since $f_X,f_Y$ are finite-to-one,   $H' \subset (\km^n)^2$ both projections of $H'$ to
$\km^n$ are  finite-to-one, and $H'$ is nonsingluar by the smoothness of $(X,Y,h)$.    
We can thus define the Jacobian $J'$ of $H'$ at any  point.
We have:
$$ \JRV(h)(q) =  \Pi(c) \inv \Pi(d) J'(qc \inv, q' d \inv) \in \RV$$

We also define $\JG(h)(q) = \sum \g' - \sum \g \in \G$ (writing $\G$ additively).  {\em Note that this depends only on the value of $h$ at $q$.}  We have:

$$\valr \JRV(h)(q) = \JG(h)(q) $$

\<{example} Jacobian of maps on $\G$.  \rm  If $\bar{X},\bar{Y} \subset \G^n$,
we saw that a definable map $\bar{f}: \bar{X} \to \bar{Y}$ lifts to $\RV$
iff it is piecewise given by an element of $GL_n(\Zz)$ composed with a translation.  Assume
$\bar{f}$ is given by a matrix $M \in GL_n(\Zz)$,   let 
 $X = \valr^{-1}(\bar{X}), Y = \valr^{-1}(\bar{Y})$, and let $f: X \to Y$ be given
 by the same matrix, but multiplicatively.  Then $X,Y,f$ are smooth, and   
 $$J(f)(x) = \Pi(y) \Pi(x)^{-1} \det M$$
 where $y = f(x)$, and $\det (M) = \pm 1$.  \>{example}

%
%We view $\G$ as discrete, and define topological notions on $\RV^m$ using
%the Zariski topology on $(\k^*)^m$, so as to be translation invariant.  Thus a definable
%subset $Y \subseteq \RV^m$ is {\em closed} if for any $c \in \RV^m$, $c Y \meet (\k^*)^m$
%is closed.  An  {\em open neighborhood} of $p$ in $Y$ is a 
%  translate $W=pU$ of a Zariski open neighborhood of $1$ in $ (\k^*)^n \meet p \inv Y$.
%By a {\em regular function} on a neighborhood of $p$ into $\RV$, we mean a 
%definable function $f: W \to \RV$
%such that $f(px) = ch(x)$,  where $c \in \RV$ and
%$h: U \to \k^*$ is a regular function.

%(Note that for any definable function $f$ into $\RV$, $\val f$ takes a unique value on a sufficiently small open neighborhood of $p$. )

%A closed definable set $Y  \subseteq \RV^m$ is a {\em smooth open} if each point
%has a neighborhood that translates to a  nonsingular variety in $(\k^*)^m$.  
%   The partial derivatives are then defined by
%$\partial_i f (y) = cp_i ^{-1}  \partial_i h (p \inv y) $.  

 \sssec{Alternative: $\G$-weighted polynomials}

 \lbl{weighting}
 We have seen that the geometry on $\valr \inv (\g)$ ($\g \in \G^n$) translates
 to the geometry on $(\k^*)^n$, but this is true for the general notions and not for specific varieties; a definable subset of $C(\g) = \val \inv(\g)$ does not correspond canonically to any
definable subset of $\val \inv (0)$.  An invariant approach is therefore useful.
  Let $\G_0 = \G(<\emptyset>)$.
$X=(X_1,\ldots,X_n)$ be variables,
 $\g = (\g_1,\ldots,\g_n) \in \G_0^n$, and let $\nu=(\nu(1),\ldots,\nu(n)) \in \Nn^n$ denote a multi-index.  By a $\g$-weighted monomial we mean an expression
$ a_\nu X^\nu$ with $a_\nu$ a definable element  of $\RV$, such that $\valr (a_\nu )+ \sum  \nu(i)  \g_i= 0 \in \G$.  Let $Mon(\g,\nu)$ be the set of $\g-$-weighted monomials of exponent $\nu$,
together with $0$.  Then $Mon(\g,\nu) \m \{0\}$ is a copy of $\valr \inv (- (a_\nu )+ \sum  \nu(i)  \g_i)$; so $Mon(\g,\nu)$ is a one-dimensional $\k$-space.  In particular addition is
defined in $Mon(\g,\nu)$.  We also have a natural multiplication
$Mon(\g,\nu) \times Mon(\g,\nu') \to Mon(\g, \nu+\nu')$.  Let 
$R[X;\g] = \oplus_{\nu \in \Nn^n} Mon(\g,\nu)$.  This is a finitely generates graded  $\k$-algebra.
It may be viewed  as an affine coordinate ring of $C[\g]$; but  the 
ring of the product $C[\g,\g']$ is $R[X,X';(\g,\g')] $, in general a bigger ring than
$R[X,\g] \tensor_\k R[X',\g']$.   Nevertheless, a Zariski closed subset of $C(\g)$
corresponds to a radical ideal of $R[X';\g]$.   In this way notions such as smoothness
may be attributed to closed or constructible subsets of any $C(\g)$ in an invariant way.

\<{defn} \lbl{rvvol} Let $X,Y \in  \Ob \RVni$ and let $h: X \to Y$ be a definable bijection.
$h$ is {\em measure preserving} if $\JG h(x) = 0$ for all $x \in X$, and
$\JRV h (x) = 1$ for   all $x \in X$ away from a set of $\RV$-dimension $<n$.  
If only the first condition holds, we say $h$ is $\G$-measure preserving.

For $X,Y \in \RVleqn$, we say $h: X \to Y$ is measure preserving
if this is true of the $\RV[n]$-component  of $h$.

 $\RVvol[n,\idot]$ (respectively  $\RVvolg[n,\idot]$)  is the subcategory of $\RVni$ with the 
same objects, and whose morphisms are the measure preserving 
(respectively $\G$-measure preserving) definable bijections.  

$\RVvol[\leq n ,\idot] = \oplus_{k <n} \RVvolg[k,\idot] \oplus \RVvol[n,\idot]$.
\>{defn}

Note that when $X,Y \in \Ob \RVni$,  a bijection $h: X \to Y$ is $\G$-measure preserving iff it  leaves invariant the
sets $S_\g = \{(a_1,\ldots,a_n): \sum_{i=1}^n \valr(a_i) = \g \}$.

\ssec{Comparing the derivatives}  Consider a definable function $F: \VF \to \VF$ 
lying above $f: \RV \to \RV$, i.e. $\rv F = f \rv$.    
The fibers of the map $\rv: \VF \to \RV$ above $\k$, for instance,  are open balls of valuative radius $0$, whereas 
the derivative is defined  on the scale of balls of radius $r$ for $r \to + \infty$.  Thus
the comparison between the derivatives of $F$ and $f$  is not tautological.  Nevertheless one obtains the expected relation almost everywhere.

While this case of the affine line would suffice (using the usual technique of partial derivatives), it is easier to place oneself in the more general context
of curves.  More precisely we consider definable sets $C$ together with 
finite-to-one maps $f: C \to \RV$.  Let $\L C$ and $\rho: \L C \to C$ be as above.

 In the following lemma, $H',h'$ denote respectively 
 the $\VF,\RV$-derivatives of functions $H,h$
 defined on objects of $\VF[1],\RV[1]$ respectively.
 
\<{prop}\lbl{rv-compat}   Let $C_i \subseteq \RV^*$ be definable sets, $f_i: C_i \to \RV$ finite-to-one
definable maps ($i=1,2$).  Let $h: C_1 \to C_2$ be a definable
bijection, and let $H:   \L C_1 \to \L C_2$ be a lifting of $h$, i.e. $\rho H = h \rho$.  Then
\begin{enumerate}
  \item  For  all but finitely many $c \in C_1$, $h$ is differentiable at $c$, 
  $H$ is differentiable at any  $x \in \L c$,
and $\rv H'(x) = h' (\rv(x))$.  
  \item For all $c \in C_1$, %for generic $x \in \L c$, 
  $H$ is differentiable at a generic $x \in \L c$, and $\val H'(x)= (\valr h')(x) = \val(f_2(h(x)))-\val(f_1(x))$.
 % In particular if $\val H'$ is constant on each $\Lc$, then for all but finitely many $x$,
\end{enumerate}
\>{prop}

\<{proof}  (1)  Let $Z'$ be the set of $x \in \L C_1$ such that $H$ is not differentiable at $x$ (a finite set) or that $\rv (H'(x)) \neq h' (\rv(x))$.  We have to show that $\rho(Z') \subseteq C$ is finite,
or equivalently that $f_1 \circ \rho (Z')$ is finite.
Otherwise there exists $c \in \rho(Z')$ with $c \notin \acl(A)$.  By 
 \lemref{1gen}, the formula $\rv(x)=f_1(c)$ generates a complete type 
 $q$ over $A(c)$; it defines a transitive open ball $b_c$ over $A(c)$.  
 Since $\rho \circ H = \rho \circ h$, we have
 %$r \inv (c) = \{c\} \times b_c$; since 
   %$h(c)=c$,  $H(r \inv (c)) = r \inv (c)$, i.e. 
   $H(c,y) = (c,H_c(y))$ for some $A(c)$-definable
  bijection  $H_c$ of $b_c$.     By \lemref{image2rd}, for some $e_0 \in \RV$, $\rv(H(u)-H(v)) = e_0 \rv(u-v)$ for
  all   $u,v \in b_c$; so $\rv( (H(u)-H(v))/(u-v)) = e_0$.  Since $H$ is differentiable almost everywhere on $b_c$
  (\lemref{diffble})
   and $b_c$ is transitive, it is differentiable at every point.
    Clearly $\rv H'(u)=e_0$, contradicting the definition of $Z'$.  
  
    (2) follows from \lemref{image2l} \>{proof}

\<{cor} \lbl{rv-partial}  Let   
$\bX \in \Ob \RV[n]$, %$X$ a definable subset of $ \L \bX$, 
$F: \L \bX \to \VF^n$ a definable function, $f: \L \bX \to \RV^n$ a definable function.
Assume $\rv F(x) = f (\rv(x))$.  Then \propref{rv-compat} applies   
for each partial derivative  of $F$.  In particular, \begin{itemize}
 \item  For  all  $c \in X$ away from a set of smaller dimension, for  all   $x \in \L c$, $F$ is differentiable at $x$,
$f$ is differentiable at $c$, and $\rv \JVF(F)(x) = \JRV(f) (x)$.  

   \item For all $c \in X$, for generic $x \in \L c$, 
  $F$ is differentiable at $x$, and $\val  \JVF(F)(x)= (\JG f)(x) $ 
\end{itemize}
\qed \>{cor}  

\<{cor}  Let $\bX,\bY \in \Ob \RV[\leq n]$, $f \in \Mor_{\RV[\leq n]} (\bX,\bY)$, 
$F \in \Mor_{\VFvol[n]}( \L \bX, \L \bY)$.  Assume $\rv F(x) = f (\rv(x))$.
Then $f \in \Mor_{\RVvol[n]} (\bX,\bY)$. \qed \>{cor}

\proof Follows from \corref{rv-partial}.  
 
% 
%Moreover, if $f: X \to Y$ is an invertible germ, then 
%$df: X \to Hom(TX,TY)$ is a well-defined germ, as is the determinant  $J(f): X \to Hom(\Omega^nTX,\Omega^nTY)$. 

% (To define $df$, let $\bX,\bY,F$ be the Zariski closures
%of $X,Y$ and the graph of $F$ respectively; let $\pr_X: F \to \bX, \pr_Y: F \to \bY$ be the
%projections; and let $df = d \pr_Y  d \pr_X \inv$.)
%
%Define an $\RV$-bundle $\Omega^n TX / \Mmm \Omega^n TX$, so that the fibers (minus zero) are $\RV$-principal homogeneous spaces. 
%Then $J(f)$ induces a map $\Jrv: \Omega^n TX / \Mmm \Omega^n TX \to \Omega^n TY / \Mmm \Omega^n TY$.  
% 
%Given a trivialization of $TX,TY$, or at least of $\Omega^nTX,\Omega^nTY$  (a definable volume form),  
% we can view $J(f): X \to \VF$, and $\Jrv =   \rv \circ J: X \to \RV$.  For the latter it suffices to 
% have  an $\RV$-volume form, i.e. a trivialization
% of $\Omega^n TX / \Mmm \Omega^n TX $.   
% When $X \subseteq \VF^n \times \RV^m$, we take a standard trivialization. 

%

\>{section}%{\RV}

\<{section}{Lifting functions from $\RV$ to $\VF$}

\<{prop}\lbl{rvlift}  Let $\T$ be an effective $\V$-minimal theory.
    Let $X \subset \RV^k$ be definable and
let  $\phi_1,\phi_2: X \to \RV^n$ be two definable maps with 
 finite fibers.  Then there exists a definable bijection 
 $F: X \times_{\phi_1,\rv}(\VF^\times)^n \to X \times_{\phi_2,\rv} (\VF^\times)^n$, 
 commuting with the natural projections to $X$.
     \>{prop}

\<{proof}  Let $A=\dcl(\emptyset) \meet (\VF \union  \G)$.   
If $b \in \dcl(\emptyset) \meet \RV$, then viewed as a ball $b$ has
a point $a \in A$; since the valuative radius of $b$ is also in $A$, 
we have $b \in \dcl(A)$.   Thus  $\phi_1,\phi_2,X$ are $\ACVF_A$-definable.  
Any $ACVF_A$-definable bijection $F$ is a fortiori $\T$-definable; so 
the Proposition for $\ACVF_A$ implies the Proposition for $\T$.  
Moreover $\ACVF_A$ is $\V$-minimal and  effective,
since any  algebaic ball of $\ACVF_A$ is $\T_A$-algebraic and hence has 
a point in $\VF(A)^{alg}$.  
Thus 
we may assume
$\T = ACVF_A$.

 The proof will be asymmetric, concentrating on $\phi_1 X$.  

We may definably partition $X$, and prove the proposition on each piece.  
   
   Consider first the case where $\phi_1: X \to U$ and $\phi_2: X \to V$ are bijections
to definable subsets $U,V \subseteq  (\k^*)^k$.
 Our task   is to lift the bijection $f=\phi_2 \phi_1 \inv$ to $\VF^n$.
   A definable subset of $\k^n$ (such as $\phi_i(X)$) is a disjoint union of smooth varieties. 
 We thus consider a definable bijection $f: U \to V$
between $\k$ varieties $U \subset \k^n$ and  $V \subset \k^n$. 
  Induction on $\dim(U)$ will allow us to remove
 a subset of $U$ of smaller  dimension.  Hence we may assume $U$
 is smooth, cut out by $h=(h_1,\ldots,h_l)$, $TU = \Ker(dh)$,  $f = (f_1,\ldots,f_n)$   where $f_i$ are regular on $U$ (defined on
an open subset of $\k^n$), and $df$ is injective on $TU$ at each
point of $U$.    Thus the common kernel of $dh_1,\ldots,dh_l,df_1,\ldots,df_n$ equals $0$.

   It follows  that at a generic point of $U$ (i.e. every point outside a proper subvariety), if $Q$ is a sufficiently generic $n \times l$
matrix of elements of $A$ (or integers), and we let $f_i' = f_i + Qh$,
then the common kernel of $df_1',\ldots,df_n'$ vanishes.  
Note that $f_i | U = f_i' | U$.  Let $W$ be a smooth variety
contained in $f(U)$ and whose complement in $f(U)$ is a constructible set of dimension smaller
than $\dim(U)$.  Replacing $U$ by $f^{-1}(W)$, we may assume
$f(U)$ is also a smooth variety.

Let $\tilde{U} = \res^{-1}(U)$. 
Lift each $f_i'$ to a polynomial $F_i$ over
$\Oo$, with definable coefficients.  This is possible by effectiveness of $\ACVF_A$.
Obtain a regular map $F$, whose Jacobian is invertible
at points of $\tilde{U}$.    We have $\res \circ F = f \circ res$.
Since $f$ is 1-1 on $U$, the invertibility of $dF$ implies that $F$
is 1-1 on $\tilde{U}$.  Moreover, by  Hensel's lemma,
$F: \rv^{-1}(U) \to \rv^{-1}(W)$ is bijective. 

  Next consider the case where
  in place of a bijection $f: U \to V$ we have a finite-to-finite
correspondence $\tf \subset U \times V$ (where 
   $U=\phi_1(X), V=\phi_2(X)$), $\tf=\{(\phi_1(x),\phi_2(x)): x \in X\}$.  We may 
   take $\tf \subset U \times V$ to be a subvariety, unramified
   and quasi-finite over $U$ and over $V$; and we can take $U,V$
   to be smooth varieties.    As before  we can lift $\tf$ to a correspondence
$\tF \subset \tU \times \tV$, such that $\tF \meet \rv^{-1}(u) \times \rv^{-1}(v)$ is a bijection $\rv^{-1}(u) \to \rv^{-1}(v)$ whenever $(u,v) \in \tf$.
It follows that  a bijection
 $X \times_{\phi_1,\rv}(\VF^\times)^n \to X \times_{\phi_2,\rv} (\VF^\times)^n$
 is given by:  
$(x,y) \mapsto (x,y')$ iff $(y,y') \in \tF$.  
 
If  $\phi_1: X \to U$ and $\phi_2: X \to V$ are bijections
to definable subsets $U,V$, each contained in a single coset of $(\k^*)^k$, 
say $U \subseteq C(\g), V \subseteq C(\g')$ for some $\g,\g' \in \G^k$
(cf.   \secref{weighting}.)  
   Let $Z=(Z_1,\ldots,Z_k)$ be
variables, $R[Z;\g]$ be the subring of
$\VF[Z]$ consisting of polynomials $\sum a_{\nu} Z^{\nu}$, with
$val(a_{\nu}) + \sum_{i=1}^k  \nu(i) \g_i = 0$, and $a_{\nu}$ a definable element of $\VF$.
There is a natural homomorphism $R'[Z;\g] \to R[Z;\g]$, where
$R[Z;\g]$ is the coordinate ring of $C(\g)$.  By effectivity, this homomorphism
is surjective.  The proof now proceeds in exactly the same way as above.

 This proves the proposition in case $\valr \phi_i(X)$ consists of one point.  
 
Next assume $\valr \phi_2$ consists of one point, and  $\valr \phi_1(X)$ is finite.  So $\phi_1(X)$ lies in
the union of finitely many cosets $(C(a): a \in E)$,
with $E$ finite.  

 For $a \in E$,
$A(a)$ remains almost $\VF,\G$-generated; since the proposition
is true for  $\phi_1 \inv C(a)$ (definable in $\T_{A(a)}$), then  
so by the one-coset case an appropriate isomorphism $F$ exists;
and the finitely many $F$ obtained this way can then be glued
together, to yield a map defined over $A$.  

The case of $\valr \phi_1,\valr \phi_2$ both finite is treated similarly.

 This proves the existence of a lifting in case $\valr \phi_i(X)$ is finite.  Now for
 the general case.
   
 \Claim{}  Let $P \subset X$ be a complete type.  Then there exists a definable
 $D$ with $P \subset D \subset X$ , and definable functions
 $\theta$ on $\valr(\phi_1(D))$ and $\theta'$ on $\valr(\phi_2(D))$ such that  for $x \in D$, 
 $\theta(\valr(\phi_1(x))) = \valr \phi_2(x)$,
  $\theta'(\valr(\phi_2 (x))) = \valr \phi_1(x)$.  
  
\proof  Let $a \in P$, $\g_i = \valr(\phi_i(a))$.  Then 
$\g_2 $ is definable over some points of $\phi_1^{-1} \valr^{-1}(\g_1)$.
But $\valr^{-1}(\g_1)$ is a coset of $\k^*$, and $\phi_1$ is finite-to-one,
so  $\phi_1^{-1} \valr^{-1}(\g_1)$ is orthogonal to $\G$. 
 Thus
$\g_2$ is algebraic over $\g_1$.  Since $\G$ is linearly ordered,
$\g_2$ is definable over $\g_1$; so $\g_2 = \theta(\g_1)$
for some definable $\theta$.  Similarly $\g_1 = \theta'(\g_2)$.
Clearly $\theta$ restricts to a bijection $ \valr \phi_1 P \to \valr \phi_2P$,
with inverse $\theta'$.  By \lemref{11nd} there exists a definable
$D$ with $\theta \phi_1 = \phi_2 , \phi_1 = \theta' \phi_2$ on 
$D$.  

Now by compactness, there exist finitely many
$(D_i,\theta_i,\theta_i')$ as in the claim with $\union _i D_i = X$.  
We may cut down the $D_i$ successively, so we may assume the union
is disjoint.  But in this case the proposition reduces to the case of each
individual $D_i$, so we may assume $X=D$.   Let $B_i = \valr \phi_i(X)$.
Given $b \in B_1$, let $X_b = (\valr \phi_1)^{-1}(b)$.  Then by the
case already considered there exists an $A(b)$-definable
$F_b: X_b \times_{\phi_1,\rv}(\VF^\times)^n \to X_b \times_{\phi_2,\rv} (\VF^\times)^n$.
Let $F = \union_{b \in B_1} F_b$.   By \lemref{collect},
$F:   X \times_{\phi_1,\rv}(\VF^\times)^n \to X \times_{\phi_2,\rv} (\VF^\times)^n$
is bijective. (see discussion in \secref{varconst}.)

                                                                                                     \>{proof}

We note a corollary:

 \lemm{red1c2}  Let $\T$ be $\V$-minimal and effective, and let
 $A$ be an almost $(\VF,\G)$-generated structure.    Then $A$ is  effective. \>{lem}
 
 \prf  By \lemref{rveffective} it suffices to show $A$ is $\rv$-effective.
 Note that if $A \subseteq \acl(\emptyset)$, then $\T$ is $\rv$-effective iff $\T_A$ is $\rv$-effective (see proof of \lemref{average} $(2 \to 3)$.)  Thus it suffices to show that if $A_0 = \acl(A_0)$, $a \in \VF \union \G$,
 and $\T'=\T_{A_0}$ is effective then so is $\T'(a)$.  The case $a \in \G$
 is included in \corref{red1c0}, so assume $a \in \VF$.  Let $P$ be the intersection of all
 $A_0$-definable balls containing $a$.  If $P$ is transitive over $A_0$, then by \lemref{rv-transitive} we have $\RV(A_0(a))=\RV(A_0)$, so $\rv$-effectivity remains true trivially.
 Otherwise $P$ is centered over $A_0$, hence has an $A_0$-definable point,
 and by translation we may assume $0 \in P$.  $a$ is then a generic point of $P$ over $A_0$.
 Let $c \in  \RV(A_0(a))$; we must show that $\rv \inv c$ is centered over $A_0(a)$.  
 By  \lemref{1gen}, if $c \in  \RV(A_0(a))$ then $c=f(d)$ for some
 $A_0$-definable function $f: \RV \to \RV$, where $d = \rv(a)$. By \lemref{rvlift}
 there exists an $A_0$-definable function $F: \VF \to \VF$ lifting $f$.  Then 
 $F(d) \in \rv \inv (c)$.  \eprf

\sssec{Base change:  summary}  \lbl{base-change}

Base change from $\T$ to $\T_A$ preserves
 $\V$-minimality, effectiveness and being resolved,
if $A$ is $\VF$-generated; $\V$-minimality and effectiveness, if $A$ is $\RV$-generated;
$\V$-minimality, if $A$ is $\G$-generated.  (Lemmas \ref{red1c2}, \ref{red1c1},
\ref{red1c0}.  The resolved case follows using  \ref{curve-selection}.)

Though the notion of a morphism $g: (X_1,\phi_1) \to (X_2,\phi_2)$ does
not depend on  $\phi_1,\phi_2$, recall that the $\RV$-Jacobian of $g$ is defined
with reference to these finite-to-one maps.

\<{lem} \lbl{rvlift-m}  Let $\T$ be $\V$-minimal and effective.  
    Let $X_i \subset \RV^{k_i}$ be definable and
let  $\phi_i: X \to \RV^n$ be   definable maps with 
 finite fibers; let 
 $g: X_1 \to X_2$ 
 %$g: (X_1,\phi_1) \to (X_2,\phi_2)$ 
 be a definable bijection.  
Assume  given in addition 
a definable function $\d: X_1 \to \RV$, such that:
\begin{enumerate}
  \item $\valr \d(x) = \JG g (x)$ for all $x \in X_1$.
  \item $\d(x) = \JRV g (x)$ for almost all $x \in X_1$ (i.e. all $x$ outside a set of dimension $<n$.)
\end{enumerate}
 Then  there exists a definable bijection 
 $G: X_1 \times_{\phi_1,\rv}(\VF^\times)^n \to X_2 \times_{\phi_2,\rv} (\VF^\times)^n$
 such that $\rho_2 \circ G = g \circ \rho_1$, where $\rho_i$ is the  natural projections to the $X_i$,
 and  such  that for 
   any  $x \in X_1 \times _{\phi_1,\rv} (VF^*)^n$,  
  $G$ is differentiable at $x$, and 
  $\rv (\JVF(G)(x)) = \d(x)$  
\>{lem}

\<{proof}    We follow closely the proof of \propref{rvlift}.  As there,
we may assume $\T = ACVF_A$, with  $A$ be an almost 
$(\VF,\G)$-generated  substructure.  

We first assume that $\valr \phi_1(X_1)$ is a single point of $\G^n$
  
  Then $X_1$ can be definably embedded into $\k^N$ for some $N$, and it follows from
  the orthogonality of 
  $\k$ and $\G$ that
  the image of $X_1$ in $\G$ under any definable map is finite. Thus  
  $\phi_2 X_2$ is contained in finitely many cosets $(C(a): a \in S)$
of $(\k^*)^n$; by partitioning $X_1$ working in $\T_{A(a)}$, we may assume $\phi_2 X_2 $ is contained in a single coset.  (cf. \lemref{collect}).   

As in \propref{rvlift}, we may assume $\phi_i X  \subseteq \k^n$, and indeed that 
$\phi_1 X = U, \phi_2 X = V$ are smooth varieties.   If $\dim(U)=n$, then 
the lift constructed in \propref{rvlift} satisfies $\rv (\JVF(G))(x) = \JRV g (x)$ for
$x \in X \times _{\phi_1,\rv} \VF^n$; thus
by assumption (2), we have $\rv (\JVF(G)))(x) = \d(x)$  for almost all $x$.  The exceptional
points have dimension $<n$, and may be partitioned into smooth varieties of dimension $<n$.
Thus we are reduced to the case $\dim(U)<n$.  We prove it by induction on $\dim(U)$. 
In this case  choose any lifting $G_0$.
We have an error term $e(x) = \rv(\JVF(G_0))(x) \inv \d(x) $.  Now $A(x)$ is
almost $\VF,\G$-generated, and so balls $\rv \inv (y)$ contain definable points; 
thus $e(x) = \rv E(x)$ for some definable $E: (X \times _{\phi_1,\rv} \VF^n) \to \VF$.
 Since $U$ is a smooth subvariety of $\k^n$
of positive codimension, some regular $h$ on $\k^n$ vanishes on $V$, 
while some partial derivative (say $h_{1}$)   vanishes only on a smaller-dimensional 
subvariety.   By induction may assume $h_1$ vanishes nowhere.  Lift $h$ to $H$; so $H_1$ lifts $h_1$.  
   Compose $G_0$ 
with a map fixing all coordinates but the first, and multiplying the first coordinate
by $E(x)H(y)/ H_1(y)$.  (Here $x = g \inv (y)$.) Where $h$ vanishes,  this has Jacobian $E(x)$; so the composition
has $\RV$- Jacobian $\d(x)$ as required.   

Now in general.  For any $\g \in \G^n$ let $X_1(\g) = \{x \in X_1: \valr \phi_1 (x) = \g \}$,
$X_2(\g) = g (X_1(\g))$.  By definition of $\JRV$ and $\JG$, $\JRV( g | X_2(\g)) = \JRV(g) | X_2(\g)$ and likewise $\JG$.  By the  case already analyzed 
  (for the sets $X_1(\g), X_2(\g)$ defined in $ACVF_{A(\g)}$)
  there exists an $A(\g)$-definable bijection 
$G_\g: X_1(\g) \times_{\phi_1,\rv} (VF^\times)^n \to X_2(\g) \times _{\phi_2,\rv} (VF^\times)^n$
with $\rv(\JVF(G_\g)(x)) = \d (x)$.  
 As in \lemref{collect} one can extend the $G_\g$ by compactness to definable
 sets containing $\g$, cover $X_1$ by finitely many such definable sets, and glue
 together to obtain a single function $G$ with the same property.

% {proof for maps on \G?}   After  a multiplicative  change of 
%variables on the target side, we may   rewrite $f$
%as $f(x',x'') =   (b_1(x'),b_2(x')+x'')$.  
%If $b_1(x')$ does not determine
%$b_2(x')$, i.e. $b_1(x')=b_1(y')$, $b_2(x') \neq b_2(y')$,
% then by choosing a generic $x''$, and letting $y'' = x'' +b_2(x')-b_2(y')$
% we have
%$f(x',x'') = (b_1(x'),b_2(x')+x'') = (b_1(y'),b_2(y')+y'') = f(y',y'')$
%contradicting injectivity.  (Note that $\valr(y'')=\valr(x'')$
%and so they have the same type.)  Thus $b_2= \phi \circ b_1$
%for some definable function $\phi$.  Now we can compose
%with the function $(z',z'') \mapsto (z',z''/\tilde{\phi}(z'))$,
%where $\tilde{\phi}$ lifts $\phi$.  (The existence of a lift $\tilde{\phi}$
%follows from the case $r=0$.)  Thus we are reduced to the case: 
%$f(x',x'') = (b_1(x'),x'')$; and indeed to $b_1$, i.e. again to the case $r=0$. 
                                                                                                     \>{proof}

{\bf Remark}  Assume $Id_X: (X,\phi_1) \to (X,\phi_2)$ has Jacobian $1$ everywhere.  
Then it is possible to find $F$ that is everywhere differentiable, of Jacobian precisely equal to $1$.  
At the  before the point where
Hensel's lemma is quoted, it is possible to multiply the function
by $J(F)^{-1}$ (not effecting the reduction, since $J(F) \in 1 + \Mm$.)
Then one obtains on each such coset a function of Jacobian $1$,
and therefore globally.

{\bf Example}  Let $\phi_2(x)=\phi_1(x)^m$.  A definable bijection 
  $$X \times_{\phi_1,\rv}(\VF^\times)^n \to X \times_{\phi_2,\rv} (\VF^\times)^n$$
is given by $(x,y) \mapsto (x,y^m)$.  (If $\rv(u)=\phi(x)^m$, there
 exists a unique $y$ with $\rv(y) = \phi(x)$ and $y^m=u$.)

\<{example} \lbl{rvlift-ex}   \propref{rvlift} need not remain valid over an $\RV$-generated base set.  \rm Let $A = \dcl(c)\}$, $c$ a transcendental point of $\k$.  
Let $f_1(y)=y, f_2(y)=1$,
 $\L(Y,f_i):= \VF \times_{\rv,f_i} Y = \{(x,y) \in \RV \times Y: \rv(x)=y \}$.
 Then $\L(Y,f),  \L(Y,f')$ are both open balls; 
%Let $Y = \{ \pm \sqrt (c) \}$.
% Let $f(y)=y$, and let $f'(y)=1$.  Let  $\L(Y,f):= \VF \times_{\rv,f} Y = \{(x,y) \in \RV \times Y: \rv(x)=y \}$, and similarly $ \L(Y,f')  = \VF \times_{\rv,f} Y$.  
  %are both unions of two open balls.
over any field $A'$ containing $A$, they are definably isomorphic, using a translation.
% on each ball.
 But these balls are not definably isomorphic over $A$.
\>{example}

\>{section}

  \<{section}{Special bijections and $\RV$-blowups}   \lbl{descent-iso}
  
  We work with a $\V$-minimal theory $\T$.  Recall the lift $\L : \RV[\leq n,\idot] \to \VFni$, with
$\rho_X: \L X \to X$.  
   Our present goal is an intrinsic description in terms of $\RV$ of the congruence relation:  $\L X \iso \L Y$.

$A$ will denote a $(\VF,\G,\RV)$-generated substructure of a model of $\T$.
Note that $\T_A$ is also $\V$-minimal 
(\corref{red1c1})
so any lemma proved for $\T$ under our assumptions can be used for any $\T_A$.

 The word ``definable'' below refers to $\T$.  The categories $\VF,\RV[*]$ defined below thus depend on $\T$; when necessary, we will denote them $\VF_\T$, etc.
When $\T$ has the form $\T=\T^0_A$ for fixed $\T^0$ but varying $A$, we write $\VF_A$, etc. 

\ssec{Special bijections}
Let  $X \subseteq \VF^{n+1} \times \RV^m$ be $\conjrv$-invariant.   Say
$$X= \{(x,y,u) \in \VF \times \VF^n \times \RV^m:  (\rv(x),\rv(y),u) \in \bar{X} \}$$   
(We allow $x$ to be any of the $n+1$ coordinates, $y$ the others.)

 Let $s(y,u)$ be a definable
function into $\VF$ with $\conjrv$-invariant domain of definition
$$\dom(s) = \{(y,u): (\rv(y),u) \in \bar{S}\}$$
and $\theta(u)$ a definable function on $pr_u(\dom(s))$ into $\RV$, 
 such
that $(s(y,u),y,u) \in X$ and $\rv(s(y,u))=\theta(u)$ for $(y,u) \in \dom(s)$.
Note that $\theta$ is uniquely defined (given $s$) if it exists.
%Note that   if $(x,y,u) \in X$ (with $x \in K, y \in \VF^n$), $(y,u) \in \dom(s)$, then $\rv(s(y,u))=\rv(x)$.  (Since $\rv(x) = \theta(u) = \rv(s(y,u))$.)  \VF
Let
$$X_1 = \{(x,y,u) \in X: (\rv(y),u) \in \bar{S}, \rv(x)=\theta(u) \}, \ \ \ \ X_2 =X  \m X_1  $$
$$X_1' = \{(x,y,u)  \in \VF\times \dom(s):  \val (x) > \valr \theta(u) \}  $$ 
and let $X' =    X_1' \du  X_2 $.
Also define $e_s: X' \to X$ to be the identity on $X_2$, and
  $$e_s(x,y,u)= (x+s(y,u),y,u)$$
   on $X_1'$.

\<{defn}\lbl{special-reduction}  $e_s:X' \to X$ is a definable bijection, called an
{\em elementary bijection}.  
        \qed
   \>{defn}

\<{lem} \lbl{special-relation}\begin{enumerate}
  \item  If $X$ is $\conjrv$-invariant, so is   $X' $.  If $X \to \VF^{n+1}$ is finite-to-one, the same is true of $X'$.
  \item   If $X_i = \L {\bar X}_i$, $X'_1 = \L{\bar X}'_1$,
then   $\bar{X}'_1$ is isomorphic to 
 $(\RVp \du \{1\}) \times \bar{S}$, while $\bar{X}_1$ is isomorphic to $\bar{S}$.   
  \item  If the projection $X \to \VF^{n+1}$ has finite fibers, then so do  the projection 
  $\dom(s) \to \VF^n$, and also the projection $\bar{S} \to \RV^n$,  $(y',u) \mapsto y'$.
    \item  $e_s$ has partial derivative matrix $I$ everywhere, hence has Jacobian $1$.
  Thus if   $F: X \to Y$ is such that $\rv \JVF F$ factors through $\rho_X$, then 
  $\rv \JVF (F \circ e_s)$ factors through $\rho_{X'}$. 
\end{enumerate}
 
  \>{lem}
 \<{proof}  (1,4) are clear.  The first isomorphism of (2)
 is obtained by dividing $x$ by $\theta(u)$, the second is
    evident.  For (3), note that if $(y,u) \in \dom(s)$ then $(s(y,u),y,u) \in X$
    so by the assumption $u \in \acl(y,s(y,u))$.  But for fixed $y$, $\{s(y,u): u \in \dom(s)\}$
    is finite, by \lemref{fmrfr}.  So in fact $(y,u) \in \dom(s)$ implies $u \in \acl(y)$.  
    Hence   $(y',u) \in \bar{S}$ implies $u \in \acl(y)$ for any $y$ with $\rv(y)=y'$,
    so (fixing such a $y$) $\{u: (y',u) \in \bar{S}\}$ is finite for any given $y'$.  \>{proof}

    A {\em special bijection} is a composition of elementary  bijections and 
 {\em auxiliary bijections} $(x_1,\ldots,x_n,u) \mapsto (x_1,\ldots,x_n,u,\rv(x_1),\ldots,\rv(x_n))$.
 
% (The auxiliary bijections are also $\conjrv$-invariant morphisms.)
 
An elementary bijection depends on the data $s$ of a partial section  of $X \to \VF^n \times \RV^m$.
Conversely, given $s$, if $\rv(s(y,u))$ depends only on $u$ we can define $\theta(u) = \rv(s(y,u))$ 
and obtain a special bijection.  If not, we can   apply an  auxiliary bijection   to $X  \subseteq \VF^{n} \times \RV^{m}$, and  obtain a set $X' \subseteq \VF \times \RV^{m+n}$, such that $\rv(x) = pr_{m+1} (u)$ for $(x,u) \in X'$.
For such a set $X'$, the condition for existence of $\theta$ is automatic 
%  $(s(y,u),y,u) \in X$ implies  $\rv(s(y,u))=\theta(u)$,
%for appropriate $\theta$ (one of the projections.)   
and we can define an elementary bijection $X'' \to X'$ based on $s$, and obtain a special 
bijection $X'' \to X$ as the composition.

  The classes of auxiliary morphisms and elementary morphisms  are all
closed under disjoint union with any identity morphism, 
and it follows that the class special morphisms under disjoint unions.  
%	 We would also close under ``glueing'', i.e.  the disjoint union of two morphisms, but this is unnecessary;

%     show that any \rvinv  morphism between \rvinv objects
% can be obtained as a composition of special bijections, their inverses, and \rvinv-bijections
% between \rvinv objects.   We begin with dimension $1$.

 \ssec{Special bijections in  one variable and families of $\RV$-valued functions.}
  
 We consider here special bijections in dimension $1$.    
 An elementary   bijection $X' \to X$ in dimension one   involves a finite set $B$ of $\rv$-balls, and a set of ``centers'' of these balls (i.e. a set $T$ containing a unique point $t(b)$ of each $b \in B$), and translates each ball so as to be centered at $0$ (while fixing the $\RV$ coordinates.)     We say that $X' \to X$ blows up the balls in $B$, with centers $T$.

Given a special bijection $h': X' \to X$, let $\fnr(X;h')$ be the set of definable functions
$X \to \RV$ of the form $H(\rho_{X'} ((h')^{-1}(x)))$, where $H$ is a  definable
function.  This is a finitely generated set of definable functions $X \to \RV$.
There will usually be no ambiguity in writing $\fnr(X,X' \to X)$ instead.

 Note that while a special bijection is an isomorphism in $\VF$,  an asymmetry exists
If $e: X' \to X$ is a special bijection, then $\fnr(X,X) \subseteq \fnr(X,X'\to X)$, usually
properly.  

What is the effect on $\fnr$ of passing from $X'$ to $X''$, where $X'' \to X'$ is a special bijection?
The auxiliary bijections have no effect.  Assume $\rv$ is already a coordinate function of
$X'$.  Consider  an elementary   
bijection $e_s: X'' \to X'$. 
 Let $B = \{(x,u) \in X': u \in \dom(s) \}$.   
 Then
the characteristic function $1_B$ lies in $\fnr(X',Id_{X'})$;  so $1_B \circ (h') \inv$
lies in $\fnr(X,h')$.  Using this, we see that $\fnr(X',e_s)$ is generated over
$\fnr(X',Id_{X'})$ by the function $B \to \RV$, $(x,u) \mapsto \rv(x-s(u))$.
(Extended by $0$ outside $B$.)    Thus if $h'' = h' \circ e_s: X'' \to X$,
$\fnr(X,h'')$ is generated over $\fnr(X,h')$ by the composition of the function 
$(x,u) \mapsto \rv(x-s(u))$
with $(h') \inv$.

Conversely,   if $B$ is a finite union of open balls whose
 characteristic function  lies in $\fnr(X,h')$, and if there exists 
 a definable set $T$ of representatives (one point $t(b)$ in each ball $b$ of $B$), and a function 
 $\phi=(\phi_1,\ldots,\phi_n)$, $\phi_i \in \fnr(X,h')$, with $\phi$ injective on $T$,
 then  one can find a special bijection $X'' \to X'$ 
with composition $h'': X'' \to X$, such that $\fnr(X,h'')$ is generated over
$\fnr(X,h')$ by $y \mapsto \rv(y-t(y))$, $y \in b \in B$.  
Namely let $\dom(s) = \phi(T)$, and 
for $u \in \dom(s)$ set 
$s(u) = h' \inv( t)$ if $t \in T$ and $\phi(t)=u$.  
In this situation, we will say that the balls in $B$ are blown up by $X'' \to X'$, with centers $T$.
Let  $\theta(u)=\rv(s(u))$.  Because $X' \to X$ 
may already have blown up some of the balls in $B$, $\fnr(X,h'')$ is generated
over $\fnr(X,h')$ by the restriction of $y \mapsto \rv(y-t(y))$ to some sub-ball of $b$,
possibly proper.  Nevertheless, we have:

\lemm{claim}  The function $y \mapsto \rv(y-t(y))$ on $B$ lies in $\fnr(X,h'')$. \>{lem}

\proof  This follows from the following more general
 
\Claim{} Let $c \in \VF, b  \in \fB$ be definable, with $c \in b$.  Let $b'$ be an $\rv$-ball
with $c \in b'$.  Then the function $\rv(x-c)$ on $b$ is 
 generated by its restriction to $b'$, $\rv$, and the characteristic function of $b$.

\proof   Let $x \in b \setminus b'$.  From $\rv(x)$ compute $\val(x)$.   If  $\val(x)< \val(c)$, $\rv(x-c)=\rv(c)$.  If 
 $\val(x)>  \val(c)$, $\rv(x-c)=\rv(x)$.  When $\val(x)=\val(c)$, but $x \notin b'$, $\rv(x-c) = \rv(x)-\rv(c)$.  Recall here that 
 $\valr^{-1}(\gamma)$ is the nonzero part of a $\k$-vector space; subtraction, for distinct elements $u,v$, can 
 therefore be defined by $u-v = u(u^{-1}v -1)$.  \qed

Thus any
special bijection can be understood as blowing up a certain finite number of balls
(in a certain sequence and with certain centers.) 
We will say that a special bijection
$X'' \to X'$ is subordinate to a given partition of $X$ if 
if each ball  blown up by $X'' \to X'$
is contained in some class of the partition.

  It will sometimes be more convenient to 
work with the sets of functions $\fnr(X,h)$ than with the special bijections $h$ themselves. 

We observe that 
any finite set of definable functions $X \to \RV$ is contained in $\fnr(X;h)$ for some $X',h$:

\<{lem} \lbl{rvfactor1}   Let $X \subseteq \VF \times \RV^*$ be \rvinv, and let 
$f: X \to (\RV \union \G)$
% (or:  $f: X \to \VF$ with finite support) 
be a definable map. Then there
exists a $\conjrv$-invariant $X' \subseteq \VF \times \RV^*$   a special bijection $h: X' \to X$,  and a definable function $t$ 
such that $ t \circ \rho_{X'} = f \circ h$.   
%
%If some special $h'':X'' \to X$ is given in advance, we can take $h$ to factor through $X'' \to X$.
Moreover, if $X = \union_{i=1}^m P_i$ is a finite partition of $X$ into sets whose
characteristic functions factor through $\rho$,  we can find $X' \to X$ subordinate
to this partition. \>{lem}

\<{proof}  Say $X \subseteq \VF \times \RV^m$; let $\pi: X \to \VF, \pi': X \to \RV^m$ be the
projections.  Applying an auxiliary bijection, we may assume $\rv( \pi(x)) = pr_m \pi'(x)$,
i.e. $\rv(\pi(x))$ agrees with one of the coordinates of $\pi'(x)$.   We now claim that there exists a finite $F' \subseteq \RV^m$, such that away from 
$\pi'^{-1}(F')$, $f$ factors through $\pi'$.  To prove this, it suffices to show that 
if $p$ is a complete type of $X$ and $\pi'_* p$  is non-algebraic (i.e. not contained in a finite definable set),   then $f|p$ factors through $\pi'$; this follows from 
\lemref{1gen}.

We can thus restrict attention to $\pi' \inv (F')$; our special bijections will be the identity
away from this.  So we may assume $\pi'(X)$ is finite.  Recall  that 
(since an auxiliary bijection has been applied)
$\rv$ is constant on each fiber of $\pi'$.  
In this case there is no problem
relativizing to each fiber of $\pi'$, and then collecting them together (\lemref{collect}), 
  we may assume
in fact that $\pi'(X)$ consists of a single point $\{\uu\}$.  In this case the partition
(since it is defined via $\rho$) will automatically be respected.
 
The rest of the proof  is similar to \lemref{tr1}.   We consider first   functions $f$ 
with finite support $F$ (i.e. $f(x)=0$ for $x \notin F$,) and prove the analog of the
statement of the lemma for them.    If $F =  \{0\} \times \{\uu\}$  
then $F = \rho^{-1}(\{(0)\} \times \{\uu\}))$ so the claim is clear.    If $F = \{(x_0,\uu)\}$,
let $s: \{\uu\}\to \VF$, $s(\uu) = x_0$.  Applying $e_s$  returns us to the previous case. 
If $F = F_0 \times \{\uu\}$ has more than one point, we    use induction on the
 number of points.     Let $s(\uu)$ be the average of  $F_0$.   Apply the special bijection $e_s$.  Then the result
 is a situation where $\rv$ is no longer constant on the fiber.  Applying a auxiliary bijection 
 to make it constant again, the fibers of $F \to \RV^{m+1}$ become  smaller.
  
The case of the characteristic function of a finite union of balls is similar (following \lemref{tr1}.)

Now consider a general function $f$.  Having disposed of the case of 
characteristic
function,
 it suffices to treat $f$ on each piece of any given partition.  Thus 
we can assume $f$ has the form of \corref{cell1.5}, $f(x) =H( \rv(x-n(x)))$.
Translating by the $n(x)$ as in the previous cases, we may assume $n(x)=0$. But 
t then again
$f$ factors through $\rho$ and $\rv$, so one additional auxiliary bijection suffices. 
                                                                                                     \>{proof}

\<{cor}\lbl{rvfactor1c}  Let $X,Y \subseteq \VF^n \times \RV^*$, and let $f: X \to Y$ be a definable bijection.  
Then there exists a special bijection $h: X' \to X$, and $t$ such that 
$\rho_Y \circ (f \circ h)= t \circ \rho_{X'}$     

It can be found subordinate to a given finite partition, factoring through $\rho_X$.
 \qed\>{cor}

We wish to obtain a symmetric version of \corref{rvfactor1c}.   We will say that bijections $f,g: X \to Y$  
{\em differ by special bijections } if there exist special bijections $h_1,h_2$ with
$h_2g=fh_1$.   We show that 
  every definable bijection
 between \rvinv objects differs  by special bijections from 
a \rvinv bijection. 

\<{lem} \lbl{rvconj1}  
Let $X \subseteq \VF \times \RV^m$, $Y \subseteq \VF \times \RV^{m'}$ be
 definable, $\conjrv$-invariant; let $F: X \to Y$ be a definable bijection.  
  Then there exist   special bijections $h_X:X' \to X$, $h_Y: Y' \to Y$, and 
  an $\conjrv$-invariant bijection
  $F': X' \to Y'$ with
$F = h_Y F' h_X^{-1}$.   I.e. 
  $F$ differs from an $\conjrv$-invariant 
    morphism  by  special bijections.    
    \>{lem}

\prf 
It suffices to find  $h_X,h_Y$ such that $\fnr(X,h_X) = F \circ \fnr(Y,h_Y)$; for then 
    we can let $F' = h_Y^{-1}F h_X$.  

Let $X = \union_{i=1}^m P_i$ be a partition as in 
\propref{image2}.  By \lemref{rvfactor1}, there exist $X_0,Y_1$ and 
special bijections $X_0 \to X$, $Y_1\to Y$, such that the characteristic functions
of the sets $P_i$ (respectively the sets $F(P_i)$) are in $\fnr(X,X_0 \to X)
(respectively \fnr(Y,Y_1 \to Y))$. 

By \corref{rvfactor1c}, one can find a special $X_1 \to X_0$ such that 
$\fnr(X, X_1 \to X)$ contains $F \circ \fnr(Y,Y_1 \to Y)$.  By another application
of the same, one can find a special bijection $Y_* \to Y_1$  subordinate to $\{F(P_i)\}$  such that  
\beq{rvconj1-1}  \fnr(Y,Y_* \to Y)  \supseteq F^{-1} \circ \fnr(X,X_1 \to X)  \eeq

Now $Y_*$ is obtained by composing a sequence $Y_*=Y_m \to \ldots \to   Y_1$ of elementary   bijections and auxiliary bijections.  
We define inductively $X_m \to \ldots \to X_2 \to X_1$, such that 

\beq{rvconj1-2}  \fnr(Y,Y_k \to Y) \circ F \subseteq \fnr(X,X_k \to X) \eeq
 
  Let $k \geq 1$.  
$Y_{k+1}$
is obtained by blowing up a finite union of balls $B$ of $Y$, with definable set $T$ of representatives  
such that some $\phi \in \fnr(Y,Y_k \to Y)$ is injective on $T$; 
and then $\fnr(Y,Y_{k+1} \to Y)$ is generated over $\fnr(Y,Y_k \to Y)$ by $\psi$, where for $y \in b \in  B$
$\psi(y)= \rv(y-t(b))$ (\lemref{claim}).    
    By the choice of the partition $\{P_i\}$, 
%\propref{image2}, 
$F \inv (B)$ is also a finite union of balls.

Now $F \inv(B)$, with center set $F \inv (T)$, can serve as data for a special bijection:
the requirement about the characteristic function of $B$ and the injective
function on $T$ being in $\fnr$ are satisfied by virtue of \lemref{claim}.
We can thus   define $X_{k+1} \to X_k$ so as to blow up $F \inv (B)$ with center
set $F \inv (T)$.  
  By \lemref{image2l}, 
$\rv(F(x)-F(x'))$ is a function of $\rv(x-x')$ (and conversely) on each of these balls,
so $\fnr(X,X_{k+1} \to X)$ is generated over $\fnr(X,X_k)$ by 
$\psi \circ F$.   Hence \eqref{rvconj1-2} remains valid for $k+1$.  

Now by \eqref{rvconj1-1}, $\fnr(X,X_1 \to X) \subseteq \fnr(Y, Y_* \to Y) \circ F$; 
  since the generators match at each stage, by induction on   $k \leq m$
\beq{rvconj1-3}   \fnr(X,X_k \to X)  \subseteq \fnr(Y, Y_m \to Y) \circ F   \eeq

By \eqref{rvconj1-2} and \eqref{rvconj1-3} for $k=m$, 
$\fnr(X,X_m \to X) = \fnr(Y, Y_* \to Y) \circ F$
  \eprf %$\fnr(Y,Y_k \to Y) \circ F \subseteq \fnr(X,X_k \to X)$

For the sake of possible future refinements, we note that the proof of \lemref{rvconj1} 
shows also:   
\<{lem} \lbl{rvconj1.1} Let $X \subseteq \VF \times \RV^m$, $Y \subseteq \VF \times \RV^{m'}$ be
 definable, $\conjrv$-invariant; let $F: X \to Y$ be a definable bijection.    If a  \propref{image2} - partition for $F$ has characteristic functions factoring through $\rho_X,\rho_Y$, and if $F$ is $\conjrv$-invariant, then
for any special bijection $h_X': X' \to X$
there exists a special bijection $h_Y': Y' \to Y'$ such that $(h_Y') \inv F h_X'$ is
$\conjrv$-invariant.  \qed  \>{lem}

 \ssec{Several variables} \lbl{special}
 
  We will show now in general  that any definable map from an    \rvinv object
  to $\RV$ factors through the inverse of a special bijection, and the standard map $\rho$.

\<{lem} \lbl{rvfactor}    Let $X \subseteq \VF^n \times \RV^m$
be $\conjrv$-invariant,  and let 
$\phi: X \to (\RV \union \G)$.   Then there
exists   a special bijection $h: X' \to X$,  and a definable function $\tau$ 
such that $ \tau \circ \rho_{X'} = \phi \circ h$.  
\>{lem}

\<{proof}  By induction on $n$.  
For $n=0$ we can take $X'=X$, since $\rho_{X} = Id_X$.  

For $n=1$ and $X \subseteq \VF$, by \lemref{rvfactor1}, there exists $\mu = \mu(X,\phi) \in \Nn$ such that the
lemma holds for some $h$ that is a composition of $\mu$ elementary   and auxiliary bijections.
It is easy to verify the semi-continuity of $\mu$ with respect to the definable topology:
 if $X_t$ is a definable family of definable sets, so that $X_b$ is $A(b)$-definable,
and 
$\mu(X_b,\phi |X_b)=m$, then there exists a definable
set $D$ with $b \in D$ and such that if $b' \in D$, then $\mu(X_{b'}, \phi | X_{b'}) \leq m$.  

Assume the lemma known for $n$ and 
suppose $X   \subseteq \VF \times Y$, with  $Y  \subseteq \VF^n \times \RV^m$.  For any
$b \in Y$, let $X_b = \{x: (x,b) \in X \} \subseteq \VF$; so $X_b$ is $A(b)$-definable.
 
Let $\mu=  \max _b \mu(X_b,\phi | X_b)$.   Consider first the case $\mu=0$.  
Then $\phi | X_b = \tau_b \circ \rho | X_b$,
for some $A(b)$-definable function $\tau_b: \RV^m \to (\RV \union \G)$.  By
stable embededness and elimination of imaginaries in $\RV \union \G$,
there exists (\secref{stab-emb})  a canonical parameter $d \in (\RV \union \G)^l$, and an $A$-definable
function $\tau$, such that $\tau_b(t) = \tau(d,t)$; and $d$ itself is definable 
from $\tau_b$, so we can write $d=\delta(b)$ for some definable $\delta: Y \to (\RV \union \G)^l$.  Using the induction hypothesis for $(Y,\delta)$ in place of $(X,\phi)$,
we find that there exists a $\conjrv$-invariant $Y' \subseteq \VF^n \times \RV^*$, a special $h_Y: Y' \to Y$,
and a definable $\tau_Y$, such that $\tau_Y \circ \rho_{Y'} = \delta \circ h_Y$.
Let $X' = X \times_Y Y'$, $h(x,y')=(x,h_Y(y'))$.    
An elementary   bijection to $Y$ determines one to $X$, where the function
$s$ does not make use of the first coordinate; so  $h: X' \to X$ is special.   In this case the lemma is proved:
  $\phi \circ h (x,y') = \phi(x,h_Y(y')) = \tau(\delta(h_Y(y')),\rho(x,y)) = \tau(\tau_Y(\rho_{Y'}(y')),\rho(x,y))$. 

Next suppose $\mu>0$.  Applying an auxiliary bijection, we may assume that 
for some definable function (in fact projection) $p$, $\rv(x) = p(u)$ for $(x,y,u) \in X$.
For each $b \in Y(M)$ (with $M$ any model of $\T_A$) there exists
  an elementary   bijection $h_b: X'_b \to X_b$, such that 
 $\mu(X'_b,\phi | X'_b) < \mu$; $h_b$ is
determined by   $s_b , \theta_b$,with $s_b \in \rv (s_b) =  \theta_b$, and
$(s_b,\theta_b) \in X$.  (The $u$-variables have been absorbed into $b$.)  By
compactness, one can take $s_b=s(b)$ and $\theta_b  = \theta'(b)$ for
some definable functions $s,\theta'$.   By the inductive hypothesis applied to
$(Y,\theta')$, as in the previous paragraph, we can  assume $\theta'(y,u) = \theta(u)$
for some definable $\theta$.  Applying the special bijection with data $(s,\theta)$
now amounts to blowing up $(s_b,\theta_b)$ uniformly over each $b$, and thus 
reduces the value of $\mu$.     \>{proof}

%{\bf Remark on the proof of \lemref{rvfactor}}   The above proof uses induction on $n$ and $\mu$;
%at each stage, the ambient space is split as $\VF \times (\VF^n \times \RV^m)$, where
%any additional $\RV$ coordinates obtained while blowing up $X$ are pushed to the right.
%As a result, even if $X \in \VFr$, during the induction one obtains sets $Y$ that
%may not be in $\VFr$, (since $Y \to \VF$ may not be finite-to-one)  and \lemref{rvfactor1} is used over any $\RV$ base.  
% However this is inessential; we could instead split the space
%as $(\VF \times \RV^l) \times (\VF^n \times \RV^m)$, in such a way that the projection
%$Y$ to $\VF^n \times \RV^m$ is also in $\VFr$.    
%
%Similarly in \lemref{rvconj1rel} below, the case: $X_i \in \VFr$ suffices.
\<{question}  \lbl{q} Is \propref{rvconj1} true in higher dimensions?  
\>{question}

\<{cor} \lbl{rvfactorcor}  Let $X \subseteq \VF^n \times \RV^m$ be definable. Then every definable function    $\phi: X \to \G$ factors through a definable function $X \to \RV^*$. \>{cor}

\proof By \lemref{tr2} we may assume $X$ is $\conjrv$-invariant; now the corollary
follows from \lemref{rvfactor}.  

(It is convenient to note this here, but it can also be proved with the methods of \S 3; the main point is
that on the generic type of a   ball with center $c$, every function into $\RV \union \G$ factors through 
$\rv(x-c)$; while on a transitive ball, every function into $\RV \union \G$ is constant.)

Consider pairs $(X',f')$ with $X',f': X' \to \VF^n$ definable, 
and such that $f'$ has $\RV$-fibers.  A bijection $g: X' \to X''$ is said to be 
{\em relatively unary} (with respect to $f',f''$) if  it commutes with $n-1$ coordinate projections,
i.e. $pr_i f'' g = pr_i f'$ for all but at most one value of $i$.

Given $X \subseteq \VF^{n} \times \RV^m$, we view it as a pair $(X,f)$ with $f$ 
the projection to $\VF^n$.   Thus for $X,Y \subseteq \VF^{n} \times \RV^*$, 
the notion: $F: X \to Y$ is relatively unary is defined.

Note that the elementary   bijections are relatively unary, as are the auxiliary bijections.

\<{lem}\lbl{decompose-unary}  Let $X,Y\subseteq \VF^n \times \RV^*$, and let $F: X \to  Y$ be a definable 
bijection.  Then $F$
can be written as the composition of relatively unary morphisms of $\VFni$.  \>{lem}

\<{proof}

We have $X$ with two finite-to-one maps    $f,g: X \to \VF^n$ (the projection, and 
the composition of $F$ with the projection $Y \to \VF^n$.)  We must decompose 
the identity $X \to X$ into a composition of 
   relatively unary maps $(X,f) \to (X,g)$.  

Begin with the case $n=2$; we are given $(X,f_1,f_2)$ and $(X,g_1,g_2)$.

\Claim{}  There exists   a definable partition of $X$ into sets $X_{ij}$ such that
$(f_i,g_j): X \to \VF^2$ is finite-to-one. 

\proof  Let $a \in X$.  We wish to show that for some $i,j$,
$a \in \acl(f_i(a),g_j(a))$.  This follows from    the exchange principle for algebraic closure in $\VF$:  if $a \in \acl(\emptyset)$, there is nothing to show.  Otherwise 
$g_j(a) \notin \acl(\emptyset)$ for some $j$; in this case either $a \in \acl(f_1(a),g_j(a))$
or $f_1(a) \in \acl(g_j(a))$, and then $a \in \acl(f_2(a),g_j(a))$.  The Claim follows by compactness.

 Let $h:X' \to X$ be a special bijection such
that  the characteristic functions of $X_{ij}$ are in $\fnr(X,X' \to X)$.  (\lemref{rvfactor}).
Since $h$ is composition of relatively unary bijections,
 we may replace $X$ by $X'$ (and $f_i,g_i$ by
$f_i \circ h, g_i \circ h$ respectively.)  So we may assume the characteristic function 
of $X_{ij}$ is in $\fnr(X,X)$, i.e. $X_{ij} \in VFr[n]$.  But then it suffices to 
treat each $X_{ij}$ separately; say $X_{11}$.  In this case the identity map on $X$
takes 
$$(X,f_1,f_2) \mapsto (X,f_1,g_1) \mapsto (X,g_2,g_1) \mapsto (X,g_2,g_1-g_2)
\mapsto (X,g_1,g_1-g_2) \mapsto (X,g_1,g_2)$$
where each step is relatively unary. 

When $n>2$, we move between $(X,f_1,\ldots,f_n)$ and $(X,g_1,\ldots,g_n)$,
by partitioning, and on a given piece  replacing each $f_i$ by some $g_j$, one at a time.    
\>{proof}

\ssec{$\RV$-blowups} % 
 
 We now define the  
 $\RV$-counterparts of the special bijections, which will be called $\RV$-blowups.  These will not be bijections; the kernel of 
 the homomorphism $\L: \SG[\RV] \to \SG[\VF]$ will be seen to be obtained by formally
 inverting $\RV$-blowups.   Let $\RVpi = \{x \in \RV: \val(x) > 0 \union \{\infty\}\} \subseteq \RVi$.
  In the $\RV[\leq 1]$-presentation, $\RVpi = [\RVp]_1 + [1]_0$.  
  (cf. \secref{altRV}). 
%  (Cf. discussion following \defref{RVcat}.)

 \<{defn}  \lbl{rvblowup} 
 \begin{enumerate}
  \item  Let $\bY =(Y,f) \in \Ob \RVi[n,\idot]$ be such that $f_n(y) \in \acl(f_1(y),\ldots,f_{n-1}(y))$, and $f_n(y) \neq \infty$.
Let $Y' = Y \times \RVpi$.  For $(y,t) \in Y'$, define $f'=(f'_1,\ldots,f'_n)$ by: $f'_i(y,t) = f_i(y)$ for $i<n$, 
$f'_n(y,t) = tf_n(y)$.  Then $\widetilde{\bY} = (Y',f')$ is an {\em elementary blowup} of $\bY$.  It comes
with the projection map $Y' \to Y$.
  \item 
Let $\bX=(X,g) \in \Ob \RVi[n,\idot]$, $X= X' \du X''$, $g' = g | X'$, $g'' = g | X''$,
 and let $\phi:   \bY \to (X',g')$ be an $\RVvol$-isomorphism.  Then the  {\em $\RV$- blowup}   $\widetilde{\bX}_{\phi}$ 
is defined to be $\widetilde{\bY} + (X'',g'') = (Y' \du X'', f' \du g'')$.  It comes with $b: Y' \du X'' \to X$,
defined to be the identity on $X''$, and the projection on $Y'$.   $X'$ is called the {\em 
blowup locus} of $b: \widetilde{\bX}_{\phi} \to \bX$.  
 
\end{enumerate}

An {\em iterated $\RV$-blowup} is obtained by finitely many iterations of $\RV$- blowups.
 
  \>{defn} 
 
Since blowups in the sense of algebraic geometry will not occur in this paper, we will
say `blowup' for $\RV$-blowup.  
 
\<{remark} \lbl{rvblowup-rem}  In the definition of an elementary blowup, $\dim_{\RV} (Y) < n$.   For such $Y$, 
$\phi: \bY \to (X',g')$ is an $\RVvol[\leq n, \idot]$-isomorphism iff it is an 
$\RVvolg$-isomorphism (\defref{rvvol}).  \>{remark}
    
\lemm{generators}  \begin{enumerate} 
  
  \item  Let $\bY'$  be an  elementary blowup of $\bY$.
  $\bY'$   is $\RVvol[n,\idot]$- isomorphic to $\bY'' = (Y'', f'')$,
  with $Y'' = \{(y,t) \in Y \times \RVi: \valr(t) > f_n(y) \}, f''(y,t)=(f_1(y),\ldots,f_{n-1}(y),t)$.

\item  An  elementary blowup
$\bY'$ of $\bY$ is $\RVi[n,\idot]$- isomorphic to $(Y \times \RVi,f')$ for any
$f'$ isogenous to $(f_1,\ldots,f_n,t)$.

\item  Up to isomorphism, the blowup depends only on blowup locus.
  In other words
if $X,X' ,g,g' $ are as in \defref{rvblowup}, and
$\phi_i: \bY_i \to (X',g')$ $(i=1,2)$ are   isomorphisms, then $\widetilde{\bX}_{\phi_1} , \widetilde{\bX}_{\phi_2}$
are $\bX$-isomorphic in $\RVvol[n,\idot]$.

\end{enumerate}  \>{lem}

\prf  (1)   The isomorphism is given by $(y,t) \mapsto (y,t f_n(y))$

(2) The identity map on $Y \times \RV$ is an $\RVi[n,\idot]$ isomorphism.

(3)  Let
$\psi_0 =  \phi_2 \inv \phi_1$, and define $\psi_1:   Y_1 \times \RVpi \to Y_2 \times \RVpi$
by  $\psi(y,t) = (\psi_0(y), t)$.   The sum of the values of the $n$ coordinates
of $\widetilde{\bY_i}$ is then $(\sum_{i<n} \valr f_i) + (\valr(t)+\valr f_n)$
in both cases.  Since by assumption $\psi_0: Y_1 \to Y_2$ is an $\RVvol$-isomorphism,
it preserves $\sum_{i \leq n} \valr f_i$ and so $\psi_1$ too is  an $\RVvolg$-isomorphism;
thus $\JRV(\psi_1) \in \k^*$ a.e. 
Let $ \theta  Y_1 \to \k^*$ be a definable
map such that $\theta =   \JRV(\psi_1)$ almost everywhere.
 Define $\psi :  Y_1 \times \RVpi \to Y_2 \times \RVpi$ by
 $\psi(y,t) = (\psi_0(y), t/\theta(y) )$.  Then one computes immediately
 that $\JRV(\psi)=1$, so $\psi$ is an $\RVvol[n,\idot]$-isomorphism, and hence so is
 $\psi \du Id_{X''} : \widetilde{\bX}_{\phi_1} \to \widetilde{\bX}_{\phi_2}$
 
%(6)  Follows from (3):  a blowup is   a direct sum of an identity map with a
% one of the form of  \lemref{generators} (3), and this characterization is $\RVvol$-isomorphism invariant
% over $\bX$.
\eprf

Here is a coordinate-free description of $\RV$-blowups;   we will not really use it in 
the subsequent development.

 \lemm{form}  \begin{enumerate}
\item  Let $\bY =(Y,g) \in \Ob \RVi[n,\idot]$, with $\dim(g(Y)) < n$; let $f: Y \to \RV^{n-1}$
be isogenous to $g$.  
  Let $h: Y \to \RV$ be definable, with $h(y) \in \acl(g(y))$ for $y \in Y$, and 
  with $\sum(g)=\sum(f)+\valr(h)$.  
Let $Y' = Y \times \RVpi$, and $f'(y,t) = (f(y), t h(y))$.   Then $\bY' = (Y', f')$ with   the projection map to $Y$
is a blowup.  
\item Let $\bY'' \to \bY$ be a  blowup  with blowup locus $Y$.  Then 
there exist $f,h$ such that with $\bY'$ as in (3), $\bY'',\bY'$ are isomorphic
over $\bY$.

\end{enumerate}
\>{lem}

\proof 
 (1)   Since $\dim_{\RV}(g(Y)) < n$,
$Id_Y: (Y,(f,h)) \to (Y,g)$ 
is an $\RVvol$-isomorphism. Use this as $\phi$ in the definition of blowup.

(2)  With notation as in \defref{rvblowup}, let  $h = g_n \circ \phi \inv, f = (g_1,\ldots,g_{n-1}) \circ \phi \inv$.

\<{defn}   For   $\fC = \RV[\leq n,\idot]$ or $\fC = \RVvol[\leq n,\idot]$, let 
 $\Isp[\leq n]$ be the set of pairs $  (\bX_1,\bX_2)  \in \Ob \fC$ such that there exist iterated blowups $b_i:  \widetilde{\bX}_i \to \bX_i$ and
  an  $\fC$- isomorphism $F:  \widetilde{\bX}_1 \to  \widetilde{\bX}_2$. \>{defn}

When $n$ is clear from the context, we will just write $\Isp$.  
\<{defn}
Let $1_0$ denote the one-element object of $\RV[0]$.  
Given a definable set $X \subseteq \RV^n$
% and a finite-to-one $f: X \to \RV^k$,$\bX = (X,f)$,
 %for each $n \geq k$,
 let  $\bX_n$ denote   $(X,Id_X) \in \RV[n]$, and $[\bX]_n$ the class   in $\SG(\RV[n])$.
   Write
 $[1]_1$ for $[\{1\}]_1$ (where $\{1\}$ is the singleton set of the identity element of $\k$.)

  \>{defn}

\lemm{diamond}  Let $\fC = \RV[\leq n,\idot]$ or $\fC = \RVvol[\leq n,\idot]$.

\begin{enumerate}  

  \item Let $f:\bX_1 \to \bX_2$ be an $\fC$- isomorphism, and let $b_1:\widetilde{\bX_1}  \to \bX_1$ be a blowup.
  Then there exists a blowup $b_2: \widetilde{\bX_2} \to \bX_2$ and a  $\fC$- isomorphism 
  $F: \widetilde{\bX_1} \to \widetilde{\bX_2}$
  with $b_2F = fb_1$.   
  \item  If $b: \widetilde{\bX} \to \bX$ is a blowup, then so are
  $b \du Id: \widetilde{\bX} \du \bZ \to \bX \du \bZ$ and  $(b \times Id): \widetilde{\bX} \times \bZ \to \bX \times \bZ$.   
  \item Let $b_i: \widetilde{\bX}_{\phi_i} \to \bX$ be a  blowup ($i=1,2$).  Then there exist  blowups $b_i': \bZ_i \to \widetilde{\bX}_{\phi_i}$
  and an isomorphism $F: \bZ_1 \to \bZ_2$ such that $b_2b_2' F = b_1b_1'$.
  \item Same as (1,2,3) for iterated blowups. 
  \item   $\Isp$ is an equivalence relation.  It induces a semi-ring congruence on
  $\SG \RV[*,\idot] $, respectively $\SG \RVvol[*,\idot]$.
  \item  As a semi-ring congruence on  $\SG \RV[*,\idot] $, $\Isp$ is generated by
  $([1]_1, [\RVp]_1 + 1_0)$.
 
  \end{enumerate}
\>{lem}
\prf (1) This reduces to the case  of elementary blowups.  If $\fC= \RVvol[n,\idot]$ then the composition
$f \circ b_1$ is already a  blowup.    % (with $\phi=f \circ \phi_{b_1}$).    
If $\fC  = \RV[\leq n,\idot]$, it is also clear using \lemref{generators} (2).

(2) From the definition of blowup.  % the iterated case, by induction.

(3)  If $b_1$ is the identity,
let $b_1' = b_2, b_2'=Id, F=Id$.  Similarly if $b_2$ is the identity.  If $X = X' \du X''$ and the statement is true above $X'$ and above $X''$, 
then by glueing   it is true also above $X$.  
 We thus reduce to the case that $b_1,b_2$ both are  blowups with blowup locus
 equal to $X$.   But then by \lemref{generators} (3), there exists an isomorphism
 $F: \widetilde{\bX}_{\phi_1} \to \widetilde{\bX}_{\phi_2}$ over $\bX$.  Let $b_1'=b_2'=Id$.

(4) For (1,2) the induction is immediate.     For (3), write $k$-blowup as shorthand for
 ``an iteration of $k$ blowups''.   We show by induction on $k_1,k'$ a more precise form:
 
\Claim{} If $\bX_1 \to \bX$ is a $k_1$-blowup, and $\bX' \to \bX$ is a $k'$ blowup,
then there exists an $k'$-blowup $\bZ_1' \to \bX_1$   a $k_1$-blowup $\bZ' \to \bX$, and an $\RVvol[n,\idot]$- isomorphism
$\bZ_1' \to \bZ'$ over $\bX$.

 If $k_1=k'=1$ this is (3).  So say  $k'>1$.    The map $\bX' \to \bX$ is a composition 
$\bX' \to \bX_2 \to \bX$, where $\bX_2 \to \bX$ is a  $k'-1$-blowup  
and $\bX' \to \bX_2$ is a blowup.  By induction there is a  $k'-1$ blowup $\bZ_1 \to \bX_1$ and a $k_1$-blowup $\bZ_2 \to \bX_2$  and an  $\RVvol[n,\idot]$- isomorphism
$ \bZ_1 \to \bZ_2$
 over $\bX$.  
 
 By induction again there is a  blowup 
and $\bZ_2' \to \bZ_2$, a $k_1$-blowup $\bZ' \to \bX'$
 an
$\RVvol[n,\idot]$- isomorphism   $\bZ' \to \bZ_2$ over $\bX_2$.   
By (1) there exists a blowup $\bZ_1' \to \bZ_1$ and an $\RVvol[n,\idot]$- isomorphism
$\bZ_1' \to \bZ_2'$, making the $\bZ_1,\bZ_2,\bZ_1',\bZ_2$-square commute.
  So $\bZ_1 \to \bX_1$ is a $k'$-blowup, $\bZ' \to \bX'$ is a $k_1$-blowup,
  and we have a composed  isomorphism $\bZ_1' \to \bZ_2' \to \bZ'$ over $\bX$.

(5)  If $(\bX_1,\bX_2) , (\bX_2,\bX_3) \in \Isp$, there are iterated blowups 
$\bX_1' \to \bX_1, \bX_2' \to \bX_2$ and an isomorphism $\bX_1' \to \bX_2'$;
and also $\bX_2'' \to \bX_2, \bX_3' \to \bX_3$ and $\bX_2'' \to \bX_3$.
Using (3) for iterated blowups,  there exist iterated blowups $\widehat{\bX_2}' \to \bX_2', \widehat{\bX_2}'' \to \bX_2''$,
and an isomorphism $\widehat{\bX_2}' \to  \widehat{\bX_2}' $.  By (1) for
iterated blowups there are iterated blowup $\widehat{\bX_1} \to \bX_1', \widehat{\bX_3} \to \bX_3$
and isomorphisms  $\widehat{\bX_1} \to \widehat{\bX_2}'$, $\widehat{\bX_2}'' \to  \widehat{\bX_3}$, with the natural diagrams commuting.  Composing we obtain 
$\widehat{\bX_1} \to  \widehat{\bX_3}$, showing that $(\bX_1,\bX_3) \in \Isp$.  Hence  $\Isp$ is an equivalence
relation. 

Isomorphic objects are $\Isp$-equivalent, so an equivalence relation on the semiring
$\SG \fC$ is induced.   If $(X_1,X_2) \in \Isp$ then 
by (2), $( X_1 \du Z, X_2 \du Z) \in \Isp$,
and $(X_1 \times Z, X_2 \times Z) \in \Isp$.  It follows
that $\Isp$   induces a congruence on the semiring $\SG \fC$.  

(6)  
 We can blow up  $1_1$ to $ \RVp_1 + 1_0 $, so $([1]_1, [\RVp]_1 + 1_0) \in \Isp$. 
 Conversely,   under the conditions of \defref{rvblowup}, let $\bY^- = [(Y,f_1,\ldots,f_{n-1})]$;
 then $[\bY] = [(Y,f_1,\ldots,f_{n-1},0)] = [\bY^-]  \times [1]_1$ by \lemref{generators},
  and we have:
$$[\widetilde{\bX}_{\bY}] =  
  [\bY]_{n-1} + [\bY]_{n-1} \times  [\RVp]_1  + [\bX'']  \cong_{\Isp} [\bY] \times [1]_1+[\bX''] = [\bX]$$
 modulo the congruence generated by $([1]_1, [\RVp]_1 + 1_0)$.

 \eprf

We now relate special bijections to blowing ups.  
 Given $\bX=(X,f), \bX' = (X',f') \in \RVni$,
say $\bX,\bX'$ are {\em strongly isomorphic} if there exists a bijection
$\phi:  X \to X'$ with $f' = \phi f$.    Strong isomorphisms are always   in $\RVvol[n,\idot]$.

Up to strong isomorphism, an elementary blowup
of $(Y,f)$ can be put in a different form:  $\widetilde(\bY) \iso (Y'',f'')$,
$Y'' = \{(z,y):  y \in Y, \valr(z) > \valr f_n(y) \}, f_i(z,y)=f_i(y)$ for $i<n$, $f_n(z,y)=z$.
The strong isomorphism $Y'' \to Y'$ is given by  $(z,y) \mapsto (y,z/f_n(y))$. 
This matches precisely the definition of special bijection, and makes evident the following lemma.

\lemm{sp-bl} \begin{enumerate}  Let $\fC = \RVi [n, \idot]$ or $\RVvol[\leq n, \idot]$.

\item  $\bX,\bY$ are strongly isomorphic over $\RV^n$ iff $\L \bX, \L \bY$ are isomorphic over
the projection to  $\VF^n$.

\item    Let $\bX,\bX'  \in \RVleqn$, and let $G: \L \bX' \to \L \bX$ be an auxiliary special  
bijection.   Then $\bX'$ is isomorphic to   $\bX$ over $\RV^n$.  

\item    Let $\bX,\bX'  \in \RVleqn$, and let $G: \L \bX' \to \L \bX$ be an elementary  
bijection.   Then $\bX'$ is strongly isomorphic  to  a blowup of $\bX$.

 \item Let $\bX,\bX'  \in \RVleqn$, and let $G: \L \bX' \to \L \bX$ be a special 
bijection.   Then $\bX'$ is strongly isomorphic to  an iterated blowup of $\bX$. 

 %  \item 

\item  Assume $\T$ is effective.   If $\bY  \to \bX$ is an $\RV$-blowup,   there exists $\bY'$ strongly isomorphic
to $\bY$ over $\bX$  and  an elementary bijection $c: \L \bY' \to \L \bY$ lying over $\bY' \to \bY$.
%If $\bY \to \bX$ is an iterated $\RV$-blowups, there 

\end{enumerate}

 \>{lem}

\prf

(1)  Clear, using \lemref{resolve-cor}.

(2)  This is a special case of (1).

(3)  Clear from the definitions.
 
 (4) Clear from (1-3).

(5)   It suffices to consider elementary blowups; we use the notation in the definition
there.  
  So   $f_n(x) \in \acl(f_1(x),\ldots,f_{n-1}(x))$ for $x \in \phi(Y)$.  By    effectiveness and \lemref{red1c2} there exists a definable function 
$s(x,y_1,\ldots,y_{n-1})$ such that if $\rv(y_i) = f_i(x)$ for $i=1,\ldots,n-1$, then 
$\rv s(x,y) = f_n(x)$.  This $s$ is the additional data needed for an elementary bijection.
%The final statement   follows from this together with \propref{rvlift}.    
\eprf

\lemm{fibrations}   Let $\bX=(X,f), \bX'=(X',f')  \in \RVleqn$, and let $h: X \to W \subseteq  \RV^*,h':  X' \to W $
be  definable maps.  Let $X_c  = h \inv (c), \bX_c = (X_c,f|X_c)$ and similarly 
$\bX'_c$.  If $(\bX_c,\bX'_c) \in \Isp (\RV_{c}[n,\idot])$ then 
$(\bX,\bX') \in \Isp$.
 \>{lem}
 
 \prf \lemref{collect} applies to 
$\RVvolg$-isomorphisms, and hence using \remref{rvblowup-rem}   
also to blowups.   It also applies to $\RVleqn$-isomorphisms; hence to $\Isp$-equivalence.   \eprf

\lemm{isp1}  If $(\bX,\bY) \in \Isp$ then $\L \bX \iso \L \bY$.
 \>{lem}
\prf Clear, since $\L   [1]_1$ is the unit open ball around $1$,   
$\L( [\RVp]_1$ is the punctured unit open ball around $0$, and $\L 1_0   = \{0\}$. 
\eprf

\ssec{The kernel of $\L$}

\<{defn}  $\VFR[k,l,\idot]$ is the set of pairs $\bX = (X,f)$, 
with  $X  \subseteq \VF^k \times \RV^*$, $f: X \to \RVi^l$, and such that $f$ factors through 
the projection $pr_{\RV}(X)$ of $X$ to the $\RV$-coordinates.  $\Isp$ is the equivalence relation on $\VFR[k,l,\idot]$:  
$$(X,Y) \in \Isp \iff  (X_a,Y_a) \in \Isp(\T_a) \text{ for each } a \in \VF^k$$
$\SG \VFR$ is the set of equivalence classes.
\>{defn}

By the usual compactness argument, if $(X,Y) \in \Isp$ then there are uniform formulas
demonstrating this.  The relative versions of Lemmas % \ref{isp1}, % \ref{isp2},
\ref{generators},  \ref{diamond} follow.

If $\bU=(U,f) \in \VFR[k,l,\idot]$, and for $u \in U$ we are uniformly given $\bV_u =(V_u,g_u) \in \VFR[k',l',\idot]$ 
we can define a sum $\sum_{u \in U} \bV_u \in \VFR[k+k',l+l',\idot]$: it is the set 
$\du_{u \in U} \V_u$, with the function $(u,v) \mapsto (f(u),g_u(v))$.  When necessary, we denote this operation
$\sum^{(k,l;k',l')}$.  The special case $k=l=0$
  is understood as the default case.

By \propref{rvconj1}, the inverse of $\L: \RV[1,\idot] \to \VF[1,\idot]$ induces an isomorphism
$I_1^1:\SG \VF[1,\idot] \to \SG \RV[1,\idot] / \Isp$.
$$I([X]) = [Y]/\Isp \iff  [\LY ] = [X]$$

Let $J$ be a  finite set of $k$ elements.  For $j \in J$, let $\pi^j: \VF^k \times \RV^*  \to \VF^{J-\{j\}} \times \RV^*$ be the projection forgetting the $j$'th 
$\VF$ coordinate.     We will write $\VF^k,\VF^{k-1}$ for $\VF^J, \VF^{J-\{j\}}$ respectively
when the identity of the indices is not important.  

Let $\bX = (X,f) \in \VFR[k,l,\idot]$.  By  assumption, $f$ factors through $\pi^j$.
We view the image $(\pi^j X,f)$ 
as an element of $\VFR[k-1,l,\idot]$.   Note that
  each fiber of $\pi^j$ is in $\VF[1,\idot]$.   

%Recall that $I^1_1$, a morphism of semigroups, lifts to a map $\VF[1,\idot] \to \RV[1,\idot]$.  
Relativizing $I^1_1$ to  $\pi^j$, we obtain a map
 %$\VFR[k,l,\idot] \to   \VFR[k-1,l+1,\idot] $ 
  %(namely $I^j(X) = \du_{y \in (\pi^j,f)(X)} ^{(k-1,l; 0,1)} I_1^1(X_y)$)
% and in particular
$$I^j = I^j_{k,l} : \VFR[k,l,\idot] \to   \SG \VFR[k-1,l+1,\idot]/\Isp$$
%(Actually only the semigroup morphism $I^1_1$ is needed to define $I^j$, not the lifting, since
%Relatively $\Isp$ equivalent sets are in particular $\Isp$ equivalent.)

\lemm{isp4} Let $\bX=(X,f), \bX' = (X',f') \in \VFR[k,l,\idot]$.  
 \begin{enumerate}
  \item $I^j$ commutes with maps into $\RV$:  if $h: \bX \to W \subseteq \RV^*$ is definable,
  $\bX_c = h \inv (c)$, then $I^j(\bX) = \sum_{c \in W} I^j(\bX_c)$
% \item If $(\bX,\bX')$ are effectively isomorphic over $\VF^k$, then $(I^j(\bX),I^j(\bX')) \in \Isp$.
  \item    If $([\bX],[\bX']) \in \Isp$ then $(I^j(\bX),I^j(\bX')) \in \Isp$.
  \item $I^j$ induces a map $\SG \VFR[k,l,\idot] /\Isp \to \SG \VFR[k-1,l+1,\idot]/\Isp$.
\end{enumerate}            \>{lem}
 
\prf (1)  This reduces to the case of $I^1_1$, where it is an immediate consequence
of the uniquness, and the fact that $\L$ commutes with maps into $\RV$ in the same sense.

(2)  All equivalences here are relative to the $k-1$ coordinates of $\VF$ other than $j$,
so we may assume $k=1$.  

For $a \in \VF$, $([\bX_a]  ,[\bX'_a]) \in \Isp(\T_a)$.  By stable embeddedness of $\RV$,
there exists $\alpha = \alpha(a) \in \RV^*$ such that $\bX_a,\bX'_a$ are $\T_\alpha$-definable and 
$([\bX]_a  ,[\bX']_a) \in \Isp(\T_\alpha)$.  Fibering over the map $\alpha$ we may assume
by (1) and \lemref{fibrations} that $\alpha$ is constant; so for some $W \in \VF[1], \bY,\bY' \in \RV[l,\idot]$,
we have
$\bX = W  \times \bY, \bX' = W \times \bY'$, and $([\bY],[\bY']) \in \Isp$.  Then
$I^j(\bX) = I^j(W) \times \bY, I^j(\bX') = I^j(W) \times \bY'$ and the conclusion is clear.   

(3) by (2). 

\eprf

 \<{lem}\lbl{fubini}  Let $\bX =(X,f), X \subseteq \VF^J \times \RV^\infty, f: X \to \RV^l$.   If $j \neq j' \in J$ then $I^j I^{j'} = I^{j'} I^j:   \SG \VFR[k,l,\idot] /\Isp \to \SG \VFR[k-2,l+2,\idot]/\Isp$.
\>{lem}

\prf  We may assume $S=\{1,2\}, j=1,j'=2$, since all is relative to $\VF^{S \m \{j,j'\}}$.  
By \lemref{isp4} (1)  it suffices to prove the statement for each fiber of a given definable
map into $\RV$.  
   
Hence we may assume  $X \subseteq \VF^2$ and  $f$ is constant; and by \lemref{cell2}, we can assume $X$ is a basic 2-cell:  
$$X= \{(x,y): x \in X_1, \rv(y-G(x))=\alpha_1\} \ \ \ \ \ \ \ \ \ \ X_1= \rv \inv (\delta_1) + c_1 $$
The case where $G$ is constant is easy, since then $X$ is a finite union of  rectangles. 
Otherwise $G$ is invertible, and by
 niceness of $G$ we can also write:
$$X =  \{(x,y): y \in X_2, \rv(x-G \inv (y))=\beta \},  \ \ \ \ \ \ \ \ \ \ X_2= \rv \inv (\delta_2) + c_2 $$ 
We immediately compute:  
$$I_2I_1(X) =(\d_1,\alpha_1), \ \ \ I_1I_2(X) =   (\a_2,\d_2)$$
%$$I_2I_1(X) = \rv \inv (\d_1) \times \rv \inv (\alpha_1), \ \ \ I_1I_2(X) =   \rv \inv (\a_2)
%\times \rv \inv(\d_2)$$
  Clearly $[(\d_1,\alpha_1)]_2 = [(\a_2,\d_2)]_2$.  

    \eprf

\<{prop} \lbl{change-of-variable-1} Let $\bX,\bY \in \RV[\leq n,\idot]$.  If $\L \bX ,\L \bY$ are
%effectively
 isomorphic, then $([X],[Y]) \in \Isp$.  \>{prop}

\prf  
define $I = I_1 \ldots I_n: \VFni = \VFR[n,0,\idot]  \to \VFR[0,n,\idot] =  \RV[\leq n,\idot]$.  Let   $V \in \VFni$.

\Claim{1.}   If $\si \in Sym(n)$   then $I = I_{\si(1)} \ldots I_{\si(n)}$.  

\proof  We may assume $\si$ just permutes two adjacent coordinates, say $2,3$ out of 
$1,2,3,4$.  Then $I=I_1I_2I_3I_4 = I_1 I_3I_2 I_4 $ by \lemref{fubini}. 

\Claim{2.}   When $F: V \to F(V)$ is a relatively unary bijection,  
we have $I(V)= I(F(V))$

\proof  By Claim 1 we may assume $F$  is relatively unary with respect to $\pr^n$.
So $F(V_a) = F(V)_a$ where $V_a, F(V)_a$ are the $\pr^n$-fibers.  By definition 
of $I^1_1$ we have $I^1_1(V_a) = I^1_1(F(V)_a) \in \RV[1,\idot](\T_a)$;
but  by definition of $I^n$ , $I_n(V)_a  = I^1_1(V_a)$.  So $I^n(V)=I^n(F(V))$
and thus $I(V)=I(F(V))$.   

\Claim{3.}  When $F: V \to F(V)$ is any definable bijection, $I(V) = I(F(V))$.  

\proof Immediate from Claim 2 and \lemref{decompose-unary}.

Now turning to  the statement of the  Proposition, 
assume $\L \bX, \L \bY$ are isomorphic.  
%since $\L \bX, \L \bY$
%are effectively isomorphic, we may assume by \propref{resolve} that there
%exists a definable isomorphism $F: \L \bX \to \L \bY$.  
We compute inductively that $\L( \bX) = [\bX]$.  
By Claim 3, 
$ [\bX] = I( \L \bX) = I(\L \bY) = [\bY]$.    \eprf

\<{notation}  Let $\L^*: \SG(\VF) \to \SG(\RV[*])/Isp$ be the inverse map to $\L$. \>{notation}

 \<{remark} \lbl{cov-2}  \rm  When $\T$ is $\rv$- effective, one can restate the conclusion of  
\propref{change-of-variable-1} as follows:  if $X,Y \in \VFni$ are \rvinv and $F:
X \to Y$ is a definable bijection, then there exist special bijections $X' \to X$
and $Y' \to Y$ and an \rvinv -definable bijection $G:X' \to Y'$.  (This follows from 
\propref{change-of-variable-1}, \propref{rvlift}, \lemref{diamond}, and \lemref{sp-bl}.
The   effectiveness hypothesis is actually unnecessary here, as will be seen in the proof
of \propref{change-of-variable-m}.  
Perhaps Question \ref{q}  can be answered simply by tracing the connection between $F$ and $G$  through the proof.
\>{remark}

 \>{section}  %descent - isomorphisms

 \<{section}{Definable sets over $\VF$ and $\RV$: the main theorems}
 
In stating the theorems we 
restrict attention to $\VF[n]$, i.e. to definable subsets of varieties, though the proof was
given more generally for $\VFni$  (definable subsets of $\VF^n \times \RV^*$.)

\<{subsection}{Definable subsets of varieties}  

Let $\T$ be $\V$-minimal.
We will look at the category of definable subsets of varieties, and
definable maps between them.  The results will be stated for $\VF[n]$; analogous statements
for $\VFni$ are true with the same proofs.

We define three variants of the sets of objects.   $\VF '' [n]$ is the category of  $\leq  n$-dimensional definable sets over $\VF$, i.e. of definable subsets of $n$-dimensional varieties.
Let $\VF [n] $ be the category of definable subsets   $X \subseteq \VF^n \times \RV^*$  such that the projection $X \to \VF^n$ has  finite fibers.  
  $\VF'[n]$ is the category of definable subsets $X$
of $V \times \RV^*$, where $V$ ranges over all $\VF(A)$-definable sets of dimension $n$,   
$m \in \Nn$,  such that the projection $X \to V$ is finite-to-one.  
   $\VF,\VF',\VF''$ are the unions over all $n$.   In all cases, the morphisms $\Mor(X,Y)$ are the definable functions $X \to Y$.

\<{lem}  \lbl{loy} The natural inclusion of $\VF[n]$ in $\VF'[ n]$ is an equivalence.   If 
$\T$ is effective, so is the inclusion of $\VF''[ n]$ in $\VF'[ n]$.   \>{lem}

\<{proof}  We will omit the index $\leq n$.  
 The inclusion is fully faithful by definition, and we have to show that
it hits every $\VF '$ -isomorphism type; in other words, that any definable
$X \subseteq (V \times \RV^m)$ is definably isomorphic to some $X' \subseteq \VF^n \times \RV^{m+l}$, for some $l$ (with $n = \dim(V)$.)  Definable isomorphisms can 
be glued on pieces, so we may assume $V$ is affine, and admits a finite-to-one
map $h: V \to \VF^m$.  By   \lemref{finite}, each fiber $h^{-1}(a)$ is $A(a)$-definably isomorphic to some $F(a) \subseteq \RV^l$.
By compactness, $F$ can be chosen uniformly definable,
$F(a) = \{y \in \RV^l: (a,y) \in F \}$ for some definable
 $F \subseteq   \VF^m \times \RV^l$; and there exists a definable isomorphism
 $\beta: V \to F$, over $\VF^m$.  Let $\a(v,t) = (\beta(v),t)$, $X' = \alpha(X)$.  
 
 Now assume $\T$ is effective.
 %$A \subseteq \acl(\VF(A) \union \G(A))$.    
 Let $X \in \Ob \VF'$;
  $X \subseteq V \times \RV^m$, $V \subseteq \VF^n$,
  such that the projection $X \to V$ has   finite fibers.  
 Then by effectivity, for any $v \in V$ (over any extension field), 
if $(v,c_1,\ldots,c_m) \in X$ then each $c_i$, viewed as a ball, has a point defined over $A(v)$.
Hence the partial map $V \times \VF^m \to X$, $(v,x_1,\ldots,x_m) \mapsto (v,\rv(x_1),\ldots,\rv(x_m))$ has an $A$-definable section; the image of this
section is a subset $S$ of $V \times \VF^m$, definably isomorphic to $X$;
and the Zariski closure $V'$ of $S$ in $V \times \VF^m$ has dimension $\leq \dim(V)$.
\>{proof}

The following definition and proposition apply both to the category of definable sets,
and to the definable sets with volume forms. 

\<{defn} \lbl{eff-iso} $X,Y$ are {\em effectively isomorphic} if 

for any effective $A$, $X,Y$ are definably isomorphic in $\T_A$.  If $\SGe(\VF)$ is the semi-ring of effective isomorphic classes
of definable sets.  $\K(\VF)$ is the correpsponding ring.  Similarly $\SGe(\VF[n])$, etc. 
\>{defn}

  Over an effective base, in particular if $\T$ is effective over
any field-generated base,  effectively isomorphic is the same as isomorphic.  But \exref{tr2m}
shows that this is not so in general.

\<{prop}   \lbl{eff}  Let $T$ be $\V$-minimal, or a finitely generated extension of a $\V$-minimal theory.
The following conditions are equivalent:  Let $X,Y \in \VF[n]$.  \begin{enumerate}
  \item $[\L^* X ] = [\L^* Y]$ in $\SG(\RV[\leq n])/\Isp[\leq n]$
  \item There exists a definable family ${\mathcal F}$ of definable bijections 
  $X \to Y$, such that for any effective structure $A$, $F(A) \neq \emptyset$.
  \item $X,Y$ are effectively isomorphic.
  \item $X,Y$ are definably isomorphic over any $A$ such that  $\VF^*(A) \to \RV(A)$
  is surjective. 
  \item For some finite $A_0 \subseteq \RV(<\emptyset>) $, 
  $X,Y$ are definably isomorphic over any $A$ such that 
  $A_0 \subseteq \rv(\VF^*(A))$.
\end{enumerate} 
\>{prop}

\prf  

(1) implies (5):  By \propref{rvlift} (\ref{rvlift-m} in the measured case), the given isomorphism
$[\L^* X ] \to [\L^* Y]$ lifts to an isomorphism $\L \L^* X \to \L \L^* Y$; since $\T_A \supseteq \ACVF_A$,
this is also a $\T_A$ isomorphism;   
 it can be composed with the isomorphisms $X \to \L \L^*X, Y \to \L \L^*Y$.

(2) implies (3), (5) implies (4) implies (3)   trivially.

 (3) implies (1),(2):  Let $E_{eff}$ be as in \propref{resolve}.  By (3), $X,Y$
 are $E_{eff}$-isomorphic.  By \propref{change-of-variable-1},
  $[\L^* X] = [\L^* Y]$ in $\SG(\RV_{E_{eff}}[*])/\Isp$.
 But $\RV(E_{eff}), \G(E_{eff}) \subseteq \dcl(\emptyset)$, so every
 $E_{eff}$-definable relation on $\RV$ is definable; i.e. $\RV_{E_{eff}}, \RV$ are the same structure.  So (1) holds.  

 Now by assumption,   there exists an 
 $E_{eff}$-definable
bijection $f': X \to Y$.  $f'$ is an $E_{eff}$-definable element of a definable family
$ {\mathcal F}$ of definable bijections 
  $X \to Y$.  Since this family has an $E_{eff}$-point, and $E_{eff}$ embeds into any effective $B$,
  it has a $B$ point too.  Thus (3) implies (2).

\eprf

\>{subsection}
\<{subsection}{Invariants of all definable maps}

 Let $[X]$ denote the class of $X$ in $\SGe(\VF[n])$.  
 
 \<{prop} \lbl{summ}    Let $\T$ be $\V$-minimal.
 There exists a canonical isomorphism of Grothendieck semigroups
 $$\ints:  \SGe(\VF[n]) \to \SG(\RV[\leq n]) / \Isp[\leq n]$$
satisfying
   $$\ints [X] = W / \Isp[\leq n] \iff [X]=[\L W] \in \SGe(\VF[n])  $$
 \>{prop}

\prf  Recall \defref{LL}.  Given $\bX=(X,f) \in \Ob \RV[k]$ we have $\L \bX \in \Ob \VF[k] \subseteq \Ob \VF[n]$.  If $\bX,\bX'$ are isomorphic,then by \propref{rvlift},
$\L \bX, \L \bX'$ are effectively isomorphic.  Direct sums are clearly respected,
so we have a semigroup homomorphism $\L: \SG(\RV[\leq n]) \to \SGe(\VF[n]) $.
It is surjective by \propref{tr2}.   By \propref{eff}, the kernel is 
precisely $\Isp[\leq n]$.  Inverting, we obtain $\ints$. 
\eprf

\<{defn} \lbl{ispd}  
Let $\SG\VF[n] / (\dim < n)$ be the Grothendieck ring of the category of 
definable subsets of $n$-dimensional varieties, and essential bijections between them.
Let $\Ispd$ be the congruence on $RV[n]$ generated by pairs $(X, X \times \RVp)$
(where $X \subseteq \RV^*$ is definable, of dimension $<n$.)  \>{defn}

\<{cor} \lbl{summ-mod-n-1}  $\ints$ induces an isomorphism $\SGe(\VF[n])/(\dim < n) \to \RV[n] / \Ispd$ 
\qed \>{cor}

  \<{cor}  \lbl{summ-c}  Let $A,B \in \RV[\leq n]$.
  Let $n'>n$, and let $A_{N},B_N$ be their images in $\RV[\leq N]$.
  If $(A_N,B_N) \in \Isp[\leq N]$ then $(A,B) \in \Isp[\leq n]$. 
  \>{cor}
 \proof  By \propref{summ}, $(A,B) \in \Isp[\  \leq n]$ iff $\L A, \L B$ are definably 
 isomorphic; this latter condition does not depend on $n$.  

   \qed
  
 Putting \propref{summ} together for all $n$ we obtain:  
 
 \<{thm}  \lbl{summ-inf}  Let $\T$ be $\V$-minimal.
  There exists a canonical isomorphism of filtered semirings
 
 $$\ints:  \SG(\VF) \to \SG(\RV[*]) / \Isp$$
 Let $[X]$ denote the class of $X$ in $\SG(\VF)$.
Then 
 $$\ints  [X] = W / \Isp  \iff [X]=[\L W] \in \SGe(\VF)  $$ 

\qed \>{thm}
 
On the other hand, using the Grothendieck group isomorphisms of \propref{summ} and passing
to the limit:  

\<{cor}   \lbl{summ-inf-K0}  Let $\T$ be $\V$-minimal.  The isomorphisms   of 
\propref{summ}  inducing an isomorphism of  Grothendieck groups:
 $$\int^{K}   :  \Ke(\VF[n]) \to \K(\RV[n])$$
The isomorphism  $\ints$ of \thmref{summ-inf}
induces  an injective  ring homomorphism 
$$\int^{\K}: \Ke(\VF) \to   \K(\RV)[J^{-1}]$$
where $J = \{1\}_1 - [\RVp]_1 \in \K(\RV)$.
\>{cor}

\<{proof}  We may work over an effective base.   With subtraction allowed, the generating relation of $\Isp$ can be read as:
$[\{1\}]_0 = \{1\}_1 - [\RVp]_1 := J$; so that the groupfication of 
$\SG(\RV[\leq n]) / \Isp[\leq n] $ is isomorphic to $\K(\RV[n])$, via the   embedding
of $\SG(\RV[n])$ as a direct factor in $\SG(\RV[\leq n])$.  Thus the groupification of the homomorphism
of \thmref{summ-inf} is a homomorphism
$$\int^{\K}: \K(\VF) \to \lim_{n \to \infty} \K(\RV[n])$$
where the direct limit system maps are   given
by $[X]_d \mapsto ([X]_{d+1} - ([X]_d \times (\RVp))) = [X_d] J$.  This direct limit
embeds into $\K(\RV)[J^{-1}]$ by mapping $X \in \K(\RV[n])$ to $XJ^{-n}$.  \>{proof}

\>{subsection}

  \ssec{Definable   volume forms: $\VF$} 
  We will now define the category 
  $\VFm[n]$ of ``$n$-dimensional $\T_A$-definable sets with definable volume 
forms, up to $\RV$- equivalence'' and the same up to $\G$-equivalence.  
We will represent the forms as functions to $\RV$, that transform in the way volume
forms do. 
 
By way of motivation:  in a local field with an absolute value, a top differential form $\omega$ induces
a measure $\abs{d \omega }$.  For a regular isomorphism $f: V \to V'$  have $\omega = h f ^* \omega' $ for a unique $h$,
and $f$ is   measure preserving between $V,\abs{\omega})$ and $(V',\abs{\omega'})$ iff
$\abs{h}=1$.   %For open subsets $V \subseteq K^n$, one conventionally takes $\omeg

%The category of smooth varieties with top forms, up to measure preserving isomorphisms, 
%can also be presented 

We do not work with an absolute value into the reals, but instead define the analogue
using the map $\rv$ or, a coarser version, the map $\val$ into $\G$.   When $\G=\Zz$,
the latter is the  is the usual practice in Denef-style motivic integration.  Using $\rv$
 leaves room for considering an absolute value on the residue
field, and iterating the integration functorially when places are composed, for instance
$\Cc((x))((y)) \to \Cc((x)) \to \Cc$.    This functoriality will be described in a sequel.  

In the definition below, the words ``almost every $y \in Y$'' will mean:  for all $y$ outside
a set of $\VF$-dimension $< \dim_\VF(Y)$.

\<{defn}   \lbl{VFmcat}
    $\Ob \VFm[n,\idot]$ consists of  pairs  $(Y,\om)$,
where $Y$ is a definable subset of $\VF^n \times \RV^*$, and   $\om: Y \to \RV$ is a definable
map.  
%of the $\RV^m$-coordinates. 
 A morphism  $(Y,\om) \to (Y',\om')$ is a definable essential 
bijection $F$ such that for almost every $y \in Y$, 
$$ \om (y)=   \om' (F(y)) \cdot \rv( \JVF F(y)) $$
(We will say:  ``$F:(Y,\om) \to (Y',\om')$ is  measure preserving''.)

  $\mu_{\Gamma} \VF[n,\idot]$ is the category of pairs $(Y,\om)$ with $\om: Y \to \G$
a definable function %of the $\RV$ coordinates; 
 A morphism  $(Y,\om) \to (Y',\om')$ is a definable essential
bijection $F: Y \to Y'$ such that for almost every $y \in Y$,
$$ \om (y)=   \om' (F(y)) + \val ( \JVF F(y)) $$ 
(``$F:(Y,\om) \to (Y',\om')$ is $\G$-measure preserving''.)

  $\VFmn, \VFgn$ are the full subcategories of $\VFm[n,\idot],\mu_{\Gamma}\VF[n,\idot]$ (respectively)
   whose objects admit a finite-to-one map to $\VF^n$.

 \>{defn}
 In this definition, let $t_1(y),\ldots,t_n(y)$ be the $\VF$-coordinates of $y \in Y$.   One can think of the form as $\om(y) dt_1 \cdot \ldots \cdot dt_n$.

Note that   $\VFvol$ of \defref{vfvol} is isomorphic to the full subcategory of $\VFm$ whose objects are pairs $(Y,1)$.

 \<{remark}  \rm When $\T$ is $\V$-minimal and effective, the data $\om$ of an object $(Y,\om)$
 of $\VFmn$  can be written as $\rv \circ \Phi$  
 for some $\Phi: Y \to \VF$.  (Write $\om = \bar{\om} \circ \rv \circ F$ for some $F$, and
 use \propref{rvlift} to lift $\bar{\om} $ to some $G$, so that $\om = \rv \circ G \circ F$.)   It is thus possible to view 
    $\om$   as the $\RV$-image (resp. $\G$-image) of a definable volume form on $Y$. 
    One could equivalently take $\om$ to be a definable section of $\Lambda^n TY / (1+\Mm)$,
    where $TY$ is the (appropriately defined) tangent bundle, $\Lambda^n$ the $n$-th
    exterior power with $n=\dim(Y)$.

    For $\VF_{\G}$ the category we take is slightly more flexible than taking varieties with 
    absolute values of volume forms, even if $\T$ is $\V$-minimal and effective, in that 
    expressions such as $\int |\sqrt{x}| dx$ are allowed.  
    %If we demand that
    %$\om$ lifts to $\VF$ and to $\RV$ we obtain an isomorphism of Grothendieck groups as well.
 \>{remark}

In either of these categories, one could restrict the objects to bounded ones: 
\<{defn} \lbl{boundedsets}
Let   $\VFmb[n]$ be the full subcategory of  $\VFmn$ whose objects are bounded definable sets, with bounded definable 
forms $\om$. Similarly one defines   $\VFmbg$.
\>{defn}

Here {\em bounded} means that there is a lower bound on the valuation of any coordinate
of any element of the set.  A similar definition applies in $\RV$ and $\RVm$.

 Note that if an object of $\VFmn$ is $\VFmn$-isomorphic to an object
of $\VFmb[n]$, it must lie in $\VFmb[n]$.

\ssec{Definable volume forms: $\RV$}

We will define a category $\RVm[n]$
of   definable subsets of $(\RV)^m$, with additional data that can be viewed as a volume form. 
Unlike $\VFm[n]$, in   $\RVm[n]$ subsets of dimension $<n$ are {\em not} ignored:   a point of $\RV^n$ corresponds
to an open polydisc of $\VF^n$, with nonzero $n$-dimensional volume.  

In particular, the Jacobian of a morphism needs to be defined at every point, not just
away from a smaller-dimensional set.   However, in accord with  \lemref{rvlift-m},
it may be modified by $\k^*$-multiplication on a smaller-dimensional set.  

% Note that it will not do to take simply the lower-dimensional definition of the Jacobian.
% E.g. if $\Aa^1 $ is cut into $G_m$ and $(0)$,
%the map $(x+a)^2$ restricts to one with Jacobian $1$ at $0$, but this is the wrong multiplier;
%should be $2a$. 

%  Thus in this case if we allowed cutting at $0$, some additional data
%would be needed.  We will instead give up on unrestricted cut-and-paste. and insist
%that the sets are open in the sense above, and that the 
%morphisms are differentiable.  

\<{defn}  \lbl{RVm}

 The  objects of $\RVm[n]$ are definable triples $(X,f,\om)$,  $X \subseteq \RV^{n+m}$,
 $f: X \to \RV^n$ finite-to-one,  
 and $\om: X \to \RV$.
 
% $\RVm[n]$ has a notion of disjoint union.  
 We   define a multiplication 
 $\RVm[n] \times \RVm[n'] \to \RVm[n+n']$ by $(X,f,\om) \times (X',f',\om') = (X \times X', f \times f', \om \cdot \om')$.  Here $\om \cdot \om'(x,x')  = \om(x)\om'(x')$.

 Given $\bX = (X,f,\om)$, we define an object $\L \bX $ of $\VF[n]$; namely
 $(\L X, \L f, \L \om)$ where $\L X = X \times_{f,\rv} (\VF^\times)^n$, $\L f (a,b) = f(a,\rv(b))$,
 $\L \om (a,b) =\om(a,\rv(b))$.  (Sometimes we will write $f,\om$ for $\L f, \L \om$.)

 A  morphism $\a: \bX = (X,f,\om) \to \bX' = (X',f',\om')$ is a   definable bijection $\a: X \to X'$, 

 such that 
 $$ \om (y)=   \om' (\a(y)) \cdot \rv( \JRV( \a) (y)) \text{\em \  for almost all } y  $$

 where ``almost all'' means: away from a set $Y$ with $\dim_{\RV}(f(Y)) <n$; and
 $$\valr  \om (y) + \sum_{i=1}^n \valr f_i(y) =       \valr \om' (\a(y)) +   \sum_{i=1}^n \valr f'_i(\a y)  \text{ \em \    for all }y  $$
 
  The objects of $\RVmg[n]$ are  triples $(X,f,\om)$, with $f: X \to \RV^n$, $\om: X \to \G$.
A morphism $\a: (X,f,\om) \to (X',f',\om')$ is a definable bijection $\a: X \to X'$  such that 
$\valr  \om (y) + \sum_{i=1}^n \valr f_i(y) =       \valr \om' (\a(y)) +   \sum_{i=1}^n \valr f'_i(\a y)  \text{ \em \    for all }y  $.   Disjoint sums and products are defined as for $\RVm$.
 
  $\RESmg[n]$ is the full subcategory of $\RVmg[n]$ with objects $(X,f,\om)$, 
  such that $\valr (X)$ is finite.  In this case, $\om$ takes finitely many values too.
 \>{defn}

$\SGe \RVm[n]$   is the Grothendieck semigroup of $\RVm[n]$ with respect to
effective isomorphism.  $\SGe \RVm$   is   the direct sum $\oplus_n \SGe \RVm[n]$;
it clearly inherits a semiring structure from Cartesian multplication, 
$(X,f,\om) \times (X',f',\om') = (X \times X', (f,f'), \om \cdot \om')$.

The morphisms of $\RVmg[n]$ are called $\G$- measure preserving.

The category $\RVvol[n,\idot]$ of \defref{rvvol} is isomorphic  to   the full subcategory whose
objects have $\om=1$.

{\bf Remark}  
The semiring $\SGe \RVvol$ is naturally a subsemiring of  $\SGe \RVm$. The latter is obtained by inverting $[\{a\}]_1$ for $a \in \RV$ and
taking the 0'th graded component.  This process is needed in order to identify integrals of 
functions in $n$ variables with volumes in $n+1$ variables.
Thus as semirings they are closely related.
But if  the dimension grading is taken into account, the subsemiring of $\RV$-volumes
  contains finer information connected to  integrability of forms.

\ssec{The kernel of $\L$ in the measured case}

The description of the kernel of $\L$ on the semigroups of definable sets with 
volume forms is essentially the same as for definable sets.  
We will now run through the proof, indicating the modficiations.  The principal
change is the introduction of a category with fewer morphisms, defined
not only with reference to $\RV$ but also to $\VF$.  For effective bases,
the category is identical to $\RVm$, so it will be invisible in the statements of the main theorems; 
but during the induction in the proof, bases will not in general be effective and
the mixed category introduced here has better properties.  

Both the introduction of the various  intermediate categories and the 
repetition of the proof  would be unnecessary if we had a positive answer to Question \ref{q}.
In this case the proof of  \lemref{rvconj1-m} would immediately lift to higher dimensions.
Indeed the   characterization 
 of the kernel of the map $\L$ on Grothendieck groups would be uniformized 
 not only for the categories we consider, but for a range of  categories carrying
 more structure.  
 
  The integer $n$ will be fixed in this subsection.
  
 \<{lem} \lbl{carryform}  Let $(X,\om) \in \Ob \VFm [n,\idot]$, $Y \in \Ob \VF [n,\idot]$, and let $F: Y \to X$ 
 be a definable bijection.  \begin{enumerate}
  \item There exists  $\psi: Y \to \RV$  such that
 $F:  (Y,\psi) \to (X,\om)$ is  measure preserving.  
  \item 
$\psi$ is essentially unique in the
 sense that if $\psi'$ meets the same condition,
 then $\psi,\psi'$ are equal away from a subset of $X$ of smaller dimension.   
  \item Dually, given $F, X,Y,\psi$
 there exists an essentially unique $\om$ such that  $F:  (Y,\psi) \to (X,\om)$ is 
  measure preserving. 
  \item 
  \lemref{decompose-unary} applies to   $\VFm[n,\idot]$ and to $\VFm[n]$
  %, as well as to $\mu \VF [n, \idot]$.
\end{enumerate}

   \>{lem}

\prf   (1,2) Let $\psi(y) = \om (\a(y)) \cdot \rv( \JRV( \a) (y))$.  By definition of $\VFm$
this works, and is the only choice ``almost everywhere''. 
(3)  follows from the case of  $F \inv$. 
(4) 
Now let $\bX,\bY \in \Ob \VFm[n]$ and let $F \in \Mor_{\VFm[n]}(X,Y)$.  We have $\bX=(X,\omega_X),\bY=(Y,\omega_Y)$
for some $X,Y \in \Ob \VF[n]$ and $\omega_X: X \to \RV, \omega_Y: Y \to \RV$.  By
  \lemref{decompose-unary} there exist $X=X_1, \ldots,X_n=Y \in \Ob \VF[n]$ and
  essentially unary $F_i: X_i \to X_{i+1}$ with $F= F_{n-1} \circ \cdots \circ F_1$.
 Let $\omega_1 = \omega_X$, and inductively let $\omega_{i+1}$ be such that
 $F_i \in \Mor_{\VFm[n]} ((X_i,\omega_i),(X_{i+1},\omega_{i+1})$.  Then 
 $F \in \Mor_{\VFm[n]}((X,\omega),(Y,\omega_n))$.  By the uniqueness it follows that
 $\omega_Y,\omega_n$ are essentially equal.  \eprf

 \<{defn}     
 Given $\bX,\bY \in \Ob \RVmn$
 call a definable bijection $h: X \to Y$ {\em liftable}  if there exists $F \in \Mor_{\VFm[n,\idot]}(\LX,\LY)$
with $\rho_Y F = h \rho_X$.

 Let $\fC= \RVlmn$ be the subcategory of $\RVmn$
consisting of all objects and liftable morphisms.  
\>{defn}

By \propref{rv-compat}, liftable morphisms must preserve the volume forms, so
$\fC$ is a subcategory of $\RVmn$.

Over an effective base, $\fC = \RVmn$ (\lemref{rvlift-m}), and the condition of existence of $s$ 
in \defref{rvblowup-m} (1) below  is equivalent to:  $f_n(y) \in \acl(f_1(y),\ldots,f_{n-1}(y))$.

 \<{defn}  \lbl{rvblowup-m} 
 \begin{enumerate}
  \item  Let $\bY =(Y,f,\om) \in \Ob \RVmn$ be such that
  % $f_n(y) \in \acl(f_1(y),\ldots,f_{n-1}(y))$,
  there exists $s: \bY \times_{f_1,\ldots,f_{n-1}} \VF^{n-1} \to \VF$
with $\rv (s (y,u_1,\ldots,u_{n-1})) = f_n(y)$.  
  Let $Y' = Y \times \RVp$.  For $(y,t) \in Y'$, define $f'=(f'_1,\ldots,f'_n)$ by: $f'_i(y,t) = f_i(y)$ for $i<n$, 
$f'_n(y,t) = tf_n(y)$.  Let $\om'(y,t)=\om(y)$.  Then $\widetilde{\bY} = (Y',f',\omega')$ is an {\em elementary blowup} of $\bY$.  It comes
with the projection map $Y' \to Y$.
  \item 
Let $\bX=(X,g,\om) \in \Ob \RVm[n,\idot]$, $X= X' \du X''$, $g' = g | X'$, $g'' = g | X''$, $\om' = \om | X',
\om'' = \om | X''$, 
 and let $\phi:   \bY \to (X',g',\om')$ be a $\RVlmn$-isomorphism.  Then the  {\em $\RV$- blowup}   $\widetilde{\bX}_{\phi}$ 
is defined to be $\widetilde{\bY} + (X'',g'',\om'') = (Y' \du X'', f' \du g'',\om' \du \om'')$.  It comes with $b: Y' \du X'' \to X$,
defined to be the identity on $X''$, and the projection on $Y'$.   $X'$ is called the {\em 
blowup locus} of $b: \widetilde{\bX}_{\phi} \to \bX$.  
 
\end{enumerate}

An {\em iterated $\RV$-blowup} is obtained by finitely many iterations of $\RV$- blowups.
 
  \>{defn}

\<{defn}    Let 
 $\Ispm[n]$ be the set of pairs $  (\bX_1,\bX_2)  \in \Ob \RVmn$ such that there exist iterated blowups $b_i:  \widetilde{\bX}_i \to \bX_i$ and
  a $\RVlmn$- isomorphism $F:  \widetilde{\bX}_1 \to  \widetilde{\bX}_2$. \>{defn}
 
When $n$ is fixed, we will simply write $\Ispm$. 
  On the other hand we will need to make explicit
 the dependence on the theory; we write $\Ispm (A)$ for the congruence $\Ispm$ of the theory
 $\T_A$.
 
When $\bX = (X,f,\om) \in \Ob \RVmn$,   $h: X \to W$ is a definable map, and $c \in W$,
define $\bX_c=(h \inv(c), f| h \inv(c), \om | h \inv(c))$.  

 Let $X_1,X_2 \in \Ob \RVmn$, and let $f_i: X_i \to Y$ be  a definable map, with $Y \subseteq \RV^*$. 
In this situation  the existence of $\mu  \RVni (<a>) $-isomorphisms between each pair of fibers
$X_1(a) , X_2(a)$ ($a \in Y$) does  not necessarily imply
 that $X_1 \iso_{\RVnm} X_2$, because of the explicit reference to dimension
  in the definition of morphisms; the dimension of the allowed exceptional
sets may accumulate over $Y$.  The definition of morphisms  for $\VFmn$  also 
allows a smaller dimensional exceptional set; but this does not create a problem 
when fibered over $W \subseteq \RV^*$, since by \lemref{vfrvdim}
 $\max _{c \in W} \dim_{\VF}(  Z_c) = \dim_{\VF}(Z)$.  Thus an $\RV$-disjoint union of $\VFmn$-isomorphisms is again a $\VFmn$-isomorphism, and it follows that the same is true for $\RVlmn$.               We thus have:   
    
\lemm{fibrations-m}   Let $\bX=(X,f,\om), \bX'=(X',f',\om)  \in \RVmn$, and let $h: X \to W \subseteq  \RV^*,h':  X' \to W $
be  definable maps.   If for each $c \in W$, $(\bX_c,\bX'_c) \in \Ispm (<c>)$ then 
$(\bX,\bX') \in \Ispm$.
 \>{lem}
 
 \prf   \lemref{collect}  applies to 
$\RVvol$-isomorphisms, and hence using \remref{rvblowup-rem}   
also to blowups.   It also applies to $\RVlmn$-isomorphisms by the discussion above,
and hence to $\Ispm$-equivalence.   \eprf

In other words, there exists a well-defined direct sum operation on $\RVmn / \Ispm$,
with respect to $\RV$-indexed systems.

 \lemm{generators-m}   (1)  Let $\bY'$  be an  elementary blowup of $\bY$.
  $\bY'$   is $\fC$- isomorphic to $\bY'' = (Y'', f'',\omega')$,
  with 
  $$Y'' = \{(y,t) \in Y \times \RVi: \valr(t) > f_n(y) \}$$
  $$ f''(y,t)=(f_1(y),\ldots,f_{n-1}(y),t,), \omega'(y,t)= \omega(y)$$

(3)  Up to isomorphism, the blowup depends only on blowup locus.
  In other words
if $X,X' ,g,g' ,\om,\om'$ are as in \defref{rvblowup-m}, and
$\phi_i: \bY_i \to (X',g',\om')$ $(i=1,2)$ are  $\RVlmn$- isomorphisms, then $\widetilde{\bX}_{\phi_1} , \widetilde{\bX}_{\phi_2}$
are $\bX$-isomorphic in $\RVlmn$.

    \>{lem}

 \prf (1) The isomorphism is given by $h((y,t))= (y,t f_n(y))$; since $f_n$ always lifts
 to a function $F_n: \L Y \to \VF$ (a coordinate projection),
  $h$ can be lifted to $H$ defined by
 $H((y,t)) = (y, t F_n(y))$.

(3)   
By assumption, 
$\phi_1,\phi_2$  lift to measure preserving maps $\Phi_i: \L \bY_i \to \L \bX'$.  
On the other hand,
by the assumption on existence of a section $s$ of $f_n$, we have measure preserving isomorphisms
$\alpha_1: \L \bY_1 \to \L \widetilde{ \bY_1}$, 
$(y,u_1,\ldots,u_n) \mapsto (y,u_1,\ldots,u_{n-1},(u_n - s)/ s)$.   Similarly
we have $\alpha_2:  \L \bY_2 \to \L   \widetilde{ \bY_2}$.  Composing we obtain
$\alpha_2 \Phi_2 \inv \Phi_1 \alpha_1 \inv: \L \widetilde{ \bY_1} \to \L \widetilde{ \bY_2}$;
it is easy to check that this is \rvinv and shows that 
$\L \widetilde{ \bY_1} , \L \widetilde{ \bY_2}$ are $\bY$-isomorphic 
in $\RVlmn$.  Taking disjoint sum with the complement $X''$ of $X'$
we obtain the result.    \eprf

{\bf Remark} There is also a parallel of \lemref{form}:   Let $\bY =(Y,g) \in \Ob \RVi[n,\idot]$, with $\dim(g(Y)) < n$; let $f: Y \to \RV^{n-1}$
be isogenous to $g$.  
  Let $h: Y \to \RV$ be definable, with $h(y) \in \acl(g(y))$ for $y \in Y$, and 
  with $\sum(g)=\sum(f)+\valr(h)$.  
Let $Y' = Y \times \RVpi$, and $f'(y,t) = (f(y), t h(y))$.   Then for appropriate $\omega'$, $\bY' = (Y', f',\omega')$ with   the projection map to $Y$
is a blowup.   This follows from \lemref{form}, and  \lemref{carryform} (3).

  {\bf Notation}  For $X \in \RVni$, $[X]=[(X,1)]$
  denotes the corresponding
  object of  $\RVmn$ with form $1$.

 \lemm{diamond-m}  \lemref{diamond} (1-5) holds for $\RVlmn$.  We also have:
 
 (6) As a semi-ring congruence on $ \SG \RVlmn$, $\Ispm$ is generated by  $([[1_\k]_1],[[\RVp]_1)$ (with the forms $1$)
  \>{lem}

 \prf  (1)-(5) go through with the same proof.  For (6),
 Let $\sim$ be  the congruence generated by this element.
By blowing up a point one sees immediately that 
 $([[1]_1],[[\RVp]_1)  \in \Ispm$, so $\sim \leq \Ispm$.     
 For the converse direction we have to show that    $(\widetilde{\bY},{\bY}) \in \sim$ whenever 
 $\widetilde{\bY}$ is a blowup of $\bY$; the elementary case suffices, since the $\RVlmn$-
 isomorphisms of \defref{rvblowup-m} (2)
 are already accounted for in the semigroup $ \SG \RVlmn$.  Now
 $\bY  = (Y,f,\omega)$ with $f_n(y) \in \RV$.
 Since $\dim(Y)<n$, we
  have $\bY \iso (Y,f',\omega')$ where $f_i'=f_i $ for $i<n$, $f'_n=1$, and
 $\omega' = f_n \omega$.   So we may assume $f_n=1$.   In this case, 
  as in
 the proof of \lemref{diamond} (6), $(\widetilde{\bY},{\bY}) \in \sim$.  \eprf

\<{defn}   
Let $J$ be a $k$-element set of natural numbers.    $\VFRm[J,l,\idot]  $
 is the set of triples $\bX = (X,f,\omega)$, 
with  $X  \subseteq \VF^J \times \RV^*$, $f: X \to \RVi^l$, $\omega: X \to \RV$, and such that $f$ and $\omega$ factor through 
the projection $pr_{\RV}(X)$ of $X$ to the $\RV$-coordinates.  $\Ispm$ is the equivalence relation on $\VFRm[J,l,\idot]$:  
$$(X,Y) \in \Ispm \iff  (X_a,Y_a) \in \Ispm(<a>) \text{ for each } a \in \VF^J$$
$\SG \VFRm$ is the set of equivalence classes.
\>{defn}
 For $j \in J$, let $\pi^j: \VF^k \times \RV^*  \to \VF^{J-\{j\}} \times \RV^*$ be the projection forgetting the $j$'th 
$\VF$ coordinate.     We will write $\VFRm[k,l,\idot], \VF^k,\VF^{k-1}$ for 
$\VFRm[J,l,\idot]$
$\VF^J, \VF^{J-\{j\}}$ respectively
when the identity of the indices is not important.

 The map $\L: \Ob  \RVmn \to \Ob \VFmn$ induces, by \lemref{rvlift-m},
 a homomorphism $\L: \SG  \RVmn \to \SG \VFmn$.   By  \propref{tr2}
 it is surjective.
 
 \lemm{sp-bl-m}  Let $\bX,\bX'  \in \RVmn$, and let $G: \L \bX' \to \L \bX$ be a special 
bijection.   Then $\bX'$ is isomorphic to  an iterated blowup of $\bX$. \>{lem}
 \prf Clear from \lemref{sp-bl}, since strong isomorphisms are also $\RVlmn$-isomorphisms.  \eprf

 \lemm{rvconj1-m}   The 
 homomorphism
 $\L: \SG  \RVm[1,\idot] \to \SG \VFm[1,\idot]$ is surjective,
with kernel equal to $\Ispm[1]$.  The image of $\SG \RVvol[1,\idot]$ is $\SG \volVF[1,\idot]$  \>{lem}

\prf   Let $\bX,\bY \in \RVm[1,\idot]$, and let $F: \L \bX \to \L \bY$ be a
definable measure preserving bijection.  We have $\bX = (X,f,\omega),
\bY=(Y,g,\omega)$ with $(X,f),(Y,g) \in \RV[1,\idot]$.  By \lemref{rvconj1}
there exist special bijections $b_X:  \L \bX' \to \L \bX, b_Y: \L \bY' \to \L \bY$
and an \rvinv definable bijection $F':\L \bX' \to \L \bY'$ such that $b_Y F' = F b_X$.  
 We used here that any \rvinv object can be written as $\L \bX'$ for some $\bX'$.
Since
$F,b_X,b_Y$ are measure preserving bijections, so is $F'$.  By \lemref{sp-bl-m}, 
  $\bX' \to \bX$ and $\bY' \to \bY$ are blowups; and $F'$ descends to a definable bijection
  between them.  This bijection is measure preserving by \lemref{rv-compat}.  Hence
  by definition $(\bX,\bY) \in \Ispm$.  \eprf

By \propref{rvconj1-m}, the inverse of $\L: \RV[1,\idot] \to \VF[1,\idot]$ induces an isomorphism
$I^{vol}_1 :\SG \volVF[1,\idot] \to \SG \RVvol[1,\idot] / \Ispm$.
$$I^{vol}_1([X]) = [Y]/\Ispm \iff  [\LY ] = [X]$$

Let $\bX = (X,f,\om) \in \VFRm[k,l,\idot]$.  By  assumption, $f,\om$ factor through $\pi^j$,
so that they can be viewed as functions on $\pi^j X$.  
We view the image $(\pi^j X,f,\om)$ 
as an element of $\VFRm[k-1,l,\idot]$.      Each fiber of $\pi^j$ is 
a subset of $\VF$; it can be viewed as an element of $\VFvol[1] \subseteq \VFm[1]
\subseteq \VFm[1,\idot]$.

\Claim{}  
Relative $\Ispm$-equivalence implies $\Ispm$-equivalence, in the following sense.  Let
$X_i  \subseteq \RV^*$ ($i=1,2$); $h_i: X_i \to   W \subseteq \RV^*$;
$f_W: W \to \RV^{l}$, $\om: W \to \RV$, and $f_i: X \to \RV^k$ be definable sets and functions.
Let $\bX_i = (X_i,(f_W \circ h_i,f_i),\om \circ h_i)$.  Let $\bX_i(w) = (X_i(w), f_i | X_i(w), \omega \circ h_i | X_i(w))$ where $X_i(w) = h_i \inv (w)$.  If $\bX_1(w), \bX_2(w) \in \Isp(<w>)$ 
for each $w \in W$, then $(\bX_1,\bX_2) \in \Ispm$.  

\prf Clear, using \lemref{fibrations-m}.  \eprf

The Claim allows us to relativize 
  $I^{vol}_1$ to  $\pi^j$.   We obtain a map
$$I^j = I^j_{k,l} : \VFRm[k,l,\idot] \to   \SG \VFRm[k-1,l+1,\idot]/\Ispm$$

\lemm{isp4-m} Let $\bX=(X,f,\om), \bX' = (X',f',\om') \in \VFRm[k,l,\idot]$.  
 \begin{enumerate}
  \item $I^j$ commutes with maps into $\RV$:  if $h: \bX \to W \subseteq \RV^*$ is definable,
  $\bX_c = h \inv (c)$, then $I^j(\bX) = \sum_{c \in W} I^j(\bX_c)$
% \item If $(\bX,\bX')$ are effectively isomorphic over $\VF^k$, then $(I^j(\bX),I^j(\bX')) \in \Ispm.
  \item    If $([\bX],[\bX']) \in \Ispm$ then $(I^j(\bX),I^j(\bX')) \in \Ispm$.
  \item $I^j$ induces a map $\SG \VFRm[k,l,\idot] /\Ispm \to \SG \VFRm[k-1,l+1,\idot]/\Ispm$.
\end{enumerate}            \>{lem}
 
\prf (1)  This reduces to the case of $I^{vol}_1$, where it is an immediate consequence
of the uniquness, and the fact that $\L$ commutes with maps into $\RV$ in the same sense.

(2)  All equivalences here are relative to the $k-1$ coordinates of $\VF$ other than $j$,
so we may assume $k=1$, and write $I$ for $I^j$.  

For $a \in \VF$, $([\bX_a]  ,[\bX'_a]) \in \Ispm(<a>)$.  By stable embeddedness of $\RV$,
there exists $\alpha = \alpha(a) \in \RV^*$ such that $\bX_a,\bX'_a$ are $<\a>$-definable 
there are $<\a>$-definable blowups $\widetilde{\bX}_a, \widetilde{\bX'}_a$ and
an $<\a>$-definable isomorphism between them, lifting to an $a$-definable isomorphism.     
Using (1) and \lemref{fibrations-m}  we may assume that $\alpha$ is constant.
So for some $W \in \Ob \VF[1], \bY,\bY' \in \RVm[l+1,\idot]$,
we have
$\bX = W  \times \bY, \bX' = W \times \bY'$, 
$\widetilde{\bY},  \widetilde{\bY'}$ are blowups of $\bY,\bY'$ respectively,
$\phi: \bY \to \bY'$ is a bijection, and for any $w \in W$ there exists a measure preserving
$F_w: \L \widetilde{\bY} \to \L \widetilde{\bY}$ lifting $\phi$.      Then
$I(\bX) = I(W) \times \bY, I(\bX') = I(W) \times \bY'$ and the 
bijection $Id_{I(W)} \times \phi$ is lifted by the measure preserving bijection
$(w,y) \mapsto (w,F_w(y))$.

(3) by (2). 

\eprf

 \<{lem}\lbl{fubini-m}  Let $\bX =(X,f,\om) \in \Ob \VFRm[J,l,\idot]$.   If $j \neq j' \in J$ then $I^j I^{j'} = I^{j'} I^j:   \SG \VFRm[J,l,\idot] /\Ispm \to \SG \VFRm[J \m \{j,j'\},l+2,\idot]/\Ispm$.
\>{lem}

\prf  We may assume $S=\{1,2\}, j=1,j'=2$, since all is relative to $\VF^{S \m \{j,j'\}}$.  
By \lemref{isp4} (1) 
and \lemref{fibrations-m}
  it suffices to prove the statement for each fiber of a given  
map into $\RV[l]$.  
Hence we may assume  $X \subseteq \VF^2$ so that  $f$ is constant; and by \lemref{cell2}, we can assume $X$ is a basic 2-cell:  
$$X= \{(x,y): x \in X_1, \rv(y-G(x))=\alpha_1\} \ \ \ \ \ \ \ \ \ \ X_1= \rv \inv (\delta_1) + c_1 $$
The case where $G$ is constant is easy, since then $X$ is a finite union of  rectangles. 
Otherwise $G$ is invertible, and by
 niceness of $G$ we can also write:
$$X =  \{(x,y): y \in X_2, \rv(x-G \inv (y))=\beta \},  \ \ \ \ \ \ \ \ \ \ X_2= \rv \inv (\delta_2) + c_2 $$ 
We immediately compute:  
$$I_2I_1(X) =(\d_1,\alpha_1), \ \ \ I_1I_2(X) =   (\a_2,\d_2)$$
and necessarily $\valr \d_1 + \valr \alpha_1 = \valr \a_2 + \valr \d_2$.  (\lemref{image2l}).
We have bijections $F_j: X \to \L I_j(X)$.  The map $F_1F_2 F_1 \inv F_2 \inv: \L I_2I_1(X) \to \L I_1I_2(X)$ lifts the   unique bijection between the singleton sets 
$\{(\d_1,\alpha_1)\} ,\{(\a_2,\d_2)\}$, and shows that  $[(\d_1,\alpha_1)]_2 = [(\a_2,\d_2)]_2$.  
    \eprf

\<{prop} \lbl{change-of-variable-m} Let $\bX,\bY \in \RVm[\leq n,\idot]$.  If $\L \bX ,\L \bY$ are
%effectively
 isomorphic, then $([X],[Y]) \in \Ispm$.  \>{prop}

\prf  Identical to the proof of \propref{change-of-variable-1}, only quoting
 \lemref{fubini-m} in place of  \lemref{fubini}, and \lemref{carryform} to enable using
  \lemref{decompose-unary}.   \eprf

\<{prop} \lbl{eff-m}  \propref{eff} is valid for $\VFm[n], \RVm[n], \Ispm[n]$. \>{prop}
\prf Same as the proof of \propref{eff},   but using
 \propref{rvlift-m} in place of \ref{rvlift} and 
   \propref{change-of-variable-m} in place of \ref{change-of-variable-1}.  \eprf

\ssec{Invariants of measure preserving maps, and some induced isomorphisms}

 \<{thm} \lbl{volumes}   Let $\T$ be $\V$-minimal.
 There exists a canonical isomorphism of Grothendieck semigroups
 $$\intm:  \SGe \VFm[n,\idot]   \to \SG(\RVm[n,\idot]) / \Ispm[n]$$
 Let $[X]$ denote the class of $X$ in $\SGe(\VFm[n])$.  Then  
   $$\intm [X] = W / \Ispm[n] \iff [X]=[\L W] \in \SGe(\VFm[n])  $$
   
   \>{thm}
   
 \prf  Given $\bX=(X,f,\om) \in \Ob \RVm[n]$ we have $\L \bX \in \Ob \VFm[n]$.

  If $\bX,\bX'$ are isomorphic, then by \lemref{rvlift-m},
$\L \bX, \L \bX'$ are effectively isomorphic.  

Direct sums are clearly respected,
so we have a semigroup homomorphism $\L: \SG(\RVm[n]) \to \SGe(\VFm[n]) $.

It is surjective by \propref{tr2}, injective by \propref{eff}.   

 Inverting, we obtain $I$.    \eprf

 Let $\Ispm'$ be the semigroup congruence on $\RVvol [n]$ generated by $((Y,f),(Y \times \RVp, f'))$ where $Y,f,f'$ are as in \defref{rvblowup}.  Let $\Ispmg$ be the congruence
on $ \SG \RVmg [n] $ generated by $([[1_\k]_1],[[\RVp]_1)$, with the constant $\G$-form $0 \in \G$.  
%We use the same letter to denote the restriction of this congruence to $\SG \RVmg \bdd$.

    Assume given a distinguished subgroup $N_1$ of the multiplicative group of the residue field $\k$.  
For example, $N_1$ may be the group of elements of norm one, with respect to some absolute value  $|,|$on $\k$.  
With this example in mind, write $|x|=1$ for $x \in N_1$.
%(cf. \secref{expansionsS})
Let $\VFma[n]$ be the subcategory of $\VF[n]$ with the same objects, and such that
$F \in \Mor_{\VFma[n]}$ iff $F \in \Mor_{\VFgn}$ and $ |JRV(F)|=1$ almost everywhere.
 Similarly define $\RVma[n]$.
    \<{thm} \lbl{volumes-2}    The isomorphism 
$\intm $ of \thmref{volumes}  induces isomorphisms:    
     
  \begin{eqnarray}{        }
 \SGe \VFvol [n]  &  \to  & \SG \RVvol[n] / \Ispm'[n]           \label{isos1}  \\
\SGe \VFvol \bdd [n]  & \to & \SG \RVvol \bdd [n] / \Ispm[n]    \label{isos2} \\
\SGe \VFm \bdd [n]  &  \to &  \SG \RVm \bdd [n] / \Ispm[n]      \label{isos4}    \\
\SGe \VFma [n] & \to & \SG \RVma[n]  / \Ispm[n]    \label{vfma}     \\
\SGe \VFmg [n] & \to & \SG \RVmg[n]  / \Ispmg[n]    \label{vfmg}  \\
%\SGe \VFmg \bdd [n]  &  \to &  \SG \RVmg \bdd [n] / \Ispmg[n]      \label{vfmgb}    
\end{eqnarray}

\>{thm}

\prf    
  Since   \propref{tr2} uses measure preserving maps,    \propref{rvlift} do not go out of the subcategory
 $\VFvol$, and $\RVvol [n]$ is a full subcategory of $\RVm [n]$, we  have \eqref{isos1}.  
 It is similarly easy to see that ``dimension $<n$'' and boundedness
 are preserved, hence \eqref{isos2}-\eqref{isos4}.
   
 We have $\SGe \VFma = \SGe \VFm / N_{\VF}$, where $N_{\VF} = \{([X,\omega],[X,g \omega]): g: X \to \RV, |g|=1 \}$.
 Similarly for $\SGe \RVma$.   Thus for \eqref{vfma} it suffices to show that
 $(\intm(\bX),\intm(\bY))  \in N_{\RV} \iff (\bX,\bY) \in N_{\VF}$.
 % or equivalently:
 %$(\L X, \L Y) \in N_{\VF}$ iff $(X,Y) \in N_{\RV}$.  
 For $X \in \Ob \VFm[n]$ or $X \in \Ob \RVm[n]$
 with $\RV$-volume form $\omega$, given $g: X \to \RV$, let $^g X$ denote the same object
 but with volume form $g \omega$.    In one direction we have to show:  
 $(\L X, \L Y) \in N_{\VF}$ if $(X,Y) \in N_{\RV}$.  This is clear since $\L (^g X) = ^g (\L X)$.
Conversely we have to show that $(\intm [^g Z] , \intm [Z]) \in N_{\RV}$.  Since $\intm$
commutes with $\RV$-sums, we may assume $g$ is constant, with value $a$.
But then $\L (^a X) = ^a (\L X)$ implies $\intm (^a Z) = ^a \intm Z$ as required.   This
gives \eqref{vfma}; \eqref{vfmg} is a special case.
\eprf

\>{section}

 \<{section}{The Grothendieck semirings of $\Gamma$}   
\lbl{exgam}

Let $T=\DOAG_A$ be the theory of divisible ordered Abelian groups $\Gamma$, 
with distinguished constants  for elements of a  subgroup $A$.  Let $\DOAG_A[*]$
be the category of all $\DOAG_A$ definable sets and bijections.  Our primary concern is
 not with $\DOAG_A$, but rather a proper subcategory $\G[*]$, having the same objects
but only piecewise integral morphisms (\defref{Gcat}.)    Our interest in $\G[*]$  
derives from this:  the morphisms of $\G[*]$
are precisely those that lift to morphisms of $\RV[*]$, and it is $\SG[\G[*]]$ that forms
a part of $\SG[\RV[*]]$.  (cf. \secref{groupext}.)
This category depends on $A$, but will nevertheless be denoted $\G[*]$  when  $A$ is fixed and understood.

We will first describe $\K(\G^{fin}[*])$, the subring of classes of finite definable sets.  
Next
  we will analyze $\K(\DOAG_A)$,   obtaining 
two   Euler 
characteristics.   This repeats earlier work
by  
Ma\v{r}\'{i}kov\'{a}.  We retain our proofs as they give a rapid
path to the Euler characteristics, but \cite{marikova} includes a complete analysis
of the semiring $\K(\DOAG_A)$, that may well be useful in future applications. 

At the level of Grothendieck rings, the categories $\G[*]_A$ and $\DOAG_A$ may be
rather close; see \lemref{oag2} and Question \ref{q2}.   But the  semiring homomorphism
$\SG(\G[*]_A) \to \K(\DOAG_A)$ is far from being an isomorphism,  and   it remains important to give a good description
of $\SG(\G[*]_A)$.      We believe
that further invariants can be
found by mapping $\SG[\G[*]]$ into the Grothendieck  semirings of other completions
of the universal theory of  ordered Abelian groups over $A$, as well as $\DOAG$, in the manner of \propref{grh0}; it is possible that all invariants 
appear in this way.  

A description of $\SG(\G[*]_A)$ would include information about the Grothendieck group
of subcategories, such as the category of bounded definable sets.   We will only
sample one bit of the information available there, in the form of a    ``volume'' map on bounded
 subsets of $\SG[\G[*]]$ into the rationals, and a discrete analogue.

 \<{defn} \lbl{Gcat} An object of  $ \G[n]$ is a finite disjoint union of  subsets of $\G^n$ defined by linear equalities and inequalities
with $\Zz$-coefficients and parameters in $A$.  
Given
$X,Y \in \Ob \G[n]$, $f \in \Mor_\G(X,Y)$ iff $f$ is a bijection, and there exists a partition $X=\union_{i=1}^n X_i$,
$M_i \in \GL_n(\Zz), a_i \in A^n$, such that for $x \in X_i$,
$$f(x) = M_i x + a_i$$    \>{defn}

 $\G[*]$ is the category of definable subsets of $\G^n$ for any $n$, with the same morphisms.  Since there are no morphisms between different dimensions, it is simply the direct sum of the categories $\G[n]$, and the Grothendieck semi-ring $\SG[\G]$ of $\G[*]$ is the graded
 direct sum of the semigroups $\SG(\G[n])$.    We will write $\K[\G]$ for the corresponding
 group.
 
 Let $\G \bdd [*]$ be the full subcategory of $\G [*]$
consisting of bounded sets, i.e. an element of  $\Ob \G \bdd [n]$ is a definable subset of $[-\g,\g]^n$
for some $\g \in \G$.

 $\G_A$ is a subcategory of $\G_{\Qq \tensor A}$ (a category with the same objects,
 but more morphisms, generated by additional translations) and this in turn is a
 subcategory of $\DOAG_{\Qq \tensor A}$.  
 
   There is therefore
 always a natural morphism from $\SG(\G_A[*])$ to the simpler semigroup  
  $\SG(\DOAG_{\Qq \tensor A})$.     We will exhibit  two independent Euler characteristics on  $\DOAG_{\Qq \tensor A}$ and show that they define an isomorphism  $\K(\DOAG_{\Qq \tensor A}) \to \Zz^2$.   Taking the  dimension 
 grading into account, this    will give rise to two families of Euler characteristics on  $\K(\G_A)$, with $\Zz[T]$-
coefficients.

 \ssec{Finite sets}       %--------------------------------------------------------
 
Let $\G^{fin}[n]$ be the full subcategory of $\G_A[n]$ consisting of finite sets.  The Grothendieck semi-ring of $\G^{fin}[*]$ embeds into the semi-rings of both   $\G_A$ and   $\RES$, within
the Grothendieck semi-ring of $\RV_A$, and we will see that $\SG (\RV_A)$ is freely generated by them over $\SG(\G^{fin}[*]$.
  We proceed to analyze $\SG(\G^{fin}[*])$ in detail.  

 Let $\eo = [0]_1  \in \SG(\G^{fin}[1])$ be the class of the singleton $\{0\}$.
 
 The unit element of $\K(\G)$ is the class of $\G^{0}$.  Note that the bijection
between $\eo$ and $\G^{0}$   is not a morphism in $\G[*]$; in fact
$1,\eo,\eo^2,\ldots$ are distinct and $\Qq$-linearly independent in
  $\K(\G)$.  The motivation for this choice of category becomes clear if one thinks
  of the lift to $\RV$:   
  the inverse image of $\eo ^n$ in $\RV$ (also denoted $\eo^n$) has dimension $n$,
  and cannot be a union of isomorphic copies of $\eo^m$ for smaller $m$.   

 Let $\K(\G^{fin})[\eo \inv]$  be the localization.
 This ring is naturally   $\Zz$-graded ring;  let 
 $\Hf$ be the 0-dimensional component.   
  
Let $\Xi_A$ be the space
of subgroups of $(\Qq \tensor A)/A$, or equivalently of subgroups of $\Qq \tensor A$ containing $A$.
  View it as a closed subspace of the Tychonoff space $2^{(\Qq \tensor A)/A}$,
via the characteristic function $1_s$ of a subgroup $s \in \Xi_A$.  Let $C(\Xi_A,\Zz)$ be the ring
of continuous functions $\Xi_A \to \Zz$ (where $\Zz$ is discrete.) 
 
    A {\em cancellation} semigroup is a semigroup where $a+b=a+c$ implies $b=c$; in other words, a subsemigroup of an Abelian group.

\<{prop}   \lbl{grh0}     
$\SG(\G^{fin}[n])$  is  a cancellation semigroup.  As a semiring,
$\SG(\G^{fin}[*])$ is generated by $\SG(\G^{fin}[1])$. 
We have:   
$$  \K(\G^{fin})[\eo \inv] = \Hf[\eo,\eo \inv]$$
$${\Hf}  \iso  C(\Xi_A,\Zz)$$    \>{prop}  

\proof           Since $\G$ is ordered,
any finite definable subset of $\G^n$ is a   union
of definable singletons.   Thus the semi-group $\SG(\G^{fin}[n])$ is freely generated by the isomorphism classes
of  singletons $a \in \G^n$, and in particular is a   cancellation semigroup.   The displayed
equality is thus clear; we proceed to prove the isomorphism.

A definable singleton of $\G^n$
has the form $(a_1,\ldots,a_n)$, where for some $N \in \Nn$,
$Na_1,\ldots,Na_n \in A$.  So $[(a_1,\ldots,a_n)] = [(a_1)]\cdots [(a_n)]$.

For any commutative ring $R$, let $Idem(R)$ be the Boolean algebra of idempotent elements of 
a commutative ring $R$
  with the operations $1,0,xy, x+y-xy$. 
Note that the elements $[(a_1,\ldots,a_n)] \eo ^{-n} \in \Hf$ belong to $Idem(\Hf)$:  
  in $\SG(\G^{fin})$:  for any $a \in \G$
we have the relation:  $[a]^2 = [a] \eo$.  Let $B$  be the Boolean subalgebra
 of
$Idem(\Hf)$ generated by the elements $[(a_1,\ldots,a_n)] \eo ^{-n}$.  For a maximal ideal $M$ of $B$, 
let $I_M$ be the ideal of $\Hf$ generated by $M$.    Note $\Hf = \Zz B$.   Hence we   have to show:
\begin{enumerate}
  \item The Stone space of $B$ is $\Xi_A$
  \item For any maximal ideal $M$ of $B$, $\Hf  / I_M \iso \Zz$ naturally.
\end{enumerate}
 
For any commutative ring $R$, a finitely generated Boolean ideal of $Idem(R)$ is generated by a single element $b$; 
if $b \neq 1$, then  $bR \neq R$ since $b(1-b)=0$.  
Thus if $M$ is a proper ideal of $Idem(R)$, then $MR$ is a proper ideal of $R$.
Applying this to $B$, viewed as a Boolean subalgebra of $Idem(\Qq \tensor \Hf)$,
we see that $I_M \meet \Zz =(0)$ for any maximal ideal $M$ of $B$.
Thus the composition $\Zz \to \Hf \to \Hf/I_M$ is injective.  On the other hand,   $\Hf$ is generated over $\Zz$ by the elements
$[a]/\eo$, and each of them equals $0$ or $1$ modulo $I_M$, so the map is surjective too.
This proves the second point.

  To prove the first, we define a map $\Phi:  \Xi_A \to Stone(B)$.  

Let $t = T/A$, $T \leq \Qq \tensor A$.  
If $[(a_1,\ldots,a_n)] = [(b_1,\ldots,b_n)]$, then some element
of $GL_n(\Zz) \ltimes A^n$ takes $(a_1,\ldots,a_n)$ to $(b_1,\ldots,b_n)$; 
in this case, if $a_i \in T$ for each $i$ then $b_i \in T$ for each $i$; so
$ \Pi_{i=1}^n 1_t(a_i +A) =  \Pi_{i=1}^n 1_t(b_i +A)$.  Thus, 
given $t \in \Xi_A$, we can
define a homomorphism $h_t: \Hf \to \Zz$ by:
$$ \frac{[(a_1,\ldots,a_n)]} {\eo^n}  \mapsto \Pi_{i=1}^n 1_t(a_i +A) $$

   Let 
$M(t) = \ker (h_t) \meet B$.  

The map $\Phi: t \mapsto M(t)$ is clearly continuous.   If $t,t'$ are distinct subgroups, let $a \in t$,
$a \notin t'$ (say); then $[a]/\eo \in M(t), [a]/\eo \notin M(t')$.  So $\Phi$ is injective.  
If $P$ is a maximal filter of $B$, let $t_P = \{a +A:  [a]/\eo \in P \}$.

\Claim{} $t_P$  is a subgroup. 

Proof:  
 
Suppose $a+A,b+A \in t_P$ and let $c=a+b$.  Then 
we have the relation 
$$[a][b] \eo= [a][b][c]$$
 in $\SG(\G^{fin})$, arising from the map
 $$(x,y,z) \mapsto (x,y,xyz)$$
Thus $([a]/\eo)([b]/\eo)(1-[c]/\eo)=0$.    As $([a]/\eo), ([b]/\eo) \in P$
we have $(1-[c]/\eo) \notin P$, so $[c]/\eo \in P$.

 Clearly $P = M(t_P)$.  Thus $\Phi$ is surjective, and so a homeomorphism. \qed

{\bf Example}  We always
have a homomorphism ${\K(\G^{fin})} \to \Zz$ (by counting points
of a finite set in the divisible hull); when $A$ is divisible, this identifies
${\K(\G^{fin})}$ with $\Zz[\tau]$.  In general, we have the surjective  morphism
$\K(\G^{fin})) \to \K(\G^{fin}_{\Qq \tensor A}) = \Zz[\tau]$.

 \<{lem}  \lbl{smallerdim}
  Let $Y$ be an $A$-definable subset of $\Gamma^n$, of dimension $<n$.  Then
  $Y$ is a finite union of $GL_n(\Zz)$-conjugates of sets $Y_i \subseteq \{c_i\} \times \G^{n-1}$,
  with $c_i \in \Qq \tensor A$.  \>{lem}
  
  \<{proof}  $Y$ can be divided into finitely many $A$-definable pieces,
each contained in some $A$-definable hyperplane of $\Gamma^n$.
So we may assume $Y$ itself is contained in some such hyperplane,
i.e. $\sum r_i y_i = c$ for some $c \in \Qq \tensor \valr(A)$.  We may
assume $r_i \in \Zz$ and $(r_1,\ldots,r_n)$ have no common divisor.
In this case   $\Zz^n / \Zz(r_1,\ldots,r_n)$ is torsion free, hence free,
so $ \Zz(r_1,\ldots,r_n)$ is a direct summand of $\Zz^n$.  Thus
after effecting a transformation of $GL_n(\Zz)$, we may assume
$(r_1,\ldots,r_n) = (1,0,\ldots,0)$, i.e. $Y$ lies in the hyperplane
$y_1 = c$.  Let $Z$ be the projection of $Y$ to the  coordinates 
 $(2,\ldots,n)$. 
 Then $Y = \{c\} \times Z$.
 \>{proof}

\ssec{Euler characteristics of $\DOAG$}

We describe two independent Euler characteristics on $A$- definable subsets of $\G$,
i.e. additive, multiplicative $\Zz[\eo]$-valued functions   invariant under all definable bijections. 
The values are in $\Zz[\eo]$ rather than $\Zz$ because $\G[*] = \oplus _n \G[n]$ is graded
by ambient dimension.      \propref{2eulers} - \lemref{euler2}
were obtained earlier in \cite{marikova}, and independently
in \cite{kageyama-fujita}.

In fact these two  Euler characteristics come from Euler characteristics of  $\DOAG_{\Qq \tensor A}$.   There, they are the only ones:

\<{prop}  \lbl{2eulers}  Let $A$ be a divisible ordered  Abelian group.   
Then $\K(\DOAG_{A}) \iso \Zz^2 $\>{prop}
 
\<{proof}   We begin by noting that there are at most two possibilities. 

  In $\DOAG$, all definable singletons are isomorphic.  The identity element
of the ring $\K(\DOAG)$ is the class of any singleton.  Thus the image of 
${\K(\G^{fin}[*])}$ in $\K(\DOAG_{A})$ is isomorphic to $\Zz$.  

%  Bounded sets reduce to finite sets in the Grothendieck group:

\Claim{}   The image of $\K(\fCgb)$ in  $\K(\DOAG_{A})$ equals
the image of ${\K(\G^{fin}[*])}$ there.  

 Translation by $a$ gives an equality of classes in $\K(\G)$, 
 $[(0,\infty) = [(a,\infty)]$ so $$[(0,a)]+[\{pt\}] = [(0,a] = 0$$
Thus bounded segments are equivalent to linear combinations of points.  
%This
%can be extended to higher dimensions using 
%  integration by parts.   Alternatively, it 
  This can be seen directly by induction on dimension and
on ambient dimension:   consider the class of a bounded set $Y \subset \G^{n+1}$.  $Y$ is a Boolean combination of
sets of the form $ \{(x,y): x \in X, f(x) < y < g(x) \}$.  This is $\DOAG_A$-isomorphic
to $Y' = \{(x,y): x \in X,  0<y<h(x) \}$, where $h(x)=g(x)-f(x)$.   Let $Z= \{(x,y): x \in X, 
y>0\}$, $Z' = \{(x,y): x \in X, y> h(x) \}$.  Then the map 
$(x,y) \mapsto (x,y+h(x))$ shows
that $[Z]=[Z']$.  On the other hand $Z'$ is the disjoint union of $Z,Y$ and a smaller-dimensional set $W$.  Thus $[Z]=[Z']=[Z']+[Y]+[W]$ so $[Y]=-[W]$, and by induction $[Y]$
lies in the image of ${\K(\G^{fin}[*])}$.

Now consider  $t=[(0,\infty)] \in \K(\G_A)$.   We have a homomorphism 
$\K(\fCgb)[t] \to \K(\G)$. 
% Using integration by parts, it can be seen to be
%surjective.    %  (Note that $\sum_{0<x} [\{x/n\}] = \sum_{0<y} [\{y\}] = t$.)   
%Here is an alternative direct argument:
To see that it is surjective, 
  again by induction it suffices to look
at sets such as $\{(x,y): x \in X, f(x)<y \}$ or $\{(x,y): x \in X, f(x)<y<g(x)\}$.
The latter is equivalent to a smaller-dimensional set, by induction, as above.
The former is equivalent to $\{(x,y): x \in X, 0<y \}$ so that it has the class
$[X] \times t$ and is thus in the image of $\K(\fCgb)[t]$.

Let $T = \{(x,y): 0<y \leq x \}$.   
The map $(x,y) \mapsto (x,y+x)$ takes $T$ to
$\{(x,y): 0<x<y \leq 2x\}$, so $2[T]  = [\{(x,y): 0< y \leq 2x \}]$.
The same map shows that $t^2-[T] = t^2-2[T]$ so $[T]=0$. 
But then $[\{(x,y): 0 < x \leq y \}] = 0$, and adding we obtain:
$0+0 = t^2 + [\{(x,x): 0<x \}] = t^2+t $.  Thus $\K(\DOAG_{A})$  
is a homomorphic image of  $\Zz[t]/ (t^2 + t) \iso \Zz^2$.  To see that the homomorphism
is bijective, it remains to exhibit a homomorphism $\K(\DOAG_{A}) \to \Zz$
with $t \mapsto 0$ and another with $t \mapsto -1$.  The two lemmas below show this, 
in a form suitable also for a dimension- graded version.

\>{proof}

 \<{lem}\lbl{euler1}  There exists a ring homomorphism $\chi_O: \K(\G) \to \Zz[\eo]$, with 
 $\chi_O( (0,\infty)) = \eo$.  It is invariant under $GL_n(\Qq)$ acting on $\G^n$.  \>{lem}
 
 \<{proof} 

Let $\RCF$ be the theory of real closed fields.  See \cite{vddries} for the existence
and definability of an  Euler characteristic map $\chi: \K(\RCF) \to \Zz$.
 For any definable $X,P,f: X \to P$ of $\RCF$, there exists $m \in \Nn$ and a definable partition
 $P= \union_{-m \leq i \leq m} P_i$, such that for any $i$, any $M \models \RCF$ and
 $b \in P_i(M)$, $\chi (X_b) =i$.  Here $X_b = f \inv (b)$, and $\chi(X_b)=i$
 iff there exists an $M$-definable partition of $X_b$ into definable cells $C_j$,
 with $\sum_j (-1)^{\dim(C_j)} = i$.

The language of $\G$ (the language of ordered Abelian groups) is contained
in the language of $\RCF$.  So if $X,P,f: X \to P$ are definable in the  language of ordered Abelian groups, they are $\RCF$-definable.  Therefore, the above result specializes, and we obtain
an   Euler characteristic map $ \chi: \K(\G_A[n]) \to \Zz$, valid
 for any divisible group $A$.  This $\chi$ is invariant under all definable bijections
 (not only the morphisms of $\G[*]$), and is additive and multiplicative.  We have 
 $\chi_O(\{0\})=1,\chi_O((a,b))=-1$ for $a<b$, and $\chi_O(0,\infty)=-1$ too
(though $(0,1)$ and $(0,\infty)$ are not definably isomorphic in 
the linear structure.)            Now let $\chi_O(X) = \chi(X) \eo^n$ for
$ X \subseteq \G^n$, and extend to $\G[*]$ by additivity.              \>{proof}

{\bf Remark}  The Euler characteristic constructed in this proof appears to depend
on an embedding of $A$ into the additive group of a model of $\RCF$.  But by the uniqueness
shown above it does not.   In fact, as pointed out to us by Van den Dries, Ealy and Ma\v{r}\'{i}kov\'{a},
an Euler characteristic with the requisite properties is defined in \cite{vddries} directly for any O-minimal 
structure; moreover the use of $\RCF$ in the lemma below can also be replaced 
by a direct inductive argument, and some simple facts about Fourier-Motzkin elimination.  

 Another Euler characteristic   can be obtained as follows:  given a definable
set $Y \subset \G^n$,
let $$\chi'(Y) =  \lim _{r \to \infty} \chi(Y \meet C_r)$$ where $C_r$ is the bounded
closed cube $[-r,r]^n$.  By O-minimality, the value of $\chi(Y \meet C_r)$
is eventually constant. 

Note that  $\chi'$  is not invariant under semi-algebraic bijections, since the bounded and unbounded open intervals are given different measures.    Still,

\<{lem} \lbl{euler2}  $\chi'$ induces a group homomorphism $\K(\G[n]) \to \Zz$; and yields 
a ring homomorphism $\K(\G[*]) \to \Zz[\eo]$.  Moreoever $\chi'$ is 
invariant under piecewise $GL_n(\Qq)  $-transformations.
 \>{lem}

\<{proof}
 $\chi'$ is clearly additive and multiplicative.  
Isomorphism   invariance can be checked as follows:  first,

\Claim{} If $X \neq \emptyset$ is defined by a finite number of weak ($\leq$)  affine  equalities and 
inequalities,  then $\chi'(X)=1$.  

\proof  It suffices to show that this is true in $(\Rr,+)$; since then it is true in any model
of the theory of divisible ordered Abelian groups.  Now we may compute the Euler characteristic $\chi$ of the bounded sets $X \meet C_r$ in $(\Rr,+,\cdot)$.  Let $p \in X$.  For large
enough $r$, $p \in X \meet C_r$ there is a definable retraction of the closed bounded set $X \meet C_r$ 
to $p$ (along lines through $p$.)  Thus $X \meet C_r$ has the same homology groups as a 
point, and so Euler characteristic $1$.    \qed

To prove the lemma we must show that
if $\phi: X \to Y$ is a definable bijection, $X, Y \subseteq \G^n$, then $\chi'(X) = \chi'(Y)$.
We use induction on $\dim(X)$.
By additivity, if $X$ is a Boolean combination of finitely many pieces, 
it suffices to prove the 
lemma for each piece.   We may therefore assume that $\phi$ is linear (rather than only
piecewise linear) on $X$.  Let $\phi'$ be a linear automorphism extending $\phi$.
Expressing
$X$ as a union of basic pieces, we may assume $X$ is defined by some inequalities
  $\sum \alpha_i x_i \leq c$, as well as some equalities and strict inequalities.  So $X$
   is convex.    We have to show that $\chi'(X) = \chi'(\phi'X)$.  Let $\bar{X}$ be the closure of $X$ (defined by the corresponding
  weak inequalities.)  Then $\bar{X} \m X$ has dimension $< \dim(X)$, so by induction
  $\chi'(\phi'(\bar{X} \m X) =  \chi'(\bar{X} \m X)$.  But $\bar{X}$ is closed and convex,
  so $\chi_{O'}(\bar{X}) = 1 = \chi_{O'}(\phi' \bar{X} )$.  Subtracting,
    $\chi'(\phi'(\bar{X}) )=  \chi'(\bar{X})$.  

Once again, using the ambient dimension grading, we can define $\chi_O' : \G[*] \to \Zz[\eo]$
with $\chi_O' (x) = \chi'(x) \eo^n$ for $x \in \G[n]$.  
\>{proof}

%We include at this point a lemma, to be used later, relevant to $\GL_n(\Zz)$-transformations.
 In the following lemma, all classes are taken in $\K(\G_A)[*]$.  Let   $e_a$ be the class in $\K(\G_A)[1]$ of 
the singleton $\{a\}$, and $\tau_a$ the class of the segment $(0,a)$.

\lemm{oag1}  Let $a \in \Qq \tensor A$, $b \in A$.
  \begin{enumerate}
 \item $\tau_a   = \tau_{a+b}$, $e_a = e_{a+b}$.
 \item  If $b<c \in A$ then $[(b,c)]=-e_0$.
 \item $e_ae_0 = e_a^2$
 \item $\tau_a(\tau_a + e_0) = 0$
 \item  If $2a \in A$ then $2 \tau_a + e_a = -e_0$, and  $e_0(e_a-e_0)=0$

\end{enumerate}

\>{lem}

 \prf  
 
(1)  $\tau_a = [(0,a)] = [(0,\infty)] - [(a,\infty)] - e_a$ and similarly $\tau_{a+b}$.
 The map $x \mapsto x+b$ shows that $[(a,\infty) = [(a+b,\infty)]$ and $e_a=e_{a+b}$, 
 hence also $\tau_a = \tau_{a+b}$.
 
(2) $[(b,c)] = [(b,\infty)] - [(c,\infty)] - e_c = -e_0$ by (1), since $c-b \in A$.

(3)  The map $(x,y) \mapsto (x,y+x)$ is an $SL_2(\Zz)$-bijection between $\{(a,0)\}$ 
 and $(a,a)$.
 
 (4)  Let $$D = \{(x,y): 0 < x <a , 0 < y \leq x \}$$ 
             $$D' = \{(x,y): 0 < y < a, 0 < x \leq y \}$$
 $$D_1 = \{(x,y): 0 < x < a, y > 0\}$$
 $$T(x,y) = (x,y+x)$$
 Then $T(D_1) = D_1 \m D$.  Since $[T(D_1)] = [D_1]$, $[D]=0$.  Similarly $[D']=0$.
 Note also   $$T((0,a) \times \{0\}) =  \{(x,x): 0<x<a \}$$
 Thus  $$0=[D]+[D'] = [(0,a)^2] + [ \{(x,x): 0<x<a \}] = \tau_a^2 + \tau_a  e_0$$
 
 (5)  Let $0<2a \in A$.  Then $[(0,a)] = [(a,2a)]$ using the map $x \mapsto 2a-x$.
 So $2 \tau_a + e_a = [(0,a) \union \{a\} \union (a,2a)] = [(0,2a)] = -e_0$ (by (2)).
 So $(-e_0-e_a)(e_0-e_a) = (2 \tau_a) (2 \tau_a + 2 e_0) = 0$ by (1).  Thus
 $e_ae_0 = e_a^2 = e_0^2$.
 \eprf

The next lemma will not be used, except as a partial indication towards the question that follows,
regarding the difference at the level of Grothendieck groups between $GL_n(\Zz)$ and
$GL_n(\Qq)$ transformations.
 Let $Ann(e_0)$ be the
annihilator ideal of $e_0$; it is a graded ideal.  
 Let $R = \K(\G_A)[*] / Ann(e_0)$, the image of $\K(\G_A)[*]$ in the localization
  $\K(\G_A)[*](e_0 \inv)$. 
In the next lemma, the classes of definable sets are taken in $R$, viewed as a subring of  $\K(\G_A)[*](e_0 \inv)$.
Let $\be_a = e_a/e_0, t_a = \tau_a / e_0$.

\lemm{oag2} Let $A' = \{a \in \Qq \tensor A:  \be_a =1  \}$.  
 
(1)  If $X \subseteq \G^n$ is definable by linear inequalities over $A$, and $T \in  GL_n(\Zz) \sd (A')^n$,
 then $[TX]=[X] \in R$.  
 
(2)  $A'$ is a subgroup of $\Qq \tensor A$.

(3)  $\be_a^2 = \be_a$, $t_a(t_a+1) = 0$

(4)  $A'$ is $2$-divisible.

\>{lem}  
 
 \prf 
 
 (1)  It suffices to show this when $T$ is a translation by an element $a \in (A')^n$.  The map
 $(x,y) \mapsto (x+y,y)$ is in $SL_{2n}(\Zz)$, hence $[X \times \{a\} ] = [TX \times \{a\}]$
 in $\K(\G_A)[2n]$.  Since $a \in (A')^n$, $[a] = e_0^n$.  So $[X] e_0^n = [TX]e_0^n$, and 
 upon dividing by $e_0^n$ the statement follows.
 
 (2) Clear from (1).   For the following clauses, note that by (1),(2), \lemref{oag1} applies with $A$ replaced by $A'$.
 
 (3) By \lemref{oag1} (3),(4) divided by $e_0^2$.
 
 (4)  By \lemref{oag1}(5) applied to $A'$, if $2a \in A'$ then $e_0(e_a-e_0) = 0$; so $\be_a -1 =0$, i.e. $a \in A'$.
  \eprf

 \<{question}\lbl{q2}  Is it true that $\K(\G_A[*])/Ann(e_0) = \K(\DOAG_A[*])/Ann(e_0)$ ?\>{question}

A positive answer would follow from an extension of (4) to odd primes, over arbitrary $A$;
by an inductive argument,
 or by integration by parts.

\ssec{Bounded sets: volume homomorphism}  

Let $\bar{A} = \Qq \tensor A$.  Recall $\G \bdd [n]$ is the category of bounded $A$-definable
subsets of $\G^n$, with piecewise $GL_n(\Zz) \sd  A $ - bijections for morphisms.
Let $Sym(\bar{A})$ be the
symmetric algebra on $A$.

\<{prop}   \lbl{volume-G}
There exists a natural ``volume'' ring  homomorphism $\K(\G \bdd [*])  \to Sym(\bar{A}) $ \>{prop}

\prf   

We first work with $\DOAG$ without parameters, defining a polynomial associated
with a family of definable sets.

Let $C(x,u) = C(x_1,\ldots,x_n;u_1,\ldots,u_m)$ be a formula of $\DOAG$.  
Write $C_b = \{x: C(x,b)\}$; this is a definable family of definable sets.
Assume the sets $C_b$ are uniformly bounded:   equivalently, as one easily sees, for some
$q \in \Nn$, for each $i$, $C(x,u)$ implies  $|x_i| \leq q \sum_j |u_j|$.
For $b \in \Rr^m$, let $v(b) = \vol C_b(\Rr^n)$.  Here $\vol$ is the Lebesgue measure.

By a {\em constructible function into $\Qq$}, we mean a   $\Qq$-linear combination of
characteristic functions of definable sets of $\DOAG$.  Let $R$ be the $\Qq$-algebra of 
constructible functions into $\Qq$.

\Claim{1} There exists a   polynomial $P_C(u)  \in R[u]$
 such that for all $b \in \Rr^m$, $ \vol C_b(\Rr^n) = P_C(b)$. 

In other words, the volume of a rational polytope is piecewise polynomial in
the parameters, with linear pieces.  
The proof of the claim is standard,  using  iterated integration.   For each $C$, fix such a polynomial $P_C$.  

At this point we re-introduce $A$.  
Any $A$-definable bounded subset of $\G^n$ has the form $C_b$ for some
$C$ as above and some $b \in \bar{A}^m$.   

\Claim{2}  If $C_b = C'_{b'}$ then $P_C(b)=P_{C'}(b')$. 

\prf  (See also below for a more algebraic proof).   Fix the formulas  $C,C'$.   Write $b = Ne$, $b' = N' e$ where $e \in \bar{A}^l$ is a vector
of $\Qq$-linearly independent elements of $\bar{A}$, and $N,N'$ are rational matrices.
Write $P_C = \sum a_\nu(u) u^\nu$ where $a_\nu$ is a constructible function into $\Qq$;
similarly $P_{C'}$.   

Note now that any formula $\psi(x_1,\ldots,x_l)$ of $\DOAG$ of dimension $l$
has a solution in $\Rr^l$ whose entries are algebraically independent.  Use this to 
find algebraically independent $\te  \in \Rr^l$ such that $C_{N \te} = C'_{N' \te}$,
and $a_\nu(N \te)=a_\nu(b), a_{\nu}(N' \te) = a'_\nu(b')$ for each multi-index $\nu$
of degree $d$.     

By definition of $P_C$ we have $P_C(N \te) = P_{C'}(N' \te)$. 
So $\sum a_\nu(b) (N \te)^\nu = \sum a'_{\nu}(b') (N' \te)^\nu$.  
 By algebraic independence, $\sum a_\nu(b) (Nv)^\nu = \sum a'_{\nu}(b') (N'v)^\nu$
 as $\Qq$-polynomials.   So $P_C(Ne) = P_{C'}(N'e)$.  \eprf

Thus we can define:  $v(C_b) = P_C(b)$.  Let us show that $v$ defines
a ring homomorphism. 

Given $C,C'$ one can find $C''$ 
such that $C''_{b,b'} = C_b \union C_{b'}$, and similarly $C'''$ with 
$C'''_{b,b'} = C_b \meet C_{b'}$.  Then $P_C+P_{C'} = P_{C''} + P_{C'''}$.
It follows that $v$ is additive.  Similarly $v$ is multiplicative, and translation invariant.
Since $| \det(M)| =1$ for 
$M \in GL_n(\Zz)$, if $\phi^M(x,u) = \phi(Mx,u)$ then $P_{\phi^M} = P_\phi$.
 
 \eprf

Van den Dries,   Ealy, and  Ma\v{r}\'{i}kov\'{a} pointed out that Claim 2 can also be reduced
to the following  statement:   if $Q \in R[u]$, $B$ is any 0-definable set of , and
 $Q$ vanishes on $B(\Rr)$, then $Q$ vanishes on $B(\G)$.   They prove it as follows: let
$\bar{B}$ be the Zariski closure of $B$;  $\bar{B}$  is clearly a finite union of linear subspaces, and by intersecting $B$ with each of these, we may assume $\bar{B}$
is linear, so it is  cut out by homogeneous linear
polynomials $Q_1,\ldots,Q_m$.   Each $Q_i$ vanishes on $B(\Rr)$ and hence on $B(\G)$.  So $Q$ lies in the (radical) ideal generated by $Q_1,\ldots,Q_m$, hence vanishes on 
$B(\G)$.    
 
\sssec{The counting homomorphism in the discrete case}

Suppose $A$ has a least positive element $1$, and 
assume given a homomorphism
$h_p: A \to \Zz_p$ for each $p$.  Then $A$ embeds into  a $\Zz$-group $\tA$, i.e.
an ordered Abelian group whose theory is the  theory $Th(\Zz)$ of $(\Zz,<,+)$. 
(We have  $\tA \meet (\Qq \tensor A)  = \{a/n \in 
 \Qq \tensor A:  (\forall p) ( n | h_p(a)) $.)
 We have
a homomorphism $[X] \mapsto [X(\tA)]$ from $\SG(\G[*]) $ to $\SG(Th(\Zz)_A)$.  On 
the other hand the polynomial formula for the number of integral points in a
polytope defined by linear equations over $\Zz$ yields a homomorphism 
$\K(Th(\Zz) \bdd [*]) \to \Qq[A]$.  By composing we obtain a homomorphism
$\K(\G \bdd [*]) \to \Qq[A]$.

{\bf Remark} Using integration by parts one can see that the homomorphism $\K(Th(\Zz) \bdd [*]) \to \Qq[A]$
above is actually an isomorphism.

 \ssec{The measured case}  

We repeat the definition of ${\mG}$
 from the introduction, along with two related categories.
  The category $\volG$ corresponds to integrable volume forms, i.e. those that can be transformed by a definable change of
variable to the standard form on a definable
  subsets of affine $n$-space.  By \lemref{seq-d}, the liftability 
condition in (2)   is equivalent to being piecewise in $GL_n(\Zz) \sd A^n$, $A^n$ being 
the group of definable points.  
%	In the definition of $\volG$, it is important to remember
%	that finite disjoint unions are allowed; otherwise in $\volG[1]$ for instance no two objects
%	are isomorphic, and direct sums cannot be formed.

 \<{defn} \lbl{Gmcat}  
 (1)   For $c= (c_1,\ldots,c_n) \in \G^n$, let $\sum (c) = \sum_{i=1}^n c_i$.

(2)  For $n \geq 0$ let  ${\mG}[n]$ be the category whose objects are pairs $(X,\om)$,
with $X \in \Ob \G[n]$ and $\om: X \to \G$ a definable map.  
 A morphism $(X,\om) \to (X',\om')$  is a definable bijection $f: X \to X'$ 
 liftable to 
a definable bijection $\valr \inv X \to \valr \inv X'$, such that
  $\sum(x) + \om(x) = \sum (x') + \om'(x')$ for $x \in X, x' = f(x)$.

 (4) Let ${\mG} \bdd [n]$ be the 
 full subcategory of ${\mG}[n]$ with objects $X \subseteq [\g,\infty)^n$
 for some $\g \in \G$.   
 
(3)  Let $\Ob {\volG}[n] $ be the set of finite disjoint unions of definable subsets of $\G^n$.
 Given
$X,Y \in \Ob {\volG}[n]$,  $f \in \Mor_{\volG[n]}(X,Y)$ iff $f \in \Mor_{\G[n]}$ and
$\sum(x) = \sum(f(x))$ for $x \in X$.  

(5) ${\mG}[*]$ is the direct sum of the ${\mG}[n]$, and similarly for the related categories.

\>{defn}

Recall the Grothendieck rings of functions from \secref{conv}.    $Fn(\G,\SG(\G))$
is a semigroup with pointwise addition.    We also have a convolution
product:   if $f$ is represented by a definable $F \subseteq \G \times \G^m$,
in the sense that $f(\g) = [F(\g)]$, 
and $g$ by a definable $G \subseteq \G \times \G^n$, let 
 $$f*g(\gamma) = [\{(\a,b,c): \a \in \G, b \in F(\a), c \in G(\g - \a)\}]$$
The coordinate $\a$ in the definition is needed in order to make the
union disjoint.  In general, it yields an element represented by a subset of $\G \times \G^{m+n+1}$ rather than $m+n$.   But
let   $Fn(\G,\SG(\G))[n]$   be the set of $[F] \in Fn(\G,\SG(\G[n]))$ such that 
$\dim(F(a)) < n$ for all but finitely many $a \in \G$.
If $f \in  Fn(\G,\SG(\G))[m]$ and $g \in Fn(\G,\SG(\G))[n]$ then $f*g \in Fn(\G,\SG(\G))[m+n]$.
Let $Fn(\G,\SG(\G))[*] = \oplus_m  Fn(\G,\SG(\G))[m]$,
a graded semiring.

\<{lem} \lbl{gamma-volume}  

(1)  $\SG(\mG) [n] \iso  Fn(\G,\SG(\G))[n]$

(2)   $\SG {\mG} \bdd[n] \iso \{f \in Fn(\G, \SG(\G \bdd))[n] :  (\exists \g_0)(\forall \g < \g_0)(f(\g)=0)\}$ 

  (3)    $\SG {\volG}[n] \iso Fn(\G,\SG(\G[n-1]))$
    \>{lem} 
    
 \prf     (1) 
Let $(X,\om) \in \Ob \mG[n]$, with $X \subseteq \G^n$ and $\om: X \to \G$.  
Let $d(x) =\om(x) + \sum(x)$.  
 For $a \in \G$, let $X_a = \{x \in X: d(x) =a \}$.  
This determines an element $F(X,\om) \in Fn(\G,\SG(\G[n]))$, namely $a \mapsto [X_a]$.
It is clear from additivity of dimension that $\dim(X_a) < n$ for all but finitely 
many $a$; so $F(X,\om) \in Fn(\G,\SG(\G))[n]$.
 If $h \in \Mor_{{\mG}[n]}(X,Y)$,
 then by definition of ${\mG}$ we have $h(X_a)=Y_a$; so $[X_a]=[Y_a]$ in $\SG(\G)[n]$.
 Conversely if for all $a \in \G$ we have $[X_a]=[Y_a]$ in $\SG(\G)[n]$, 
then $\valr \inv (X_a), \valr \inv(Y_a)$ are $a$-definably isomorphic.   
By \lemref{collect}
 there exists a definable $H: \valr \inv (X) \to \valr \inv (Y)$ 
  such that for any $x \in \valr \inv (X)$,  $H(x) = h_a(x) $ where $a= \sum \valr (x)$. 
 Clearly $H$ descends to $\bar{H}: X \to Y$; by construction $\bar{H}$ lifts to $\RV$,
and preserves $\sum + \om$, so ${\bar H}   \in \Mor_{{\mG}[n]}(X,Y)$.  We have thus shown
that $[X] \mapsto [F(X)]$ is injective.  It is clearly a   semiring homomorphism.

For surjectivity, let $g \in Fn(\G,\SG(\G))[n]$ be represented by $G \subseteq \G \times \G^n$. 
It suffices to consider either   $g$ with singleton support $\{\g_0\}$, or $g$ such that 
$\dim(G(a)) < n$ for all $a  \in \G$.  In the first case, $g = F(X,\om)$ where
$X = G(\g_0)$ and $\om(x)  = \g_0 - \sum(x)$.    In the second:
after effecting a partition and a permuation of the variables, 
we may assume $G(a) \subseteq \G^{n-1} \times \{\psi(a)\}$ for some definable function 
$\psi(a)$.  With another partition of $\G$, we may assume   $g$ is supported on $S \subseteq \G$, i.e. 
$g(x) = 0$ for $x \notin S$, and $\psi$ 
  is either injective or constant on $S$.    In fact we may assume $\psi$ is
  injective on $S$:   if $\psi$ is constant on $S$, let $G' = 
  \{(a,(b_1,\ldots,b_{n-1},b_n+a)):  (a,(b_1,\ldots,b_n)) \in G, a \in S \}$.  Then 
    $G'$ also represents $g$, and for $G'$ the function $\psi$ is injective.  
Now   let $X = \union_{a \in S} G(a)$, and let $\omega(x) = - \sum(x) + \psi^{-1} (x_n)$.      Then
$F(X,\om) =g$.

(2) follows from (1) by restricting the isomorphism.  

(3) is proved in a similar manner to (1) though more simply and we omit the details.  The key point is that $GL_n(\Zz)$ acts transitively on $\Pp^n(\Qq)$; this can be seen as a consequence
of the fact that finitely generated torsion free Abelian groups are free.  More specifically the
co-vector $(1,\ldots,1)$ is $GL_n(\Zz)$-conjugate to $(1,0,\ldots,0)$.  Thus
the catgegory $\volG[n]$ is equivalent to the one defined using the weighting $x_1$
in place of $\sum(x_i)$.  For this category the assertion is clear.

    \eprf

This lemma reduces the study of $\SG({\mG})$ to that of $\SG(\G)$.

 \>{section}  %The Grothendieck ring of $\G$}

\<{section}{The Grothendieck semirings of $\RV$}
\ssec{Decomposition  to $\G,\RES$.}

Recall that $\RV$ is a structure with an exact sequence
$$0 \to \k^* \to RV \to_\valr \G \to 0$$
We study here the Grothendieck semiring of $\RV$ in a theory $\T_{\RV}$ satisfying 
the assumptions of  \lemref{seq-d}.   The intended case is the structure induced from
 $ACVF_A$  for some $\RV,\G$-generated base structure $A$.)

We show that the Grothendieck ring of $\RV$ decomposes
into a tensor product of those of $\RES$, and of $\G$.  

The category $\G[*]$ was described in \secref{exgam}.  We used $GL_n(\Zz)$ rather
than $GL_n(\Qq)$ morphisms.  The reason is given by:  

\<{lem} \lbl{fcg}   The 
morphisms of $\G[n]$ are precisely those definable maps that lift to morphisms of $\RV[n]$.
The map $X \mapsto \valr^{-1}(X)$ therefore induces a functor $\G[n] \to \RV[n]$,
yielding an embedding of Grothendieck semirings $\SG[\G[n]] \to \SG[\RV[n]]$. 
  \>{lem}
 
\<{proof}  Any morphism of $\G[*]$ obviously lifts to $\RV$, since $GL_n(\Zz)$
acts on $C^n$ for any group $C$.   The converse is a consequence of \lemref{rvfn}.  \>{proof}

We also have an inclusion morphism $\SG(\RES) \to \SG(\RV)$.

Observe that  $\sggf$ forms a part of both $\SG(\RES[*])$ and $\SG(\G[*])$:  the embedding of $\SG(\G[*])$ into $\SG(\RV[*])$ takes $\sggf$
 to a subring of $\SG(\RES[*])$, namely the subring generated by the pullbacks $\valr(\g)$,
 $\g \in \G$ a definable point.  
 
 Given two semirings $R_1,R_2$ and a homomorphism $f_i: S \to R_i$, define $R_1 \tensor_S R_2$ by the universal property for triples $(R,h_1,h_2)$, with $R$ a semi-ring and 
 $h_i: R_i \to R$ a semiring homomorphism, satisfying $h_1f_1 = h_2f_2$.  
 
 We have a natural map $\SG(\RES) \tensor \SG(\G[*]) \to \SG(\RV)$,
 $[X] \tensor [Y] \mapsto [X \times \valr \inv (Y)]$.
  By the universal property it induces
 a map on  $\SG(\RES) \tensor_{\sggf}  \SG(\G[*])$.  A typical element of the image
 is represented by a definable set of the form $\du (X_i \times \valr \inv (Y_i))$,
 with $X_i \subseteq \RES^*, Y_i \subseteq \G^*$.

 \<{prop} \lbl{tensor} The natural map 
$\SG(\RES) \tensor_{\sggf}  \SG(\G[*]) \to \SG(\RV)$
is an isomorphism. \>{prop}

\prf  Surjectivity is \corref{seq-c}.   We will prove injectivity.  In this proof, $X \tensor Y$ will
always denote an element of $\SG(\RES) \tensor_{\sggf}  \SG(\G[*])$.

\Claim{1} Any element of $\SG(\G[*])$ can be expressed as $\sum_{j=1}^l [Y_j] \times \{p_j\}$,
for some $Y_j \subseteq \G^{m_j}$, $\dim(Y_j) = m_j$, and $p_j \in \G^{l_j}$.  

\prf  Let $Y \subseteq \G^m$ be definable.  If $\dim(Y)<m$, then $Y$ can be partitioned
into finitely many sets $Y_j$, each of which lies in some definable affine hypersurface
$\sum_{i=1}^m \alpha_i x_i = c$, with $\alpha_i \in \Qq$, not all $0$.  In other words
$x \mapsto \alpha \cdot x $ is constant on $Y_j$, where $\alpha = (\alpha_1,\ldots,\alpha_m)$.
We may assume
that each $\alpha_i \in \Zz$ and that they are relatively prime.
Then $(\alpha)$ is the first row of a matrix $M \in GL_m(\Zz)$.
The map $x \mapsto Mx$ takes $Y_j$
to a set of the form $Y_j' \times \{c\}$, $Y_j' \subseteq \G^{m-1}$.  Since $[MY_j]=[Y_j]$ in
$\SG(\G[*])$,    the Claim follows by induction.
     \eprf

\Claim{2}  Any element of  $\SG(\RES) \tensor_{\sggf}  \SG(\G[*])$ can be represented
as $\sum_{i=1}^k X_i  \tensor \valr \inv Y_i$, where $X_i  \subseteq \RES^{n_i}$ and $Y_i  
\subseteq \G^{m_i}$ are definable sets, and $m_i = \dim Y_i$.  

\prf   By definition  of  $\SG(\RES) \tensor_{\sggf}  \SG(\G[*])$ and by Claim 1, any element is a sum of tensors $X \tensor \valr \inv (Y \times \{p\})$; using the $\tensor_{\sggf} $-relation, 
$X \tensor \valr \inv (Y \times \{p\}) = (X \times \valr \inv(p)) \tensor Y$.     \eprf

Now let $X_i,X_i' \subseteq \RES^*, Y_i,Y_i' \subseteq \G^*$ be definable sets, and let   
$$F: \du (X_i \times \valr \inv (Y_i)) \to \du (X'_{i'} \times \valr \inv (Y'_{i'}))$$
 be  a definable isomorphism.  Let $m$ be the maximal dimension $m$ of any $Y_i$ or $Y'_{i'}$.
Assume (by Claim (2)):

 (*)    for each  $i'$,
 $Y'_{i'} \subseteq \G^{\dim(Y'_{i'})}$ and similarly for the $Y_i$. 

\Claim{3}    

Let $P$ be a complete type of $Y_i$ of dimension $m$,
and $Q$ a complete type of $X_i$.   
 Then $F(Q \times \valr \inv P) = Q' \times \valr \inv P'$
where $Q'$ is a complete type of some $X'_{i'}$, and $P'$ a complete type type
of $Y'_{i'}$.    

Moreover there exist definable sets $\tP,\tQ,\tP', \tQ '$ containing
$P,Q,P',Q'$ respectively, such that 

\begin{enumerate}
  \item $F$ restricts to a bijection $\tQ \times \valr \inv \tP  \to \tQ' \times \valr \inv \tP '$
  \item    there exist definable bijections $f: \tP \to \tP'$ and
$g: \tQ \to \tQ'$.  
  \item For any $x \in \tQ, y \in \tP$, $F$ restricts to a bijection
$\{x \} \times \valr \inv (y) \to \{f(x)\} \times \valr \inv (g(y))$.
  \end{enumerate}  

\prf  By \lemref{orthplus}, $\valr \inv (P)$ is a complete type; by the same lemma,
 $Q  \times \valr \inv (P)$   is complete; hence so is $F(Q  \times \valr \inv (P))$.
   We have $F(Q \times \valr \inv (P)) \subseteq (X'_{i'} \times \valr \inv (Y'_{i'})) $
for some $i'$.  Let $Q' = pr_1( F(Q \times \valr \inv (P))), V' = pr_2(F(Q \times \valr \inv (P)))$, $P'=\valr(V') \subseteq Y'_{i'}$.  
where $pr_1: X'_{i'} \times \valr \inv (Y'_{i'}) \to X_i \subseteq \RES, pr_2:  X'_{i'} \times \valr \inv (Y'_{i'})
\to  \valr \inv (Y'_{i'})$ are the projections.  Then $Q',V',P'$ are complete types.  We have $m = \dim(P') \geq \dim(Y'_{i'})$ so by maximality of $m$, equality holds.  We thus
have $P' \subseteq \G^{\dim(P')}$.  
By \lemref{orthplus}, $Q'  \times \valr \inv (P')$ is also complete type.
So $F(Q \times P) = Q' \times \valr \inv P'$.                     
 
By one more use of \lemref{orthplus}, the function $f_y: x \mapsto pr_1 F(x,y)$, whose
graph is a subset of the stable set  $Q \times Q'$, cannot depend on $y \in P$.
Thus $f_y=f$, i.e. $F(x,y) = (f(x),pr_2 F(x,y))$.

 Since $Q \times \valr \inv(y)$ is stable, $\valr pr_2 F$ must be constant on it; so 
 $\valr pr_2 F(x,y) = g(y)$ on $P \times Q$.  This shows that (3) of the   ``moreover''
 holds on $P \times Q$.  By compactness, it holds on some definable  $ \tQ \times \tP$
 (and we may take $f$ injective on $\tQ$, and $g$ on $\tP$.)  
 Let $\tQ' =f(\tQ), \tP' = g(\tP)$.  Then (1,2) hold also.    
      \eprf

\Claim{4} Assume (*) holds.  Then there exist finitely many definable $Y_i^j$  ($j=0,\ldots,N_i$)
and $X_i^j$ such that $\dim(Y_i^0)<m$, and the conclusion of Claim 3 holds on 
each $X_i^j \times \valr \inv Y_i^j$ for $j \geq 1$.    Moreover we may take the $Y_i^j, X_i^j$
pairwise disjoint.  
 
\prf  This follows from Claim 3 by compactness; the disjointness can be achieved by
noting that if Claim 3 (3) holds for $\tP,\tQ$ then it holds for their definable subsets too.
  \eprf

 We now show 
  that if $ \du (X_i \times \valr \inv (Y_i))$ and$ \du (X'_{i'} \times \valr \inv (Y'_{i'}))$
  are definably isomorphic then 
   $\sum_{i' }[X_{i'}]\tensor [Y_{i'}] = \sum_i [X_i \tensor Y_i] $.
   We use induction on  
  the maximal dimension $m$ of any $Y_i$ or $Y'_{i'}$,   and also
  on  the number of indices $i$ such that $\dim(Y_i) =m$.  Say $\dim(Y_1)=m$.  

  By Claim 2,  
 without changing $\sum_{i'} X'_{i'} \tensor \valr \inv (Y'_{i'}))$ as an element
 of  $\SG(\RES) \tensor_{\sggf}  \SG(\G[*])$,  we can arrange that 
$ \dim(Y_{i'}) = m_{i'}$, i.e.  (*)  holds. 
   So Claims 3,4 apply.   
   
  The  $Y_1^j$ for $j \geq 1$ may be removed from $Y_1$,
 if their images are correspondingly excised from the appropriate $Y'_j$, 
 since   $[\tQ] \tensor _{\sggf} [\tP] = [f(\tQ)] \tensor _{\sggf} [g(\tP)]$.  What is left in $Y_1$ has
 $\G$-dimension $<m$, and so by induction the classes are equal.

The injectivity and the Proposition follow.   

\eprf

For applications to $\VF$, we need a version of \propref{tensor} keeping track of dimensions.
 Below, the tensor product is in the category of graded semirings.
 
\<{cor}\lbl{tensor-graded} 
  The natural map 
$\SG(\RES[*]) \tensor_{\sggf}  \SG(\G[*]) \to \SG(\RV[*])$
is an isomorphism.   
 \>{cor}

\prf   
For each $n$ we have a surjective   homomorphism
 $$\oplus_{k=1}^n \SG(\RES[k]) \tensor \SG(\G[n-k])  \to \SG(\RV[n])$$
 $\SG \RV[n]$ can be identified with a subset of the semiring $\SG \RV$, namely
$\{[X]: \dim(X) \leq n \}$.
The proof of \propref{tensor} shows that the kernel is generated
by relations of the form 
$$(X \times \valr \inv (Y) ) \tensor Z = X \tensor ( Y  \tensor Z)$$
when  $Y \in \sggf$ and $\dim(X)+ \dim(\valr \inv (Y))+ \dim(\valr \inv(Z))=n$.  These relations are taken into account
in the ring $   \SG(\RES[*]) \tensor_{\sggf}  \SG(\G[*]) $, so that the natural map  
 $ \SG(\RES[*]) \tensor_{\sggf}  \SG(\G[*]) \to \SG(\RV[*])$ is injective, hence an isomorphism. \eprf

%

%If  $R$ is a  ring and $A$ is a subgroup of $R$,   and $B$ of $R'$,
%a  group homomorphism $f: A \to B$ 
%  is said to be a {\em partial ring homomorphism} if
%$f(ab)=f(a)f(b)$ whenver $a,b,ab \in A$.   If a homomorphism  $g: B \to A$ is also
%given, and  $g \circ f = Id_B$
%we say that $f$ is a retraction of partial rings.

Recall the classes $e_a = [\{a\}])_1$ in $\K(\G[1])$, defined for $a \in \G(<\emptyset>)$
They are in $\sggf$, hence identified
with classes in $\K(\RES[1])$, namely $e_a = [\valr \inv (a) ]$.   When denoting
classes of varieties $V$ over the residue field, we will write $[V]$ for $[V(\k)]$, when no
confusion can arise.

\<{defn} \lbl{!K}  Let $I!$ be the ideal of $\K(RES[*])$ generated by all differences
$e_a -e_0$, where $a \in \G(<\emptyset>)$.   Let $\Ks(\RES[*]) = \K(\RES[*])/ I!$.  \>{defn}

By \lemref{oag1} (3), the natural homomorphism  $\K(\RES[*])$ into the localization of
$\K(\RV[*])$ by all classes $e_a$ factors through $\Ks(\RES[*])$.

Since $I!$ is a homogeneous ideal,  $\Ks(\RES[*])$ is a graded ring.

   The  theorem that follows, when combined with the canonical isomorphisms
 $\K(\VF[n]) \to \K( \RV[\leq n])/\Isp$ and  $\K(\VF) \to  \K(\RV[*])/ \Isp $, %\oplus_{n \geq 0} \K( \RV[n]) / \Isp$,  
yields homomorphisms
$$\Xint{{\rm R}}: \K( \VF) \to  \Ks( \RES)[[\Aa_1(\k)] \inv ]$$
$$\Xint{{\rm R}}': \K( \VF) \to  \Ks( \RES) $$

 \<{thm}\lbl{retract2}  \begin{enumerate}
  \item 
There exists a group homomorphism  
$$\Ee_n:  \K( \RV[\leq n])/\Isp  \to \Ks( \RES  [n])$$
 with:
  $$[\RVp]_1 \mapsto -[\Aa^{n-1} \times G_m]_n$$ and
$$[X]_k \mapsto [X \times \Aa^{n-k}]_n$$
for $X \in \RES[k]$
   \item There exists a   ring homomorphism  $\Ee  :   \K(\RV[*])/\Isp   \to  \Ks( \RES)[[\Aa_1] \inv ]$ with $\Ee([X]_k) = [X]_k / \Aa^k $ for $X \in \RES[k]$.
  \item  There exists a group homomorphism 
$${\mathcal E}_n': \K (\RV[\leq n])/\Isp \to \Ks (\RES [n])$$
with  $[\RVp]_1 \mapsto 0$, and
$[X]_k \mapsto [X]_n$ for $X \in \RES[k]$.

   \item There exists a   ring homomorphism  $\Ee ' :   \K(\RV[*])/\Isp   \to  \Ks( \RES)  $ with $\Ee([X]_k) = [X]_k   $ for $X \in \RES[k]$.
 
\end{enumerate}
 \>{thm}

\prf  
(1) 
We first define a homomorphism $\chi[m]:  \K ( \RV[m]) \to \Ks(\RES[m])$.
 By \corref{tensor-graded},  
$$\K (\RV[m] )= \oplus_{l=1}^m \K (\RES[m-l]) \tensor_{\sggf} \K ( \G[l]  )$$
Let $\chi_0=Id_{\K \RES[m]}$.  For $l \geq 1$
recall the homomorphism ${\chi}: \K ( \G[l]) \to \Zz$ of \lemref{euler1}. 
It induces
$\chi_l: \K (\RES[k]) \tensor_{\sggf}  \K ( \G[l]  ) \to \Ks(\RES[k]) $
by $  a \tensor b \mapsto \chi(b) \cdot [G_m]^l \cdot a$.

Define  a  group homomorphism
$$ \chi[m]: \K ( \RV[m]  )   \to \K (\RES[m])  , \ \ \ \chi[m]=\oplus_l \chi_l$$ 
 We have $$\chi[m_1+m_2](ab) = \chi[m_1](a)\chi[m_2](b)$$ when $a \in \K(\RV[m_1]), b \in \K(\RV[m_2])$.   This can be checked
on homogeneous elements with respect to the grading $ \oplus_{l} \SG (\RES[m-l]) \tensor \SG ( \G[l]  )$.

We compute:  $\chi[1]([\RVp]_1) = \chi_1(1 \tensor [\G^{>0}]_1)  = -[G_m]  \in \K(\RES[1])$.

Next define a group homomorphism $\beta_m: \Ks( \RES[m]) \to \Ks(\RES[n])$
 by $\beta_m([X]) = [X \times \Aa^{n-m}]$.  

Define $\gamma: \oplus_{m \leq n} \K( \RV[m]) \to \Ks(\RES[n])$ by $\gamma = \sum_m \beta_m \circ \chi[m]$.

Then $\gamma$ is a group homomorphism, and
 $\gamma(a) \gamma(b)  = \gamma(ab) \times [\Aa^n]$
for $a \in \K(\RV[m_1]), b \in \K(\RV[m_2])$, $m_1+m_2 \leq n$.
  Again this is easy to verify on homogeneous
elements.

Finally we compute $\gamma$ on the standard generator $J=[\RVp]_1 + [1]_0 - [1]_1$ of $\Isp$.
Since $\chi[1]([\RVp]_1) =  -[G_m] $, we have  
$$\gamma([\RVp]_1) = \beta_1( -[G_m]) =
-[G_m \times \Aa^{n-1}]_1$$
  On the other hand 
  $$\gamma([1]_0) = \beta_0([1]_0) = [\Aa^n]_n$$
   $$\gamma([1]_1) = \beta_1([1]_1) = [\Aa^{n-1}]_n$$
   So  
$\gamma(J) = [\Aa^{n-1}]_{n-1} \times (-[G_m]_1 + [\Aa^1]_1 - [1]_1) = 0$.  
 A  homomorphism $\K ( \RV[\leq n]) /\Isp \to \!K ( \RES [n])$ is thus induced.

(2)  For $a \in \K(\RV[m])$, let $\Ee(a) = \beta_m(a) / [\Aa^m]$.  For any large enough $n$
we have $\Ee(a) = \Ee_n(a) / [\Aa^n]$.  
The formulas in (1) prove
that $\Ee$ is a ring homomorphism.  

(3,4)   The proof   is similar, using $\chi'$ from \lemref{euler2} 
 in place of  $\chi$ of \lemref{euler1}, and the identity in place of $\beta_m$.

\eprf

\<{cor}\lbl{embkrv}  The natural morphism $\K (\RES[n]) \to  \K( \RV[\leq n] ) / \Isp$ has kernel contained in $!I$.  \qed
\>{cor}

\<{lem}  \lbl{ls0}  Let $\T=\ACVF_{F((t))}$ or $\T=\ACVFR_{F((t))}$, $F$ a field of characteristic $0$, with $\val(F) = (0)$, $\val(F((t)))=\Zz$, 
and $\val(t) = 1 \in \Zz$.    Then there exists a retraction 
$\rho_t: \SG(\RES) \to \SG(Var_F)$.  It induces a retraction $\Ks(\RES) \to \K(\Var_F) $\>{lem}

\prf      
Let $t_n \in F((t))^{alg}$ be such that $t_1=t$ and $t_{nm}^n  = t_m$.   For $\a = m/n \in \Qq$,
with $m \in \Zz$, $n \in \Nn$, let $t_\a = t_n^m$.  So $\a \to t_a$ is a homomorphism
$\Qq \to G_m(F((t))^{alg})$. 

Let   $V(\a) = \valr \inv (\a)$.  Let $\bt _\a = \rv(t_\a)$.  Then $\bt_a \in V(\a)$.  

Let $X \in \RES[n]$.  Then for some $\a_1,\ldots,\a_n \in \Qq$ we have
 $X \subseteq \Pi_{i=1}^n V(\a_i)$, where $V(\a_i) = \valr \inv (\a_i)$.  Define
 $f(x_1,\ldots,x_n) = (x_1/\bt_{\a_1},\ldots,x_n / \bt_{\a_n})$.  Then $f$ is $F((t^{1/m}))$-definable for some $m$, but not in general definable.  Nevertheless, $F(X) =: Y$ is 
 definable.  This is because the Galois group $G=Aut(F^a((t^{1/m}))/F^a((t)))$ extends
 to a group of valued field automorphisms $Aut(\k((t^{1/m}))/\k((t)))$ fixing the entire residue field $\k$; while $Y \subseteq \k$; thus $G$ fixes $Y$ pointwise and hence setwise.

The map $X \mapsto Y$ of definable sets described above clearly respects disjoint unions. 
It also respects definable bijections:  if $h: X \to X'$ is a definable bijection, $Y=f(X),Y'=F(Y')$,
then $f h f \inv$ is an $F((t^{1/\infty}))$-definable bijection $Y \to Y'$; by the Galois argument
above, it is in fact definable.

The definable subsets of $\k$ are just the $F$-constructible sets.
Thus we have an induced homomorphism $\rho_t: \SG(\RES) \to \SG(Var_F)$; it is clearly
the identity on $\SG(\RES)$.  It induces a homomorphism $\K(\RES) \to \K(\Var_F)$. 

Finally $\rho_t( \valr \inv(\a)) = [G_m]$ for any $\a \in \Qq$; so a homomorphism
on $\Ks(\RES)$ is induced.
\eprf 

This example can be generalized as follows.  Let $L$ be a valued field with residue field $F$ of
  characteristic $0$, $\T=\ACVF_L$ or $\ACVFR_L$.  Let $A= \res(L)$, $\bA = \Qq \tensor A$,
and let $t: \bA \to G_m(L^a)$ be a monomorphism, with $t(A) \subseteq G_m(L)$.
Then there exists a retraction $\rho_t: \SG(\RES) \to \SG(Var_F)$.

From \thmref{retract2} and \lemref{ls0} we obtain the example discussed in the introduction:

\<{prop} \lbl{ls1}  Let $\T=\ACVFR_{F((t))}$, $F$ a field of characteristic $0$, with $\val(F) = (0)$
and $\val(t) = 1 \in \Zz$.   Then there exists a ring homomorphism
$\Ee_t:  \K(\VF) \to \K(\Var_F)[[\Aa^1] \inv]$, with $[\Mm] \mapsto  -  [G_m]/[G_a]$,
$\L([X]_k) \mapsto [X]_k / [\Aa^k]$ 
for $X \in \Var_F[k]$.  There is also a ring homomorphism
$\Ee_t':  \K(\VF) \to \K(\Var_F)$ with $\L([X]_k) \mapsto [X]_k$.
  \>{prop}

\ssec{Decomposition of $\RVm$}

An analogous decomposition is valid for the measured Grothendieck
semiring $\RVmg$ (\defref{RVm}).

\lemm{mgtorvm}  There exists a homomorphism $\SG  {\mG}[n] \to \SG   \RVmg [n]$
 with $[(X,\om)] \mapsto [(\valr \inv (X), Id, \om \circ \valr)]$.  \>{lem}
 
\prf 
We have to show that an $\mG[n]$-isomorphism $X \to Y$ lifts
to an $\RVmg [n]$-isomorphism. This follows immediately from the definitions.
\eprf
 
    Recall $\RESmg$ from \defref{RVm}.  
  Along the lines of \lemref{gamma-volume},
 we can also describe $\SG \RESmg[n]$ as the semigroup of functions 
 with finite support
 $\G \to \SG(\RES [n])$. 
  We also have
  the inclusion $\SG \RESmg [n] \to \SG \RVmg [n]$, $[(X,f)] \mapsto [(X,f,1)]$.   
   
  Let ${\mG}^{fin}[n]$ be full subcategory
  of ${\mG} [n]$ whose objects are finite.  We have a homomorphism
  $\sggfm[n] \to  \RESmg [n]$, $(X,\om) \mapsto (\valr \inv(X),Id,\om \circ \valr)$.  
    As before  we obtain 
  a homomorphism $\SG \RESmg [*]  \tensor_{\sggfm} \SG({\mG}[*]) \to \SG(\RV[*])$.  
  
%:

Let  $\RESvolg$ be the full subcategory of 
  $\RVvolg$  whose objects are in $\RES$; this is the same as $\RV$ except that
  morphisms must respect $\sum \valr$.  
  Let  $\volG^{fin}$ be the subcategory of  finite objects of $\volG$.   
  
\<{prop} \lbl{tensor-m} 

(1)  The natural map 
$\SG(\RESmg[*]) \tensor_{\sggfm}  \SG({\mG}[*]) \to \SG(\RVmg[*])$
is an isomorphism.

(2)  So is $\SG(\RESvolg[*]) \tensor_{\SG(\volG^{fin}[*])}  \SG(\volG[*]) \to \SG(\RVvolg[*])$.

(3)  The decompositions of this section preserve
the subsemirings of bounded sets.  
  
\>{prop}

\prf  We first prove  surjectivity in (1).  By the surjectivity in \corref{tensor-graded}, it suffices to consider
a class $c=[(X \times \valr \inv (Y), f, \om)]$ 
with $X \in \RES[k], Y  \subseteq \G^l$, $f(x,y)=(f_0(x),y)$, and $\om: X \times (\val r \inv (Y)) \to \RV$.  In fact as in \propref{tensor} we may take $\dim(Y) = l$, and inductively we may assume
 that any   class $[(X' \times Y', f', \om')]$
with $\dim(Y')<l$ is in the image.  Since we may remove a subset of $Y$ of smaller dimension, applying
 \lemref{orthplus}  to 
$\om: X \times \valr \inv (Y) \to \G$, 
  we may assume $\om(x,y) = \om'(\g)$ when $\valr(y)=\g$.  Now 
 $c = [(X,f_0,1)] \tensor [(Y,\om')]$.
 
The proof of surjectivity in (2) is similar.  %, without $\om$.  

 The proof of   injectivity in (1,2) is the same as of \propref{tensor} and \corref{tensor-graded}.  (3) is clear by inspection of the homomorphisms.  \eprf

 \def\D{\Delta}
 We now deduce \thmref{d+}.   For a finite extension $L$ of $\Qq_p$, write
  $\vol_L(U)$ for $\vol_L(U(L))$.
 Let $r$ be the ramification degree, i.e. $\val(L^*) = (1/r) \Zz$.  Let $Q=q^r$.
 The normalization is such that $\Mm$ has volume $1$; so an open ball of 
 valuative radius $\g$ has volume $q^{r\g} = Q^\g$.  Thus the volume of $\valr \inv (\g)$ 
 is $(q-1) Q^\g$.   Also the norm satisfies  $|y| = Q^{\val(y)}$.
 
  \proof[Proof of \thmref{d+}]  For ${a} \in \G^k$ let 
 $Z({a}) =  \{x \in \Oo_L^n: \val(f_1(x))={a}_1 \ldots \val(f_k(x))={a}_k \} $.  Then 
 
$$ \int_{\Oo_L^n} |f|^s = \sum_{{a} \in (\G^{\geq 0})^k}  {Q}^{  s \cdot {a}} \vol_L(Z({a})) $$
According to \propref{tr2} and \propref{tensor-m} we can write
$$ Z({a})  \sim   \du_{i=1}^\nu \L \bX_{i}   \times \L \D_{i}({a})  $$
where  $\D_{i} $ is a definable subset of $ \G^{k+n_2({i})} $, $h^{i}: \D_{i} \to \G^k$ the projection to the first $k$ coordinates, 
$\D_{i}(a)= \{d \in \G^{n_2(i)}:  h^{i}(d) = a \} $,
  $\bX_{i} = (X_{i},f_{i}) \in \RES[n_1({i})] $, and $\sim$ denotes equivalence up to an admissible
transformation.  Thus 

$$ \vol_L(Z({a})) = \vol_L( \du_{{i}=1}^\nu \L \bX_{i}   \times \L \D_{i}({a}) )  = \sum_{{i}=1}^\nu \vol_L(\L \bX_{i}) \vol_L(\L \D_{i}({a}))  $$

  If $b=(b_1, \ldots, b_{k+n_2({i}) }) \in \D_i$,
let $h^i_0(b)$ be the sum of the last $n_2(i)$ coordinates.
 
Since $\valr$ takes only finitely many values on a definable subset of $\RES$, we may
assume $\sum \valr (f(x)) = \g(i)$ is constant on $x \in X_i$. 
 %$f(X_{i}) \subseteq   V_{\g({i},1)} \times \ldots \times V_{\g({i},n_1({i}))}$.
%  Let $\gamma({i}) = \g({i},1)+ \ldots + \g({i},{n_1({i})})$.  
  Then 
 $\vol_L (\L X_{i} (L)) =  Q^{  \gamma({i})} |X_{i}(L)|$.  Thus

\beq{d3}  \int_{\Oo_L^n} |f|^s =\sum_{i} |X_{i}(L)| Q^{\g({i})} \sum _{{a} \in (\G^{\geq 0})^k} {Q}^{s \cdot {a}} \vol_L ( \L \D_{i}({a})) \eeq

Now $\vol_L( \L \D_{i}({a}) ) = \sum_{{b} \in \D_{i}, h({b})={a}} (q-1)^{n_2({i})} Q^{h_0(b)}$.  So

\beq{d4}  \sum _{{a} \in (\G^{\geq 0})^k} {Q}^{s \cdot a } \vol_L ( \L \D_{i}({a})) = 
\sum_{{b} \in \D_{i}}  {Q}^{s_1h^i_1(b)+ \ldots + s_kh^i_k(b)} (q-1)^{n_2({i})} Q^{h_0(b)}
% = \sum_{b \in \D_{i}} Q^{h_0(b)+  s_1h_1(b)+ \ldots + s_kh_k(b)} (q-1)^{n_2({i})}  
= (q-1)^{n_2(i)} \ev_{h^i,s,Q}(\Delta_i) \eeq
The theorem follows from Equations (\ref{d3}),(\ref{d4}).  \qed

Let $A$ be the set of definable points of $\G$.   Recall
that for $X \subseteq \RV$, $[X]_1$ denotes the class $[(X,Id_X)] \in \RV[1]$ of $X$ 
with the identity map to $\RV$,
and the constant form $1$.    
For $a \in A$, 
let $\tilde{e}_a = [(\valr \inv(0),Id,a)] \in \RES[1]$, $f_a = [\{1\}_\k,Id,a] \in \RES[1]$
 where $a$ in the third coordinate
is the constant function with value $a$.  If $a$ lifts to a definable point $d$ of $\RV$,
multiplication by $d$ shows that $\tilde{e}_a = [\valr \inv(a),Id,0], f_a = [\{d\},Id,0]$.    
Note $\tilde{e}_a\tilde{e}_b=\tilde{e}_{a+b}\tilde{e}_0$; and
   $\tilde{e}_0 = [G_m]$.   
Let $\tau_a \in \RES[1]$ be the class of $(\valr \inv ((a,\infty)), Id,0)$.
  The generating relation of $\Ispmg$ is thus 
$(\tau_0,f_0)$ (\lemref{diamond-m} (6)).     Let $\fh$ be the class of
$[(\RVp,Id,x \inv)]$.

Let $!I^0_\mu$ be the ideal of   $ \K( \RESmg [*])$ generated by the relations
$\tilde{e}_{a+b} = [(\valr \inv(a),Id,b)]$, where $a,b \in A$, $b$ denoting the constant function $b$.
Let $!I_\mu$ be the ideal generated by $!I^0_\mu$ 
as well as the element $[\Aa_1]_1$.

\<{thm}\lbl{retract-m}    
There exist  two graded ring homomorphisms   
$$\inte , \inte': \Ke( \VFmg[*]) = \K(\RVmg[*])/\Ispm   \to \K( \RESmg [*]) / {!I_\mu}$$
 such that the composition $K(\RESmg[*]) \to   \K(\RVmg[*])/\Ispm   \to \K( \RESmg [*]) / {!I_\mu}$ equals the natural projection $\pi: \K( \RESmg [*])  \to \K( \RESmg [*]) / {!I_\mu}$;
 with   $$\inte \fh = -[\{0_\k\}]_1, \ \  \inte'  \fh  = 0$$ 
 \>{thm}

 \prf  The identification $\Ke(\VFmg[*]) = \K(\RVmg[*]) / \Ispm$ is given by 
  \thmref{volumes}, and we work with $ \K(\RVmg[*]) / \Ispm$.
 
   According to \propref{tensor-m}, we can identify 
          $$\K( \RVmg[*]) = \K ( \RESmg[*] ) \tensor_{\sggfm}  \K(\mG[*])$$
          
We first construct  two homomorphisms of
 graded rings $R,R': \K( \RVmg[*]) \to  \K( \RESmg [*])/{!I_\mu}$.  
 This 
amounts to finding  graded ring homomorphisms 
$\K(\mG[*]) \to \K( \RESmg [*])/{!I_\mu}$,   agreeing with $\pi$ on the graded ring $\sggfm$.  
It will be simpler to work with $R,R'$ together, i.e. construct 
$$R'' = (R,R'):   \K(\mG[n]) \to  (\K ( \RESmg [n])/{!I_\mu})^2$$ 
Recall from   \lemref{gamma-volume}    the isomorphism

            $$\phi: \K(\mG[n]) \to  Fn(\G,\K(\G))[n]$$
            
 Let $\chi'': \K (\G[n]) \to \Zz^2$ be the Euler characteristic of \propref{2eulers};
so that $\chi'' = (\chi,\chi')$, cf. Lemmas \ref{euler1},\ref{euler2}.
 We obtain by composition a map
 $E_n''=(E_n,E_n'):   Fn(\G,\K(\G[n])) \to Fn(\G, \Zz)^2$.   Here $Fn(\G,\Zz)$ is the group of  
 functions $g: \G \to \Zz$ such that  $g(\G)$
 is finite and $g \inv (z)$ is a definable subset of $\G$ (a finite union of definable intervals and points.)  Thus $Fn(\G,\Zz)$ is freely generated as an Abelian group by 
  $\{p_a,q_a,r\}$ where $r$ is the constant function $1$, and for $a \in A$,  $p_a$,$q_a$ are the characteristic functions of $\{a \}, \{(a,\infty) \}$, respectively.   
  Define $\psi_n: Fn(\G,\Zz) \to \K( \RESmg [*]) $:
 $$\psi_m(p_a) = [G_m]^{n-1} \tilde{e}_a = [G_m]^n f_a, \  \ \psi_n(q_a) = -[G_m]^{n} f_a , \ \ \psi_n(r)=0$$  

 For $u \in \K(\mG[n]) $, let $R''(u) = \psi_n(E''_n(\phi(u)))$. 
 
\Claim{}  $R'': \K(\mG[*]) \to \K( \RESmg [m])^2$ is a graded ring homomorphism.

\prf:   We have already seen
that $\phi$ is a ring homorphism, so it remains to show this for $\psi_* \circ E''_*$.  Now 
by \propref{2eulers}, $\chi''(Y)=\chi''(Y')$ iff $[Y]=[Y']$ in the Grothendieck group of
$\DOAG$.  Hence given a families $Y_t,Y_{t'}$ of pairwise disjoint sets with
$\chi''(Y_t)=\chi''(Y'_t)$, by \lemref{collect} we have $\chi''(\union_t Y_t) = \chi''(\union_t Y'_t)$.
From this and the definition of multiplication in $Fn(\G,\K(\G))[*]$, and the multiplicativity
of $E_n''$,  it follows that if $E_n''(f)=E_n''(f')$ and $E_m''(g)=E_m''(g')$ then
$E_{n+m}''(fg) = E_{n+m}''(fg)$.   In other words, $E''_*$ is a graded homomorphism
from  into $ (Fn(\G, \Zz)^2,\star)$  for some uniquely determined multplication $\star$ 
on $Fn(\G, \Zz)^2$.  Clearly $(a,b)\star(c,d)=(a {*_1} c,b {*_2} d)$ for two operations 
$*_1,*_2$ on $Fn(\G,\Zz)$.  

Now we can compute these operations explicitly on the generators:

$$p_a * p_b = p_{a+b}, \ p_a * q_b = q_{a+b}, \ q_a * q_b = -q_{a+b}$$
for both $*_1$ and $*_2$, and
$$r {*_1} \tilde{e}_a = r, \   r {*_1} q_a = -r, \    r {*_1} r =r$$
$$r {*_2} \tilde{e}_a = -r, \   r {*_2} q_a = 0, \    r {*_2} r = -r$$

composing with $\psi$, we see that $R''$ is indeed a graded ring homomorphism.   \eprf

Let $R,R'$ be the components of $R''$.

\Claim{}  $R,R',\pi$ agree on $\sggfm$.     $R(\tau_0)=R'(\tau_0)=-\tilde{e}_0$.

This is a direct computation.   It folllows that    $R,R'$ induce homomorphisms
$\K( \RVmg[*]) \to  \K( \RESmg [*])/ $.   Since $\tilde{e}_0 + f_0 = [(\Aa_1,Id,0)]$, modulo $!I_\mu$
    both $R,R'$ equalize $\Ispmg$, and  hence induce 
homomorphisms on $\K(\RVmg[*])/\Ispmg \to  \K( \RESmg [*]) / {!I_\mu}$.
 \eprf
 
{\bf Remark}

  The construction is heavily, perhaps completely constrained.  The value of $\psi_m(p_a)$ is determined 
 by the tensor relation over $\sggfm$.  The value of $\psi_m(q_a)$ is determined by
 the relation $\Isp$.  The choice $\psi(r)=0$ is not forced, but the multiplicative
 relation shows that either $r$ or $-r$ is idempotent, so one has a product of two rings,
 with $\psi(r)=0$ and with $\psi(r)= \pm 1$.  In the latter case we obtain the isomorphisms 
 of \thmref{retract2}.  Thus the only choice involved is to factor the fibers of an element of
 $Fn(\G,\K(\G))[n]$  through
 $\chi''$, i.e. through $\K(\DOAG)$.  It is possible that $\K(\G[n]) = \K(\DOAG[n])$
 (cf. Question \ref{q2}).  In this case, $\inte,\inte',\Xint{{\rm R}},\Xint{{\rm R}}'$ are injective as a quadruple, and determine $\K(\VFm[*])$ completely, at least when localized
 by the volume of a unit ball.

\>{section} %\>{section}  %The Grothendieck ring of $\RV$}

\<{section}{Integration with an additive character}

\lbl{additivecharS}

Let $\Om = \VF / \Mm$.   Let $\psi: \VF \to \Om$ be the canonical map.

Motivation:  for any $p$,  $\Om(\Qq_p)$ can be identified with the $p$-th power roots of unity via an additive character on $\Qq_p$.  For other local fields, the universal $\psi$ we use
is tantamount to integration with respect to all additive characters of conductor $\Mm$ at once.  Thus $\Om$
is our motivic analog of the roots of unity, and the natural map $\VF \to \VF / \Mm$
an analogue of a generic additive character.

Throughout this paper we have been able to avoid subtractions and work with semi-groups,
but here it appears to be   essential to work with a group  or at least a cancellation semigroup.   The reason is that we will 
introduce, as the essential feature of integration with an additive character, an 
identification of the integral of a function   $f$ with $f+g$ if $g$ is $\Oo$-invariant.  
This corresponds to the rule that the sum over a subgroup of a nontrivial character vanishes.
Now for any $h: \Om \to \SG(\VFm)$, it is easy to construct $h': \Om \to \SG(\VFm)$ such that
$h+h'$ is $\Oo$-invariant.  Thus if $f+h=f'+h$ for some $h$, then $f=f+h+h' = f'+h+h' = f'$.
So cancellation appears to come of itself.

If we allow all definable sets and volume forms, a great deal of collapsing is caused
by the cancellation rule.  We thus use the classical remedy and  work with bounded sets and volume forms.   The setting is flexible and can be compatible with stricter notions of 
boundedness.  
This is only a partial remedy in the case of higher dimensional local fields,
cf. \exref{hdlf}.
 
The theory can be carried out for any of the settings we considered.
Let $\RR$ be one of these groups or rings, 
 with $\DD$ the corresponding data.  For instance:  $\DD$ the set of pairs
 $(X,\phi)$ with $X$ a bounded definable subset of $ \VF^n \times \RV^*$,
 and $\phi: X \to \RV$ a bounded definable function;  $\RR$ the corresponding
 Grothendieck ring.   Similarly we can take $\G$-volumes, or pure isomorphism
 invariants without volume forms.  In this last case there is no point restricting to bounded sets.
 As we saw, two Euler characteristics into the Grothendieck group of varieties over
 $\RES$ do survive.
 
 In each case, we we think of  $\RR$ as a Grothendieck ring of associated
 $\RV$-data, modulo a canonical ideal.   
 %We denote the $\RV$-version of the data
 %by $\DD^\RV$.  
 %(So we have a canonical map $\L: \DD^\RV \to \DD$, inducing an isomorphism
 %of the Grothendieck rings.)
 
   Everything can be graded by dimension, but for the moment 
   we have no   need to keep track of it,  so in the volume case we can take the direct sum over all $n$ or fix one $n$ and omit it
 from the notation.    

The corresponding group for the theory $\T_A$ or $\T_{<a>}$ will be denoted
$\RR_A,\RR_a$, etc.  When $V$ is a definable set, we let $\DD_V$, $\RR_V$
denote the corresponding objects over $V$.  For instance in the case of bounded
$\RV$-volumes, $\DD_V$ is the set of pairs $(X \subseteq V \times W, \phi: X \to \RV^*)$  
such that for any $a \in V$, $(X_a,\phi | X_a)$ with $X_a$ bounded.

If $\RR$ is our definable analog of the real numbers (as recipients of values of $p$-adic integration), the group ring $\CC=\RR[\Om]$ will take the role of the complex numbers.
We have a canonical group homomorphism $(\VF,+) \to \Om \subseteq G_m(\CC)$, corresponding
to a generic additive character.

Integration with an additive character can be presented in two ways: in terms
of definable functions $f: X \to \Om$ (Riemann style), where we wish to evaluate expressions
such as  $ \int_X  f(x)  \phi (x)$; classically $f$ usually has the form $\psi( h(x))$, where $h$ is 
a regular function and $\psi$ is the additive character.  Or we can treat  definable functions
$F: \Om \to \RR$ (Lebesgue style), and evaluation $\int_{\om \in \Om} F(\om)$.  We will work with the latter.  Given this,
to reconstruct a Riemann style integral, given $f: X \to \Om$, and an
$\RR$-valued volume form $\phi$ on $X$,   let 
$$F(\om) = \int _{f^{-1}(\om)} \phi(x) $$
Then we can define 
$$ \int_X  f(x)  \phi (x) = \int_{\om \in \Om} \om  F(\om)$$

It thus suffices to define the  
integral of a definable function on $\Om$.  Such a function
can be interpreted as an $\Mm$-invariant function on $\VF$.  We impose one rule
(cancellation):  
the integral of a function that is constant on each $\Oo$-class equals zero.  
The integral is  a homomorphism on
the group of $\Mm$-invariant functions $\VF \to \RR$, vanishing on the $\Oo$-invariant ones.
 We give a full description of the quotient group, showing that the universal
homomorphism of this type factors through a similar group on the residue field.

 Recall the group $Fn(V,\RR)$ of \secref{functions}.   We will not need to refer to the 
 dimension grading explicitly.

If $V$ is a definable group, $V$ acts on 
on $Fn(V,\RR)$ by translation.   

\<{defn}  For a definable subgroup $W$ of $V$, let 
$Fn(V,\RR)^W$ be the set of $W$-invariant elements of $Fn(V,\RR)$:  
they are represented by a definable
$X$, such that if $t \in W$
and $a \in V$ then $X[a],X[a+t]$ represent the same class
in $\K(\VFm_{a,t})[n]$.    \>{defn}

\<{lem} \lbl{reps}  An element
of $Fn(\VF,\RR)^{\Mm}$ can be represented by an $\Mm$-invariant $X \subseteq (\VF \times *)$.  \>{lem}

\prf  Let $Y \in \DD^{RV}_\VF$ represent an element of $Fn(\VF,\RR)^{\Mm}$.   So each fiber
$Y_a \in \DD^{RV}$.  
By  \lemref{resolve-cor}, for $\ba \in \VF/\Mm$ one can find $Y'_{\ba} \in \DD^{RV}$ 
such that for some $a \in \VF$ with $a+ \Mm = \ba$, $Y_a = Y'_{\ba}$.  
As in \lemref{collect} there exists $Y' \in \RR_{\VF/\Mm}$ such that 
 $Y'_{\ba}$ to be the fiber of $Y' $ over $\ba$. Pulling back to $\VF$ gives the required
$\Mm$-invariant
representative.  \eprf

Since the equivalence is defined in terms of    effective isomorphism, \defref{eff-iso},
it is clear  that two elements of $\DD_\Omega$ are equivalent 
iff the corresponding pullbacks to $  Fn(\VF,\RR)^{\Mm}$ are equivalent.

The groups  $Fn(\VF,\RR)^{\Mm}$ and $Fn(\VF/\Mm, \RR)$
can thus be identified.

 Note that effective isomorphism
agrees with pointwise isomorphism for   $Fn(\VF,\RR)^{\Mm}$, but not for $Fn(\VF/\Mm, \RR)$.
 
The group we seek to 
describe is   $\AA = \AA_{\T} = Fn(\VF,\RR)^{\Mm}/Fn(\VF,\RR)^\Oo$.  
The quotient corresponds to the cancellation rule discussed earlier.

Let $Fn(\k,\RR)$ be the Grothendieck group of functions $\k \to \RR$,  with 
addition induced from $\RR$.

Let $\CC=\RR[\Om]$ be
the ring
of definable functions $\Om \to \RR$ with finite support, convolution product.

{\bf Remark}  $\CC$ embeds into the Galois-invariant elements of the 
 abstract group ring 
 $ \RR_ \tT [\Om_\tT] $, where 
   $\tT = \T_{\acl(\emptyset)}$.

  The additive group $\k = \Oo/\Mm$
is a subgroup of  $\Om = \VF / \Mm$, and so acts on $\Om$ by translation.  It also
acts naturally on $ Fn(\k,\RR)$.   This gives two actions on
$Fn(\k,\CC) = Fn(\k,\RR)[\Om]$.   Let $Fn(\k,\CC)_\k$  denote the coinvariants
with respect to the anti-diagonal action, i.e. the largest quotient on which the two
actions coincide.
 %Let $\IAk$ be the smallest subgroup equating these two actions,
%$ Write:
 %$$ Fn(\k,\CC)_\k = Fn(\k,\CC) / \IAk$$

 In general,  the upper index denotes invariants, the lower index co-invariants.

 $Fn(\VF,\RR)$ is the ring of  definable functions
 from $\VF$ to $\RR$.     $Fn(\k,\RR)$ is the ring of  definable
 functions from $\k$ to $\RR$.   $Fn(\k,\CC)$ is the ring of  
 definable functions from $\k$ to $\CC$; equivalently, it is the set of Galois-invariant 
 elements of   the group ring $Fn(\k,\RR)[\Omega]$. 
   
   The action of $\k$ on $Fn(\k,\CC)$ is by translation
 on $\k$, and negative translation on $\Omega$ and hence on $\CC$.
 The term (Const) 
refers to the image of the constant functions of $Fn(\k,\CC)$ in $Fn(\k,\CC)_\k$
(it is isomorphic to  $(\CC/\k)$.)

 \<{thm} \lbl{additivechar}  There exists a canonical isomorphism
$ Fn(\k,\CC)_\k  / (Const)  \isomto  Fn(\VF,\RR)^\Mm / Fn(\VF,\RR)^\Oo $
  \>{thm}

\prf    
Let $\AA_{fin} $ be the subring of $Fn(\VF,\RR)^\Mm$ consisting of functions
represented by elements of $Fn(\VF,\DD)^{\Mm}$ 
 whose support projects to a finite subset
of $\VF / \Oo$.

A definable function on $\k$ can be viewed as an $\Mm$-invariant 
function on $\Oo$;  this gives

\beq{add1}   Fn(\k,\RR)  \isomto  Fn(\Oo,\RR)^{\Mm}  \eeq

On the other hand we can define a homomorphism 

\beq{add2}  
 Fn(\Oo,\RR)^{\Mm} [\Om] \to  \AA_{fin}   \ \ : \ \ \ 
    \sum_{\om \in W}  a(\om) \om   \mapsto \sum_{\om \in W} a(\om)_{\om}   \eeq

where $W$ is a finite $A$-definable subset of $\Om$, $a: W \to  Fn(\Oo,\RR)^{\Mm}$
  is
an $A$-definable function,
(so that  $\sum_{a \in W} a(\om) \om$ is a typical element of the group
ring $  Fn(\Oo,\RR)^{\Mm} [\Om] $), 
and $b_{\om}$ is the translation of $b$ by 
$\om$, i.e. $b_{\om}(x) = b(x-\om)$.

  \eqref{add2} is surjective:   let $f \in \AA_{fin}$ be represented
  by $F$, with support $Z$, a finite union of translates of $\Oo$.  
By \lemref{red1c1}  there exists
a finite definable set $W$, meeting each ball of $Z$   in a unique point.  
Define $a: W \to Fn(\Oo,\RR)^{\Mm}$ by 
$$a(\om) = (f | \om+\Oo)_{-\om}$$
  Then \eqref{add2} maps $\sum a(\om) \om$ to $f$.
 
   The kernel of \eqref{add2} is the equalizer of the two actions of $\k$.  Composing with \eqref{add1} we obtain an isomorphism $(Fn(\k,\RR)[\Omega])_\k \isomto \AA_{fin}$,
   or equivalently

\beq{add3}  Fn(\k,\CC)_\k \isomto \AA_{fin} \eeq
 
 The last ingredient is the homomorphism

\beq{add5} \AA_{fin} \to \AA  \eeq
We need to show that it is surjective, and to describe the kernel.

Using the representation $\DD$ of elements of $\RR$ by $\RV$-data, 
An element of $\AA$  is  
 is represented by an  $\Mm$-invariant definable  $W \subset \VF \times \RV^*$.

  By \lemref{Om-rv}, for each coset $C$ of $\Oo$ in $\VF$ apart
from a finite number, $W  \meet (C \times  \RV^{n+l}) $ is
invariant under translation of the first coordinate by elements of $\Oo$. 
   Thus $W$ is the disjoint sum of an $\Oo$-invariant set $W'$
   and a set $W'' \subset \VF \times 
   \RV^{*}$ projecting to a finite union $Z$ of cosets of $\Oo$ in $\VF$, i.e.
   representing a function in $\AA_{fin}$.  
   
 Clearly $W' \times_{\RV^n} \VF^n$ lies in $Fn(\VF,\RR)^\Oo$.   
 
 Thus \eqref{add5} is surjective; the kernel is $\AA_{fin} ^{\Oo}$.  Composing
 \eqref{add3},\eqref{add5} we obtain an isomorphism 
 
 $$ \AA \isomto (Fn(\k,\RR)[\Omega])_\k / (Const)$$
 
 Using the identification $Fn(\k,\RR)[\Omega] \iso Fn(\k,\CC)$, the theorem follows.
 
     \eprf

Note that $Fn(\k,\CC)^\k  \iso \CC$, via $Fn(\k,\CC) \iso Fn(\k \times \Omega, \RR)_{fin}$.
%**********

\ssec{Definable distributions}

 \def\RRdf{{\mathcal R}_{df}}  \def\CCdf{{\CC}_{df}}
  \def\dd{{\mathfrak d}}   
$\RR$ is graded by dimension ($\VF$-presentation) or  ambient dimension ($\RV$-presentation.)
Write $\RR=\oplus_{n \geq 0} \RR[n]$.

Let  $\RRdf$ be the dimension-free version:
first form the localization $\RR[ [0]_1^{-1}] $, where $[0]_1$ is the class of the point
$1 \in \RV$, as an element of $\RV[1]$.  Equivalently, $[0]_1^n$ is the volume of the
open  $n$-dimensional polydisc $\Oo^n$.  Let $\RRdf$ be the zero-dimensional component
of this localization.  %Then $\RR \leq \oplus_{n \in \Nn}  \RRdf [0]_1^n$.
Similarly define $\CCdf$; so that  $\CCdf = \RRdf[\Omega]$.   We can also define
$\SG(\DD)_{df}$, and check that the groupification is $\RRdf$.

Given $a=(a_1,\ldots,a_n) \in \VF^n$ and $\g=(\g_1,\ldots,\g_n) \in \G^n$, let $B(a,\g) = \Pi_{i=1}^n B(a_i,\g_i)$, where $B(a_i,\g_i) = \{c \in \VF: \val(c-a_i) > \g_i \}$.  Call $B(a,\g)$
an open poly-disc of dimensions $\g$.  If $\g \in \G$, let
$B(a,\g) = B(a,(\g,\ldots,\g))$ (the open cube of side $\g$.)

Note that $[B(0,\g)]$ is invertible in $\RRdf$, in each dimension.  In particular
in dimension $1$, $[B(0,\g)][B(0,-\g)] = [0]_1^2$.   Note also: 
$[B(a,\g)] \oeq{a} [B(0,\g)]$. 

We proceed to define integrals of definable functions.

Let $U$ be a bounded definable subset of $\VF^n$.
A definable function $f: U \to \SG(\DD)_{df}$ has the form $[0]_1^{-m} F$, where $F: U \to \SG \DD[m]$
is a definable function, represented by some $\bar{F} \in \DD[m+n]_U$.     In case $\bar{F}$
can be taken bounded, define
$$\int_U f = [0]_1^{-m+n} [F]_{n+m} $$
We say that $f$ is {\em boundedly represented} in this case.

  In particular $\vol(U) = \int_U 1 = [0]_1^{-m} [U]_m$
is treated as a pure number now, without dimension units.
(Check independence of the choices.)

This extends by linearity to $\int_U f$ for $f: U \to \RRdf$, provided $f$ can be expressed
as the difference of two boundedly represented functions $U \to \SG(\DD)_{df}$.  

We now note that averaging twice, with appropriate weighting, is the same as doing it once.
The function $\g'$ in the lemmas below corresponds to a partition of $U$ into cubes;
$\g'(u)$ is the side of the cube around $u \in U$.  

\<{lem} \lbl{double-int}  Let $U$ be a bounded open subset of $\VF^n$, $f$ a boundedly represented function on $U$.   Let $\g': U \to \G$ be a definable
 function such that if $u \in U$ and $u' \in B(u,\g'(u))$ then $u' \in U$ and  $\g'(u') = \g'(u)$.
   Then
 $$\int_U f = \int_U [\vol(B(u,\g'(u)))^{-1} \int_{B(u,\g'(u))} f]$$
\>{lem}

\prf  Let $f = [0]_1^{-m} F$, where $F: U \to \SG \DD[m]$ is bounded.  We
have $\vol(B(u,\g')) = [0]_1 ^{-n} [\g'(u)]^{n}$ so 
$$\vol(B(u,\g'))^{-1}  =  [0]_1 ^{n} [\g'(u)]^{-n} = [0]_1^{-n} [-\g'(u)]^{n}$$
Thus, multiplying by $[0]_1^{3n+m}$,  we have to show
$$ [0]_1 ^{2n} [F]  = [-\g'(u)]^{n}]  [\{(u,u',z):  u  \in U, u' \in B(u,\g'(u)), (u',z) \in F \}]$$
Now $u' \in B(u,\g'(u))$ iff $u \in B(u',\g'(u'))$.  Applying the measure preserving
bijection $(u,u',z) \mapsto (u-u',u',z')$ we see that the  
$ [\{(u,u',z):  u  \in U, u' \in B(u,\g'(u)), (u',z) \in F \}] = [\g]_1^n [\{(u',z): (u',z) \in F ]$
so the equality is clear.
\eprf

We now define the integral of definable functions into $\CCdf$.  By definition, such
a function is a finite sum of products $fg$, with $f \in Fn(U,\RRdf)$ and $g \in Fn(U,\Om)$.
Define
$$ \int_U fg = \int_{\om \in \Om} \om \int_{g \inv (\om)} f $$
and extend by linearity.

Note that this is defined as soon as $g$ is boundedly represented.  (Again, check independence
of choices.)

\<{defn}  A definable distribution  on an open $U \subseteq \VF^n$ is a definable function 
$\dd: U \times \G \to \CCdf$, such that $\dd(a,\g) = \dd(a',\g)$ if $B(a,\g)=B(a',\g)$,   and whenever $\g'> \g$
in each coordinate,
$$\dd(b,\g) = \int _{u \in B(b,\g)} \vol(B(0,\g'))^{-1} \dd(u,\g')$$

\>{defn}

As in \lemref{reps}, the invariance condition means that  $\dd$ can be viewed as a function on open polydiscs,
and we will view it this way below.

If $\dd$ takes values in $\RRdf$, we say it is $\RRdf$-valued.
By definition, $\dd$ can be written as a finite sum $\sum \omega_i \dd_i$, 
where $\dd_i$ is an $\RRdf$-valued function;     in fact $\dd_i$ is 
an $\RRdf$-valued distribution.  

We wish to strengthen the definition of a distribution so as to apply to sub-polydiscs
of variable size.  For this we need a preliminary lemma.

 \lemm{dstrict}  Let $U=B(a,\g)$ be a polydisc.  Let $\g': B(a,\g) \to \G$ be a definable
 function such that $\g'(u') = \g'(u)$ for $u' \in B(u,\g'(u))$.  Then $\g'$ is bounded on $U$.
 \>{lem}
 
 \prf suppose for contradiction that
$\g'$ is not bounded   on $B(a,\g)$, i.e. 
$$(\forall \d \in \G)(\exists u \in B(a,\g))(\g'(u)>\d)$$
 This will not change if we add a generic element of $\G$ to the base,
so we may assume $\G(\dcl(\emptyset)) \neq (0)$.  By \lemref{resolve}, there exists a resolved structure
with the same $\RV$-part as $<\emptyset>$; hence we may
assume $\T$ is resolved.  By \secref{base-change}
any $\VF$-generated structure
is resolved.  By \lemref{curve-selection}, for any $M \models \T$ and $c \in \VF(M)$, $\acl(c)$ is an elementary submodel of $M$.   Consider $c$ with $\val(c) \models p_0$,
where $p_0$ is the generic type at $\infty$ of elements of $\G$,
i.e. $p_0 | A = \{x > \delta:  \delta \in \G(A) \}$.  Since 
$$\acl(c) \models (\forall \d \in \G)(\exists u \in B(a,\g))(\g'(u)>\d)$$
 there exists
$e \in \acl(c)$ with $e \in B(a,\g)$ and $\g'(e) > \val(c)$.  By \lemref{limitexists2},
there exists $e_0 \in \acl(\emptyset)$ such that $(c,e) \to (0,e_0)$.  In particular
$e_0 \in B(a,\g)$. But then since $e \to e_0$ and $\g'(e_0) \in \G(\acl(\emptyset))$
we have $e \in B(e_0,\g'(e_0))$.  So $\g'(e) = \g'(e_0)$.   But then $\g'(e_0) > \val(c)$,
contradicting the choice of $c$.     \eprf

\lemm{dist-strict} \begin{enumerate}
 
  \item Let  $\dd: U \times \G \to \CCdf$ be a definable distribution.  Let
  $\g': U \to \G$ be a definable function with $\g'(u)>\gamma$, such that
$\g'(u') = \g'(u)$ for $u' \in B(u,\g'(u))$.  Then 
\beq{strict}  \dd(b,\g) = \int _{u \in B(b,\g)} \vol(B(0,\g'(u))^{-1} \dd(u,\g'(u))  \eeq
  
  \item Let  $\dd_1,\dd_2$ be definable distributions on $U$ such that for any $x \in U$, for all large enough $\g \in \G$,
for any $y \in B(x,\g)$ and any $\g'>\g$, $\dd_1(B(y,\g')) = \dd_2(B(y,\g'))$. 
 Then $\dd_1=\dd_2$.
\end{enumerate}  \>{lem}

\prf (1)  To prove  \eqref{strict}, fix $b,\g$.   We may assume $U = B(b,\g)$.    Using 
\lemref{dstrict}, pick a constant $\g''$ with 
$\g'' > \g'(u)$ for all $u \in B(b,\g)$.  Use the definition of a distribution with respect to $\g''$ to
compute both $\dd(B(b,\g))$ and  for each $u$ $\dd(u,\g'(u))$, and compare the 
integrals using \lemref{double-int}.

 (3) Define $\g'(u)$ to be the smallest $\g'$ such that for all $\g'' > \g'$ and
all $y \in B(u,\g)$, $\dd_1(B(y,\g'')) = \dd_2(B(y,\g''))$.  It is clear that
$\g'(u') = \g'(u)$ for $u' \in B(u,\g'(u))$.     \eqref{strict} gives the same integral formula
for $\dd_1(b,\g)$ and $\dd_2(b,\g)$. 
 \eprf

Let $\dd$ be a definable distribution, and $U$ an arbitrary bounded open set.  We can define $\dd(U)$
as follows.  For any $x \in U$, let 
$\rho(x,U)$ be the smallest $\rho \in \G$ such that $B(x,\rho) \subseteq U$.
Let $B(x,U) = B(x,\rho(x,U))$; this is the largest open cube around $x$ contained in $U$.
Note that two such cubes $B(x,U),B(x',U)$ are disjoint or equal.  Define

$$\dd(U) = \int_{x \in U} \vol(B(x,U))^{-1} \dd(x,\rho(x,U))$$

More generally,  if $h$ is a locally constant function on $\VF^n$ into $\RRdf$ with bounded support,
we can define

\beq{test} \dd(h) = \int_{x \in \VF^n} h(x) [B(x,h)]^{-1} \dd(x,\rho(x,U))  \eeq
where now $B(x,h) = B(x, \rho(x,U))$ is the largest open cube around $x$ on which $h$ is constant.

%If $\dd$ is a distribution on $\VF^n$, 
%Let $U$ an open subset of $\VF^n$,   $g$ a definable
%function $U \to \CCdf$.  I   $g$ is integrable, 

\<{prop}\lbl{locallyL1}  Let $\dd$ be a definable distribution.  Then there exists  a definable 
open set $G \subseteq \VF^n$ whose complement $Z$ has dimension $<n$,  
and a definable function
$g: G \to \CCdf$ such that for any polydisc $U \subseteq G$
$$ \dd(U) = \int_U g$$
 \>{prop}

\prf  Since $\dd$ is a finite sum of $\RRdf$-valued distributions, we may assume
it is $\RRdf$-valued.  
Given $a \in \VF^n$, we have a function $\a_a: \G \to \RRdf$ defined by
$\a_a(\rho) = \dd(B(a,\rho))$.   Using the $\RV$-description of $\RR$, 
and the stable embeddedness of $\RV \union \G$, we see that $\a_a$ has a canonical
code $c(a) \in (\RV \union \G)^*$.   

Let $G$ be the union of all polydiscs $W$ such that $c$ is constant on $W$.  Let
$Z = \VF^n \m G$.  By \lemref{opensets}, $\dim(Z)<n$.

\Claim{} Let $W$ be a polydisc such that $c$ is constant on $W$.  Then for
some   $r \in \RRdf$, for any polydisc $U=B(a,\rho) \subseteq W$
$  \dd(a,\rho) = r \vol(U)$.

.

\prf  Since $c$ is constant on $W$, for some function $\d$,
all $\rho$ and all $b \in W$ with $B(w,\rho) \subseteq W$, we have
$\dd(B(w,\rho)) = \d(\rho)$.  By definition of a distribution we have, for any $a \in W$
$$\d(\rho) \vol(B(a,\rho'))  \oeq{a} \vol B(a,\rho)  \d(\rho')$$
Now $\vol (B(a,\rho)) \oeq{a}  \vol B(0,\rho) $.
So $\d(\rho) \vol(B(0,\rho')) \oeq{a} \vol B(0,\rho) \d(\rho')$.  Since this holds
for any $a \in W$, by \propref{resolve} we have
$$\d(\rho) \vol( B(0,\rho'))) =  \vol B(0,\rho) \d(\rho')$$
  So $\d(\rho)/\vol B(0,\rho) =r$
is constant. The Claim follows.\eprf
The Proposition follows too using \lemref{dist-strict}.
\eprf

\ssec{Fourier transform} Let $\psi$ be the tautological projection $K \to K/\Mm = \Om$.

Let $g: \VF^n \to \CCdf$ be a definable function, bounded on bounded subsets of $\VF^n$.   Define a function
$\FF(g)$ by
$$\FF(g)(U) = \int_{y \in \VF} g(y) (\int_{x \in U} \psi(x \cdot y))$$
This makes sense since for a given $U$, $(\int_{x \in U} \psi(x \cdot y))$ vanishes for $y$ outside a certain 
polydisc (with sides inverse to $U$.)  Moreover,

\lemm{fourier}  $\FF(g)$ is a definable distribution.  

\>{lem}
\prf  This follows from Fubini, \lemref{double-int}, and chasing the definitions.  \eprf

\<{cor}\lbl{lL1}  Fix integers $n,d$.  For all local fields
$L$ of sufficiently large residue characteristic, for any polynomial $G \in L[X_1,\ldots,X_n]$ of degree $\leq d$, there exists a proper variety $V_G$ of $L^n$ such that $\FF(|G|)$ agrees with a locally constant function outside
 $V_G$. \>{cor}

\prf  By \lemref{fourier} and \lemref{locallyL1}.  \eprf

See \cite{bernstein} for the real case.

\>{section}  %distributions + additivie

\<{section}{Expansions and rational points over Henselian fields} \lbl{exp-rat}

  We have worked everywhere with the geometry of  algebraically closed valued fields,
  or more generally of $\T$, but at a geometric level; all objects and morphisms can
  be lifted to the algebraic closure, and quantifiers are interpreted there.
  
 For many purposes, we believe this is the right framework.   It includes  for
 instance    Igusa integrals $\int_{x \in X(F)} |f(x)|^s$, and we will show in a sequel
 how to interpret in it some constructions of representation theory.   See also \cite{kazhdan-int}.  
 
In other situations, however, one wishes to integrate definable sets over Henselian
fields rather than
only constructible sets; and to have a change of variable formula for definable
maps, as obtained by Denef-Loeser and Cluckers-Loeser (cf.   \cite{cluckersloeser}).  It turns
out that our formalism lends itself immediately to this generalization; we explain in this section how to recover it.  The point is that an arbitrary definable set is an $\RV$-union of
constructible ones, and the integration theory commutes with $\RV$-unions.

   We will consider $F$ that admit quantifier elimination 
in a language $\L^+$ obtained from the language of $\T$ {\em by adding relations to $\RV$
only.}  For example, if $F= Th( \Cc((X)))$, $F$ has quantifier elimination in a language
expanded with names $D_n$ for subgroups of $\G$ (with $D_n(F) = n \G(F)$.)  
  
There are two steps in moving from $F^{alg}$ to $F$.  We will try to clarify the situation by
taking them one at a time.  The
two steps are to restrict the { points} to a smaller set (the $F$-rational points),
and the enlarge the { language} to a larger one (with enough relation symbols for $F$-quantifier elimination.)  We will take these steps in the reverse order.   In \secref{expansionsS}
we show how to extend the results of this paper to expansions of the language in the $\RV$
sorts, and in \secref{rational} how to pass to sets of rational points over a Hensel field.

 The reader who wishes
to restrict attention 
to constructible integrals (still taking rational points) may skip \secref{expansionsS}, taking $\T^+ = \T$ 
in \secref{rational}.  In this case one still has  a change of variable formula
for constructible change of variable, but not for definable change of variable.  An advantage
is that the 
the target ring correspondingly involves the Grothendieck group of constructible sets and 
maps rather than definable ones, which sometimes has more faithful information;
cf. \exref{hdlf}.   

\<{subsection}{Expansions of the $\RV$ sort}  \lbl{expansionsS}
 Let $\T$ be $\V$-minimal.
 
 Let $\T^+$ be an  expansion of $\T$ obtained by adding relations to $\RV$. 
 We assume that every $M \models \T$ embeds into the restriction to the language of $\T$ of
 some $N \models \T^+$ (as $\T$ is complete, this is actually automatic.)
  By adding some more basic relations, without changing the class of definable relations, we may assume
 $\T^+$ eliminates $\RV$-quantifiers.  As $\T$ eliminates field quantifiers, and $\T^+$
 has no new atomic formulas with $\VF$ variables, $\T^+$ eliminates $\VF$- quantifiers too,
 hence all quantifiers.  
 
For instance $\T+$ may include a name for  a subfield of the residue field 
(say, pseudo-finite) or the angular coefficients the  the Denef-Pas language (where $\RV$ is split).
 Write   $+$-definable for $\T^+$-definable; similarly $tp_+$ will denote the type
in $\T^+$< etc.  The unqualified 
    words formula, type, definable closure will  refer to 
quantifier-free formulas of $\T$.   

 \<{lem}  \lbl{comp+types}
 Let $M \models \T^+$.   Let $A$ be a substructure of $M$, $c \in M$, 
 $B= A(c) \meet \RV$.
 \begin{enumerate}
  \item    $tp(c/ A \union B) \union {\T^+}_{A \union B}$ implies ${tp_+}(c /  A \union B)$.
  \item  Assume $c$ is ${\T^+}_A$-definable.  Then $c \in \dcl(A,b)$ for some
  $b \in A(c) \meet \RV$% \meet \dcl_+(A)$.  
  \end{enumerate}

  \>{lem}
 
 \<{proof}  
 (1) This follows immediately from the quantifier elimination for ${\T^+}$.  Indeed
let $\phi(x) \in {tp_+}(c/A \union B)$.  Then $\phi$ is a Boolean combination of atomic
formulas, and it is sufficient to consider the case of $\phi$ atomic, or the negation of an atomic formula.
Now since any basic function $\VF^n \to \VF$ is already in the language
 of $\T$,
 every basic function of the language of ${\T^+}$ denoting a function $VF^n \to \RV$ factors through a 
 $\T$-definable function   into $\RV$.  Hence the same is true for all terms (compositions of basic functions).  
 And any basic relation is either the equality relation on $\VF$, or else a relation between
 variables of $\RV$.  If $\phi$ is an equality or inequality between $f(x),g(x)$, it is already
 in $tp(c/A)$.  Now suppose $\phi$ is a relation $R(f_1(x),\ldots,f_n(x))$ between elements of $\RV$.
 Since $B(c) \meet \RV \subseteq B$, the formula $f_i(x)=b_i$ lies in   $tp(c/A \union B)$  for some $b_i \in B$.
 On the other hand $R(b_1,\ldots,b_n)$  is part of ${\T^+}_{B}$.  
 These formulas together  imply $R(f_1(x),\ldots,f_n(x))$.
 
 (2)    We must show that $c \in \dcl(A \union B)$.
 Let $p=tp(c/A \union B)$.  By (1), $p$ generates a complete type of ${\T^+}_{A \union B}$.  Since this is the type of $c$ and $c$ is ${\T^+}_A$-definable,  and since any model of $\T$ embeds into a model of $\T^+$, 
 $p$ has a unique solution solution in any model of $\T$.  Thus $c \in \dcl(A \union B)$.   
 \>{proof}

We will now see that  any ${\T^+}$-definable bijection decomposes
into $\T$-bijections, and bijections of the form $x \mapsto (x,j(g(x)))$ where
$g$ is a $\T$-definable map into $\RV^m$ and $j$ is a ${\T^+}$-definable map on $\RV$.

\<{cor} \lbl{+formulas}
 \begin{enumerate}
  \item Let $P$ be a ${\T^+}$-definable set.  There exist   $\T$-definable  $f: \tP \to \RV^*$
 and a ${\T^+}$-definable $Q \subseteq \RV^*$
  such that $P = f \inv Q$.   
  
  \item Let $P_1,P_2$ be ${\T^+}$-definable sets,   and let $F: P_1 \to P_2$ be a 
  ${\T^+}$-definable bijection.   Then there exist $g_i: \tP_i \to R_i \subseteq \RV^m$, $R \subseteq \RV^m$,   $h_i: R \to R_i$,  and a bijection   
  $H: \tP_1 \times_{g_1,h_1} R \to \tP_2 \times_{g_2,h_2} R$ over $R$,  all $\T$-definable,
 and ${\T^+}$-definable  $Q_i \subseteq R_i$, $Q \subseteq R$, and $j_i: Q_i \to Q$
 such that $P_i = g_i \inv Q_i$, $h_ij_i = Id_{Q_i}$, and for $x \in P_1$
$$(\diamond) \ \ \ \ \  j_1 g_1 (x) = j_2 g_2 (F(x)) =: j(x), \text{ and }
  H(x,j(x)) = (F(x),j(x))$$
 
 Moreover if $P_i \subseteq \VF^n \times \RV^m$ projects finite-to-one to $\VF^n$, then  $R \to R_i$
 is finite-to-one.
  \end{enumerate}
 \>{cor} 
 
 \<{proof}  (1)  
 Let ${\mathcal F}$ be the family of all $\T$-definable functions 
 $f: W \to \RV^m$ 
 where $W$ is a definable set.
 
\Claim{} 
 If $tp(c)=tp(d)$ and  $f(c)=f(d)$ for all $f \in   {\mathcal F}$ with $c,d \in \dom(f)$, then $c \in P \iff d \in P$. 
 \proof  We have 
 $tp(c,f(c)) = tp(d,f(d)) = tp(d,f(c))$, so $tp(c/f(c))=tp(d/f(c))$ for all $f \in {\mathcal F}$ with $c \in \dom(f)$,  and
 thus $tp(c/  B) = tp(d/ B)$, where $B=A(c) \meet \RV$.  It follows that 
 ${tp_+}(c) = {tp_+}(d)$, and in particular $c \in P \iff d \in P$.
 
   By compactness,
 there are   $(f_i,W_i)_{i=1}^m  \in {\mathcal F}$  
 such that if $c \in W_i \iff d \in W_i$ and  $f_i(c)=f_i(d)$ whenever $c,d \in W_i$, 
 then  $c \in P \iff d \in P$.  Let $\tP = \union_i W_i$, and extend $f_i$ to $\tP$ 
 by $f_i(x) = \infty$ if $x \notin W_i$.  Let $f(x) = (f_1(x),\ldots,f_m(x))$.  
 % $\{Q_j\}$ be the equivalence classes of the equivalence
 %relation:  $x \in W_i \iff y \in W_i$; for fixed $j$, define $f$ on $W_j$ by
% $f(x)= (f_i(x))_{x \in W_i}$.  
 Letting $\tP = \union_i W_i$ and
 $ Q = f(P)$, (1) follows.  
 
 For (2), consider first a ${\T^+}$- type $p={tp_+}(c_1)$, $c_1 \in P_1$.  Let $c_2 = F(c_1)$.
 Using \lemref{fg}, there exists $g^p_i \in {\mathcal F}$ such that $e_i = g^p_i(c_i)$
 generates $\dcl(c_i) \meet \RV$.   It follows as in \lemref{comp+types} (1)  that $e_i$
 generates ${\dcl_+}(c_i) \meet \RV$.   Let $e$ generate $\dcl(c_1,c_2) \meet \RV$;
 we have $e_i = h^p_i (e)$ for appropriate $\T$-definable $h^p_i$.
 Note ${\dcl_+}(c_1)={\dcl_+}(c_2)$, and so $e \in {\dcl_+}(c_i)$.   Now quantifier elimination for ${\T^+}$
 implies the stable embeddedness of $\RV$, in the same way as for $\ACVF$
(cf. \secref{stab-emb}.)
By 
  \lemref{eilem-c} 
 ${tp_+}(c_i / e_i)$ implies ${tp_+}(c_i /  \RV)$; in particular since $e \in {\dcl_+}(c_i)$
  $e=j_i^p(e_i)$ 
 for some ${\T^+}$-definable $j_i^p$.
 By  \lemref{comp+types} (2) over $\dcl(c_1)$, $c_2 \in \dcl(c_1,e)$; similarly
 $c_1 \in \dcl(c_2,e)$; so there exists a $\T$-definable invertible $H^p$ with $H^p(c_1,e)=(c_2,e)$.   The   equations $(\diamond)$ have been shown
 to hold on $p$.   Now $g_i$ extends to a $\T$-definable function 
 $g_i: \tP_i \to R_i$.  
  By compactness  $(\diamond)$ holds on some definable neighborhood
 of $p$; and by (1) this neighborhood can be taken to have the form $g_1 \inv Q_1$
 for some $Q_1$.   Finitely many such neighborhoods cover $P_1$, and the data
 can be sewed together as in (1).  We thus find $\tP_1,R,R_1,R_2,g_1,g_2,h_1,h_2,H,Q_1,j_1,j_2$
 such that $h_ij_i(x)=x$ and  $(\diamond)$ holds on $g_1 \inv Q_1 = P_1$.  Let $Q_2 = h_2j_1Q_1$; it follows that
 $P_2 = F(P_1) = g_2 \inv Q_2$.  
 
To prove the last point, since $c_2 \in \dcl(c_1,e)$ we have (\lemref{fmrfr}) $c_2 \in \acl (c_1)$.
But $e \in \dcl(c_1,c_2)$ so $e \in \acl(\dcl(c_1))$; and as $e \in \RV^m$ for some $m$,
$e \in \acl(\dcl(e_1))$.      \>{proof}

  Let $\VF^+$ be the category of $+$- definable
 subsets of varieties over $\VF \meet \dcl(\emptyset)$, and $+$-definable maps.  
 Define effective isomorphism as in \defref{eff-iso}; let $\SG^{eff}$ denote the Grothendieck
 group of effective isomorphism classes, and let $[X]$ be the class of $X$.
 
 Let $\RV^+[*]$ be the category of
 pairs $(Y,f)$, where $Y$ is a 
 $+$-definable subset  of   $X$ for some    $(X,f) \in \Ob \RV[*]$ (\defref{RVcat}).  A morphism
 $(Y,f) \to (Y',f')$ is a definable bijection $h: Y  \to Y'$ such that $f'(h(y)) \in \acl(f(y))$ for $y \in Y$.

Let $\SG(\RV^+[*])$ be the Grothendieck semi-group of isomrphisms classes of $\RV^+$; let
$\Isp$ be the congruence generated by $(J,1_1)$ where 
$J = \{1\}_0 + [\RVp]_1 $.

\<{prop} \lbl{expansions}  There exists a  canonical surjective homomorphism of Grothendieck semigroups
$$\ints: \SG(\VF^+[*]) \to \SG(\RV^+[*])/ \Isp$$
%inducing a Grothendieck group homomorphism
%$$\int^{\K}: \K(\VF^+) \to \K(\RV^+)/ \Isp$$
determined by:
  $$\ints [X] = [W] / \Isp  \iff    [X] = [\L W] $$
 \>{prop} 

\prf   %To show that the formula on the second line defines a homomorphism, 
We have
to show:

 (i) any element of  $\SG(\VF^+)$ is effectively isomorphic to one of the form $[\L W]$;
 
(ii) if $[\L W_1] = [\L W_2]$  then $([W_1],[W_2]) \in \Isp$.

(i)  By \corref{+formulas} (1), a typical element of $\SG(\VF^+)$ is represented by 
$P = f \inv Q$, where $Q \subseteq \RV^*$is $\T^+$-definable, $f: \tP \to \RV^*$ is 
$\T$-definable.   For any $a \in \RV^*$, $f \inv (a)$ is $\T_a$-definable,
and $[f \inv (a)] = [\L C_a]$ where $[C_a] = [\int f \inv(a)]$.  Since $\L$ commutes
with $\RV$-disjoint unions, it follows that $[P] = [\L W]$ where $W = \du_{a \in Q} C_a$.

(ii)  Assume $[\L W_1] = [\L W_2]$.  By \propref{resolve}, the base can be enlarged so as to 
be made effective, without change to $\RV$; thus to show that $([W_1],[W_2]) \in \Isp$
we may assume $\L W_1,\L W_2$ are isomorphic.  Let $f: \L W_1 \to \L W_2$ be an isomorphism.  Let $P_i = \L W_i$ and let $\tP_i,R_i,g_i,h_i,R,H,Q,Q_i,j_i$ be as in \corref{+formulas} (2).

Since $P_i = g_i \inv Q_i = \L W_i$,  the maximal \rvinv subset of $\tP_i$  contains $P_i$, so 
we may assume $\tP_i$ is \rvinv; in other words   $\tP_i = \L \tW_i$ for some $\T$-definable $\tW_i \in \RV[*,\idot]$ containing $W_i$.    

By \lemref{rvfactor}, there exists a special bijection  $\si: \L  \tW_i^* \to \L \tW_i$
such that $g_i \circ \si$ factors through $\rho$, i.e. for some $e_i: {\tW}_i^* \to R_i$
we have $g_i \circ \si = e_i \circ \rho$
on $\L \tW_i$.  Let ${W_i}^*$ be the pullback of $W_i$ to ${\tW}_i^*$, so that
$\si(\L {W_i^*}) = \L W_i = P_i$.  
 Then $([W_i],[{W_i}^*]) \in \Isp$, so it suffices to show that
$({W}_1^*,{W}_2^*)  \in \Isp$.  Since $P_i = g_i \inv Q_i$, we have
 ${W}_i^* =  e_i \inv Q_i$.

For $c \in R$, let $\tP_i(c) = \si \inv g_i \inv (h_i(c))$, $\tW_i(c) = e_i \inv (h_i(c))$.  
Then $\tP_i(c) = \L \tW_i(c)$.  Now 
$H$ induces  a bijection $\tP_1(c) \to \tP_2(c)$.  Thus by \propref{change-of-variable-1}, 
$(\tW_1(c),\tW_2(c)) \in \Isp$.  In particular this is true for $c \in Q$; now $h_i:Q \to Q_i$
is a bijection, and $W_i^* = \du_{c \in Q} \tW_i(c)$.  Thus
$([W_1^*],[W_2^*]) \in \Isp$.  \eprf

 {\bf Remark} 
Since the structure on $\RV$ in $\T^+$ is arbitrary, we cannot expect
the homomorphism of \corref{expansions} to be injective.  We could make it so tautologically
by modifying the category $\RV^+$, taking only {\em liftable} morphisms, i.e. those that
lift to $\VF$; we then obtain an isomorphism.  In specific cases it may be possible to check that
all morphisms are liftable.

 \>{subsection}

\<{subsection}{Transitivity}
 \def\ac{{\rm ac}}
 \def\ACVFdp{{\ACVF^{\rm DP}}}
 
{\em Motivation.}  Consider a tower of valued fields, such as $\Cc \leq \CC((s)) \leq \CC((s))((t))$.  Given a definable
 set over $\CC((s))((t))$, we can integrate with respect to the $t$-valuation, obtaining data over 
 $\CC((s))$ and the value group.  The $\CC((s))$ can then be integrated with respect to the $s$-valuation.
 On the other hand, we can consider directly the $\Zz^2$-valued valuation of $\CC((s))((t))$, and integrate
 so as to obtain an answer involving the Grothendieck group of varieties over $\Cc$.    We develop below the language for comparing these answers, and show that they coincide.  

 For simplicity we accept here a
 Denef-Pas splitting, i.e. we expand $\RV$ so as to split the sequence $\k^* \to \RV^* \to \G$. 
 Then $\rv$ splits into two maps, $\ac: \VF^* \to \k^*$ and $\val: \VF^* \to \G$.  This expansion of $\ACVF(0,0)$
 is denoted $\ACVFdp$.  Note that this falls under the framework of \secref{expansionsS}, as will the furher
 expansions below.
 
 Consider two expansions of $\ACVFdp$:  (1) expand the residue field to 
 have the structure of a valued field (itself a model of $\ACVFdp$.)  (2)  expand the value group 
 to be a lexicographically ordered product of two ordered Abelian groups.  
 Then (1),(2) yield bi-interpretable
 theories.  In more detail:  
 
 First expansion:   Rename the $\VF$ sort as $\VF_{21}$,     the residue field as $\VF_{1}$, and the value group $\G_{1}$.  
      $\VF_{1}$ carries
 a field structure; expand it to a model of $\ACVFdp$, with residue field $F_0$ and value group $\G_0$.
 Let $\ac_{21},\val_{21}$ have their natural meanings.
  
Second expansion:    Rename the $\VF$-sort as $\VF_{20}$,    the residue field as $F_0$ and the value
group as $\G_{20}$.    Add
 a  predicate  $\G_{0}$ for a proper convex  subgroup of $\G_{20}$, and a predicate $\G_1$ for  a complementary subgroup, so that $\G_{20}$ is identified with the lexicographically ordered $\G_{0} \times \G_{1}$. 
 
 \lemm{same}  The two theories described above are bi-interpretable.  A model of (1) can canonically be made into a model of (2) with the same class of definable relations; and vice versa.  \>{lem}
 
 \prf  Given (1), let $\VF_{20}=\VF_{21}$ as fields.  Define
\beq{same0} \ac_{20} = \ac_{10} \circ \ac_{21} \eeq
  Let 
 $\G_{20} = \G_1 \times \G_0$, and define $\val_{20}: \VF_{21}^* \to \G_{20}$ by
\beq{same1} \val_{20} (x) = (\val_{21}(x), \val_{10}(\ac_{21}(x)))   \eeq
  
 Conversely given (2), let $\VF_{21}=\VF_{20}$ as fields; 
$$\Oo_{21} = \{x \in \VF_{21}: (\exists t \in \G_0) (\val_{20}(x) \geq t) \}$$ 
$$\Mm_{21} =  \{x \in \VF_{21}: (\forall t \in \G_0) (\val_{20}(x) > t) \}$$ 
 $$\VF_{1} = \Oo_{21}/\Mm_{21}$$
 Let $\VF_{21}$ have the valued field structure with residue field $\VF_1$; note that
 the value group $\VF_{21}^*/ \Oo_{21}^*$ can be identified with $\G_1$.   Note that $\ker \ac_{20} \supset 1+\Mm_{21}$, so that
 factors through $\VF_1^*$, and define $\ac_{10},\ac_{21}$ so as to make \eqref{same0} hold.
 Then define $\val_{21},\val_{10}$ so that \eqref{same1} holds.
%Observe 
%and let $\Oo_{10}$ be the image of $\Oo_{20}$ in  $\VF_1$
  \eprf 
 
 %Note that $\VF_{21}, \VF_{20}$ have different valued field structures.  

Let $\VF^+[*]$ denote the category of definable subsets of $\VF_{21}$, equivalently $\VF_{20}$,
in the expansions (1) or (2). According to \propref{expansions}, and \lemref{tensor0},
we have canonical
maps $\SG(\VF^+[*]) \to \SG(\RV^+_1[*] ) /\Isp$ and $\SG(\VF^+[*]) \to \SG(\RV^+_2[*] ) /\Isp$,
where $\RV^+_i[*]$ denotes the expansion of $\RV$ according to (1),(2) respectively.  

By \propref{summ} we have canonical maps:

\beq{transi1} \SG(\VF^+[*]) \to \SG(\VF_1[*]) \tensor \SG(\G_{21}[*]) / \Isp \to (\SG(F_0) \tensor \SG(\G_{10})) \tensor \SG(\G_{21}) /\Isp_{1}  \eeq

for a certain congruence $\Isp_1$.  And on the other hand:  
 
\beq{transi2} \SG(\VF^+[*]) \to \SG(F_0[*]) \tensor \SG(\G_{20}[*]) / \Isp = \SG(F_0[*]) \tensor  (\SG(\G_{10}[*]) \tensor \SG(\G_{21}[*])) /\Isp_2  \eeq
 
 For an appropriate $\Isp_2$.   The tensor products here are over $\Zz$, in each dimension separately.

 Using transitivity of the tensor product we identify $(\SG(F_0) \tensor \SG(\G_{10})) \tensor \SG(\G_{21}) $
 with  $\SG(F_0[*]) \tensor ( \SG(\G_{10}[*]) \tensor \SG(\G_{21}[*]))$.  Then

\<{thm}\lbl{transi}  $\Isp_1,\Isp_2$ are equal and the maps of \eqref{transi1},\eqref{transi2} coincide.
\>{thm}

 \prf  It suffices to show in the opposite direction that the compositions of maps induced by $\L$

 \beq{transi3} (\SG(F_0[*]) \tensor \SG(\G_{10}[*]) \tensor \SG(\G_{21}[*])  \to \SG(\VF_1[*]) \tensor \SG(\G_{21}[*]) \to 
 \SG(\VF^+[*])  \eeq

 \beq{transi4} \SG(F_0[*]) \tensor  \SG(\G_{10}[*]) \tensor \SG(\G_{21}[*])  \to \SG(F_0[*]) \tensor \SG(\G_{20}[*])   
 \to  \SG(\VF^+[*]) \eeq
 
coincide.  
 But this reduces by $\RV$-additivity to the case of points, and by multplicativity to the individual 
 factors $F_0,\G_{21}, \G_{10}$, yielding to an obvious computation in each case.

 \eprf
 
 %:measure statement
 
 \>{subsection}

 \<{subsection}{Rational points over a Henselian subfield:  constructible sets and morphisms} 
  \lbl{rational}
 
 Let $\T$ be $\V$-minimal, and $\T^+$ an expansion of $\T$ in the $\RV$ sorts.
 
 Let $F$ be an effective substructure of a model of $\T$.  Thus
$F = (F_\VF,F_\RV)$, with $F_\VF$ a field, and $\rv(F_\VF) = F_\RV$; and ,
 $F$ closed under definable functions of $\T$.   
 For example, if 
  $\T = \T^+ = ACVF(0,0)$, this is the case   iff $F_\VF$   is a Henselian field
   and $F_\RV = F / \Mm(F)$;   any Hensel field of
   residue characteristic $0$ can be viewed in this way.  See \exref{hensel-aux}.
     
  By a {\em $+$- constructible} subset of $F^n$, we mean a set of the form $X(F) = X \meet F^n$,
  with $X$ a quantifier-free formula of $\T^+$.   Let  $\VF^+(F)$ be the 
  category  of such sets, and
  $+$-constructible functions between them.  The Grothendieck semiring $\SG \VF^+ (F)$ is thus the quotient of $\SG \VF$
  by the semiring congruence
  $$I_F = \{ ([X],[Y]):  X,Y \in \Ob \VF^+, X(F)=Y(F) \}$$
  (one can verify this is an ideal; in fact if $X(F) = Y(F)$ and $X \iso X'$,
  then there exists $Y' \iso Y$ with $X'(F) = Y'(F)$.)

   Similarly we can define $I_F^{RV}$ and
  form
  $\SG \RV (F) \iso (\SG \RV)/I_F^{RV}$.   
   As usual, let $\Isp$ denote the 
  congruence generated by $([1]_0+  [\RVp]_1, [1]_1)$, and $I_F^{RV}+\Isp$
  their sum.

 \Claim{} If $([X],[X']) \in I_F$ then $(\ints [X],\ints [X']) \in I_F^{RV} + \Isp$.  
 
 \prf
We may assume,   changing $X$ within the $\VF$-isomorphism class $[X]$, that 
$X(F)=X'(F)$.  Then $X(F) = (X \union X')(F) = X'(F)$, and it suffices to show:
$(\ints [X],\ints [X \union X']) , (\ints [X'], \ints [X \union X'])  \in I_F^{RV}$.  Thus we may
assume $X \subseteq X'$.  Let $Z=  X' \m X$.  Then $Z(F) = \emptyset$, and it suffices
to show:  $(\ints (Z), \emptyset) \in I_F^{RV}$.  Now $\ints(Z) = [Y]$ for some $Y$
with $Z$ definably isomorphic to $\L Y$.  So $\L Y (F) = \emptyset$, hence $Y(F) = \emptyset$.
Thus $([Y], \emptyset) \in I_F^{RV}$, as required.  \>{proof} 

As an immediate consequence we have:

\<{prop} \lbl{rationalpoints}  Assume $ F \leq M   \models \T$, with $F$ closed under definable
functions of $\T$.  
 The homomorphism $\ints$ of \thmref{summ-inf} induces a homomorphism
 $$\int_F: \SG \VF^+ (F) \to \SG \RV^+ (F) / \Isp$$   \qed

\>{prop}
 \>{subsection}

\<{subsection}{Quantifier elimination for Hensel fields}   
  
Let $\T$ be a $\V$-minimal theory in a language $L_\T$, with sorts $(\VF,\RV)$
(cf. \secref{language}.)
  Assume $\T$ admits
quantifier-elimination, and moreover that any definable function is given by a basic
function symbol.  This can be achieved by an   expansion-by-definition of the language.

Let $\T_h = (\T)_{\forall} \union \{(\forall y \in \RV)(\exists x \in \VF)(\rv(x)=y)$.
% (\exists x \in \G)(x>0)\}$.  

A model of $\T_h$ is thus the same as a substructure $A$ of a model of $\T$, 
such that $\RV(A) = \rv (\VF(A))$.  % and $\G(A) \neq (0)$.  

\lemm{lang}  Any formula of $L_\T$ is $\T$-equivalent to a Boolean combination of
formulas in $\VF$-variables alone, and formula $\psi(t(x),u)$ where $t$ is a sequence
of terms for functions $\VF^n \to \RV$,   $u$ is a sequence of $\RV$-variables, and $\psi$
is a formula of $\RV$ variables only.  \>{lem}

\prf  This follows from stable embeddedness of $\RV$,   \corref{rvei},  \lemref{eilem}
and the fact  (\lemref{rvfactorcor}) that definable functions into $\G$
factor through definable functions into $\RV$.  \eprf 

\<{example} \lbl{hensel-aux}  If $\T=\ACVF(0,0)$, then $\T_h$ is an expansion-by-definition of 
the  theory of  
Hensel fields of residue characteristic zero.  \>{example}
 
\prf   We must show that a Henselian valued field is definably closed in its algebraic closure,
in the two sorts $\VF,\RV$.

 Let $K  \models T_{Hensel}$, $K \leq M  \models \ACVF$.   Let $X \subseteq \VF^k \times \RV^l$,
$Y \subseteq \VF^{k'} \times \RV^{l'}$ be $\ACVF_K$-definable sets, and
 $F:X \to Y $   an $\ACVF_K$-definable bijection.  We have to show that
  $F(X \meet K^k \times \RV(K)^l) = Y \meet K^{k'} \times \RV(K)^{l'}$.    

   $K^{alg}$ is an elementary submodel of $M$; we may assume $K^{alg}=M$.    
By one
of the characterizations of Henselianity, the valuation on $K$ extends uniquely to $K^{alg}$.  Hence
every   field automorphism of $M$ over $K$ is a valued field automorphism.  So $K$ is the fixed field
of $Aut(M/K)$ (in the sense of valued fields), and hence $K=\dcl(K)$.  Since $\ACVF_K$ is effective,
any definable point of $\RV$ lifts to a definable point of $\VF$; so $ \dcl(K) \meet \RV = \RV_K$.  Thus
$K$ is definably closed in $M$ in both sorts. 
 \eprf

% Let $L_{VF}$ be the language of valued fields in the sorts $(\VF,\RV)$, \secref{language}.
 Let $L  \supset L_{\T}$;   assume $L  \m L_{\T}$ consists of relations and functions on $\RV$ only.
% Let $T_{Hensel}$ be the $L$- theory whose models are the pairs $(K, \RV_K)$ with $K$ a Hensel field  of residue characteristic $0$, and $\RV_K = K^* / \Mm_K$.  
If $A \leq M \models \T$, let $L_{\T}(A)$ be the languages enriched with constants
for each element of $A$; let $\T_{h}(A) = \T_A \union \T_{h}$, where $\T_A$ is the set of 
quantifier-free valued field formulas true of $A$.

 \<{prop} \lbl{henselQE}   $\T_{h}$ admits elimination of field quantifiers. 
 \>{prop}
 
 \prf   Let $A$ be as above.  Let $\Phi_A$ be the set of $L(A)$-   formulas with no $\VF$-quantifiers.
  
 \Claim{}  Let $\phi(x,y) \in \Phi_A$ with $x$ a free $\VF$-variable.    Then 
 $(\exists x)\phi(x,y)$ is $\T_{h}(A)$-equivalent to  a formula in $\Phi_A$.

\prf By the usual methods of compactness and absorbing the $y$-variables into $A$, it suffices to prove
this when $x$ is the only variable.  Assume first that   $\phi(x)$ is an $L_{\T}(A)$-formula.
 By
\lemref{tr1}, there exists an ACVF-definable bijection between the definable set defined by $\phi(x)$,
and a definable set of the form $\L \phi'(x',u)$, where $\phi'$ is an $L_{\T}(A)$-formula in $\RV$-variables
only (including a distinguished variable $x'$ on which $\L$ acts.) By definition of $\T_h$,  in any model of
$\T_{h}$, $\phi$ has a solution iff $\L \phi'(x',u)$ has a solution.  But clearly $\L \phi'(x',u)$
has a solution iff $\phi'(x',u)$ does.  Thus $\T_{h}(A) \models (\exists x) \phi(x) \iff (\exists x',u)\phi'(x',u)$. 
 
Now let $\phi(x)$ be an arbitrary  $\Phi_A$ formula.
Let $\Psi$ be the set of formulas of $L(A)$ involving $\RV$-variables only.  
Let $\Theta$ be the set of conjunctions of formulas of $L_{\T}(A)$ in $\VF$-variables only, and of formulas
of the form $\psi(t(x))$, where $\psi \in \Psi$ and $t$ is a term of $L_{\T}(A)$.The set of disjunctions of formulas in $\Theta$ is then closed under Boolean combinations, and under
existential $\RV$-quantification.     By \lemref{lang} it includes all
  $L_\T$-formulas, up to equivalence; and also all formulas in $\RV$-variables only.
 Thus $\phi(x)$ is a disjunction of formulas  in $ \Theta$, and we may assume
$\phi(x) \in \Theta$.  Say $\phi= \phi_0(x) \wedge \psi(t(x))$, with $\phi_0 \in L_{\T}(A)$ and $\psi \in \Psi$.
By Claim 1, for some   formula $\rho(y)$ of $\Phi_A$, we have $T_{h}(A) \models \rho(y) \iff (\exists x)(t(x)=y \& \phi_0(x))$.  Hence 
$(\exists x) \phi(x)  \iff (\exists y)(\psi(y) \& \rho(y))$.   \eprf

Quantifier elimination now follows by induction.  \eprf

{\bf Remark}  Since only field quantifiers are mentioned, this immediately extends
to expansions in the field sort.

In particular, one can split the sequence $0 \to \k^* \to \RV \to \G \to 0$ if one wishes.  
This yields the  quantifier-elimination \cite{pas} in the Denef-Pas
language.  

The results of Ax-Kochen and Ershov, and the large literature that developed around them,
appeared to require methods of ``quasi-convergent sequences''.  It is thus curious
that they can also be obtained directly from Robinson's earlier and purely ``algebraic''
quantifier elimination for ACVF.   Note that in the case of ACVF, there is no need
to expand the language to obtain QE; and then \lemref{lang} requires no proof beyond inspection of the language.

  \>{subsection} %Quantifier elimination for Hensel fields
\<{subsection}{Rational points:  definable sets and morphisms}

In this subsection  we will work   
with   completions $T$ of $\T_h \union  \{ (\exists x \in \G)(x>0)\}$. .  These are theories of valued fields of residue characteristic $0$, possibly expanded, not necessarily algebraically closed.   The language of $T$ is thus the language of $\T^+$.  
The words formula, type, definable closure will  refer to 
quantifier-free formulas of $\T^+$.   Definable closure, types
with respect to $T$ are referred to explicitly as $\Tdcl$, $\Ttp$ etc.   

Let $F \models T$.  Since $F \models \T_{\forall}$, $F$ embeds into a model
$M'$ of $\T^+$.  Since $\G(F) \neq (0)$, by
 \propref{resolve} and \lemref{curve-selection}, there exists $F' \subseteq M'$ containing $F$, 
with $\G(F')=\G(F)$, and $M=\acl(F')$ an elementary submodel of $M'$.  Hence
$F$ embeds into a model $M$ of $\T^+$ with  $\G(F)$ cofinal in $\G(M)$.

 \<{lem}  \lbl{comp-types}
 Let $F \models T$,
 $F \leq M \models \T^+$, $\G(F)$ cofinal in $\G(M)$.
 Let $A$ be a substructure of $M$, $c \in F$, $B= A(c) \meet \RV \meet F$, 
 \begin{enumerate}
  \item    $tp(c/ B) \union T_{B}$ implies $\Ttp(c /  B)$.
  \item  Assume $c$ is $T_A$-definable.  Then $c \in \dcl(A,b)$ for some
  $b \in B$.  
  \end{enumerate}

  \>{lem}
 
 \<{proof}  
 (1) This follows immediately from the quantifier elimination for $T$ and \lemref{comp+types} (1).   
 
 (2)  
We have $B  \subseteq \Tdcl(A) \meet \RV$.
 We must show that $c \in \dcl(A \union B)$.
 Let $p=tp(c/A \union B)$.  By (1), $p$ generates a complete type of $T_{A \union B}$.  Since this is the type of $c$ and $c$ is $T_A$-definable,  some formula $P$ in the language of $T_{A \union B}$ 
 with $P \in p$ has a unique solution in $F$.  Now the values of $F$ are cofinal in the value group of
 $F^a$; so $P$ cannot contain any ball around $c$  (any such ball would have an additional point of $F$, obtained by adding to $c$ some element of
 large valuation.)  Let $P'$ be the set of isolated elements of $P$; then $P'$ 
 is finite (as is the case for every definable $P$), $\T_A$-definable, and $c \in P'$.  
 By \lemref{finite}, there exists an  $\T_A$-definable bijection $f: P' \to Q$ with 
 $Q \subseteq \RV^n$.  Then $f(c) \in \dcl_{T}(A)=B$, and $c = f \inv (f(c)) \in \dcl(A \union B)$. \>{proof}

\<{cor} \lbl{hensel} Two definably isomorphic
definable subsets of $F$ have the same class in $\SG \VF^+(F)$.  \>{cor}

\prf 
$T$-definable bijections are restrictions of ${\T^+}$-definable bijections.  Hence 
\corref{+formulas} is true with $T$ replacing $\T^+$.     \eprf

 Thus   \propref{rationalpoints}
includes a change-of-variable formalism for definable bijections.
 
 \ssec{Some specializations}

\sssec{Tim Mellor's Euler characteristic} \lbl{tim}  Consider the theory $\RCVF$ of real closed valued fields.
Let $\RV_{\RCVF}$, $\RES_{\RCVF}$, $\VAL_{\RCVF}$    denote the categories of definable sets
and maps that lift to bijections of $\RCVF$ (on $\RV$ and on   the residue field, value group
respectively.  We do not need to use the sorts of $\RES$ other than the residue 
field here, say all structures $A$ of interest have $\G_A$ divisible.)    From \propref{rationalpoints} and \corref{hensel}
we obtain an isomorphism:
$\K(\RCVF) \to \K(\RV_{\RCVF}) / ([0]_1 - [\RV^{>0}]_1 - [0]_0) $.

The residue field is a model of the theory $\RCF$ of real closed fields;  
$\K(\RCF)= \Zz$ via the Euler characteristic (cf. \cite{vddries}).  Since the ambient dimension grading is respected
here,   $\K(\RES_{\RCVF}) = \Zz[t]$.

   The value 
group is a model of $\DOAG$, and moreover, any definable bijection on $\G[n]$ for
fixed $n$ lifts to 
$\RV$ and indeed to $\RCVF$.  This is because the multiplicative group of positive
elements is uniquely divisible, and so $SL_n(\Qq)$ acts  the $n$'th power of this group.
By \propref{2eulers}, $\K(\DOAG)[n] = \Zz^2$ for each $n \geq 1$,
and $\K(\VAL_{\RCVF}) =  \Zz[s]^{(2)}: = \{(f,g) \in \Zz[s]: f(0)=g(0) \}$.  

Thus $\K(\RV_{\RCVF} = \Zz[t] \tensor \Zz[s]^{(2)} \leq \Zz[t,s]^2$; and $J$ is identified
with the class $(1,1) - (0,-s) - (t,t)$.  Thus we obtain two homomorphisms 
$\K(\RV_{\RCVF})/J  \to \Zz[s]$ (one mapping $t \mapsto 1$, the other with $t \mapsto 1-s$;
and as a pair they are injective.  

Equivalently, we have found two ring homomorphisms $\chi,\chi': \K(\RCVF) \to \Zz[t]$.  One of these
 was found in \cite{mellor}.

\sssec{Cluckers-Haskell}  \lbl{cl-ha} Take the theory of the $p$-adics.  By 
   \propref{rationalpoints} and \corref{hensel}
we obtain an isomorphism:
$\K(\pCF) \to \K(\RV_{\pCF}) / \Isp $.  However $\RV_{\pCF}$ is a finite extension of
$\Zz$, and evidently $\K(\Zz)=0$, since $[[0,\infty)]=[[1,\infty)]$.  Thus $\K(\pCF)= 0$.

 \ssec{Higher dimensional local fields}  

We have seen that the Grothendieck group of definable sets with volume forms loses a great deal
of information compared to the semi-group.  Over fields with discrete value groups, restricting 
to bounded sets is helpful; in this the Grothendieck group retains information about volumes.  In 
case of higher dimensional local fields, with  value group is $A=\Zz^n$, 
simple boundedness is insufficient to save from collapse.  We show that using a simple-minded 
notion of boundedness is only partly helpful, and loses much of the volume information
(all but one $\Zz$ factor.)  

 \def\Km{{\K_\mu}}
\<{example} \lbl{hdlf}  Let $\Km \bdd (Th(\Cc((s_1))((s_2)))[n])$ be the 
Grothendieck ring of definable bounded sets and measure preserving maps in
$\Cc((s_1))((s_2))$ (with $\val(s_1) << \val(s_2))$.   Let $Q^t$   denote the class of the 
thin annulus  of radius
$t$.   
 In particular
$Q^0$ is the volume of the units of the valuation ring.   
Then in  $\Km \bdd(Th(\Cc((s_1))((s_2)))[2] )$ we have for example
$(Q^0)^2 = 0$.   To see this directly let
 $$Y = \{(x,y): \val(x)=0, \val(y)=0\}, \ X = \{(x,y): 0< 2 \val(x) < \val(s_2), \val(x)+\val(y)= 0 \}$$
 Then $X$ is bounded.  Let $f(x,y) = (x/s_1, s_1y)$.  Then $f$ is a measure preserving
 bijection $X \to X' =  \{(x,y): 0< 2 (\val(x)+\val(s_1)) < \val(s_2), \val(x)+\val(y)= 0 \}$
 But in $\Cc((s_1))((s_2))$, $2 \val(x) < \val(s_2)$ iff $2 (\val(x)+\val(s_1)) < \val(s_2)$
 so $X'(\Cc((s_1))((s_2))) = X(\Cc((s_1))((s_2)) \union Y(\Cc((s_1))((s_2)))$.
 \>{example}
 
\remm{11}   
%$$\phi(x_1,x_2,x_3)=(1-[{{x_1-x_2}\over{2}}])(1-[{{x_2-x_3}\over{2}}])(1-[{{x_1-x_3}\over{2}}])$$
%or similarly 
$(2[[0,y/2]]-[[0,y]])(2[[0,y/2)] - [[0,y)])$,
is a class of the Grothendieck group of $\G$
%(depending on $x=(x_1,x_2,x_3) \in \G^3$)
  that vanishes
identically in the $\Zz$-evaluation, but not in the
$\Zz^2$-evaluation. 
 
% In the non-archimedean case, there is always a nontrivial homomorphism to an 
% archimedean group $B$, yielding a homomorphism $\K(\G[n])$ modulo the $<n$-dimensional ideal, to $\K(\G_B[n])$

 \>{subsection}  %valued fields

\>{section}

  \<{section}{The Grothendieck group of algebraic varieties}

Let $X,Y$ be smooth nonsingular curves in $\Pp^3$, or in some other smooth projective
variety $Z$, and assume $Z \m X, Z \m Y$ are biregularly isomorphic.  Say $X,Y,Z$
are defined over $\Qq$.  Then 
 for almost all $p$, $|X(\Ff_p)| = |Y(\Ff_p)|$, as one may see 
  by counting points of $Z$, $Z \m X$ and subtracting.  It follows
from Weil's Riemann hypothesis for curves that $X,Y$ have the same genus, from
Faltings that $X,Y$ are isomorphic if the genus is 2 or more, and from Tate that $X,Y$
are isogenous if the genus is one.  It was this observation that led Kontsevich and Gromov to ask if $X,Y$ must actually be isomorphic.    We show that this is the case below.  
\footnote{This already follows from \cite{larsen-luntz}, who use different methnods.}

 \<{thm} \lbl{complement}  Let $X,Y$ be two smooth $d$-dimensional subvarieties 
  of a smooth projective $n$-dimensional variety $V$, and assume $V \m X, V \m Y$
 are biregularly isomorphic.   Then 
 $X,Y$ are stably birational, i.e. $X \times \Aa^{n-d},Y \times \Aa^{n-d}$
 are birationally equivalent.   If $X,Y$ contain no rational curves, then      $X,Y$ are  birationally equivalent.   \>{thm}

While we do not obtain a complete characterization in dimensions $>1$, the results and method 
of proof do
 show that the answer lies in synthetic geometry and is not cohomological in nature.

 Let $\VarK$ be the category of algebraic varieties over a field ${\Kz}$ of characteristic $0$.

  Let $[X]$ denote the class of a variety  $X$ in the Grothendieck semigroup $\SG(\VarK)$.
 We allow varieties to be disconnected.  As all varieties will be over  the same field $\Kz$,
 we will write $\Var$ for $\VarK$.  Let $\SG \Var_n$ be the Grothendieck semigroup 
 of varieties of dimension $\leq n$.

For the proof, we
 view ${\Kz}$ as a trivially valued subfield of a model of  $ACVF(0,0)$.  
 We work with the theory $\ACVF_{\Kz}$, so that ``definable'' means $\Kz$-definable
with quantifier-free $\ACVF$-formulas.

Note that $\RES = \k^*$ in $\ACVF_{\Kz}$; the only definable point of $\G$ is $0$,
so the only definable coset of $\k^*$ is $\k^*$ itself.  

 The residue
map is   an isomorphism on ${\Kz}$ onto a subfield ${\Kz}_{\RES}$ of the residue field  $\k$.
In particular, any smooth variety $V$ over ${\Kz}$ lifts canonically to a smooth scheme $V_\Oo = V  \tensor_{\Kz} \Oo$ over $\Oo$,
with generic fiber $V_\VF = V_\Oo \tensor _{\Oo} \VF$ and
special fiber $V_\Oo \tensor_{\Oo} \k = V \tensor_{\Kz} \k$.  
We have a reduction homomorphism $\rho_V: V(\Oo) \to V(\k)$.  We will write 
$V(\Oo) ,V(\VF)$ for $V_\Oo(\Oo), V_\VF (\VF)$.  

Given $k \leq n$ and a definable subset $X$ of $\RV^*$ of dimension $\leq k$,   let
$[X]_k$ be the class of $X$ in $\SG \RV [k] \subseteq  \SG \RV[\leq n]$.
Thus if $\dim(X)=d$ we have $n-d+1$ classes $[X]_k$, $d \leq k \leq n$, in 
different direct  factors of $\SG \RV[\leq n]$.  We also use $[X]_k$ to denote
the image of this class in $\SG \RV[\leq n] / \Isp$.   This abuse of notation is not excessive since
for $n \leq N$,  $\SG \RV[\leq n] /\Isp$
embeds in $\SG \RV [\leq N] / \Isp$ (\lemref{summ-c}).

Let $SD_d$ be the image of $\SG \RV [ \leq d]$ in $\SG \RV [\leq N] / \Isp$.  Let $WD^n_d$
be the subsemigroup of $\RV[n]$ generated by $\{[X]: \dim(X) \leq d \}$,
and use the same letter to denote the image in $\RV[\leq N] / \Isp$.  
Let $FD^n = SD_{n-1}+ WD^n_{n-1}$.    
We write $a \sim b \, (FD^n_d)$ for:
$(\exists u,v \in FD^n_d)(a+u=b+v)$. More generally, for any subsemigroup $S'$ of 
a semigroup   $S$, write $a \sim b \, (S')$ for 
$(\exists u,v \in S)(a+u=b+v)$. 

We write $\K( \RV[\leq n]) / \Isp $ for the groupification of $\SG ( \RV[\leq n] )/ \Isp $

 \lemm{comp2}  Let $V$ be a smooth projective $\k$- variety of dimension $n$, $X$ a definable subset of $V(\k)$.   
 Then  $$\ints [\rho_V \inv (X)] =[X]_n$$
\>{lem} 

\prf  Let $\bX = (X,f)$ where $f: X \to \RV^n$ is a finite-to-one map.
We have to show that $[\L \bX] = [\rho_V \inv (X)]$ in $\SG(\VF[n])$, i.e. that $\L \bX, \rho_V \inv(X)$ are definably isomorphic.  By \lemref{collect} this reduces to the case that $X$ is a point $p$.  
Find an open affine neighborhood $U$ of $V$ such that $\rho_V \inv (p) \subseteq U(\Oo)$,
and $U$  admits an \'etale map $g: V \to \Aa^n$ over $\k$.
Now $U(\Oo) \iso \Oo^n \times_{\res,g} U(\k)$.  This reduces the lemma to the case of affine space, where
it follows from the definition of $\L$.  \eprf
 
 \lemm{varclass}   Let $X$ be a a $\Kz$-variety of dimension $\leq d$.
  \begin{enumerate}
  \item  $\ints   ( X(\VF))  \in SD_d = \SG( \RV[\leq d]) / \Isp$
    \item If $X$ is a smooth complete variety   of dimension $d$  then $\ints X(\VF) = [X]_d$.
 \item If $X$ is a   variety   of dimension $d$  then 
  $\ints X(\VF)  \sim [X]_d \  (FD^d)$.
% $\ints (X(\VF)) = [X]_d + [Y]_d + [Z] \in  \SG \RV[\leq n] / \Isp$ where $\dim(Y) < d$ and $Z \in \SG \RV[<d] / \Isp$
 % \item If $X$ is a projective variety   of dimension $d$  then .... 

\end{enumerate} 
  \>{lem}

\prf 

(1) is obvious, since $\dim(X(\VF)) \leq d$.  

(2)  By Grothendieck's valuative criterion for properness,  $X(\VF) = X(\Oo)$.  We thus have
a map  $\rho_V: X(\VF)= X(\Oo) \to X(\k)$.  For $\a \in X(\k)$, let $X_\a(\VF) = \rho_V \inv (\a)$.  

Since $X$ is smooth of dimension $d$ it is covered by Zariski open
neighborhoods $U$ admitting  an \'etale map $f_U: U \to \Aa^d$, defined over $\Kz$;
let ${\mathcal S}$ be a finite  family of such pairs $(U,f_U)$, with 
$\union_{(U,f_U) \in {\mathcal S} }U  = X$.  
We may choose a definable finite-to-one $f: X \to \Aa^d$, defined over $\Kz$,
such that for any $x \in X$, for some pair $(U,f_U) \in {\mathcal S}$, $f(x)=f_U(x)$.  
We have $\L([X]_d) = \L(X,f) = \VF^d \times_{\rv,f} X(\k)$.   We have to show that $\L(X,f)$ is definably isomorphic to $X(\VF)$.
By \lemref{collect} it suffices to show that for each $\a \in X(\k)$, 
 $\VF^d \times_{\rv,f} \{a\} $
is $\a$-definably isomorphic to $X_\a(\VF)$.   Now   $\VF^d \times_{\rv,f} \{a\} \oeq{\a} \rv \inv (f(\a))$.   We have $f(\a)=f_U(\a)$ for some $(U,f) \in {\mathcal S}$ with $\a \in U$.  
Since $f_U$ is \'etale,  it  induces a bijective map 
 $U_\a(\VF) \to \rv \inv (f(\a))$.   But
 $X_\a(\VF) = U_\a(\VF)$, so the required isomorphism is proved.  
 
(3) If $X,Y$ are birationally equivalent, then $[X]_d \sim [Y]_d \  (WD^d_{<d})$,
while $X(\VF),Y(\VF)$ differ by   $\VF$-definable sets of dimension $<d$,
so $\ints (X(\VF)) \sim \ints (Y(\VF)) \  (SD_d)$.  Using the resolution of singularities
in the form:  every variety is birationally equivalent to a smooth nonsingular one, we are done by (2).  With a more complicated induction we should be able to dispense with
this use of Hironaka's theorem. 
\eprf

\lemm{comp3}   Let $V$ be a smooth projective  $\Kz$-variety,  $X,Y$ closed subvarieties, 
Let $F: V \m X \to V \m Y$    a biregular isomorphism.  Let $V_{\Oo},V_{\VF},V_{\k},F_{\VF}$,
etc. be the objects obtained by base change.  Then
$F_{\VF}$ induces a bijection  $V(\VF) \m X(\VF) \to V(\VF) \m Y(\VF)$, and
$$F_{\VF} (\rho_V \inv(X) \m X(\VF)) = \rho_V \inv(Y) \m Y(\VF)$$ \>{lem}

\prf  The first statement follows from the Lefschetz principle since $\VF$ is algebraically closed.

Since $V$ is projective, $V(\VF) = V(\Oo)$, and one can define for $v \in V$
the valuative distance $d(v,X)$, namely the greatest $\a \in \G$ such that the image of $x$ in $V(\Oo/\a)$
lies in $X(\Oo/\a)$.   

Let $\bF$ be the Zariski closure in $V^2$ of the graph of $F$. Then $\bF \meet
(V \m X) \times (V \m Y)$ is the graph of $F$.  In fact, in any algebraically closed
field $L$, we have
\beq{comp3s} \text{ if } a \in V(L) \m X(L) \text{ and } (a,b) \in \bF(L)  \text{ then } b \in V(L) \m Y(L) \eeq
and conversely. 

Suppose  for the sake of contradiction that in some $M \models ACVF_{\Kz}$ there exist
$a \in \rho_V \inv(X)$, $b \notin \rho_V \inv(Y)$, $(a,b) \in \bF$.   So $d(a,X) = \a >0$,
$d(b,Y) = 0$.  
Let $$C = \{\g \in \G: (\forall n \in \Nn) n  \g  <   \a \}$$
We may
assume by compactness that $C(M) \neq \emptyset$.  Let 
$$I = \{y \in \Oo(M):  \val(y) \notin C \}$$
so that $I$ is a prime ideal   of $\Oo(M)$.
Let $L$ be the field of fractions of $\Oo(M)/I$.  Let $\bar{a},\bar{b}$ be the images of $a,b$ in $L$.
Then $(\bar{a},\bar{b}) \in \bar{F}$, and $\bar{a} \in X$, $\bar{b} \notin Y$; contradicting 
\eqref{comp3s}
\eprf

\prf[Proof of \thmref{complement}]
  By \lemref{comp3}, there exists a definable bijection $ \rho_V \inv(X) \m X \to \rho_V \inv(Y) \m Y$. Applying $\ints: \K (\VF[n]) \to \K (\RV[\leq n] ) / \Isp$, and using   Lemmas \ref{comp2} and \ref{varclass}, we have 
    $[X]_n - [X]_d = [Y]_n - [Y]_d $.        Applying the first retraction 
    $ \K (\RV[\leq n] ) / \Isp \to \K(\RES[n])$
    of \thmref{retract2}, we obtain
$$[X_n] - [X \times \Aa^{n-d}]_n   =   [Y]_n - [Y \times \Aa^{n-d}]_n$$ 
in  $!\K(\RES[n]) =  \K(\Var_n)$.
  So
 $$  [X \times \Aa^{n-d} \du Y ]_n + [Z] = [Y \times \Aa^{n-d} \du X]_n + [Z]$$
 for some $Z$ with $\dim(Z) \leq n$, where now the equality is of classes in $\SG \Var_n$.
Counting birational equivalence classes of varieties of dimension $n$, 
we see that $X \times \Aa^{n-d}, Y \times \Aa^{n-d}$ must be  birationally equivalent. 
The last sentence follows from the lemma below.      \eprf

 \lemm{norational}  Let $X,Y$ be varieties containing no rational curve.
 Let $U$ be a   variety such that there exists a surjective morphism $\Aa^m \to U$.
 If $X \times U, Y \times U$ are birationally equivalent, then so are $X,Y$.  \>{lem}

\prf    For any variety $W$, let ${\mathcal F}(W)$ be the set of all rational maps
$g: \Aa^1   \to W$.  Write $\dom(g)$ for the maximal subset of $\Aa^1$ where
$g$ is regular; so $\dom(g)$ is cofinite in $\Aa^1$.   Let 
$R_W = \{(g(t),g(t')) \in W^2:    g \in {\mathcal F}(W),  t,t'  \in \dom(g) \}$.  Let $E_W$
be the equivalence relation generated by $R_W$, on points in the algebraic closure.
$R_W,E_W$ may not be constructible in general, but in the case we are concerned with
they are:    
 
\Claim{}  
Let $W \subseteq X \times U$ be a Zariski dense open set.  Let $\pi: W \to X$
be the projection.  Then $\pi(w)=\pi(w')$ iff $(w,w') \in E_W$ iff $(w,w') \in R_W$.  

\prf  If $g \in {\mathcal F}(U)$, then $\pi \circ g: \dom(g)  \to X$ is a regular map; hence
by assumption on $X$ it is constant.  It follows that if $(w,w') \in R_U$ then $\pi(w)=\pi(w')$,
and hence if $(w,w') \in E_U$ then $\pi(w)=\pi(w')$.  Conversely, assume $w',w'' \in W$ and
$\pi(w')=\pi(w'')$;
then $w'=(x,u'),w'=(x,u'')$ for some $x \in X$, $u',u'' \in U$.  Let $U_x = \{u \in U: (x,u) \in 
W\}$.  Since $W$ is open, $U_x$ is open in $U$.  Let $h: \Aa^m \to U$ be a surjective
morphism; let $h(v')=u',h(v'')=u''$. The line through $v',v''$ intersects $h \inv (U_x)$
in a nonempty open set.  This gives a regular map $f$ from the affine  lines, minus finitely
many points, into $U$, passing through $u',u''$.  So $t \mapsto (x,f(t))$ gives
a rational map from $\Aa^1$ to $W$, passing through $(w',w'')$; and so $(w',w'') \in R_U$
and certainly in $E_U$.  \eprf

Using the Claim, we   prove the lemma.   
Let $W_X \subseteq X \times U$ , $W_Y \subseteq Y \times U$ be Zariski dense open, 
and $F:W_X \to W_Y$ a biregular  isomorphism.  Then  $F$ takes $E_{W_X}$ to $E_{W_Y}$.
Moving now to the  category of constructible sets and maps,
quotients by constructible equivalence relations exist, and $W_X/ E_{W_X}$ is isomorphic
as a constructible set to $W_Y/ E_{W_Y}$.  Let 
 $\pi_X: W_X \to X, \pi_Y: W_Y \to Y$ be the 
 projections.  By the Claim, $W_X/ E_{W_X}  = \pi_X(W_X) =:X'$.  Similarly $W_Y/ E_{W_Y}  = \pi_Y(W_Y)=:Y'$.  Now since $W_X,W_Y$ are Zariski dense, so are $X',Y'$.  Thus
 $X,Y$ contain isomorphic Zariski dense constructible sets, so they are birationally
 equivalent.    \eprf
 
{\bf Remark}  The condition on $X,Y$ may be weakened to the statement that they contain 
no rational curve through a generic point; i.e. that 
there exist proper subvarieties $(X_i: i \in I)$ defined over $\Kz$, such that 
for any field $L \supset \Kz$, any rational curve on $X \times_\Kz L$ is contained
in some $X_i \times _\Kz L$.   
  
\>{section}

\end{document}